 \documentclass[11pt]{article}
 \usepackage{fullpage}
\date{}

\usepackage{tikz}
\usetikzlibrary{calc}
\usetikzlibrary{matrix}
\usepackage{amsmath,amssymb,amsfonts,amsthm,subfigure}
\usepackage{color}                    
\usepackage{hyperref}                 

\usepackage{graphicx}
\allowdisplaybreaks

\DeclareMathOperator{\tr}{tr}

\newcommand{\eps}{{\epsilon}}
\newcommand{\bey}{\begin{eqnarray}}
\newcommand{\eey}{\end{eqnarray}}

\newcommand{\beq}{\begin{equation}}
\newcommand{\eeq}{\end{equation}}
\theoremstyle{plain}

\theoremstyle{definition}

\theoremstyle{remark}

\title{A Study on Phase-Field Models for Brittle Fracture}

\author{Fei Zhang%
\thanks{College of Petroleum Engineering,  China University of Petroleum -- Beijing, 18 Fuxue Road, Changping,
Beijing 102249, China ({\it fzhang\_cup@outlook.com})},
\and Weizhang Huang%
\thanks{Department of Mathematics, the University of Kansas, Lawrence, Kansas 66045, U.S.A.
({\it whuang@ku.edu})},
\and Xianping Li%
\thanks{Department of Mathematics and Statistics, University of Missouri -- Kansas City,
5120 Rockhill Road, Kansas City, Missouri 64110, U.S.A. ({\it lixianp@umkc.edu})},
\and Shicheng Zhang%
\thanks{College of Petroleum Engineering,  China University of Petroleum -- Beijing, 18 Fuxue Road, Changping,
Beijing 102249, China ({\it zhangsc@cup.edu.cn})}
}

\begin{document}
\vskip 1cm
\maketitle
\begin{abstract}
In the phase-field modeling of brittle fracture, anisotropic constitutive assumptions for the degradation of stored elastic energy due to fracture are crucial to preventing cracking in compression and obtaining physically sound numerical solutions. Three energy decomposition models, the spectral decomposition, the volumetric-deviatoric split, and
an improved volumetric-deviatoric split,  and their effects on the performance of the phase-field modeling are studied. 
Meanwhile, anisotropic degradation of stiffness may lead to a small energy remaining on crack surfaces, which violates crack boundary conditions and can cause unphysical crack openings and propagation. A simple yet effective treatment for this is proposed: define a critically damaged zone with a threshold parameter and then degrade both the active and passive energies in the zone. A dynamic mesh adaptation finite element method is employed for the numerical solution of the corresponding elasticity system.  Four examples, including two benchmark ones, one with complex crack systems, and one based on an experimental setting, are considered. Numerical results show that the spectral decomposition and improved volumetric-deviatoric split models, together with the improvement treatment of crack boundary conditions, can lead to crack propagation results that are comparable with the existing computational and experimental results. It is also shown that the numerical results are not very sensitive to the parameter defining the critically damaged zone. 
\end{abstract}




\section {Introduction}

In recent years, the phase-field model for brittle fracture based on the variational approach of Francfort and Marigo \cite{FM98} has become a commonly used numerical simulation technique for engineering designs because it can handle complex cracks and crack initiation and propagation more easily than other methods. The basic idea of the phase-field modeling is to describe cracks by a continuous scalar field variable $d$, which is used to indicate whether the material is damaged or not. This variable $d$ depends on a parameter $l$ describing the actual width of the smeared cracks and has a value of 0 or close to 0 near the cracks and 1 away from the cracks.
There are three major advantages of the phase-field modeling for brittle fracture over other methods. Firstly, the behavior of the crack is completely determined by a coupled system of partial differential equations (PDEs) based on the energy functional. Therefore, additional calculations such as stress-intensity factors are not required to determine the crack initiation and propagation. Secondly, complex fracture networks can be easily handled since crack merging and branching do not require explicitly keeping track of fracture interfaces. Thirdly, sharp but smooth interfaces (instead of discontinuities) can be introduced into the displacement field.

Since it was first proposed by Bourdin et al. \cite{BFM00, FM98}, the phase-field modeling for brittle fracture has attracted considerable attention and significant progress has been made; e.g., see \cite{AMM09, Bor12, BHLV14, KM10, MVB15, MWH10, VMBV14}. However challenges still exist. In the phase-field modeling, constitutive assumptions for the degradation of energy due to fracture can be categorized into two groups, isotropic models
and anisotropic models. In the former group, the degradation function acts on the whole stored bulk energy, which means that energy is released due to fracture in both tension and compression. Thus, crack propagation may also arise under compressive load state, which is physically unrealistic. On the other hand, in order to overcome this unphysical feature, the elastically stored energy is decomposed into active and passive parts where only the former
is degraded by the phase field. Two commonly used energy decomposition models have been proposed in the past. Miehe et al. \cite{MWH10} introduced a fully anisotropic constitutive model for the degradation of energy based on the spectral decomposition of strains with the assumption that crack evolution is induced by the positive principal strains. The other model is the unilateral contact model proposed by Amor et al. \cite{AMM09} that splits the strain into volumetric and deviatoric parts, with the expansive volumetric part and the total deviatoric part being degraded. Since the choice of the energy splitting controls the energy contribution in the damage evolution, different splitting models can significantly affect numerical approximations in the phase-field modeling of cracks. In this work we shall study these models plus an improved version of the unilateral contact model. 

In the phase-field modeling, a pre-existing crack is used to initialize crack which can be modeled as a discrete crack in the geometry or an induced crack in the phase-field. For the former, it has been successfully applied in phase-field models for single initial crack. However, it is difficult to hand complex crack boundary conditions since the location of the initial crack is mesh-dependent. For the latter, an initial strain-history field is introduced to define the location of the induced crack. One of the major advantages of this treatment is that initial cracks can be placed anywhere in the domain without referring to the mesh, which makes it possible to deal with complex initial cracks.
The induced crack model was first proposed by Borden et al. \cite{Bor12} and significant improvements have been made \cite{Bor12,MVB15,SBS15}. 
However, there is a small energy remaining in the totally damaged zone due to the anisotropic degradation of stiffness. For the induced crack model, stress remains in the interior of the initial crack and increases with external loads before the crack begins to propagate. This violates the vanishing stress condition on the crack surface and often results in unphysical crack propagation.
May et al. \cite{MVB15} have observed that with the induced crack setting for a single notched shear test, numerical results are very different from those with discrete crack boundary conditions. To overcome this problem, Strol and Seelig \cite{SS15,SS16} proposed a novel treatment of crack boundary conditions in which crack orientation is taken into account so that both the positive normal stress on the crack surfaces and the shear stress along the frictionless crack surface vanish. However, the establishment of this constitutive assumption is not based on the variational approach, which makes the phase-field model more complicated to implement.

Another important issue is that the phase-field modeling approximates the original discrete problem as $l \to 0$ under the condition that $h \ll l$ or at least $h < l$, where $h$ denotes the mesh spacing. A very fine mesh is needed to fulfill the condition if only a uniform mesh is used in the computation, which increases the computational cost significantly. Moreover, cracks can propagate under continuous load. Thus, it is natural to use a dynamic mesh adaptation strategy to improve the efficiency of the simulation.  In this work we use the moving mesh PDE (MMPDE) method \cite{BHR09,HRR94a, HRR94b, HR11} to concentrate mesh elements around evolving cracks. The reader is referred to Zhang et al. \cite{ZHLZ17} for the study of the application of the MMPDE method to the phase-field modeling of brittle fracture.

The objective of this paper is twofold. The first is to investigate two commonly used energy decomposition
models, the spectral decomposition model \cite{MHW10} and the volumetric-deviatoric
split model \cite{AMM09}, and their impacts on the performance of the phase-field modeling for brittle fracture.
The former degrades the energy related to the tensile strain component while the latter degrades
both the expansive volumetric strain energy and the total deviatoric strain energy.
In the latter case, the compressive deviatoric strain energy is also accounted for contributing
to the damage process, which may cause unphysical crack propagation.
An improved volumetric-deviatoric model is proposed and studied to avoid this difficulty.

The second goal is to study the treatment of crack boundary conditions for the induced crack model.
As mention before, the induced crack model leads to a small remaining stress in the damaged area,
which can cause unrealistic crack boundary conditions such as normal stress remaining on the crack surface.
To overcome this difficulty, we propose a simple yet effective treatment of crack boundary conditions: 
define a critically damaged zone with $0 < d < d_{cr}$, where $d_{cr}$ is a positive parameter, and then
degrade both the active and passive components of the energy in this zone. We shall consider
four two-dimensional examples to verify the treatment. The first two are classical benchmark problems,
single edge notched tension and shear tests. The third example is designed to demonstrate
the ability of the models to handle complex cracks. The last example is also a single edge notched shear
test but with the physical parameters and domain geometry chosen based on an experiment setting.
Numerical results show that the spectral decomposition and improved volumetric-deviatoric models
together with the improved treatment of crack boundary conditions lead to correct crack propagation
for all of the examples.

The present work is organized as follows. Section~\ref{SEC:phase-field} is devoted to the description of the phase-field modeling with three energy decomposition models and the improved treatment of crack boundary conditions.  A moving mesh finite element method based on the MMPDE moving mesh method is described in Section~\ref{SEC:mmfem}
for the elasticity system. Numerical results for four two-dimensional examples are presented in Section~\ref{SEC:numerics}. Finally, Section~\ref{SEC:conclusions} contains conclusions.

\section {Phase-field models for brittle fracture}
\label{SEC:phase-field}
\subsection{Variational approach to elastic models with cracks}

We consider small strain isotropic elasticity models with cracks.
Let $\Omega$ be a two-dimensional bounded domain filled with elastic material and having
the boundary $\partial \Omega = \partial \Omega_t \cup \partial \Omega_u$, where the surface
traction $\overline t$ is specified on $\partial \Omega_t$ and the displacement $\overline u$
is given on $\partial \Omega_u$. Denote the displacement by $u$. Then the strain tensor is given by
\[
\eps = \frac{1}{2} \left (\nabla u + (\nabla u)^T\right ),
\]
where $\nabla {u}$ is the displacement gradient tensor.
For isotropic material, the strain energy density
is given by Hooke's law as
\begin{equation}
\Psi_{e}(\eps) = \frac{\lambda}{2} \left (\tr({\eps})\right )^2 + \mu \tr({\eps}^2) ,
\label{energy-el}
\end{equation}
where $\lambda$ and $\mu$ are the Lam\'{e} constants and $\tr(\cdot)$ denotes the trace of a tensor.

For an elastic body with a given crack $\Gamma$, we take the approach of
Francfort and Marigo \cite{FM98} to define the total energy as
\begin{equation}
\mathcal{W}(\eps,\Gamma) = \mathcal{W}_{e}(\eps,\Gamma)
+ \mathcal{W}_{c}(\Gamma)
\equiv \int_{\Omega \setminus \Gamma} \Psi_{e}(\eps) \, d \Omega + \int_{\Gamma} g_c \, dS,
\label{energy-total}
\end{equation}
where $\mathcal{W}_{e}(\eps,\Gamma)$ represents the energy stored in the bulk of the elastic body, $\mathcal{W}_{c}(\Gamma)$ is the energy input required to create the crack according to the Griffith criterion, and $g_c$ is the fracture energy density that is the amount of energy needed to
create a unit area of fracture surface.

\begin{figure} 
\centering 
\subfigure[Regularized crack $\Gamma_l(d)$ by the phase-field approximation]{\label{fig:subfig:Original_Crack}
\includegraphics[width=0.45\linewidth]{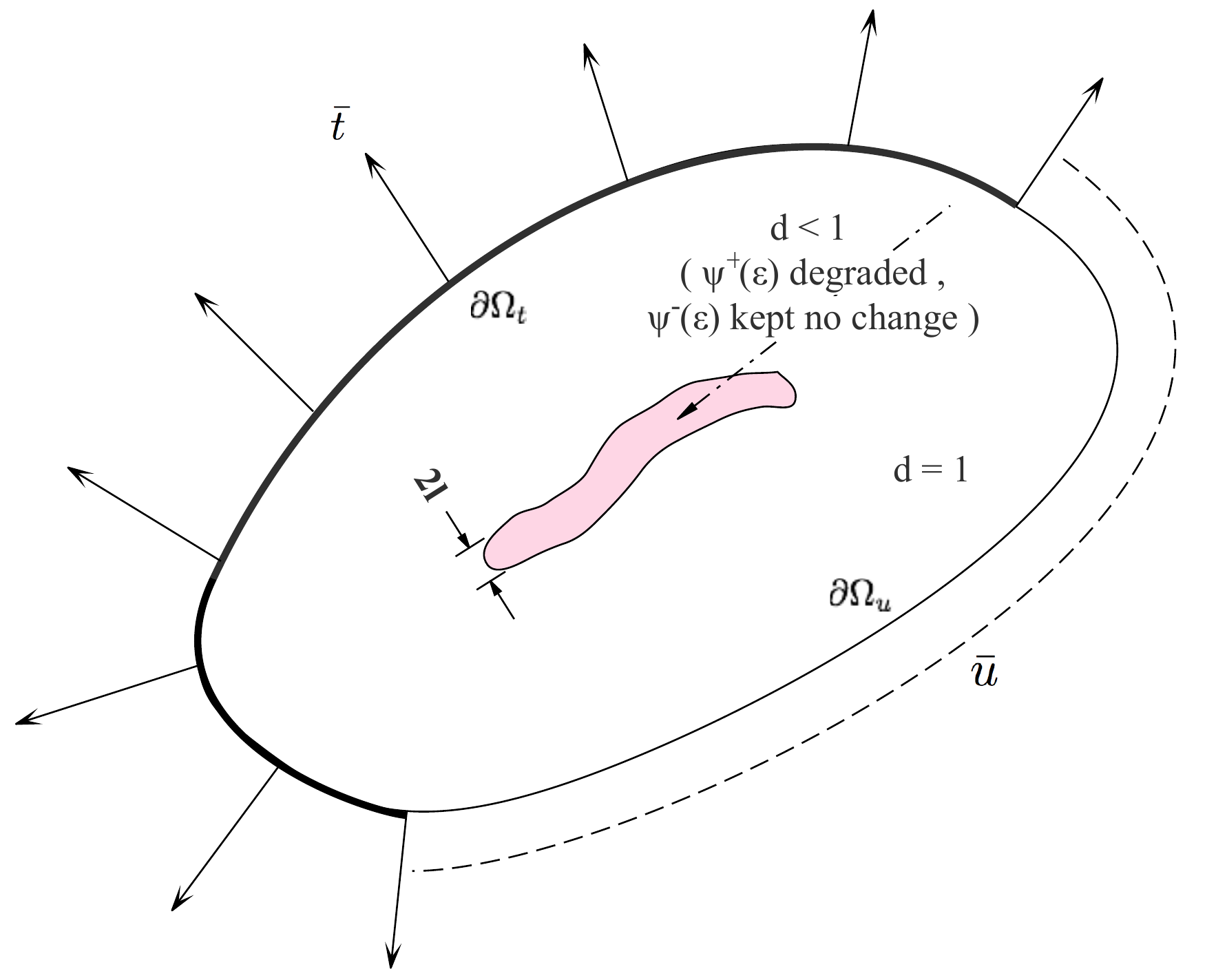}}
\hspace{10mm}
\subfigure[The critically damaged region indicated by the critical value $d_{cr}$]{\label{fig:subfig:Modified_Crack}
\includegraphics[width=0.45\linewidth]{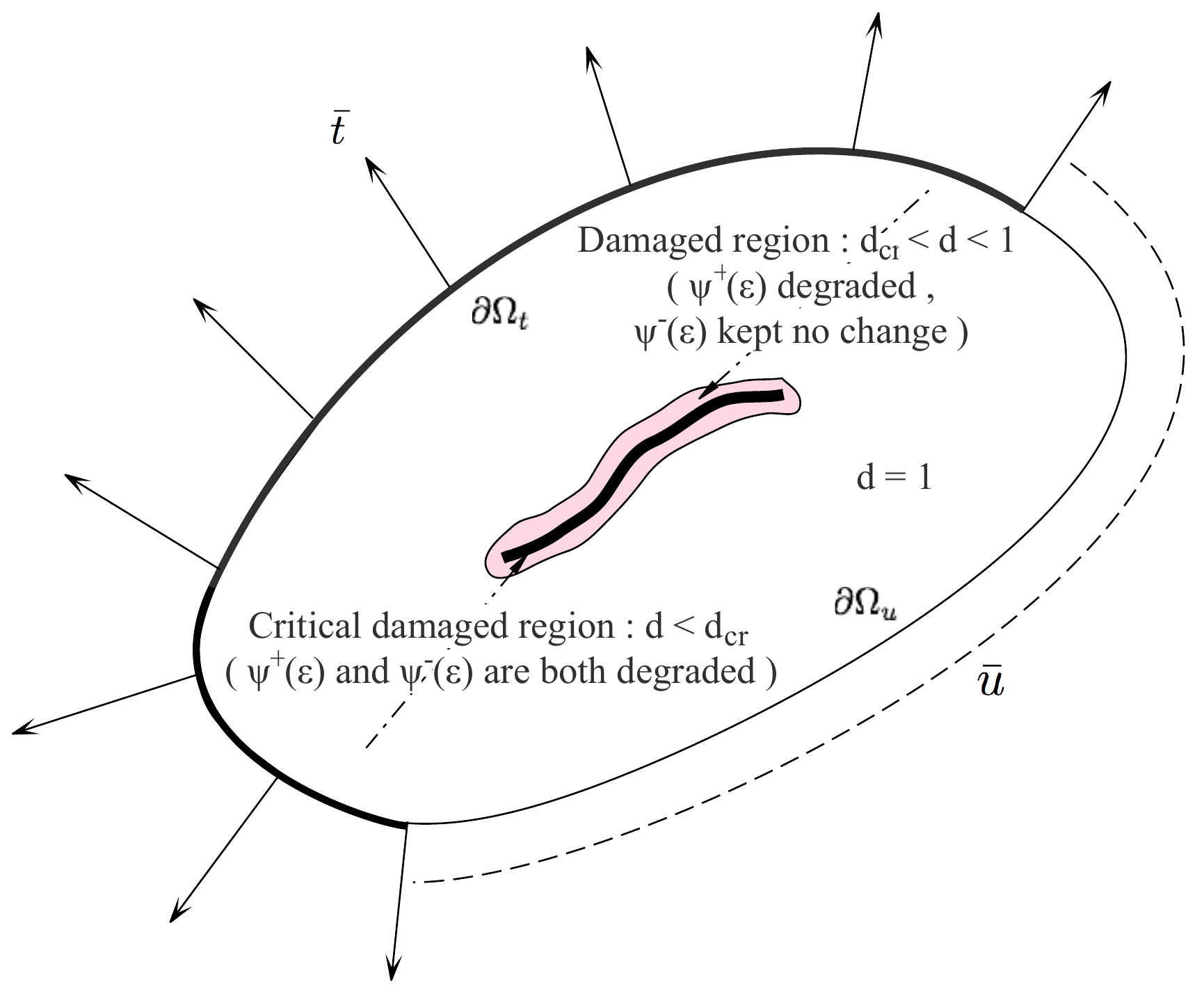}} 
\caption{A sketch of the phase-field modeling of brittle crack.}
\label{sketch_crack}
\end{figure}

The model (\ref{energy-total}) treats the crack $\Gamma$ as a discontinuous object, which has been a challenge
in numerical simulation. Here we use the phase-filed approach with which $\Gamma$ is smeared into
a continuous object. Specifically, a phase-field variable $d(x,t)$ is used to represent $\Gamma$ along with
the parameter $l$ that describes the width of the smooth approximation of the crack. 
The value of $d(x,t)$ is zero or close to zero near the crack and one away from
the crack (see Fig. \ref{fig:subfig:Original_Crack}). 
The fracture energy $\mathcal{W}_{c}(\Gamma)$ is approximated
by the smeared total fracture energy \cite{BFM00} as
\begin{equation}
\label{energy-c}
\mathcal{W}_{c}^{l}  = \int_\Omega \frac{g_c}{4 l} \left ((d-1)^2 + 4 l^2 |\nabla d|^2 \right ) d \Omega.
\end{equation}

The elastic energy needs to be modified to reflect the loss of material stiffness in the damaged zone.
A degradation function $g(d)$ is commonly used for this purpose,
which is required to satisfy the following property
\[
\begin{cases}
g(0) = 0: &\quad \text{Damage occurred for $d = 0$ and this part should vanish}; \\
g(1) = 1: &\quad \text{No damage occurs for $d = 1$}; \\
g'(0) = 0: &\quad \text{No more changes after the fully broken state}; \\
g'(1) \ne 0: &\quad \text{The damage has to be initiated at the onset.}
\end{cases}
\]
In our computation, we use a common choice $g(d) = d^2$.
A straightforward attempt for the modification of the elastic energy is to degrade the total stored energy
in the damaged region, i.e., 
\begin{equation}
\Psi_{e}^{l}(\eps,d) = g(d) \Psi_{e}(\eps) .
\label{split-0}
\end{equation}
Unfortunately, this can result in crack opening in compressed regions, which is physically unrealistic.
To avoid this, it is common to decompose the total stored energy into two components as
\begin{equation}
\Psi_{e}(\eps) = \Psi_{e,act}(\eps) + \Psi_{e,pas}(\eps),
\label{split-1}
\end{equation}
where $\Psi_{e, act}(\eps)$ and $\Psi_{e, pas}(\eps)$ represent the active and passive components, respectively,
with only the former contributing to the damage process that results in fracture.
(Two commonly used decomposition models will be discussed in the next subsection.)
Then, the damage model reads as
\begin{equation}
\Psi_{e}^{l}(\eps, d) =  g(d) \Psi_{e,act}(\eps) +  \Psi_{e,pas}(\eps),
\label{energy-el-d}
\end{equation}
and the corresponding total energy is
\begin{equation}
\mathcal{W}^l =  \mathcal{W}_e^{l} + \mathcal{W}_c^l = 
 \int_{\Omega} \left( (d^2+k_l)\Psi_{e,act}(\eps) + \Psi_{e,pas}(\eps) + 
 \frac{g_c}{4l}\left ((d-1)^2+4l^2 |\nabla d|^2 \right ) \right) d \Omega,
\label{energy-d}
\end{equation}
where $k_l \ll l$ is a regularization parameter added to avoid degeneracy. 
It is emphasized that with this setting, {\em only the active component of the elastic energy
is degraded in damaged regions.}

Letting
\[
W = (d^2+k_l)\Psi_{e,act}(\eps) + \Psi_{e,pas}(\eps) + \frac{g_c}{4l}\left ((d-1)^2+4l^2 |\nabla d|^2 \right ),
 \]
we obtain the variation of the energy as
\[
\delta \mathcal{W}^l =
\int_{\Omega} \frac{\partial W}{\partial d} \delta d \; d \Omega +
\int_{\Omega} \frac{\partial W}{\partial \nabla d} \cdot \nabla \delta d \; d \Omega +
\int_{\Omega} \frac{\partial W}{\partial \epsilon} : \epsilon(\delta u) \; d \Omega ,
\]
where $A:B$ is the inner product of tensors $A$ and $B$, i.e., 
$A:B = \sum\limits_{i,j} A_{i,j} B_{i,j}$.
Define the function spaces
\begin{gather*}
V_{u} = \left \{ \; \varphi \; | \; \varphi \in H^1(\Omega),\;  \varphi = \overline u \;
\text{on} \; \partial \Omega_u \right \},\\
V_{u} ^0= \left \{ \; \varphi \; | \; \varphi \in H^1(\Omega),\; \varphi = 0 \; \text{on} \;
\partial \Omega_u \right \},
\end{gather*}
where $H^1(\Omega)$ is a Sobolev space defined as
\[
H^1(\Omega) = \left \{ \; \varphi \; | \;  \int_{\Omega} \varphi^2 \, d \Omega < +\infty,\;
\int_{\Omega} |\nabla \varphi |^2 \, d \Omega < + \infty \right \}.
\]
Take $V_d = H^1(\Omega)$. 
The weak formulation for the phase-field model is to find $d \in V_d$ and $u \in V_{u}$ such that 
\begin{align}
& \int_\Omega \left (
 \left (2 d \mathcal{H} + \frac{g_c (d-1) }{2 l} \right ) \delta d
 +  2 g_c l \nabla d \cdot \nabla \delta d \right) d \Omega = 0 ,
 \quad \forall\; \delta d \in V_d 
\label{weak-d}
\\
& \int_\Omega \sigma : \epsilon(\delta u) d \Omega = 
\int_{\Omega_t} \overline{t} \cdot \delta {u} \; dS +
\int_\Omega f \cdot \delta {u} \; d \Omega ,
\quad \forall\; \delta u \in V_u^0
\label{weak-u}
\end{align}
where $\mathcal{H} = \Psi_{e,act}(\epsilon)$ and $\sigma$ is the Cauchy stress defined as
\begin{equation}
\sigma \equiv \frac{\partial W}{\partial \eps} = 
(d^2+k_l) \frac{\partial \Psi_{e,act}}{\partial \eps} +
\frac{\partial \Psi_{e,pas}}{\partial \eps} .
\label{sigma}
\end{equation}
To ensure that cracks can only grow (i.e., crack irreversibility),
we replace $\mathcal{H} = \Psi_{e,act}(\eps)$ in (\ref{weak-d}) by
\begin{equation}
\label{eqn-H}
\mathcal{H} = \max_{s \le t} \Psi_{e,act}(\eps)(s),
\end{equation}
where $t$ is the quasi-time corresponding to the load increments.

\subsection {Models for energy decomposition}

In this subsection we discuss two commonly used models and an improved one for energy decomposition.

\label{SEC:Various decomposition}
\subsubsection {Spectral decomposition}
\label{SEC:spectral}

We first consider a commonly-used model proposed by Miehe et al. \cite{MHW10}
based on the spectral decomposition of the strain tensor.
To this end, we define the positive-negative decomposition of a scalar function $f$ as
\[
f = f^+ + f^-,\quad f^+ = \frac{f+|f|}{2}, \quad f^- = \frac{f-|f|}{2}.
\]
For the strain tensor, we have
\[
\eps = \eps^+ + \eps^-,
\quad \eps^+ = Q \text{diag}(\lambda_1^+, ..., \lambda_n^+) Q^T,
\quad \eps^- = Q \text{diag}(\lambda_1^-, ..., \lambda_n^-) Q^T,
\]
where $Q \text{diag}(\lambda_1, ..., \lambda_n) Q^T$ is the eigen-decomposition
(or the spectral decomposition) of $\eps$. Notice that $\eps^+$ represents the tensile strain component
that contributes to the damage process resulting in crack initiation and propagation while
$\eps^-$ represents the compression strain component that does not contribute to the damage process.
Based on this, the active and passive components of the elastic energy are given by
\[
\Psi_{e,act}(\eps) = \frac{\lambda}{2} \left ( (\tr(\epsilon))^+\right )^2 + \mu \tr((\epsilon^+)^2),
\quad
\Psi_{e,pas}(\eps) = \frac{\lambda}{2} \left ( (\tr(\epsilon))^-\right )^2 + \mu \tr((\epsilon^-)^2).
\]
The elastic energy density in \eqref{energy-el-d} and the Cauchy stress in \eqref{sigma} can be written as
\[
\Psi_{e}^{l}(\eps,d) = (d^2 + k_l) 
\left( \lambda \frac{((\tr(\eps))^+)^2}{2} + \mu \tr \left( (\eps^+)^2 \right) \right) +
\left( \lambda \frac{((\tr(\eps))^-)^2}{2} + \mu \tr \left( (\eps^-)^2 \right) \right),
\]
\[
\sigma = (d^2+k_l)
\left (\lambda (\tr(\epsilon))^+ I + 2 \mu \epsilon^+ \right ) +
\left (\lambda (\tr(\epsilon))^- I + 2 \mu \epsilon^- \right ).
\]
In this model, {\em only the energy component related to the tensile strain component
is degraded in damaged regions}.

\subsubsection{Volumetric-deviatoric split}
\label{SEC:v-d}

Another commonly-used model is proposed by Amor et al. \cite{AMM09} based on
the volumetric-deviatoric split (v-d split),
\begin{equation}
\label{v-d decomposition}
\epsilon = \epsilon_{S} + \epsilon_{D},
\quad
\epsilon_{S} = \frac{1}{m} \tr(\epsilon) I,
\quad
\epsilon_{D} = \epsilon - \frac{1}{m} \tr(\epsilon) I,
\end{equation}
where $m$ is the spatial dimension, $\epsilon_{S}$ is the \textit{spherical} component 
(called the volumetric strain tensor) related to the volume change,
and $\epsilon_{D}$ is the \textit{deviatoric} component (called the strain deviator tensor)
related to distortion. In this model, {\em the expansive volumetric strain energy
and the total deviatoric strain energy are released by the creation of new cracks
whereas the compressive volumetric strain energy is not}. The active and passive components
of the elastic energy are given by
\[
\Psi_{e,act}(\eps) = \kappa_0\frac{((\tr(\eps))^+)^2}{2} + \mu \tr \left( (\eps_{D})^2 \right),
\qquad
\Psi_{e,pas}(\eps) = \kappa_0\frac{((\tr(\eps))^-)^2}{2},
\]
where $\kappa_0 = \lambda + 2 \mu/m$ is the bulk modulus of the material.
The energy density and Cauchy stress associated with this model are 
\[
\Psi_{e}^{l}(\eps,d) = (d^2 + k_l) 
\left( \kappa_0\frac{((\tr(\eps))^+)^2}{2} + \mu \tr \left( (\eps_{D})^2 \right) \right) +
\kappa_0 \frac{\left((\tr(\eps))^-\right)^2}{2},
\]
\[
\sigma = (d^2 + k_l)
\left( \kappa_0 \left( \tr(\eps) \right)^+ I + 2 \mu \eps_D \right) +
\kappa_0 \left( \tr(\eps) \right)^- I.
\]

\subsubsection {Improved volumetric-deviatoric split}
\label{SEC:improved}

According to unilateral contact conditions, crack can only open in the regions where the material tends
to expand. For the volumetric-deviatoric split model described in the previous subsection,
both the positive part of the volumetric component and the total deviatoric component of the elastic energy
are degraded in damaged regions. However, the three principal strains of the deviatoric component of the strain tensor
can be negative, which indicates that the compressive deviatoric strain energy is also accounted for
contributing to the damage process in the v-d split model.
To avoid this unphysical feature, we propose to apply the spectral decomposition to the strain deviatoric tensor
and call the resulting model as the improved v-d split model. Specifically, 
the spectral decomposition of $\eps_D$ is
\[
\eps_D = \eps_D^+ + \eps_D^-,
\qquad
\eps_{D}^+ = Q \text{diag}(\lambda_1^+, ..., \lambda_m^+) Q^T,
\qquad
\eps_{D}^- = Q \text{diag}(\lambda_1^-, ..., \lambda_m^-) Q^T,
\]
provided that $Q \text{diag}(\lambda_1, ..., \lambda_m) Q^T$ is the eigen-decomposition of $\eps_D$.
{\em The degradation applies only to the expansive volumetric and deviatoric strain energies.}
Then,  the active and passive part of the elastic energy are given by
\[
\Psi_{e,act}(\eps) = \kappa_0\frac{((\tr(\eps))^+)^2}{2} + \mu \tr \left( (\eps_{D}^+)^2 \right),
\qquad
\Psi_{e,pas}(\eps) = \kappa_0\frac{((\tr(\eps))^-)^2}{2} + \mu \tr \left( (\eps_{D}^-)^2 \right),
\]
and the corresponding energy density and Cauchy stress are 
\[
\Psi_{e}^{l}(\eps,d) = (d^2 + k_l) 
\left( \kappa_0\frac{((\tr(\eps))^+)^2}{2} + \mu \tr \left( (\eps_{D}^+)^2 \right) \right) +
\left( \kappa_0\frac{((\tr(\eps))^-)^2}{2} + \mu \tr \left( (\eps_{D}^-)^2 \right) \right),
\]
\[
\sigma = (d^2 + k_l) 
\left( \kappa_0 \left( \tr(\eps) \right)^+ I + 2 \mu \eps_D^+ \right) +
\left( \kappa_0 \left( \tr(\eps) \right)^- I + 2 \mu \eps_D^- \right).
\]

\subsection{Improved treatment of crack boundary conditions}
\label{SEC:ItCBC}

\begin{figure} 
\centering 
\subfigure[Original model: $\Psi_{e}^{l}$]{\label{fig:subfig:Original_Energy}
\includegraphics[width=0.4\linewidth]{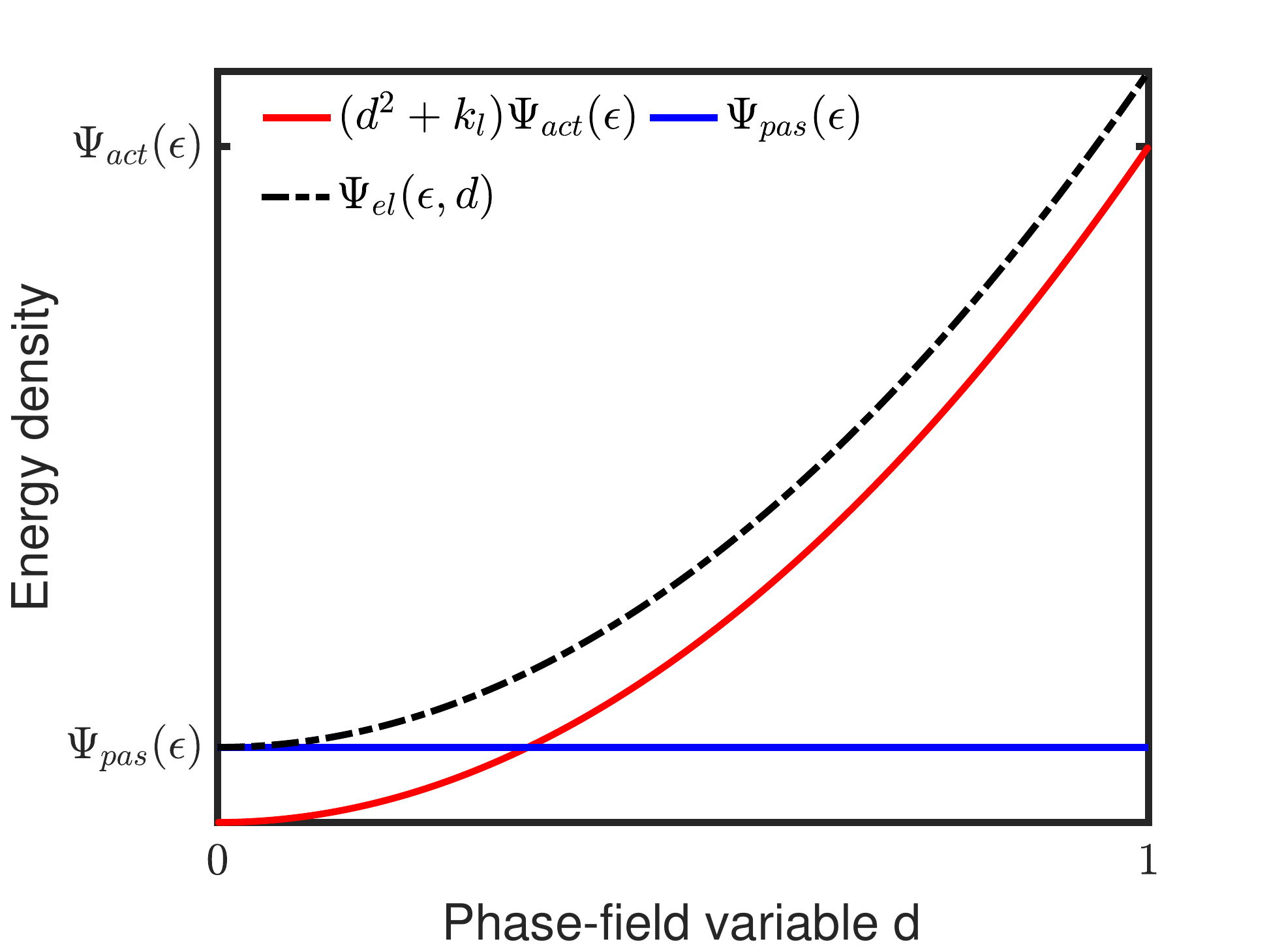}}
\subfigure[Modified model: $\Psi_{e,m}^{l}$]{\label{fig:subfig:Modified_Energy}
\includegraphics[width=0.4\linewidth]{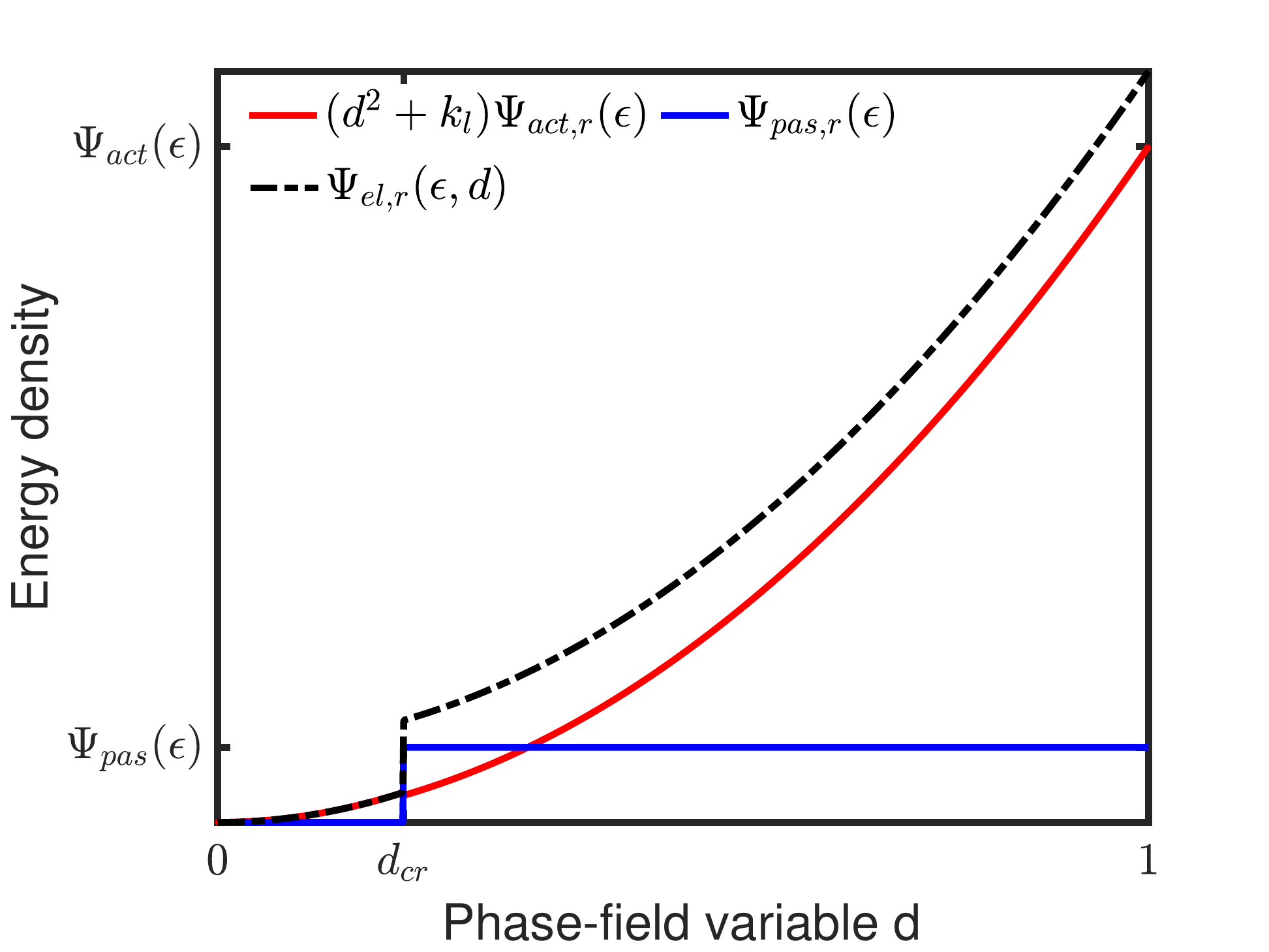}}
\subfigure[Comparison of the $\Psi_{e}^{l}$ and $\Psi_{e,m}^{l}$]{\label{fig:subfig:Comparison_Energy}
\includegraphics[width=0.4\linewidth]{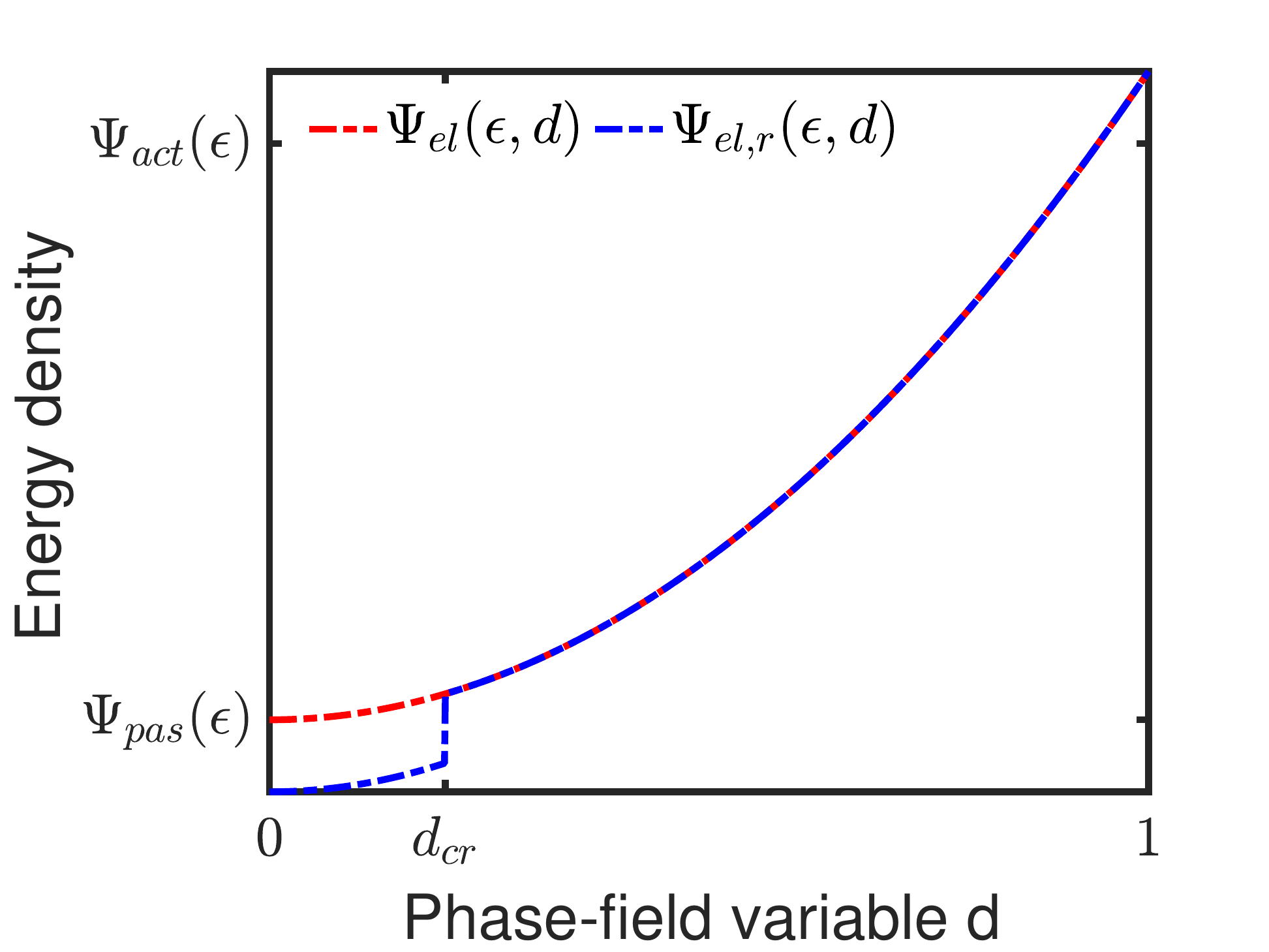}}
\caption{Energy stored in the bulk of the elastic body.}
\end{figure}

\begin{figure}
\centering 
\subfigure[Tension test]{\label{fig:subfig:Tension}
\includegraphics[width=0.4\linewidth]{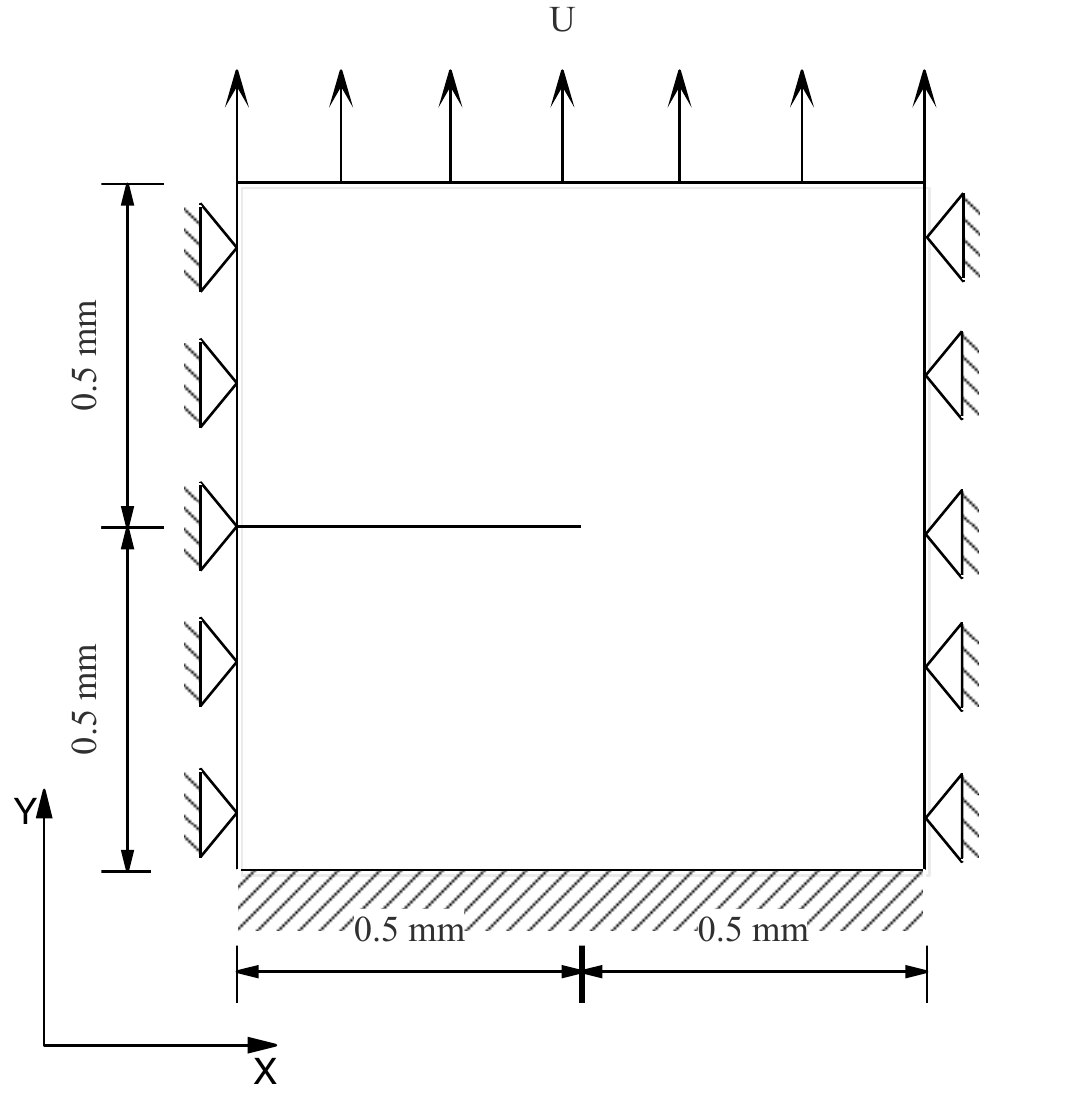}}
\subfigure[Shear test]{\label{fig:subfig:Shear}
\includegraphics[width=0.4\linewidth]{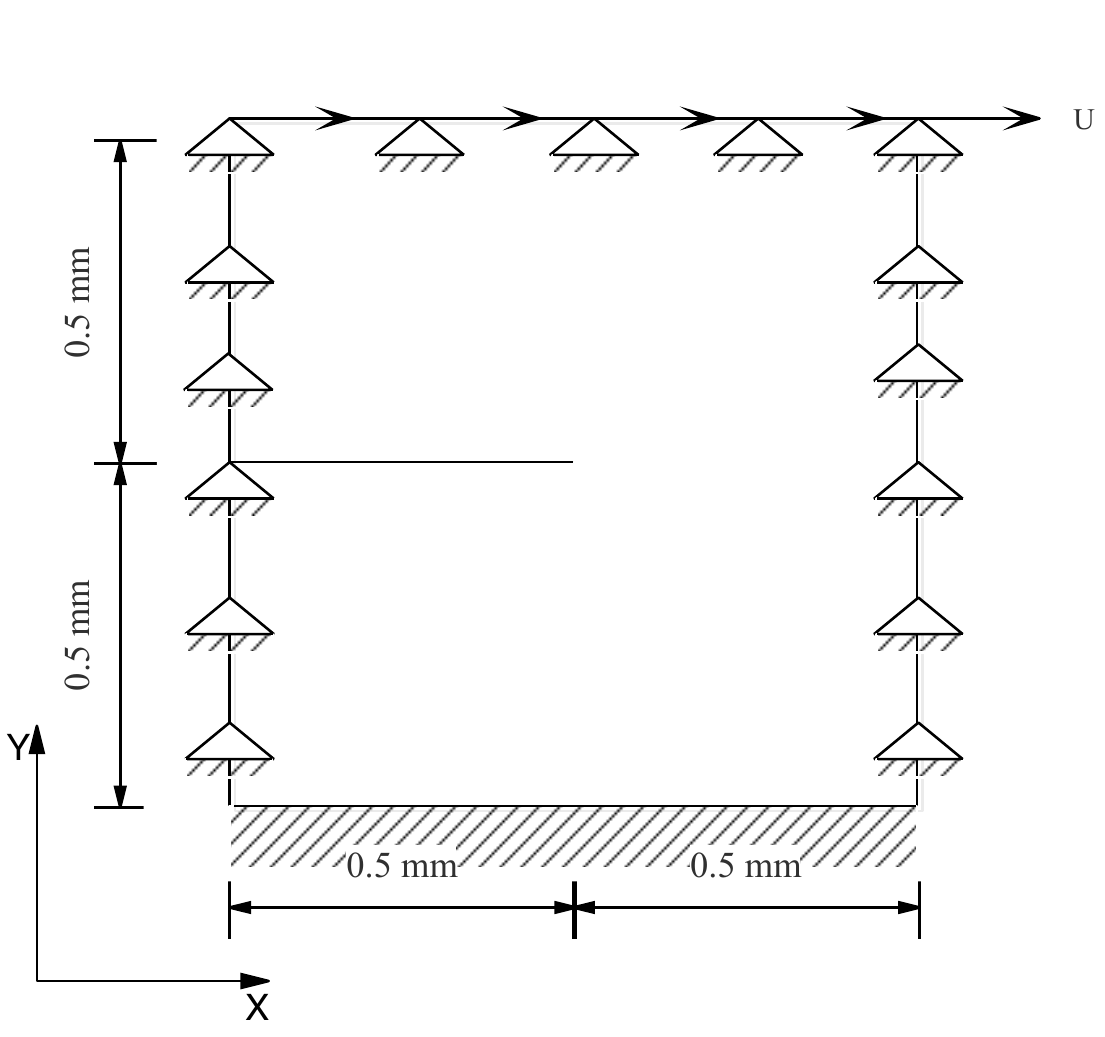}}
\caption{Domain and boundary conditions for single edge notched tests,
(a) tension test for Example 1 and (b) shear test for Example 2.}
\label{fig:sketch}
\end{figure}

\begin{figure}
\centering 
\subfigure[Mesh on deformed domain]{\label{fig:subfig:Original_Mesh}
\includegraphics[width=0.3\linewidth]{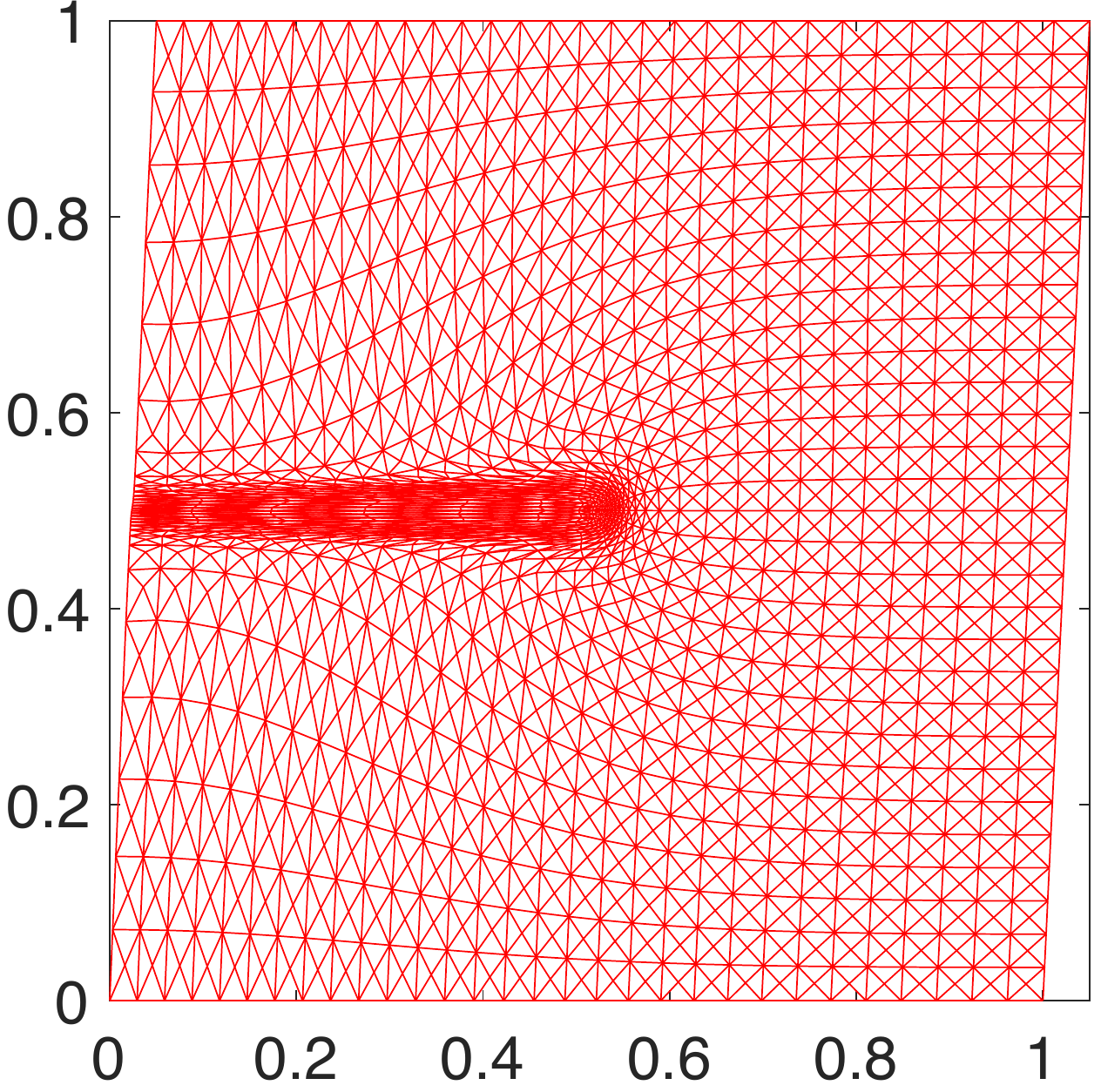}}
\subfigure[Von Mises stress distribution on the deformed domain]{\label{fig:subfig:Original_Sigma}
\includegraphics[width=0.35\linewidth]{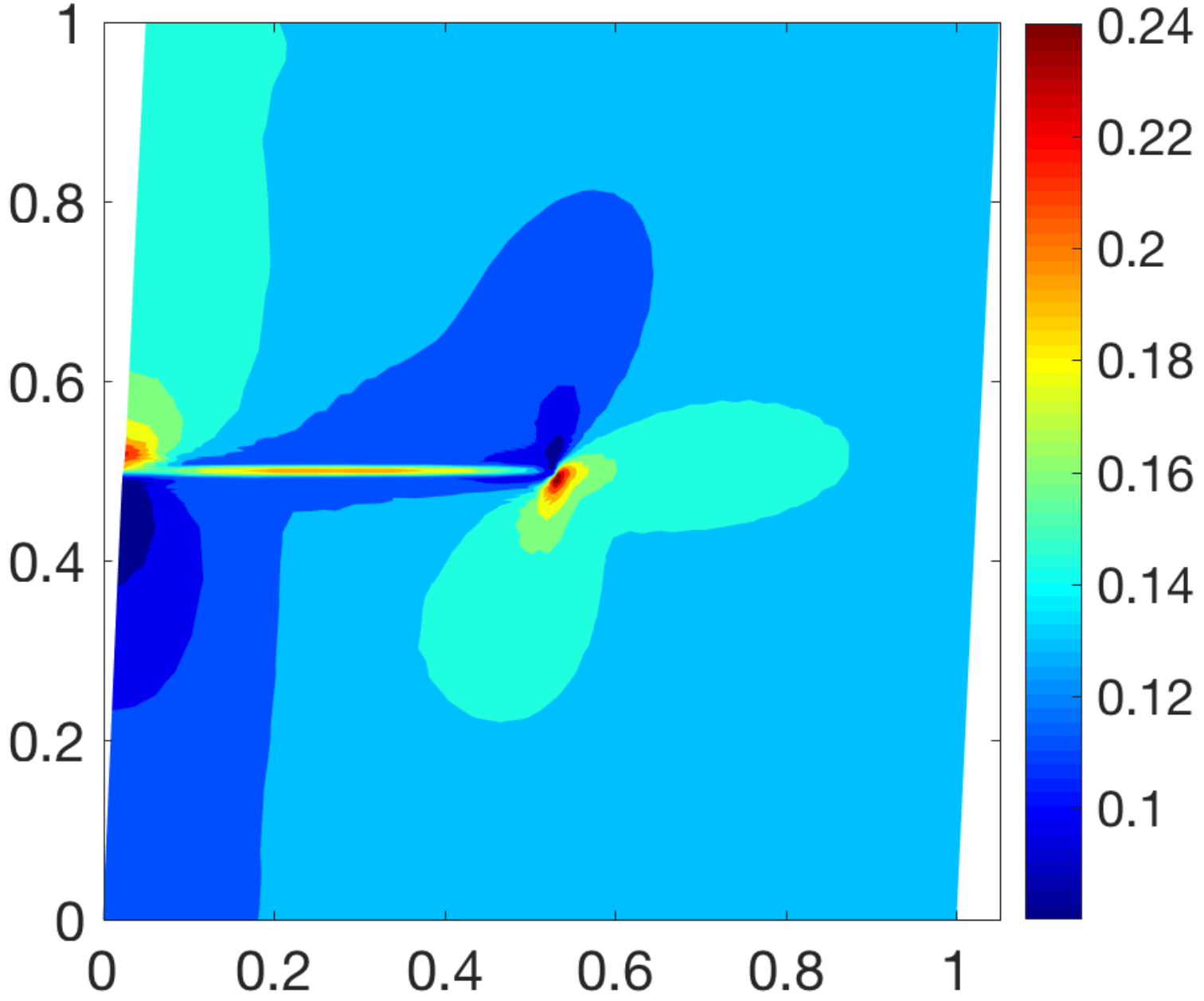}}
\caption{The single edge notched shear test for the energy spectral decomposition model.
The von Mises stress is defined as $(\sigma_x^2 + \sigma_y^2 - \sigma_x \sigma_y + 3 \tau_{xy}^2)^{\frac 1 2}$.}
\label{fig:original_formulation}
\end{figure}

We recall that on any crack $\Gamma$, the material is totally damaged, the phase-field variable is zero
(i.e., $d(x,t) = 0$), and the stress vanishes, viz., 
\begin{equation}
\sigma \cdot n |_{\Gamma} = 0,
\label{CBC-1}
\end{equation}
where $n$ is the unit outward normal to $\Gamma$.
However, there is no guarantee that this is satisfied in the above described three models
of energy decomposition. Indeed, we recall that $\Psi_{e,pas}(\eps)$ is not degraded.
By close examination, we can see that it does not vanish on $\Gamma$ in general for all of the models
described above;
see Fig.~\ref{fig:subfig:Original_Energy} for a sketch for this remaining energy in the totally damaged zone.
This can also be explained from (\ref{sigma}). On $\Gamma$, we have $d = 0$ and 
\[
\sigma = k_l \frac{\partial \Psi_{e,act}}{\partial \eps} + \frac{\partial \Psi_{e,pas}}{\partial \eps}
\approx \frac{\partial \Psi_{e,pas}}{\partial \eps},
\]
where we have used $k_l \ll l \ll 1$.  The fact that $\partial \Psi_{e,pas}/\partial \eps$ does not vanish
on $\Gamma$ in general implies that it is not guaranteed that (\ref{CBC-1}) be satisfied.
The violation of the boundary condition can lead to unphysical crack propagation.
To see the impacts, we consider a single edge notched shear test, with the domain and boundary conditions shown
in Fig. \ref{fig:subfig:Shear}.  In this case, a relative displacement is expected and no stress should
remain on the crack surface. However, the results obtained with the spectral decomposition model
show a stiff response (Fig. \ref{fig:subfig:Original_Sigma}) and unrealistic displacement under loading
(Fig. \ref{fig:subfig:Original_Mesh}). 

To avoid this problem, we propose here a modification which introduces a critically damaged zone
with $d < d_{cr}$ and degrades the total stored energy in the zone; see Fig. \ref{fig:subfig:Modified_Crack}.
Here, $d_{cr} \in (0,1)$ is a threshold. When $d_{cr} = 0$, the modified model reduces back to the original
model. On the other hand, when $d_{cr} = 1$, the total energy is degraded for the entire damaged region.
For any $d_{cr} \in (0,1)$, {\em the total energy is degraded in the critically damaged zone with $0\le d < d_{cr}$
and only the active component of the energy is degraded in the damaged zone with $d_{cr} < d \le 1$.}
More discussion on the choice of $d_{cr}$ and its effects will be given in the numerical result section
\S\ref{SEC:numerics}.  Mathematically, the modified damage model reads as
\begin{equation}
\Psi_{e,m}^{l}(\eps, d) = \left( d^2 + k_l\right) \Psi_{e,act,m}^{l}(\eps) + \Psi_{e,pas,m}^{l}(\eps),
\label{ItCBC-1}
\end{equation}
where 
\begin{align*}
& \Psi_{e,act,m}^{l}(\eps)  = 
\begin{cases}
\Psi_{e,act}^{l}(\eps) + \Psi_{e,pas}^{l}(\eps) = \Psi_{e}^{l}, & \quad \text{ for } \quad 0 \le d \le d_{cr} 
\\
\Psi_{e,act}^{l}(\eps), &\quad \text{ for }  \quad d_{cr} < d \le 1
\end{cases}
\\
& \Psi_{e,pas,m}^{l}(\eps)  = 
\begin{cases}
0, & \quad \text{ for } \quad 0 \le d \le d_{cr} 
\\
\Psi_{e,pas}^{l}(\eps), &\quad \text{ for }  \quad d_{cr} < d \le 1 .
\end{cases}
\end{align*}
By construction, $\Psi_{e,m}^{l}(\eps, d) = 0$ and $\sigma = 0$ on $\Gamma$
(assuming that $k_l = 0$) and thus the crack boundary condition (\ref{CBC-1}) is satisfied.
For simplicity, we will call this modification as ItCBC (Improved treatment of Crack Boundary Conditions).

\begin{figure} 
\centering 
\subfigure[Mesh on deformed domain]{\label{fig:subfig:Modified_Mesh}
\includegraphics[width=0.3\linewidth]{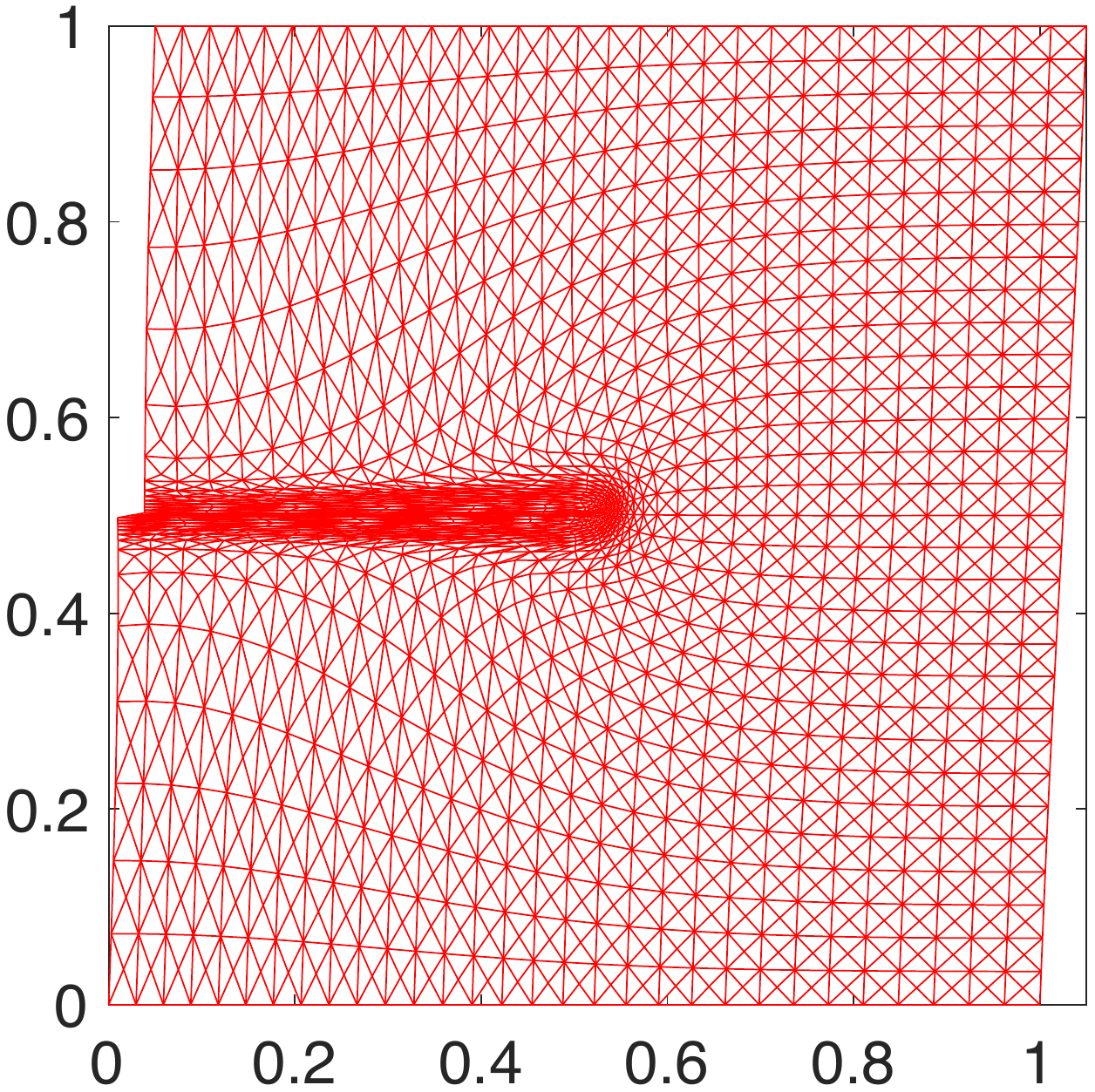}}
\subfigure[Von Mises stress distribution on the deformed domain]{\label{fig:subfig:Modified_Sigma}
\includegraphics[width=0.35\linewidth]{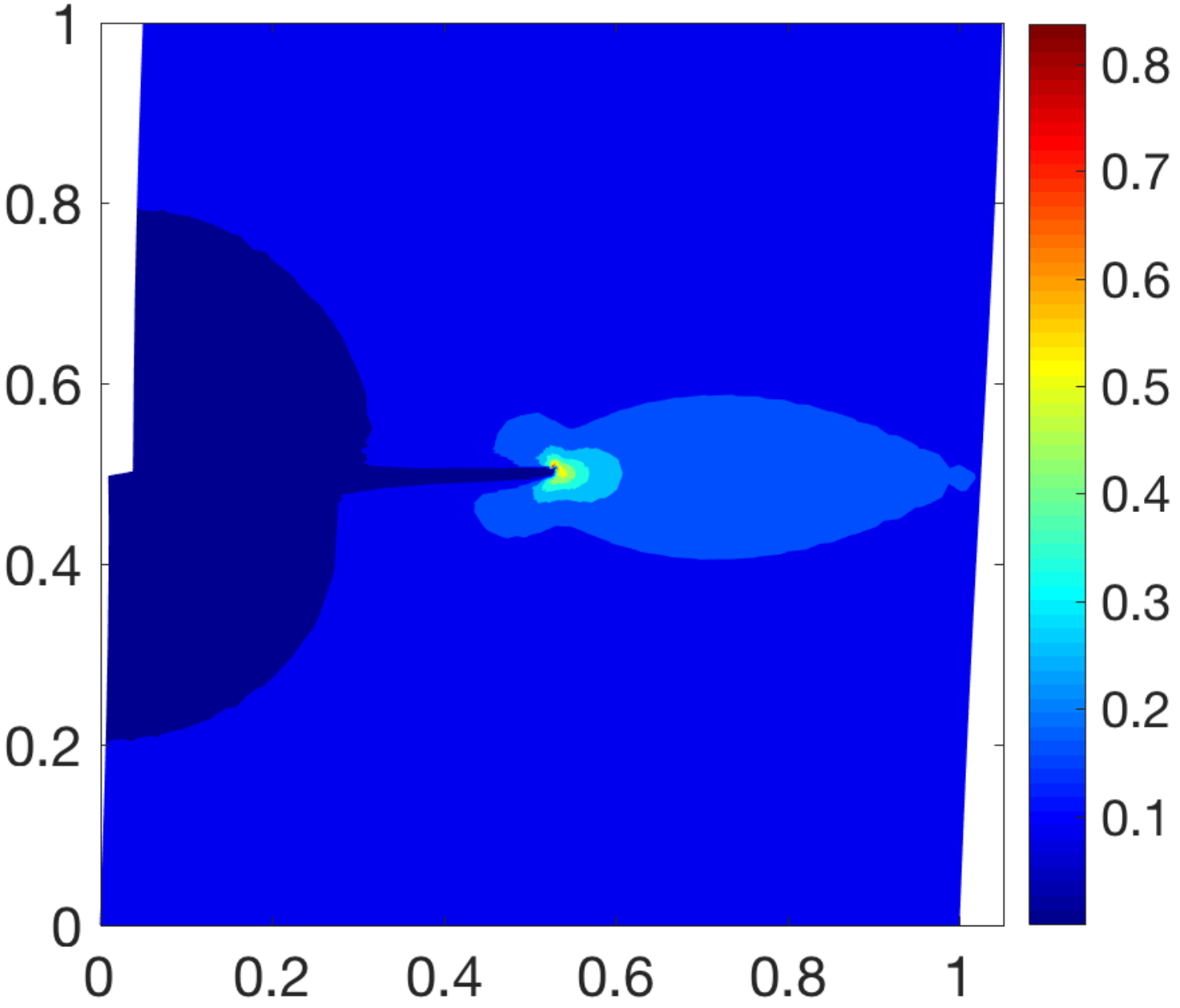}}
\caption{The single edge notched shear test with ItCBC with $d_{cr} = 0.2$ for
the energy spectral decomposition model.}
\label{fig:modified_formulation}
\end{figure}

For the single edge notched shear test considered earlier, the results with ItCBC with $d_{cr} = 0.2$
for the spectral decomposition model are shown in Fig. \ref{fig:modified_formulation}.
It can be seen that ItCBC effectively reduces the stress to zero on the crack and
the results agree well with the expected response as shown in the deformed domain.

\section {A moving mesh finite element method}
\label{SEC:mmfem}

In this section we briefly describe a moving mesh finite element method for solving
the phase-field problem (\ref{weak-d}) and (\ref{weak-u}).
It was first considered for the phase-field modeling of brittle fracture by
Zhang et al. in \cite{ZHLZ17}. The reader is referred to the reference for
more detailed discussion of the method.

\subsection{Finite element discretization and solution procedure}
\label{SEC:procedure}

Let $\mathcal{T}_h$ be a simplicial mesh for the domain $\Omega$ and denote $N$ and $N_v$ as the number of its elements and vertices, respectively. The function spaces $V_d$, $V_u$ and $V_u^0$ are approximated by
\begin{align*}
& V_d^h =  \left \{ 
\varphi_h\; | \; \varphi_h \in C^0(\overline{\Omega}); \;
\varphi_h |_K \in P_1(K),\;  \forall K \in \mathcal{T}_h 
\right \} \subset V_d ,\\
& V_u^h =  \left \{ 
\varphi_h \; | \; \varphi_h \in C^0(\overline{\Omega}) ; \; \varphi_h|_{\partial \Omega_u} = \overline{u}; \; 
\varphi_h |_K \in P_1(K),\;  \forall K \in \mathcal{T}_h 
\right \} \subset V_u , \\
& V_u^{0,h} =  \left \{ 
\varphi_h \; | \; \varphi_h \in C^0(\overline{\Omega}); \; \varphi_h|_{\partial \Omega_u} = 0; \;
\varphi_h |_K \in P_1(K),\;  \forall K \in \mathcal{T}_h
\right \} \subset V_u^0,
\end{align*}
where $P_1(K)$ is the set of polynomials of degree less than or equal to 1 defined on $K$. 
For the phase-field problem (\ref{weak-d}) and (\ref{weak-u}), the linear finite element approximation is to find $d_h \in V_d^h$ and $u_h \in V_u^h$ such that
\begin{align}
& \int_\Omega \left (
 \left (2 d_h \mathcal{H} + \frac{g_c (d_h-1) }{2 l} \right ) \varphi_h
 +  2 g_c l \nabla d_h \cdot \nabla \varphi_h \right) d \Omega = 0 ,
 \quad \forall\; \varphi_h \in V_d^h 
\label{discrete-d}
\\
& \int_\Omega \sigma(u_h) : \epsilon(\varphi_h) d \Omega = 
\int_{\Omega_t} \overline{t} \cdot \varphi_h \; dS +
\int_\Omega f \cdot \varphi_h \; d \Omega ,
\quad \forall\; \varphi_h \in V_u^{0,h}.
\label{discrete-u}
\end{align}

In our computation, we solve (\ref{discrete-d}) for $d_h$ and (\ref{discrete-u}) for $u_h$ alternately.
The procedure leads to smaller and easier systems to solve since $d_h$ and $u_h$ are decoupled
and (\ref{discrete-d}) is linear about $d_h$.
We recall that the energy density is decomposed into active and passive components 
with only the former contributing to the evolution of cracks. This decomposition results in a non-smooth
elastic energy, increases nonlinearity of the displacement system, and makes Newton's iteration
often difficult to converge. Three regularization methods have been proposed by Zhang~et~al.~\cite{ZHLZ17}
to smooth positive and negative eigenvalue functions via a switching technique (sonic-point regularization)
or convolution with a smoothed delta function. In this work, we use the sonic-point regularization method
with which the positive and negative eigenvalue functions are replaced by
\[
\lambda_{\alpha}^+ = \frac{\lambda + \sqrt{\lambda^2 + \alpha^2}}{2} , \quad
\lambda_{\alpha}^- = \frac{\lambda - \sqrt{\lambda^2 + \alpha^2}}{2},
\]
where $\alpha > 0$ is the regularization parameter. It is shown in \cite{ZHLZ17} that this
regularization can effectively make Newton's iteration convergent.

We consider the problem in a quasi-static condition, with the quasi-time $t$ being introduced to 
represent the load increments. The solution procedure from $t^n$ to $t^{n+1}$ is described
as follows.

\begin{enumerate}
\item [(i)] Suppose the mesh $\mathcal{T}_h^n$ at time $t^n$ and the history field $\mathcal{H}_h^n$ in \eqref{eqn-H} (defined on $\mathcal{T}_h^n$) are given.
\item [(ii)] Compute the phase-field variable $d_h^{n+1}$ and new mesh $\mathcal{T}_h^{n+1}$ as follows.
	\begin{itemize}
	\item Let $\mathcal{T}_h^{n+1,1} = \mathcal{T}_h^n$;
	\item For $k = 1:kk$
		\begin{itemize}
		\item[-] Compute $\mathcal{H}$ on $\mathcal{T}_h^{n+1,k}$ using linear interpolation of $\mathcal{H}_h^n$;
		\item[-] Compute $d_h^{n+1,k}$ using (\ref{discrete-d}) and $\mathcal{H}$ on $\mathcal{T}_h^{n+1,k}$;
		\item[-] If $k < kk$, compute the new mesh $\mathcal{T}_h^{n+1,k+1}$ by the MMPDE moving mesh method
		based on $\mathcal{T}_h^{n+1,k}$ and $d_h^{n+1,k}$; see Section~\ref{SEC:MMPDE}.
		\end{itemize}
	\item Let $\mathcal{T}_h^{n+1} = \mathcal{T}_h^{n+1,kk}$ and $d_h^{n+1} = d_h^{n+1,kk}$.
	\end{itemize}
\item [(iii)] Compute the displacement field $u_h^{n+1}$ by solving the nonlinear system (\ref{discrete-u})
based on $d_h^{n+1}$ and $\mathcal{T}_h^{n+1}$. Newton's iteration is used.

\item [(iv)] Compute $\Psi_{e,act,h}^{l, n+1} (\epsilon(u_h^{n+1}))$ and set
$\mathcal{H}_h^{n+1} = \max\{\Psi_{e,act,h}^{l, n+1} (\epsilon(u_h^{n+1})), \mathcal{\tilde{H}}_h^n\}$,
where $\mathcal{\tilde{H}}_h^n$ is the linear interpolation of $\mathcal{H}_h^n$ from
the old mesh $\mathcal{T}_h^n$ to the new mesh $\mathcal{T}_h^{n+1}$.
\end{enumerate}

The parameter $kk$ in the above procedure determines the adaptivity of the mesh $\mathcal{T}_h^{n+1}$ to the phase-field variable $d_h^{n+1}$. Our experience shows that $kk=5$ is sufficient for the mesh to be well adaptive
to $d_h^{n+1}$.

\subsection{The MMPDE moving mesh method}
\label{SEC:MMPDE}

We use the MMPDE moving mesh method \cite{HRR94a,HR11} for generating the new mesh $\mathcal{T}_h^{n+1}$. The method is based on the $\mathbb{M}$-uniform mesh approach with which a nonuniform adaptive mesh
is viewed as a uniform one in the metric specified by a tensor $\mathbb{M}$. The metric tensor $\mathbb{M}$
is assumed to be symmetric and uniformly positive definite on $\Omega$ and determines the shape, size, and orientation of the mesh elements through the so-called equidistribution and alignment conditions
(see (\ref{eq-1}) and (\ref{ali-1}) below). 

Denote $H(d_h)$ as a recovered Hessian of $d_h$ and assume its eigen-decomposition is given by $H(d_h) = Q \text{diag}(\lambda_1,\lambda_2) Q^T$. Then we choose the metric tensor in our computation as
\begin{equation}
\mathbb{M} = \det(I+|H(d_h)|)^{-\frac{1}{6}} (I+|H(d_h)|),
\label{M-1}
\end{equation}
where $|H(d_h)| = Q \text{diag}(|\lambda_1|,|\lambda_2|) Q^T$.
Since (\ref{M-1}) is based on the Hessian of the phase-field variable $d$, the mesh elements are concentrated around the crack where the curvature of $d$ is large.
The form (\ref{M-1}) is known to be optimal for the $L^2$ norm of linear interpolation error (e.g., see \cite{HR11}).

Denote $x_1, ...., x_{N_v}$ as the coordinates of the vertices of $\mathcal{T}_h$.
Let $\hat{\mathcal{T}}_{c,h} =  \{\hat{\xi}_1,...,\hat{\xi}_{N_v} \}$ be the reference computational mesh
which is taken as the initial physical mesh.
We also denote $\mathcal{T}_{c,h} = \{\xi_1,...,\xi_{N_v} \}$ as an intermediate computational mesh. We assume that $\mathcal{T}_h$, $\hat{\mathcal{T}}_{c,h}$, and $\mathcal{T}_{c,h}$ have the same number of elements and vertices and the same connectivity. Then for any element $K \in \mathcal{T}_h$, there exists a corresponding element $K_c \in \mathcal{T}_{c,h}$.
Let $F_K$ be the affine mapping from $K_c$ to $K$ and $F'_K$ be its Jacobian matrix. Denote the vertices of $K$ and $K_c$ as $x_0^K, x_1^K, x_2^K$ and $\xi_0^K, \xi_1^K, \xi_2^K$, respectively. 
We define the edge matrices of $K$ and $K_c$ as
\[
E_K = [x_1^K-x_0^K,x_2^K-x_0^K], \quad \hat{E}_K = [\xi_1^K-\xi_0^K, \xi_2^K-\xi_0^K].
\]
It is easy to show that
\[
F'_K = E_K \hat{E_K}^{-1}, \quad (F'_K)^{-1} = \hat{E_K} E_K^{-1}.
\]

It is known (e.g., see \cite{Hua06,HR11}) that an $\mathbb{M}$-uniform mesh $\mathcal{T}_h$
approximately satisfies the equidistribution and alignment conditions 
\begin{align}
& |K| \sqrt{\det(\mathbb{M}_K)} = \frac{|\Omega_h|\; |K_c|}{|\Omega_c|} , \quad \forall{K} \in \mathcal{T}_h
\label{eq-1}
\\
& \frac{1}{2} \tr\left( (F'_K)^T \mathbb{M}_K F'_K \right) = \det\left( (F'_K)^T \mathbb{M}_K F'_K \right)^{\frac{1}{2}}, \quad K \in \mathcal{T}_h
\label{ali-1}
\end{align}
where $|K|$ and $|K_c|$ denote the area of $K$ and $K_c$, respectively, $\mathbb{M}_K$ is the average of $\mathbb{M}$ over $K$, $\det(\cdot)$ denotes the determinant of a matrix, and
\[
|\Omega_h| = \sum_{K \in \mathcal{T}_h} |K| \sqrt{\det(\mathbb{M}_K)},
\quad |\Omega_c| = \sum_{K_c \in \mathcal{T}_{c,h}} |K_c|. 
\]
An energy function that combines the above two conditions has been proposed by Huang \cite{Hua01} as
\begin{equation}
I_h(\mathcal{T}_h; \mathcal{T}_{c,h}) = \sum_{K \in \mathcal{T}_h} |K|
G(\mathbb{J}_K, \det(\mathbb{J}_K), \mathbb{M}_K),
\end{equation}
where $\mathbb{J}_K = (F'_K)^{-1}$ and
\begin{align*}
G(\mathbb{J}_K, \det(\mathbb{J}_K), \mathbb{M}_K) & = 
\frac{1}{3} \sqrt{\det(\mathbb{M}_K)} \left( \tr(\mathbb{J}_K \mathbb{M}_K \mathbb{J}^T_K) \right)^{\frac 3 2}
+  \frac{2^{\frac 3 2}}{3} \sqrt{\det(\mathbb{M}_K)}
\left( \frac{\det(\mathbb{J}_K)}{\sqrt{\det(\mathbb{M}_K)}} \right)^{\frac 3 2} .
\end{align*}
In principle, a new physical mesh $\mathcal{T}_h^{n+1}$ can be found by minimizing $I_h$
with respect to $\mathcal{T}_{h}$ for a given $\mathcal{T}_{c,h}$ (for example, taken as $\hat{\mathcal{T}}_{c,h}$).
However, this minimization can be difficult and costly since $I_h$ is generally not convex.
Here, we use the $\xi$-formulation of the MMPDE moving mesh method. That is, we first
take $\mathcal{T}_h = \mathcal{T}_h^n$ and define the moving mesh equation
as a gradient system of $I_h$ with the coordinates of the computational vertices, i.e., 
\begin{equation}
\frac{d\xi_j}{dt} = -\frac{P_j}{\tau} \left( \frac{\partial I_h}{\partial \xi_j} \right)^T, \quad j = 1,...,N_v
\label{MMPDE1}
\end{equation}
where $\tau > 0$ is a parameter used to adjust the time of the mesh movement
and $P_j = \det(\mathbb{M}(x_j))^{\frac{1}{4}}$ is chosen to make (\ref{MMPDE1}) invariant under the scaling transformation of $\mathbb{M}$. Using the analytical formulas for the derivatives \cite{HK15a}, we can rewrite
\eqref{MMPDE1} into
\begin{equation}
\frac{d\xi_j}{dt} = \frac{P_j}{\tau}\sum_{K\in\omega_j} |K| {v}_{j_K}^K, 
\label{MMPDE2}
\end{equation}
where $\omega_j$ is the element patch associated with the $j$-th vertex, $j_K$ is the corresponding local index of the vertex on $K$, and ${v}_{j_K}^K$ is the local velocity of the vertex that is given by
\begin{displaymath}
\left[ \begin{array}{ccc}
({v}_1^K)^T \\
({v}_2^K)^T
\end{array} \right] = 
-E_K^{-1} \frac{\partial G}{\partial \mathbb{J}} - \frac{\partial G}{\partial \det(\mathbb{J})} \frac{\det(\hat{E}_K)}{\det(E_K)} \hat{E}_K^{-1},
\quad {v}_0^K = - \sum_{i=1}^2 {v}_j^K.
\end{displaymath}
The derivatives of the function $G$ in the above equation are given by
\begin{align*}
& \frac{\partial G}{\partial \mathbb{J}} = \sqrt{\det(\mathbb{M})}\left( \tr(\mathbb{J}\mathbb{M}^{-1}\mathbb{J}^T) \right)^{\frac 1 2}\mathbb{M}^{-1}\mathbb{J}^T, \\
& \frac{\partial G}{\partial \det(\mathbb{J})} = \sqrt{2} \det (\mathbb{M})^{-\frac{1}{4}} \det(\mathbb{J})^{\frac 1 2}.
\end{align*}

Note that the velocities for the boundary vertices need to be modified appropriately so that they stay
on the boundary. After that, \eqref{MMPDE2}  is integrated using
the Matlab\textsuperscript \textregistered\, ODE solver {\em ode15s}
from $t^n$ to $t^{n+1}$ with $\hat{\mathcal{T}}_{c,h}$ as the initial mesh.
The new computational mesh is denoted as $\mathcal{T}_{c,h}^{n+1}$. This mesh and the physical mesh
$\mathcal{T}_h^{n}$ form a piecewise linear correspondence, i.e.,
$\mathcal{T}_h^{n} = \Phi_h (\mathcal{T}_{c,h}^{n+1})$.
The new physical mesh is then defined as
$\mathcal{T}_h^{n+1} = \Phi_h (\hat{\mathcal{T}}_{c,h})$, which can be computed
using linear interpolation.

It is worth pointing out that an $x$-formulation of the MMPDE moving mesh method
can also be obtained by taking $\mathcal{T}_{c,h} = \hat{\mathcal{T}}_{c,h}$ and
using the gradient system of $I_h$ with respect to the coordinates of the physical vertices.
Although more complicated to implement (e.g., see \cite{HK15a}) than the $\xi$-formulation,
the $x$-formulation has the advantage that its generated mesh is theoretically and numerically
guaranteed to be nonsingular if it is so initially (cf. \cite{HK15b}). On the other hand,
a similar theoretical result has not yet been proven for the $\xi$-formulation although
numerical experiment has shown that it also produces no mesh tangling or crossing.

\section{Numerical results}
\label{SEC:numerics}

In this section we present numerical results for four examples. The first two are classical benchmark problems,
a single edge notched tension test and a shear test. The third example is designed to demonstrate the
ability of our method to handle complex cracks. The last example is also a single edge notched shear test but with
the physical parameters and domain geometry chosen based on an experiment setting.
The numerical results obtained with the three energy
decomposition models and their modified ones are presented and compared.
The effects of the critical damage threshold $d_{cr}$ and the improved treatment of crack boundary conditions
on the numerical solution are discussed. Unless stated otherwise, the following
choices of the parameters are used: $\alpha = 10^{-3}$ in the sonic-point regularization method,
$N = 6,400$ for the size of an adaptive mesh, and $l = 0.0075$~mm for the width of smeared cracks. 

\subsection{Example 1. A single edge notched tension test}

We first consider a single edge notched tension test from Miehe et al. \cite{MHW10}. The geometry and boundary conditions are show in Fig. \ref{fig:subfig:Tension}. The bottom edge of the domain is fixed. The top edge is fixed along the $x$-direction while along the $y$-direction a uniform displacement $U$ is increased with time to drive the crack propagation. The following material properties are used in our computation: $\lambda = 121.15$~kN/mm$^{2}$,
$\mu = 80.77$~kN/mm$^{2}$, and $g_c = 2.7 \times 10^{-3}$~kN/mm. 
Due to the brutal character of the crack propagation, we choose two displacement increments for the computation,
$\Delta U = 10^{-5}$~mm for the first 500 time steps and $\Delta U = 10^{-6}$~mm afterwards. For comparison purpose, we compute the surface load vector on the top edge as 
\[
F = (F_x,F_y)\equiv \int_{\text{top edge}} \sigma(\eps) \cdot n \, d l ,
\]
where $n$ is the unit outward normal to the top edge. We are interested in $F_y$ for the tension test
and $F_x$ for the shear test.

\begin{figure} 
\centering 
\subfigure[spectral decomposition]{\label{fig:subfig:TMD_U800}
\includegraphics[width=0.25\linewidth]{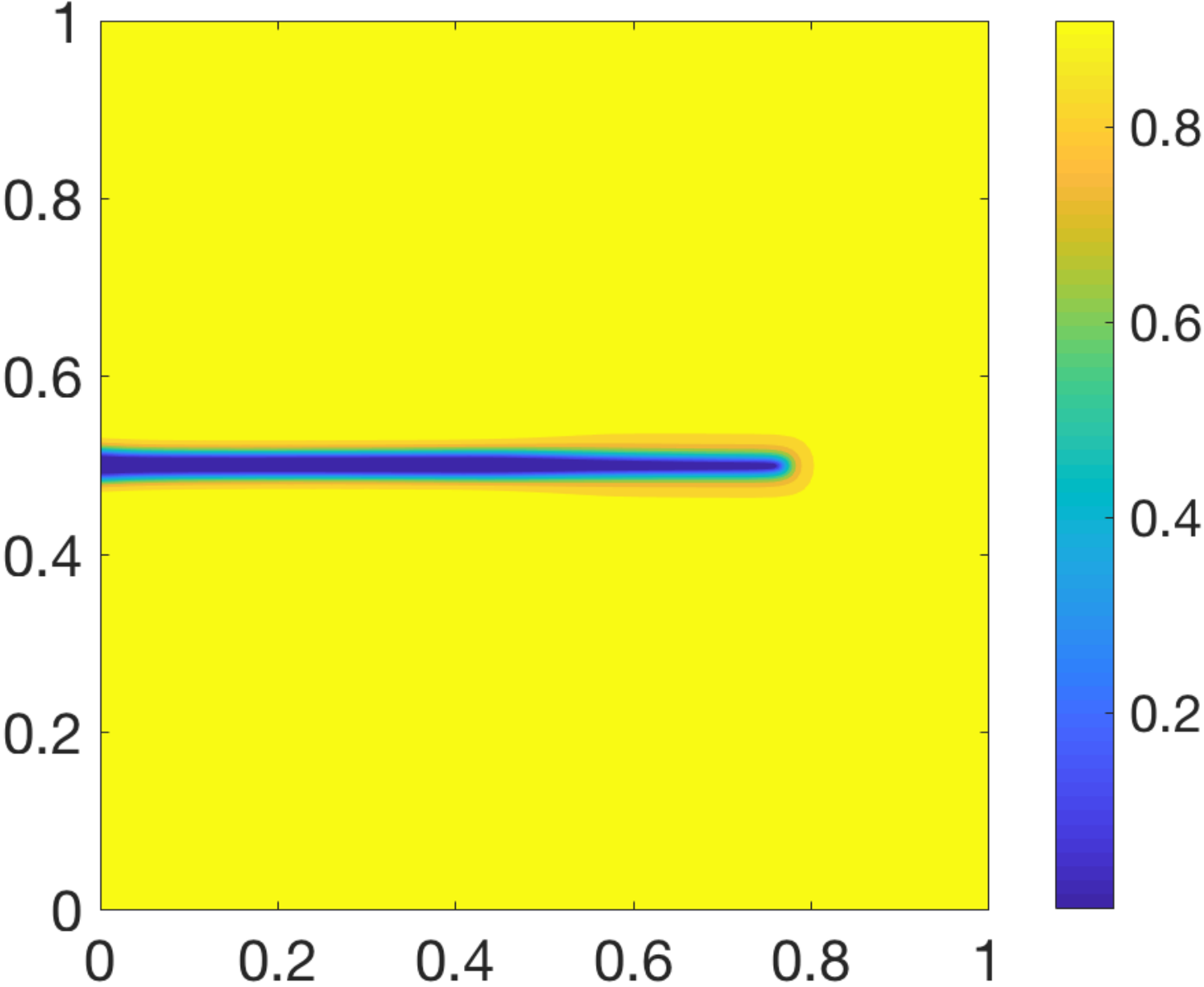}}
\subfigure[v-d split]{\label{fig:subfig:TAD_U800}
\includegraphics[width=0.25\linewidth]{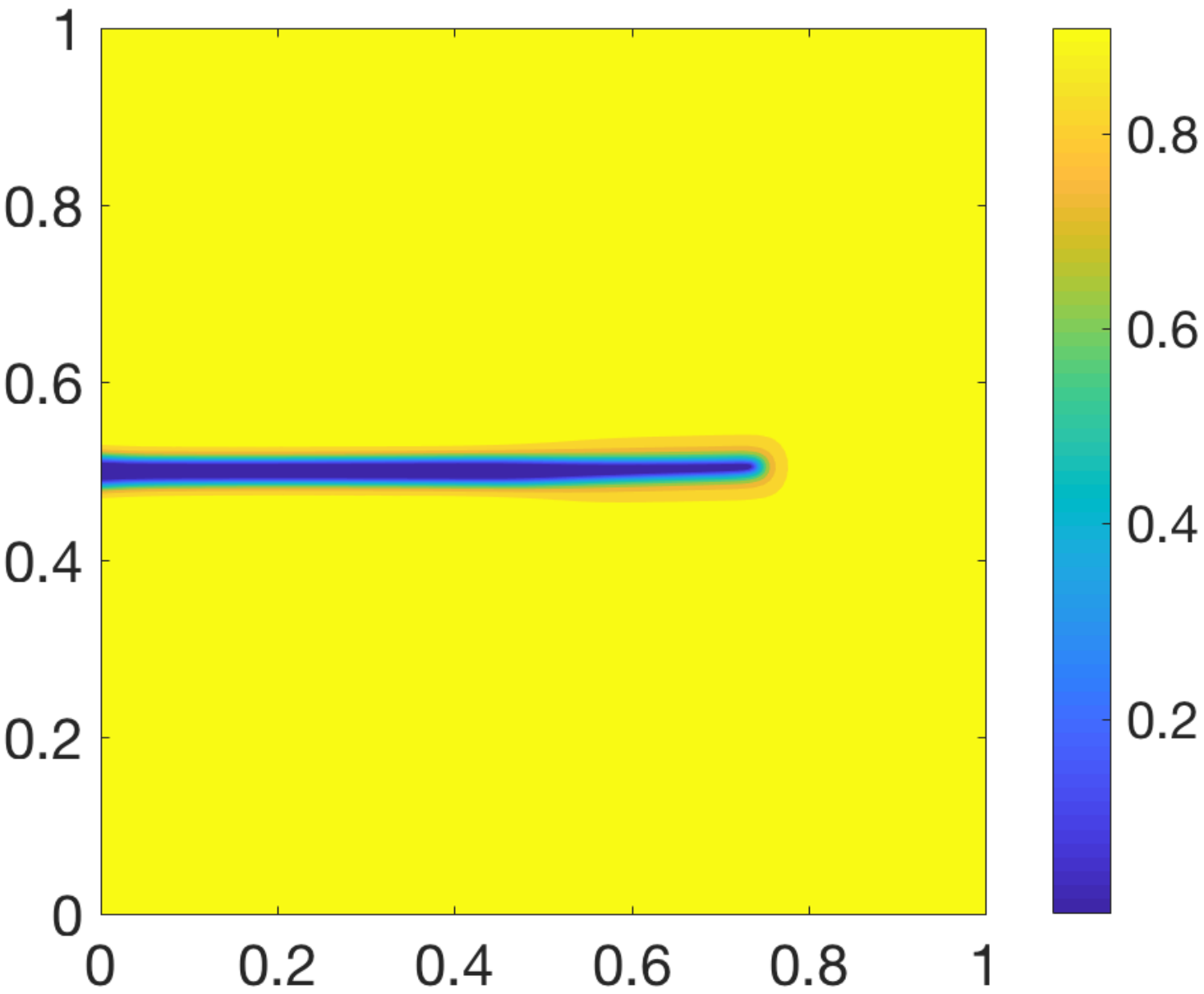}}
\subfigure[improved v-d split]{\label{fig:subfig:TND_U800}
\includegraphics[width=0.25\linewidth]{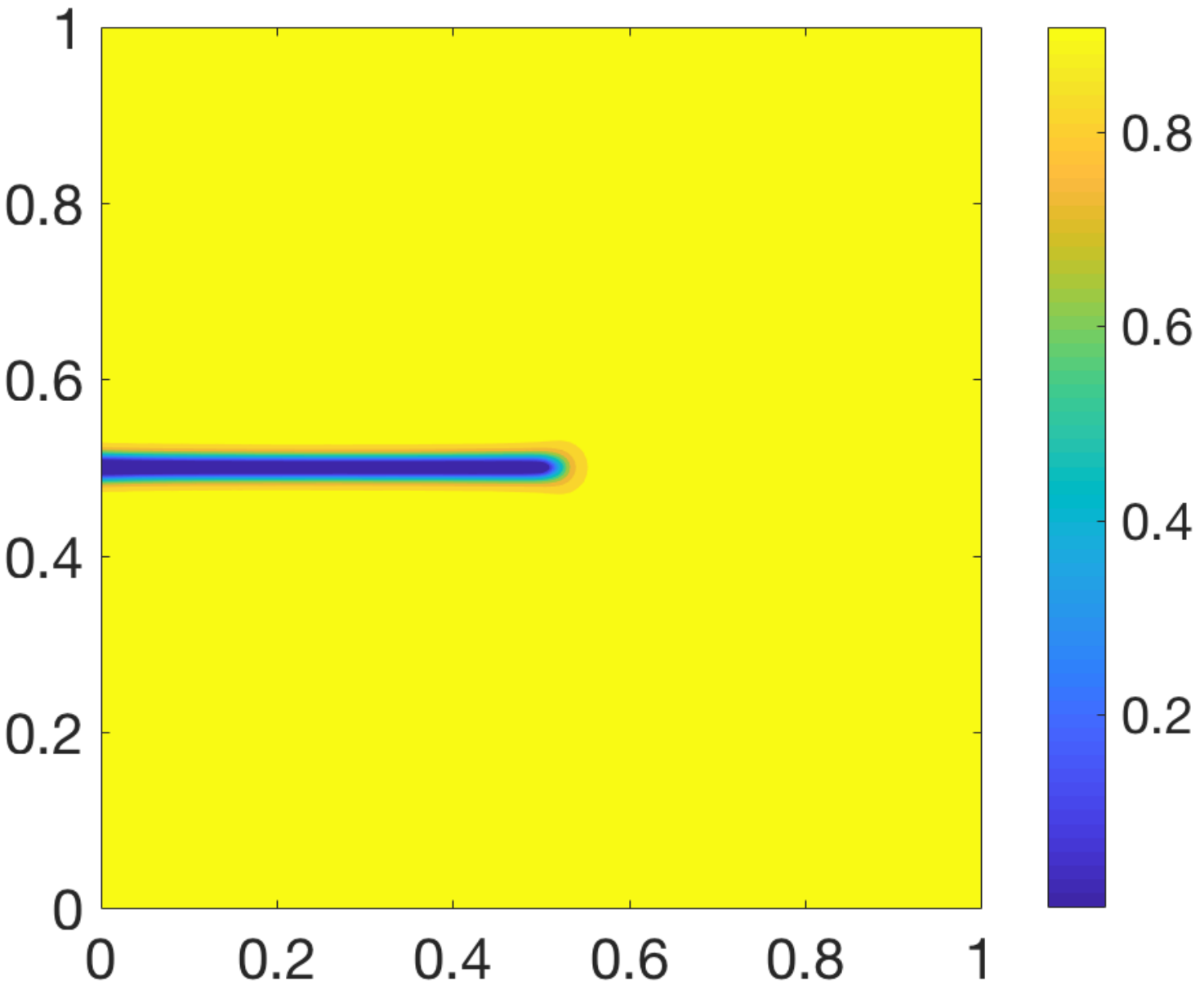}}
\vfill
\subfigure[spectral decomposition]{\label{fig:subfig:TMM_U800}
\includegraphics[width=0.25\linewidth]{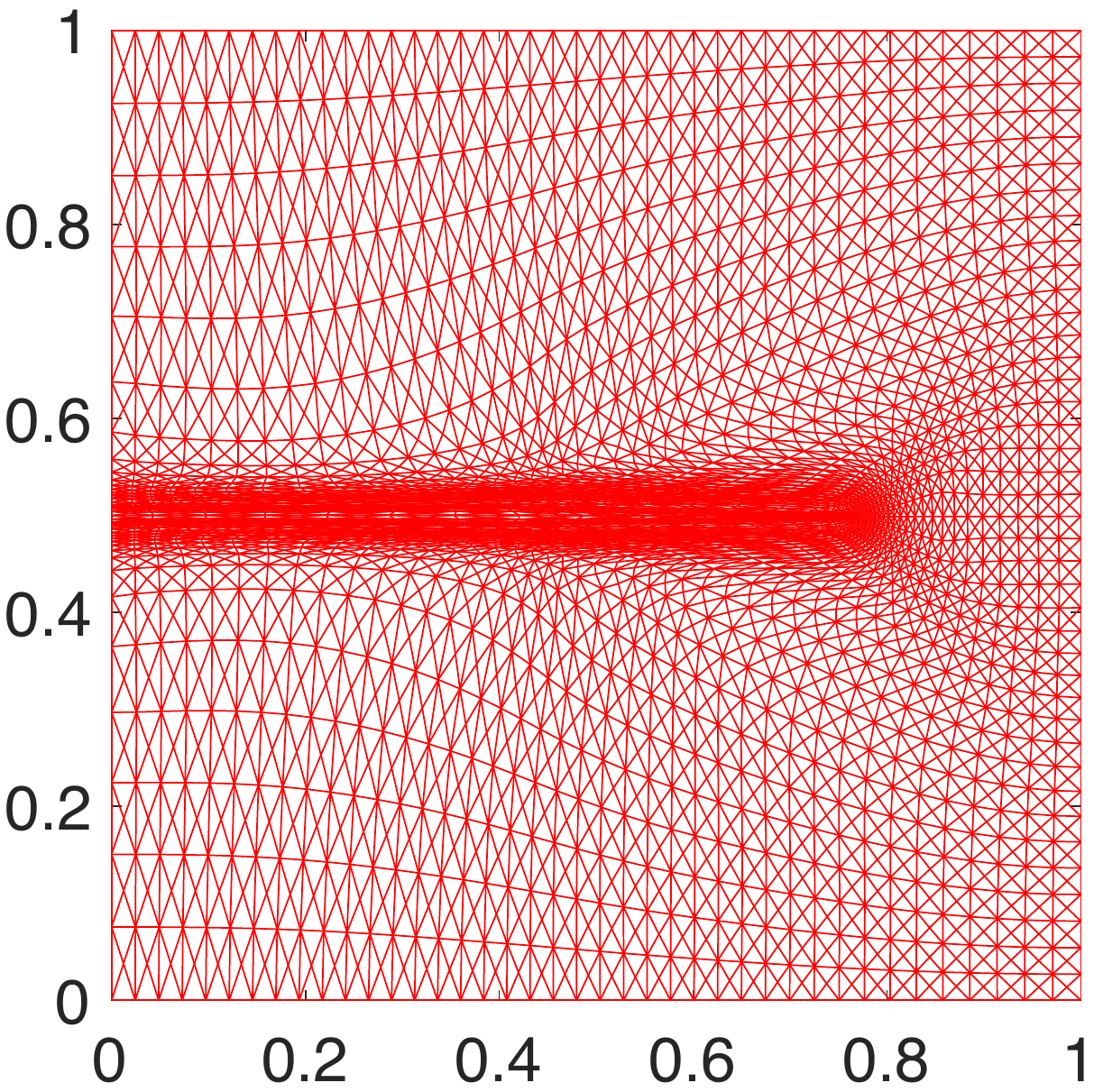}}
\subfigure[v-d split]{\label{fig:subfig:TAM_U800}
\includegraphics[width=0.25\linewidth]{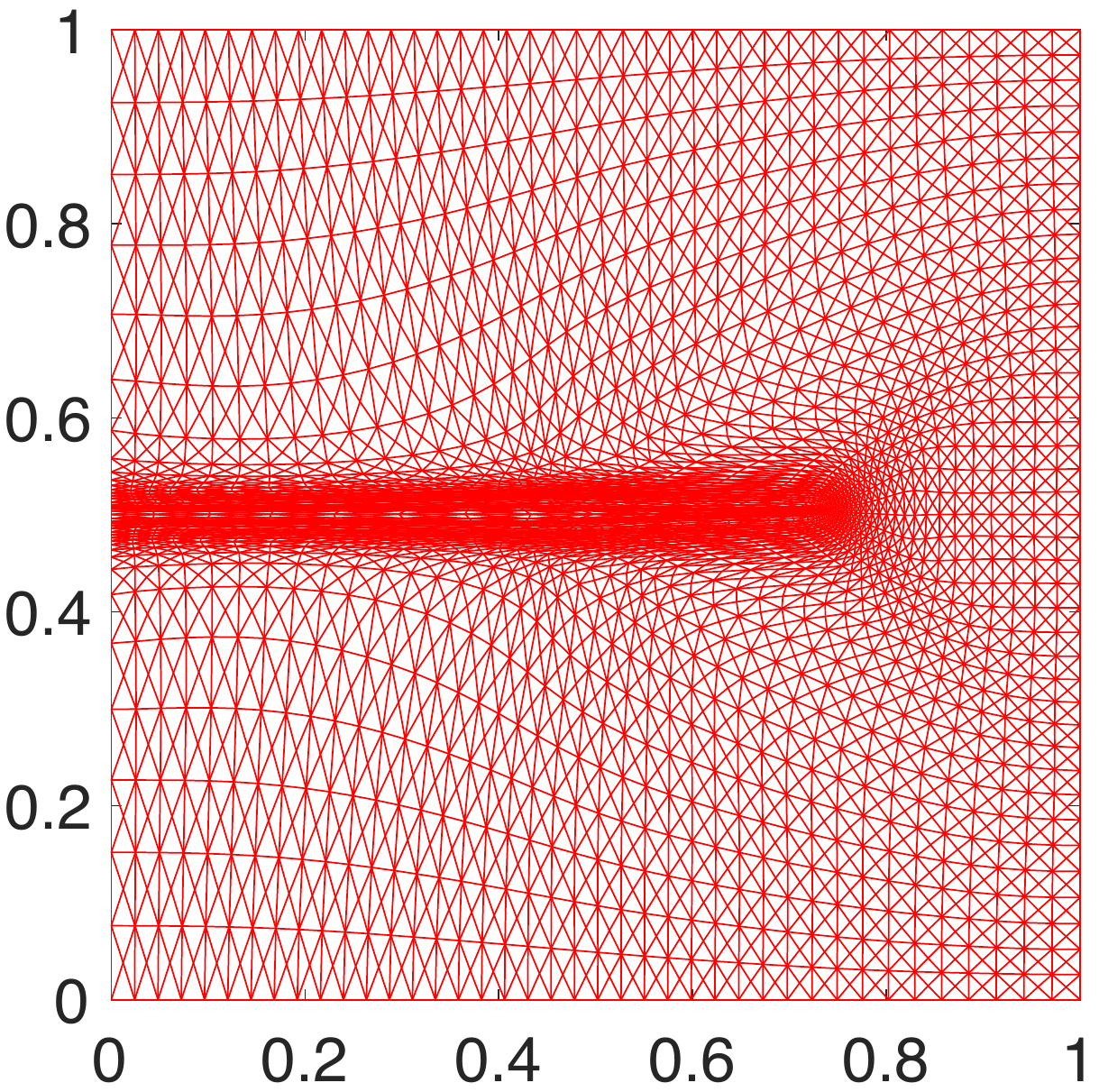}}
\subfigure[improved v-d split]{\label{fig:subfig:TNM_U800}
\includegraphics[width=0.25\linewidth]{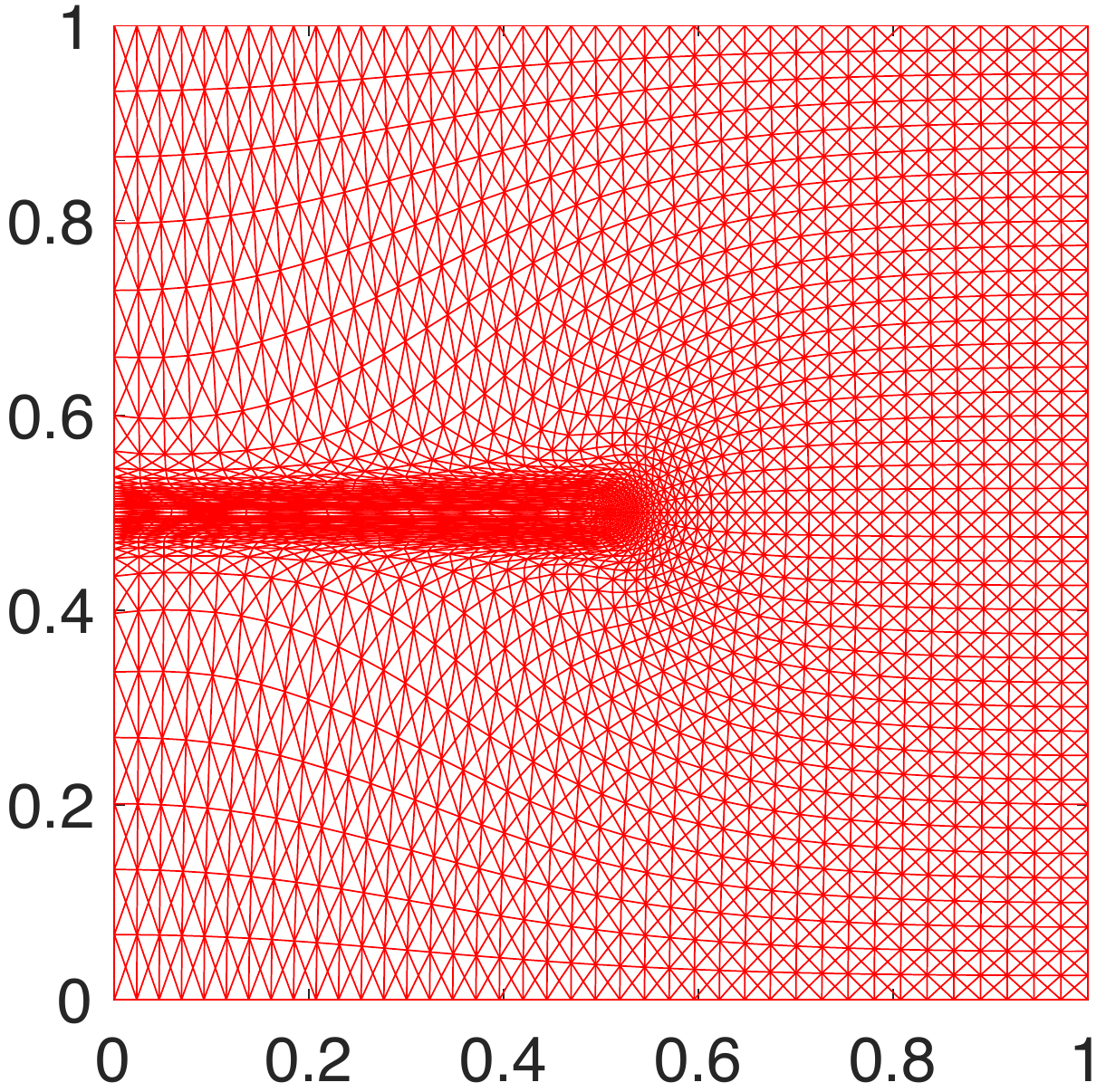}}
\vfill
\subfigure[spectral decomposition]{\label{fig:subfig:TMS_U800}
\includegraphics[width=0.25\linewidth]{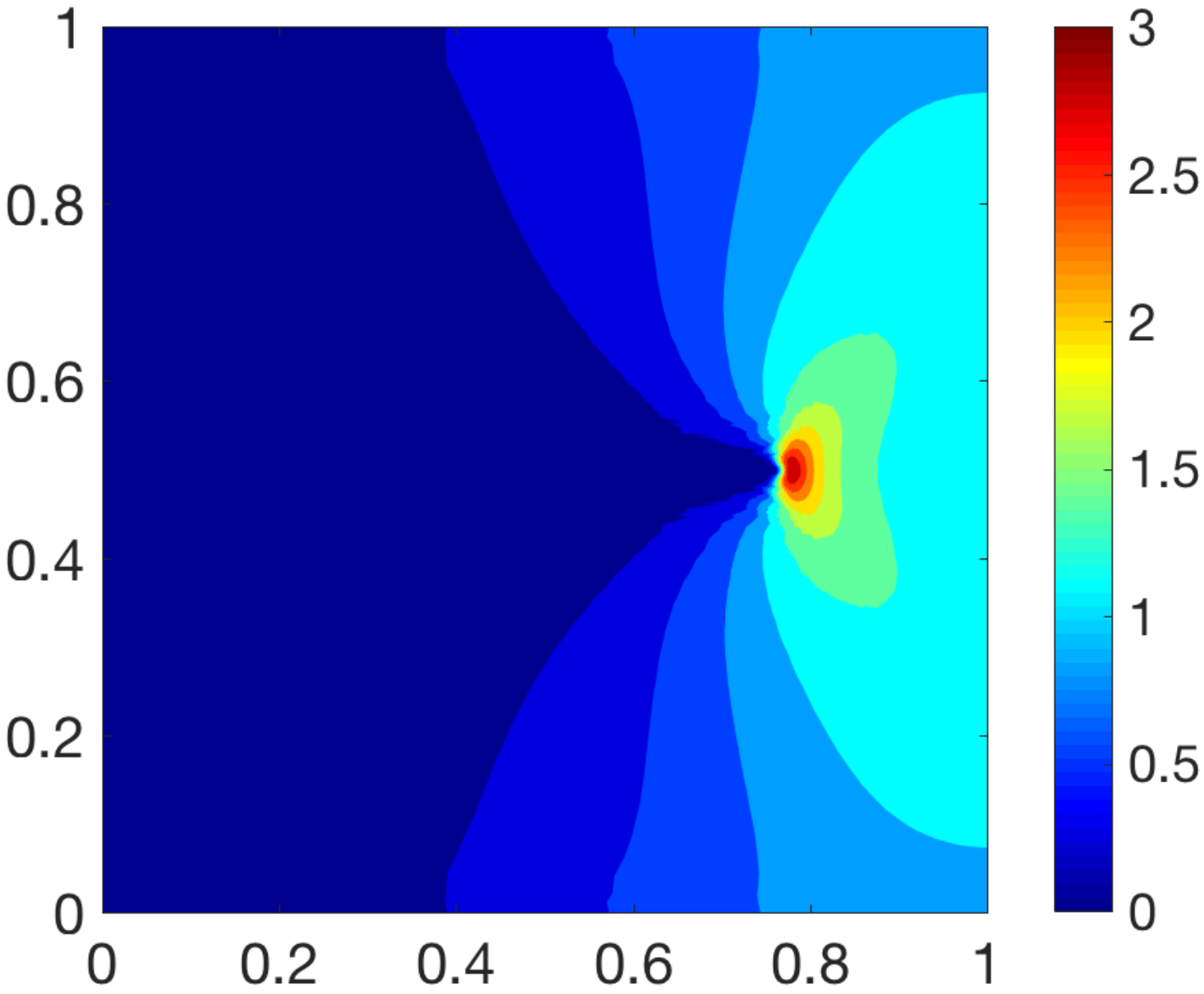}}
\subfigure[v-d split]{\label{fig:subfig:TAS_U800}
\includegraphics[width=0.25\linewidth]{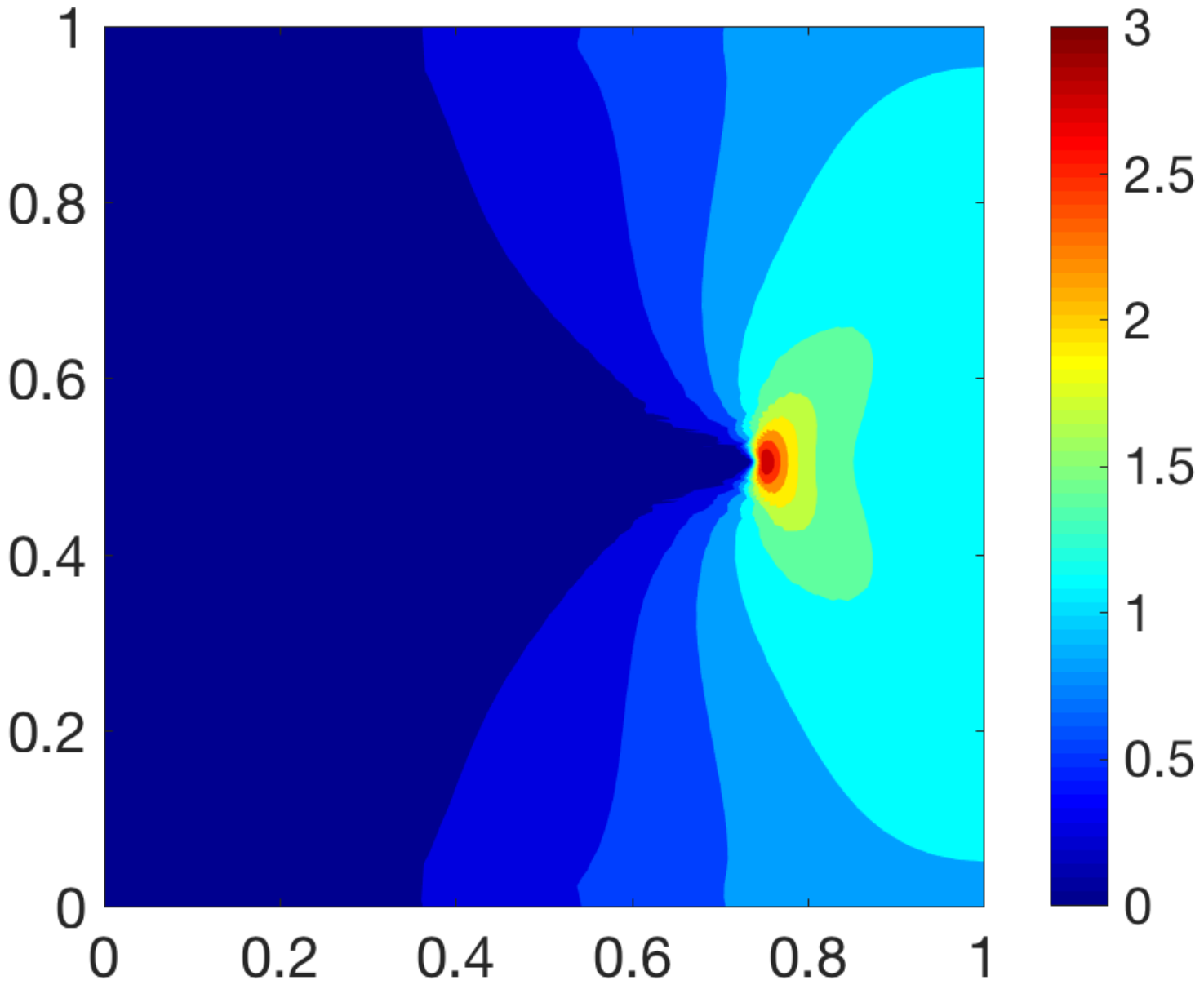}}
\subfigure[improved v-d split]{\label{fig:subfig:TNS_U800}
\includegraphics[width=0.25\linewidth]{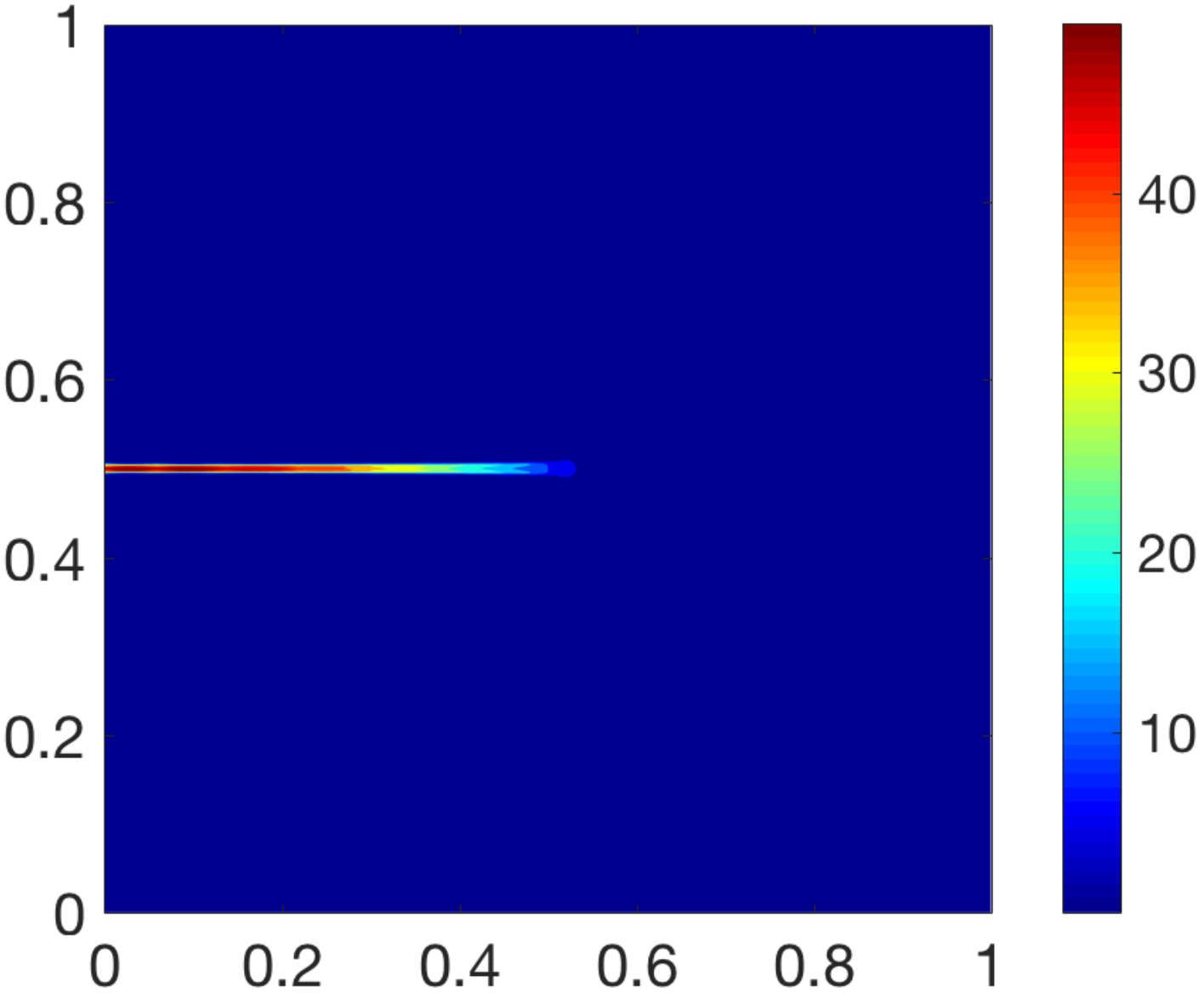}}
\caption{Example 1. Meshes and contours of the phase-field and von Mises stress distribution
are plotted at $U = 5.3 \times 10^{-3}$~mm.  Three energy decomposition models are used.}
\label{fig:tension contours}
\end{figure}

\begin{figure} 
\centering 
\includegraphics[width=0.4\linewidth]{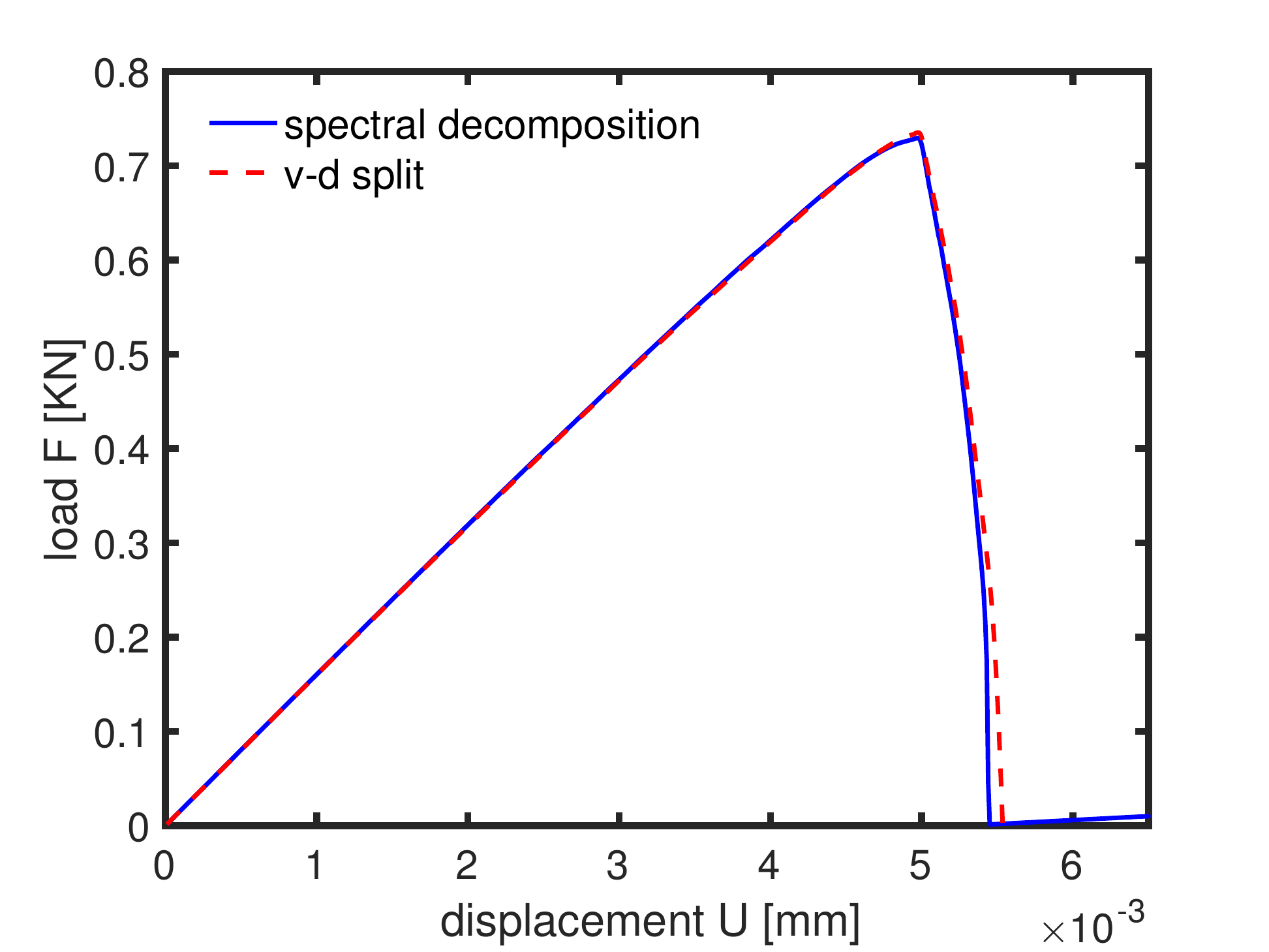}
\caption{Example 1. The load-deflection curves are obtained for spectral and v-d split decomposition models.}
\label{fig:tension curves}
\end{figure}

We now investigate the effects of different energy decomposition models and the crack boundary conditions
on the numerical results. An initial triangular mesh is constructed from a rectangular mesh by subdividing each rectangle into four triangles along diagonal directions. Typical adaptive meshes and contours of the phase field
and von Mises stress distribution using three decomposition models are shown in Fig.~\ref{fig:tension contours}.
The load-deflection curves are shown in Fig.~\ref{fig:tension curves}. 
For the spectral decomposition and v-d split method, the resulting crack path and load-deflection curves
agree well with the results obtained by Miehe at al. \cite{MWH10} using pre-refined mesh around the regions of the crack and its expected propagation path. Moreover, the stress is always concentrated in the region around the crack tip during crack evolution as shown in Figs.~\ref{fig:subfig:TMS_U800} and \ref{fig:subfig:TAS_U800}. 
Inconsistent results are obtained with for the improved v-d split model. The stress growth and concentration occur
in the totally damaged area as shown in Fig. \ref{fig:subfig:TNS_U800}. The reason of this unphysical phenomenon is that for the tension test, this split model leads to a small remaining energy in the damaged zone. 

\begin{figure} 
\centering 
\subfigure[spectral decomposition]{\label{fig:subfig:TMD_U800D04}
\includegraphics[width=0.25\linewidth]{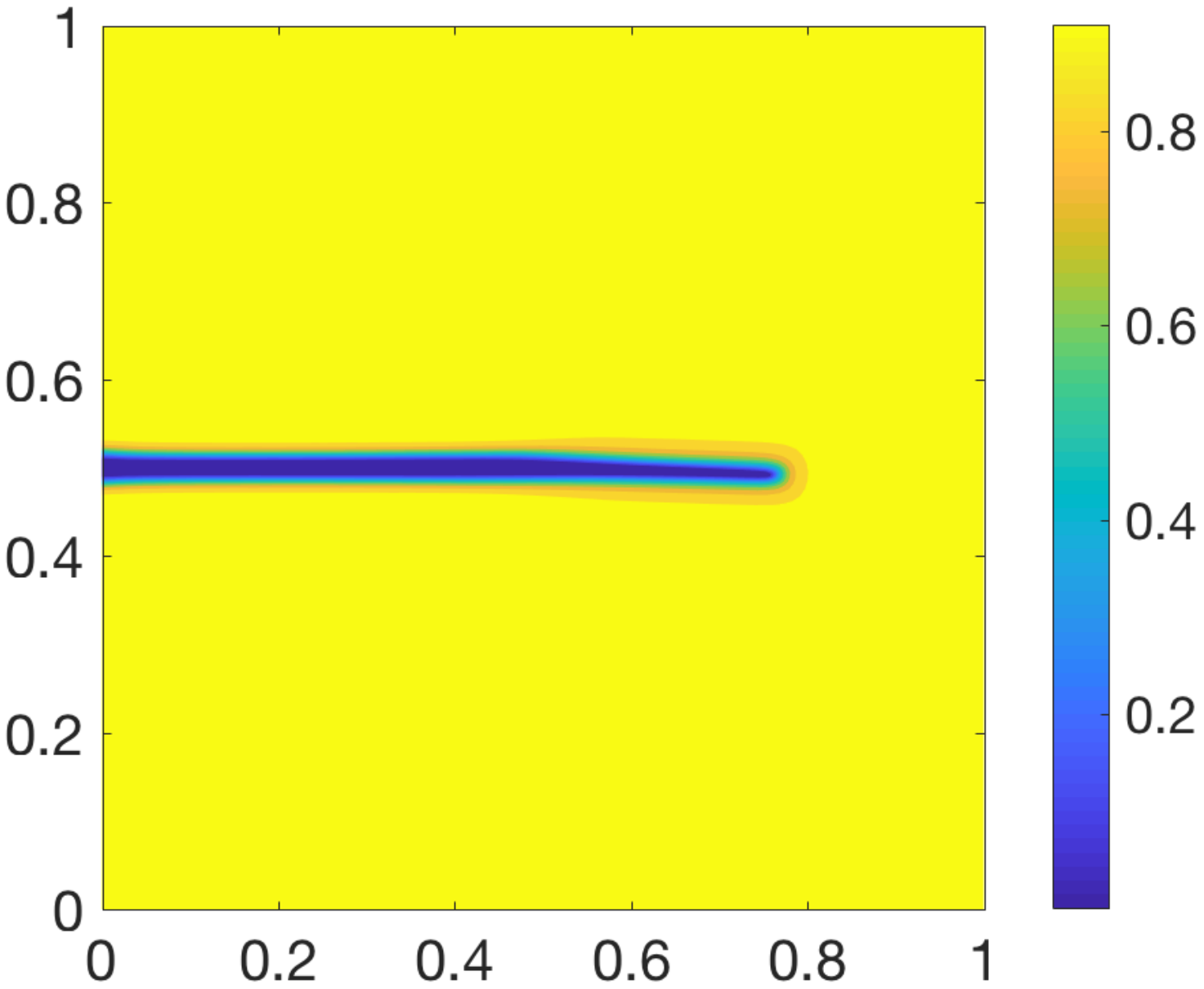}}
\subfigure[v-d split]{\label{fig:subfig:TAD_U800D04}
\includegraphics[width=0.25\linewidth]{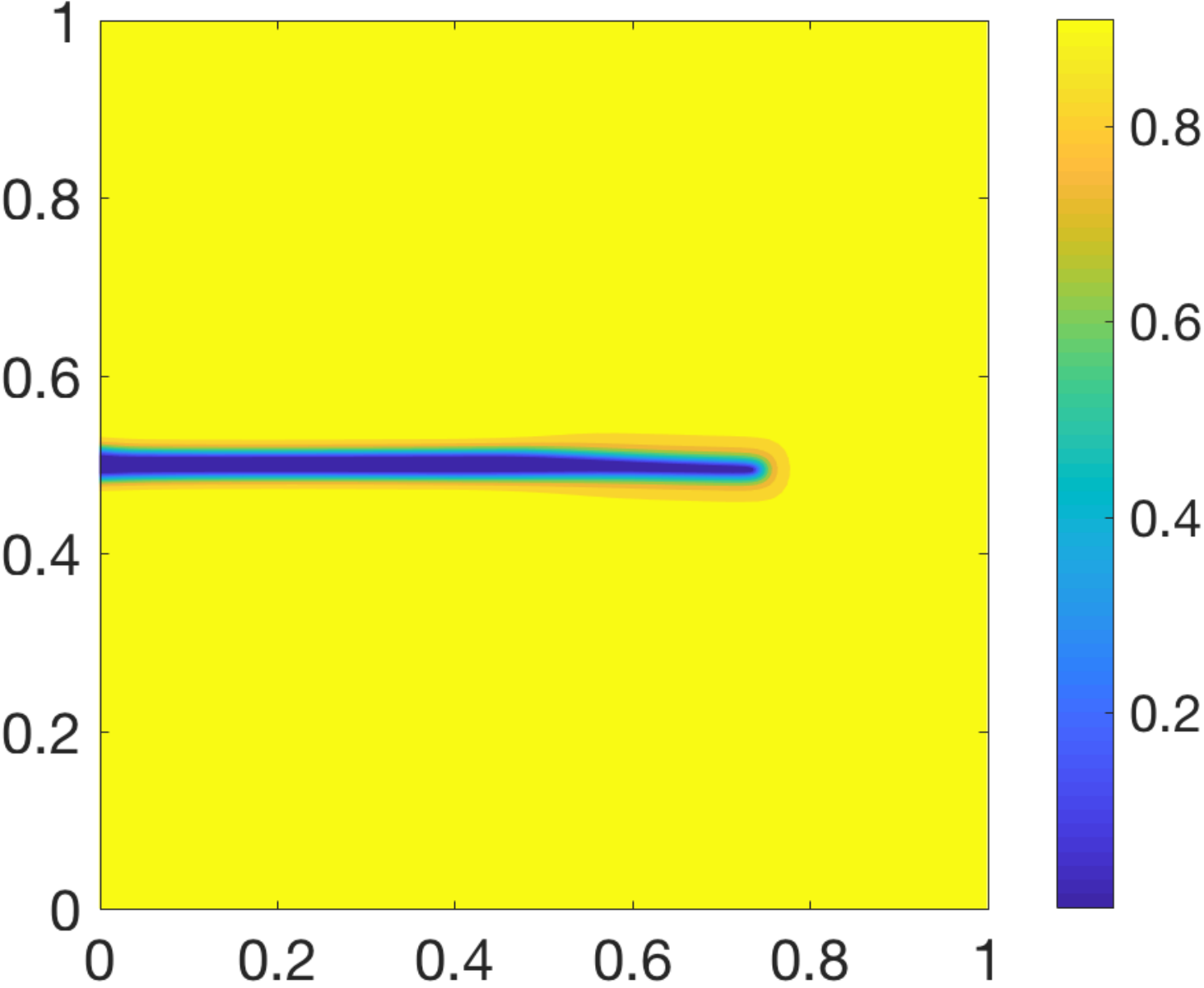}}
\subfigure[improved v-d split]{\label{fig:subfig:TND_U800D04}
\includegraphics[width=0.25\linewidth]{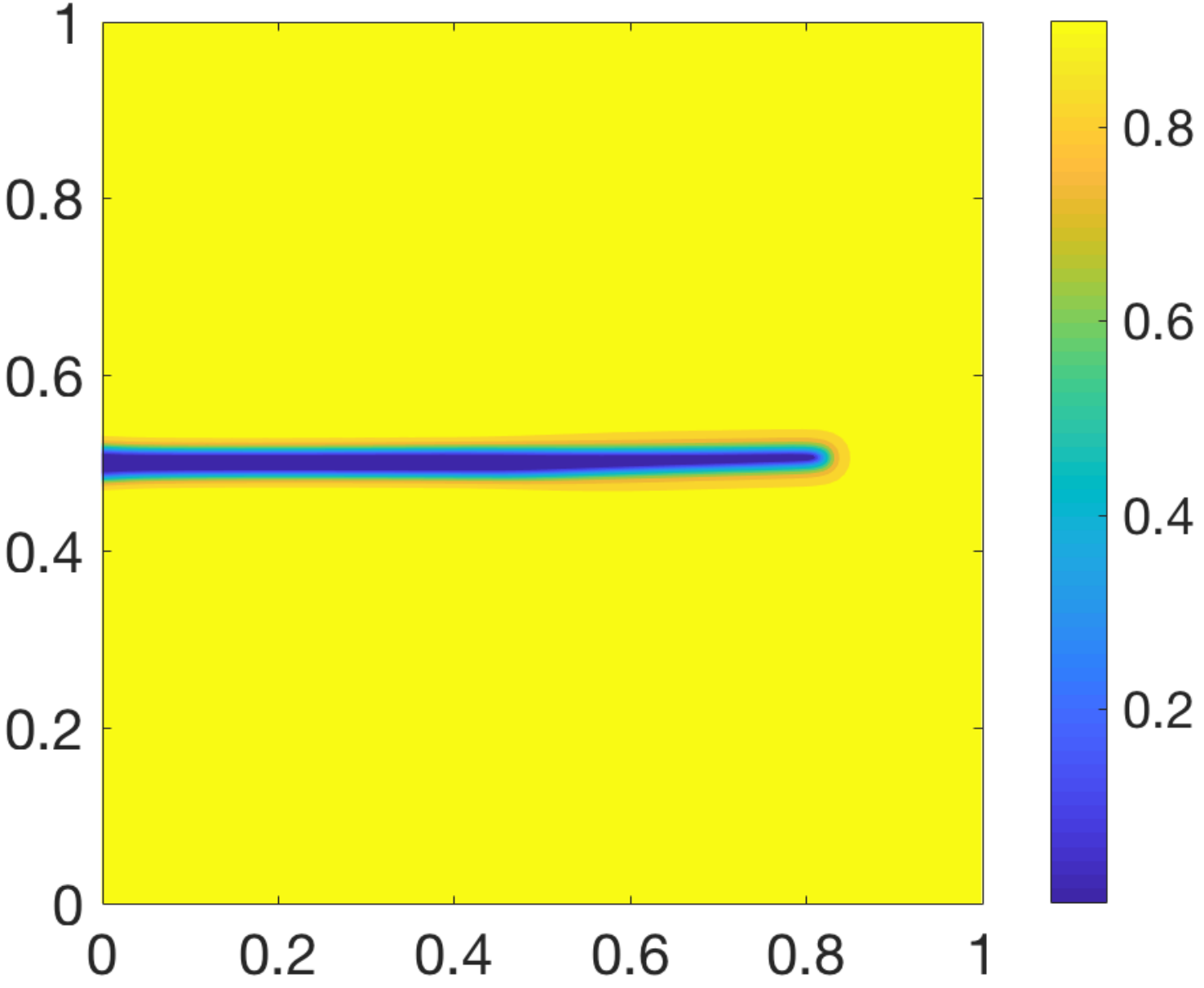}}
\vfill
\subfigure[spectral decomposition]{\label{fig:subfig:TMM_U800D04}
\includegraphics[width=0.25\linewidth]{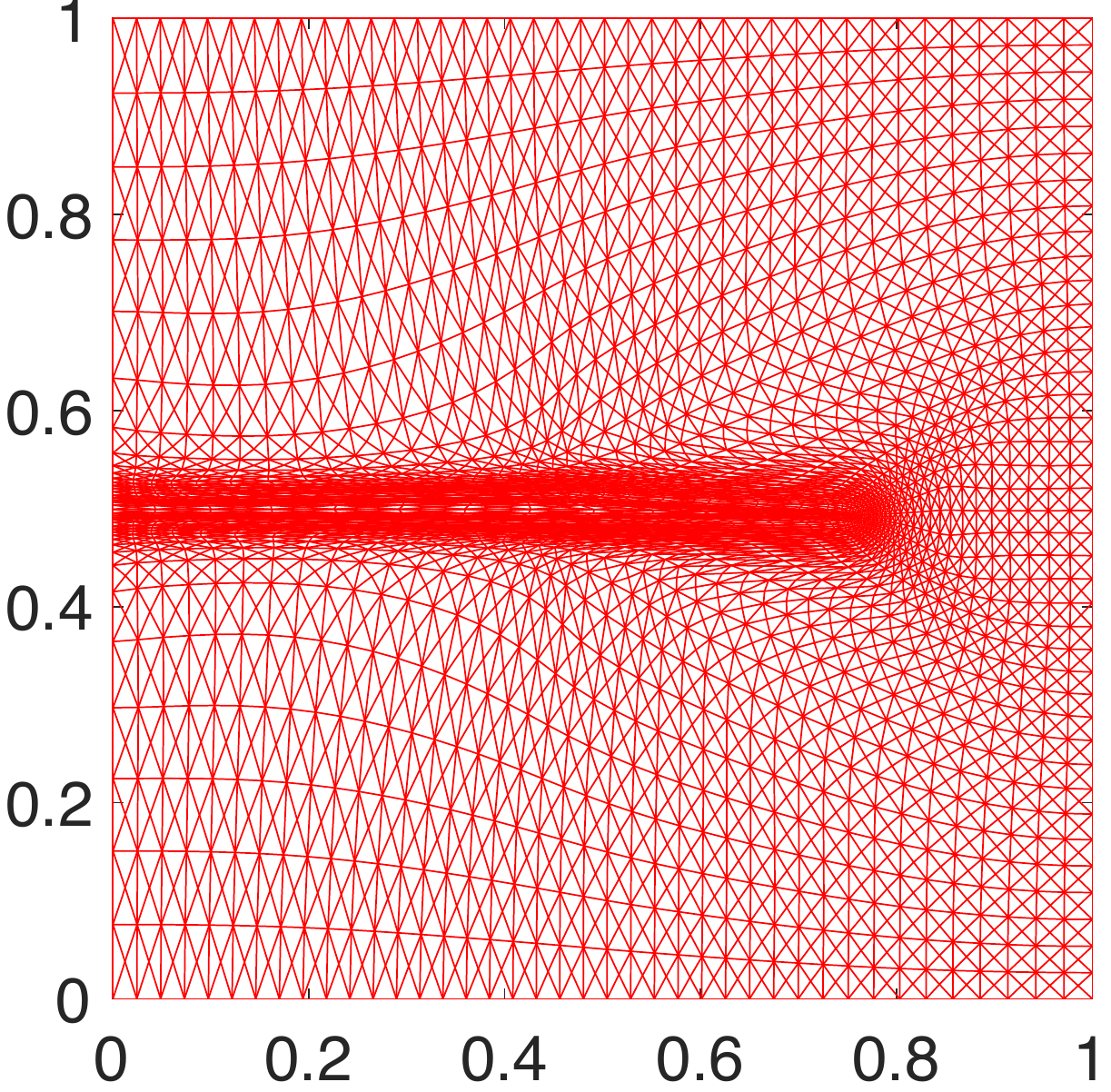}}
\subfigure[v-d split]{\label{fig:subfig:TAM_U800D04}
\includegraphics[width=0.25\linewidth]{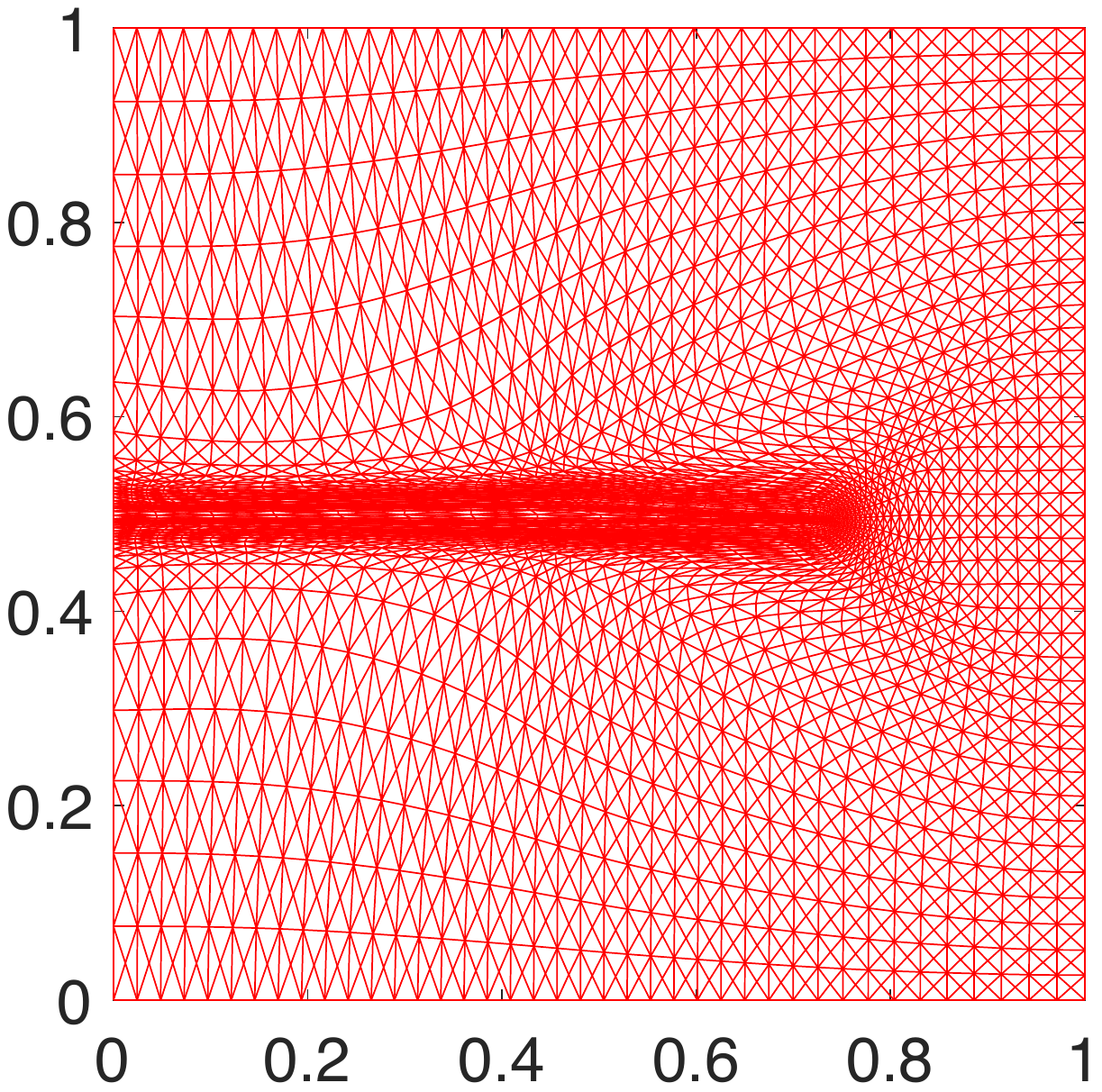}}
\subfigure[improved v-d split]{\label{fig:subfig:TNM_U800D04}
\includegraphics[width=0.25\linewidth]{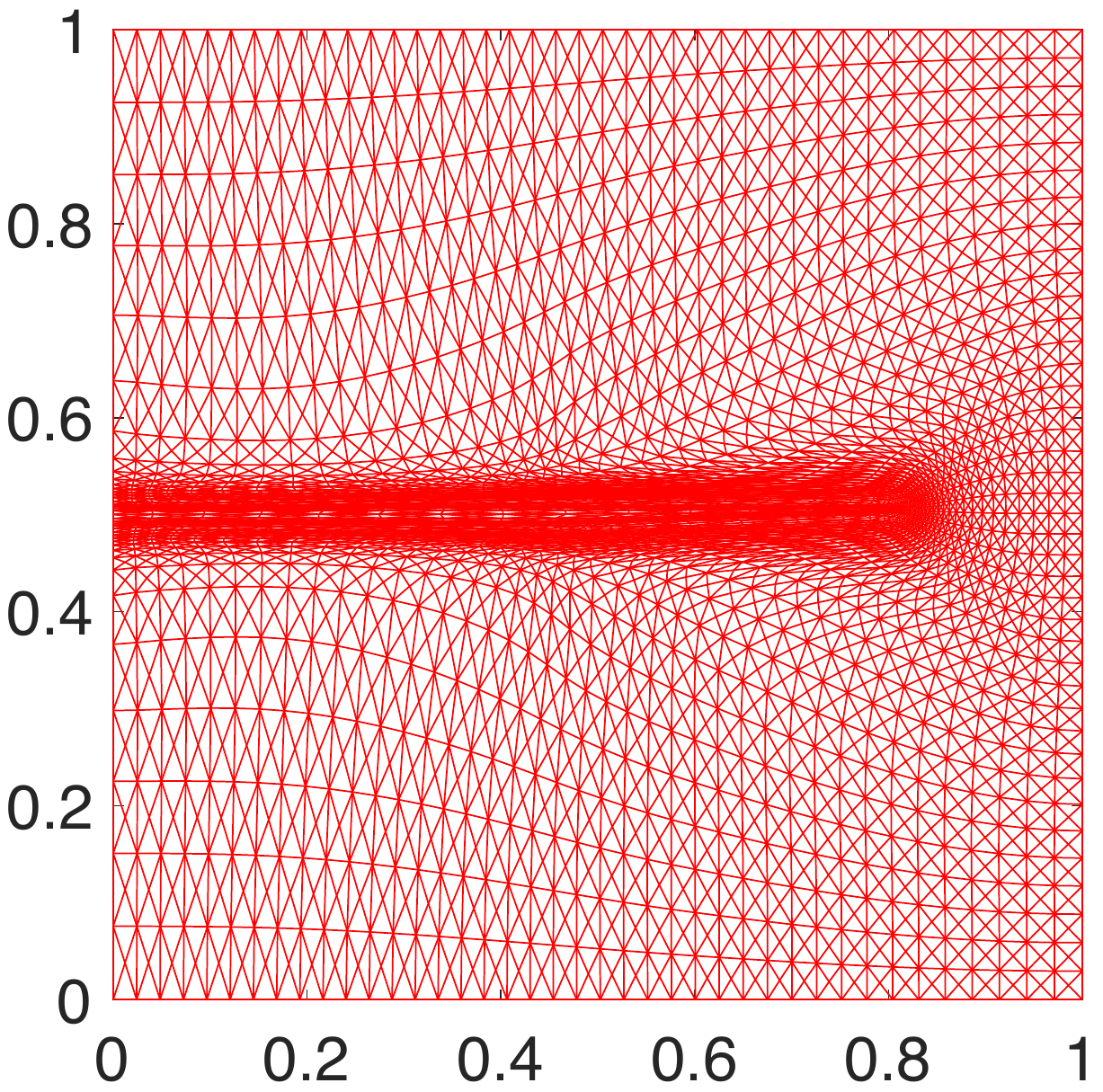}}
\vfill
\subfigure[spectral decomposition]{\label{fig:subfig:TMS_U800D04}
\includegraphics[width=0.25\linewidth]{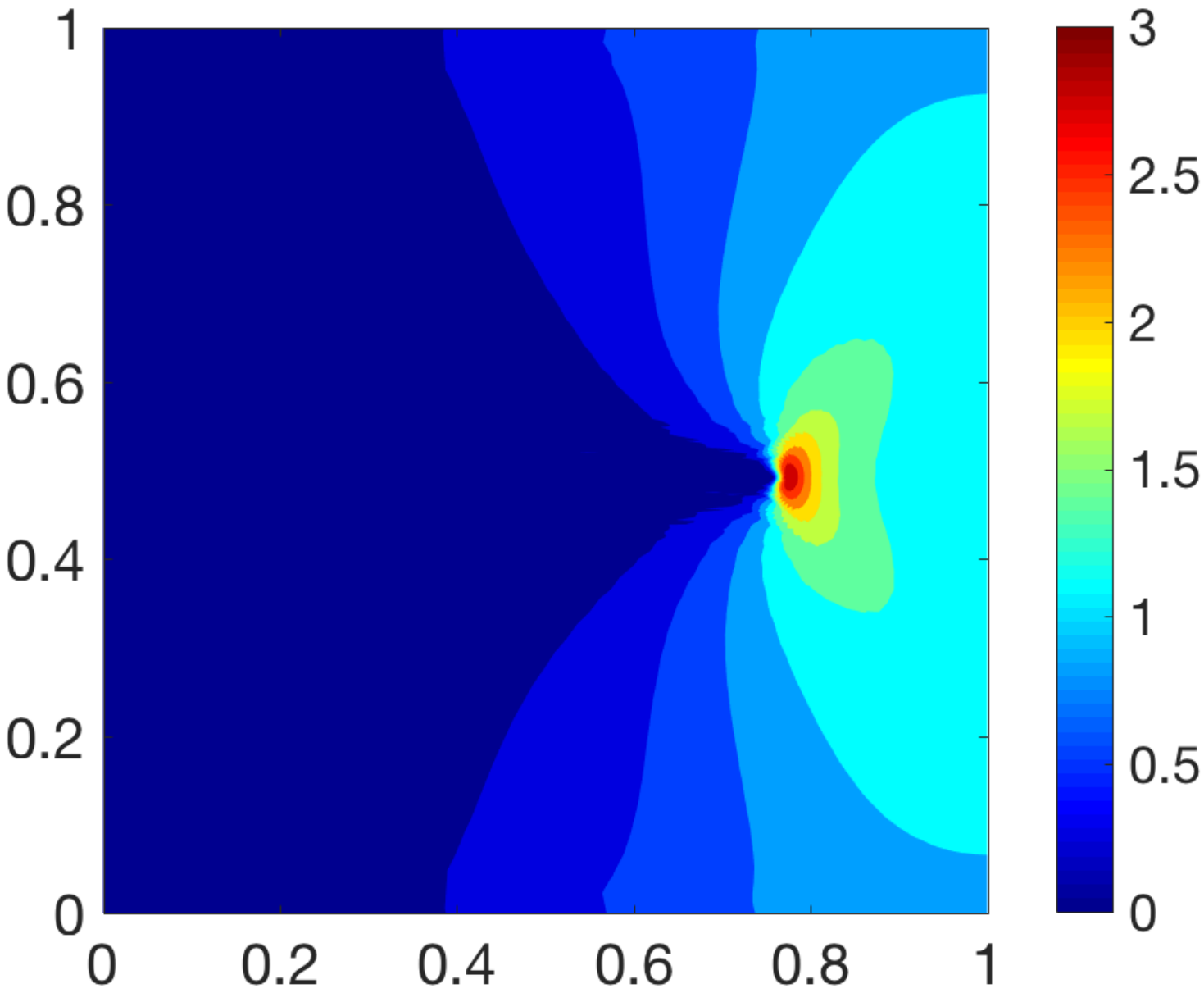}}
\subfigure[v-d split]{\label{fig:subfig:TAS_U800D04}
\includegraphics[width=0.25\linewidth]{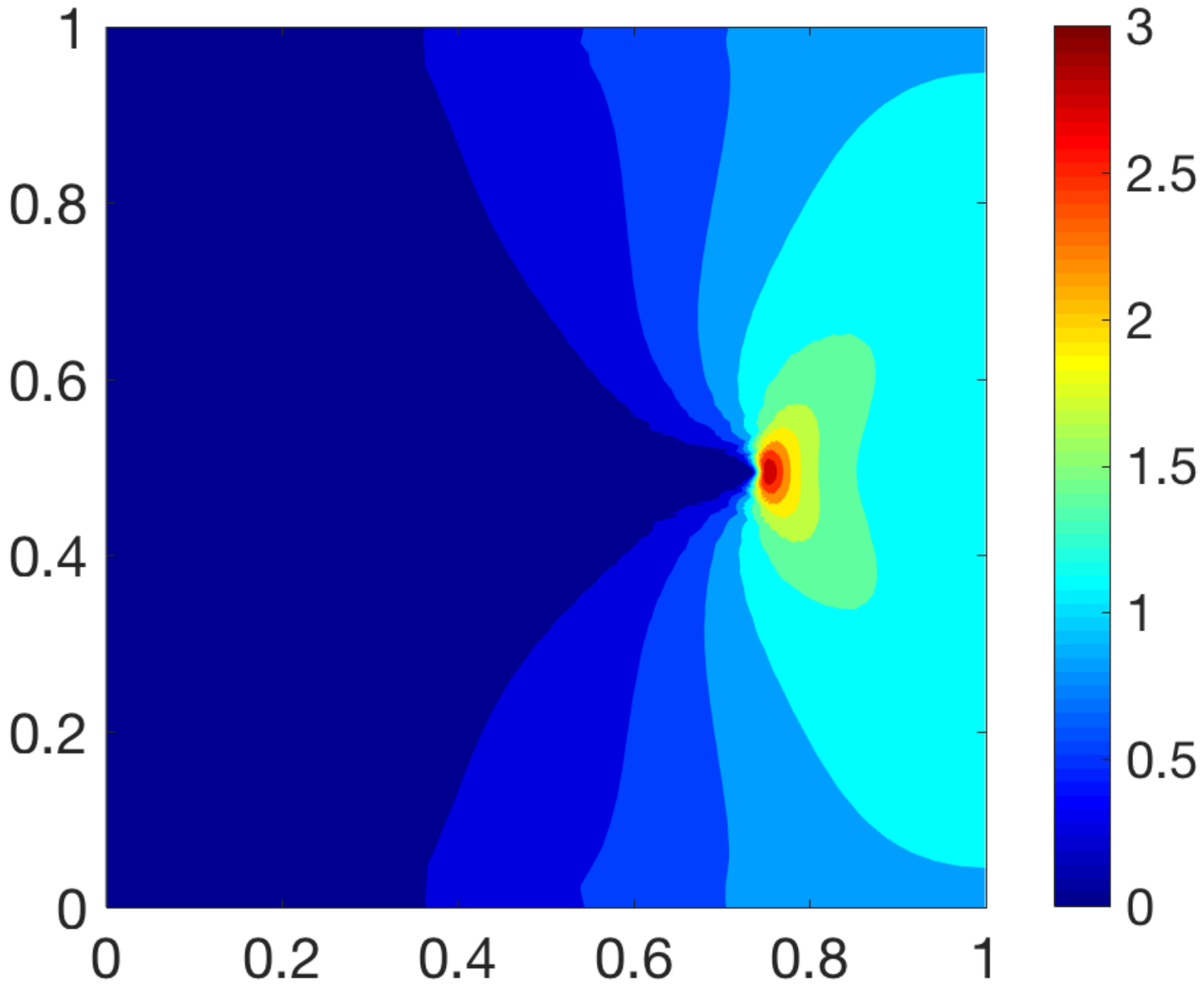}}
\subfigure[improved v-d split]{\label{fig:subfig:TNS_U800D04}
\includegraphics[width=0.25\linewidth]{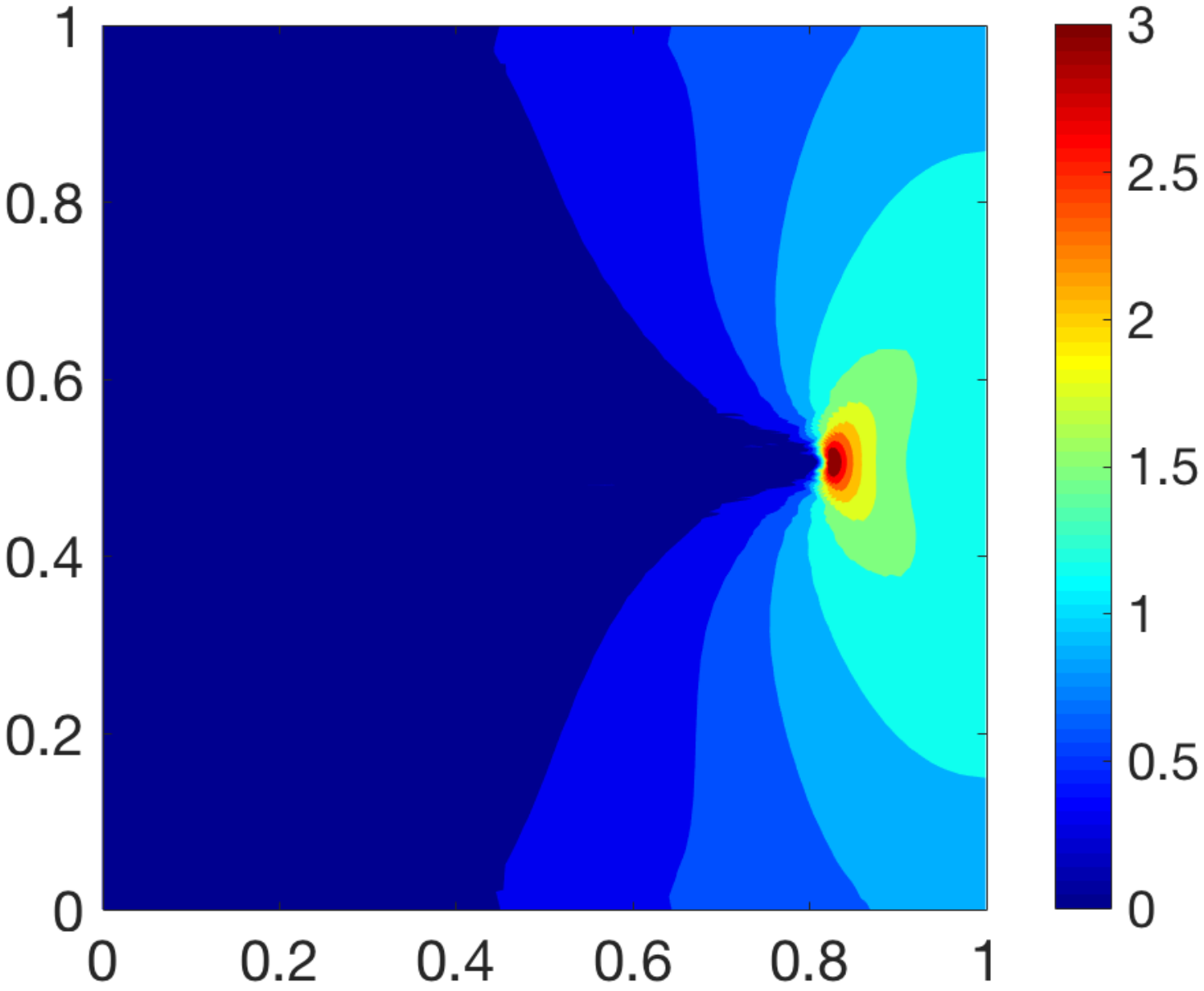}}
\caption{Example 1. Meshes and contours of the phase-field and von Mises stress distribution
are plotted at $U = 5.3 \times 10^{-3}$~mm.  Three decomposition models with ItCBC
($d_{cr} = 0.4$) are used.}
\label{fig:tension contours with dcr}
\end{figure}

Next, we examine the effects of ItCBC (cf. \S\ref{SEC:ItCBC}) and the choice of the threshold $d_{cr}$.
Using ItCBC with $d_{cr} = 0.4$, typical adaptive meshes and contours of the phase-field variable and
von Mises stress at $U = 5.3 \times 10^{-3}$~mm are shown in Fig. \ref{fig:tension contours with dcr}. 
One can see that the effects of ItCBC for both the spectral decomposition and v-d split models are small.
On the other hand, the effects of ItCBC on the improved v-d split model are significant. Indeed,
the von Mises stress distribution for the latter is now comparable to those obtained with
the spectral decomposition and v-d split models; see Fig. \ref{fig:subfig:TNS_U800D04}. 
 
\begin{figure} 
\centering 
\subfigure[spectral decomposition]{\label{fig:subfig:Tld_MieheDiffDcr}
\includegraphics[width=0.4\linewidth]{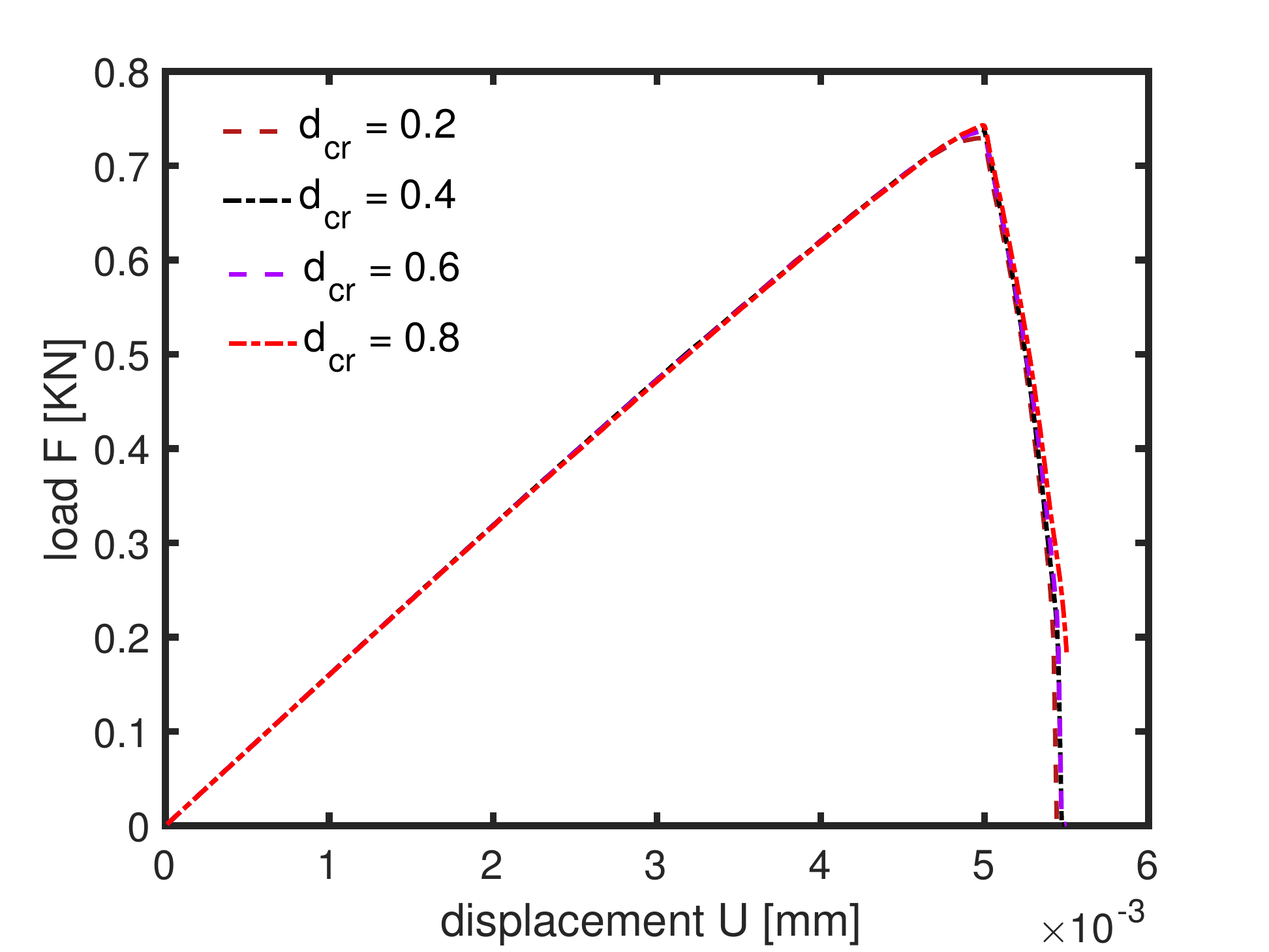}}
\subfigure[v-d split]{\label{fig:subfig:Tld_AmorDiffDcr}
\includegraphics[width=0.4\linewidth]{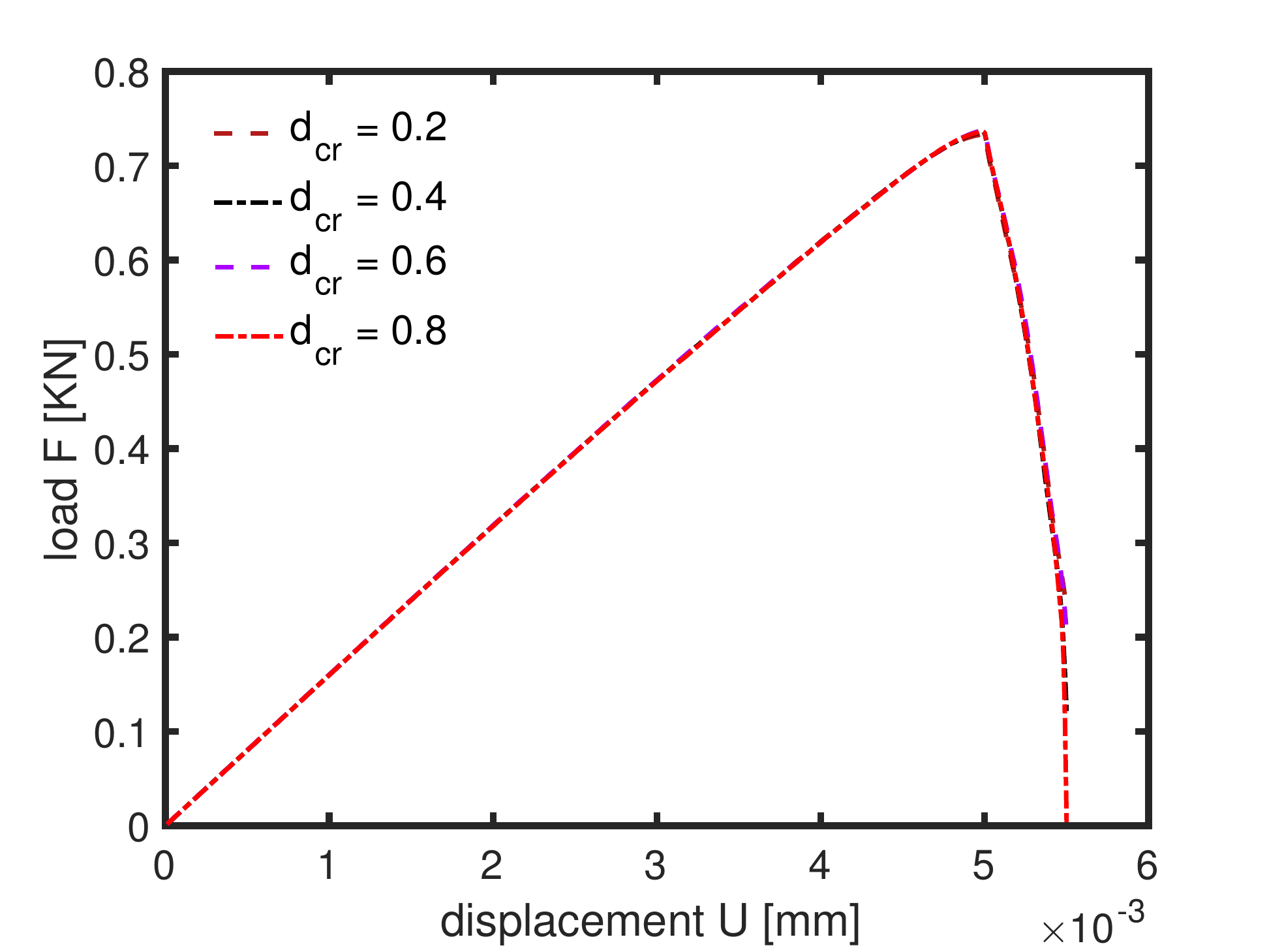}}
\subfigure[improved v-d split]{\label{fig:subfig:Tld_NewDiffDcr}
\includegraphics[width=0.4\linewidth]{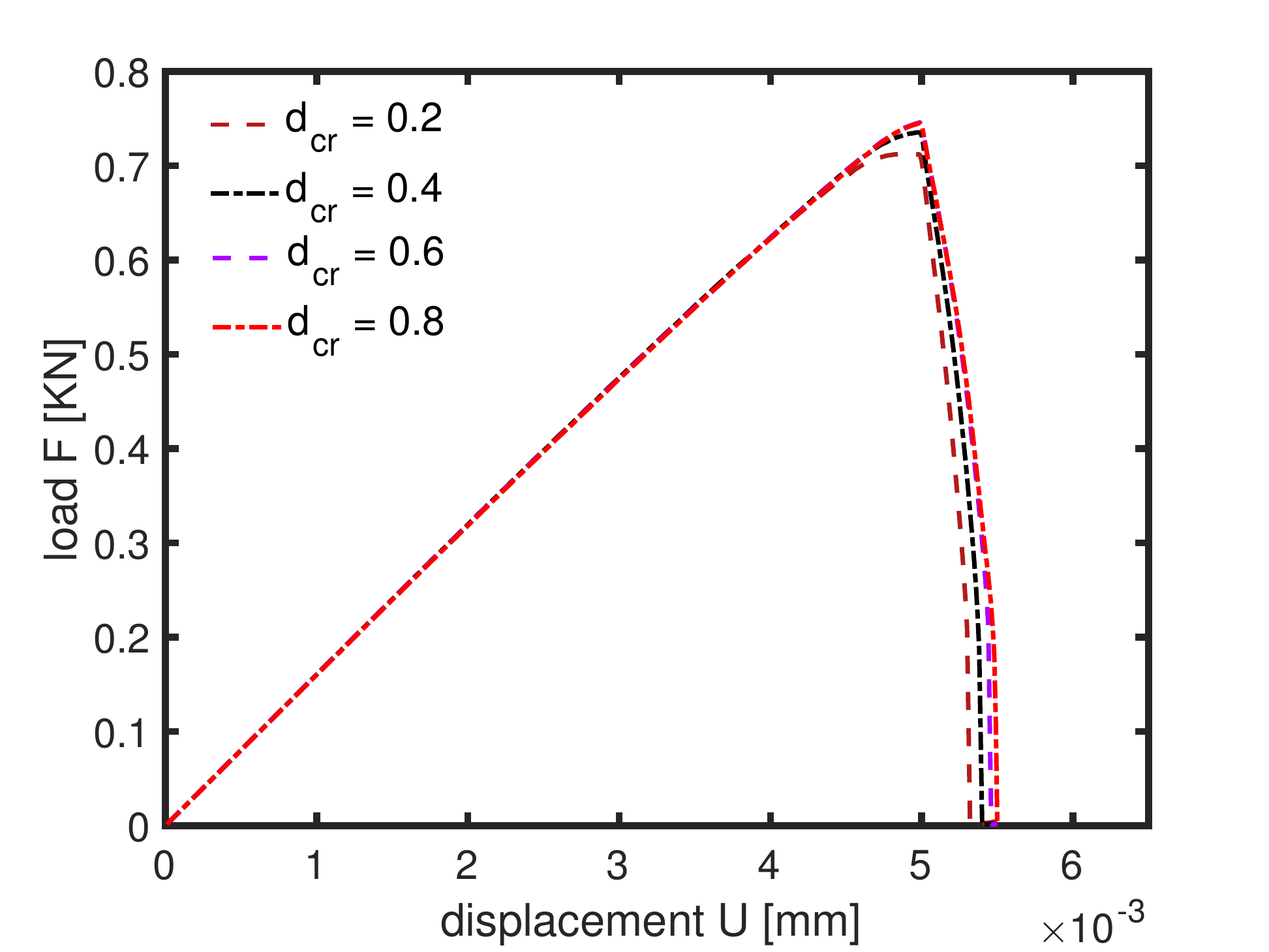}}
\caption{Example 1. The load-deflection curves are obtained using three decomposition models with
or without ItCBC.
(a) spectral decomposition with various $d_{cr}$;
(b) v-d split with various $d_{cr}$;
(c) improved v-d split with various $d_{cr}$.
}
\label{fig:tension effects of the dcr}
\end{figure}

\begin{figure} 
\centering 
\subfigure[$d_{cr} = 0.2$]{\label{fig:subfig:Tld_DiffSplit02}
\includegraphics[width=0.4\linewidth]{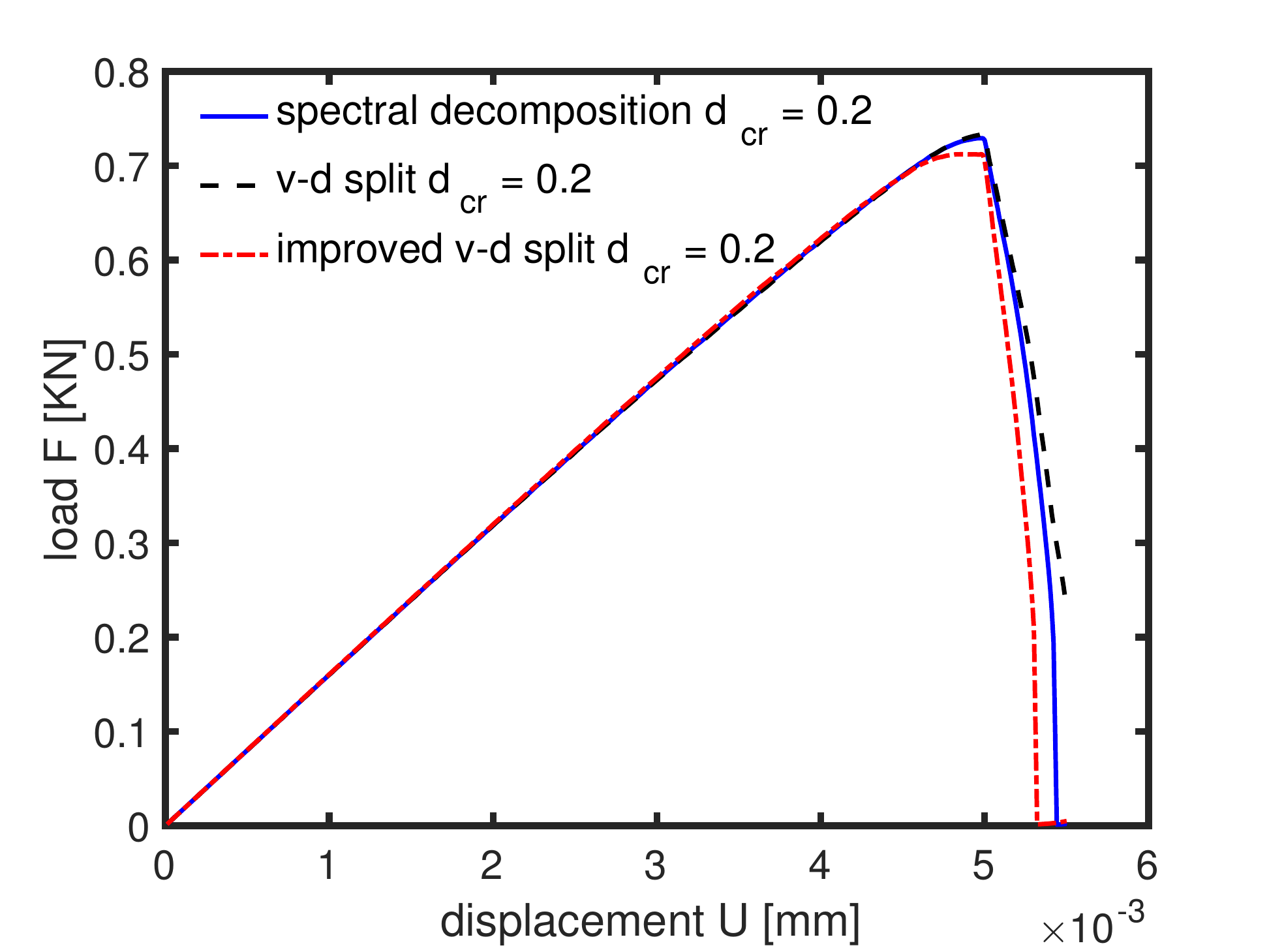}}
\subfigure[$d_{cr} = 0.4$]{\label{fig:subfig:Tld_DiffSplit04}
\includegraphics[width=0.4\linewidth]{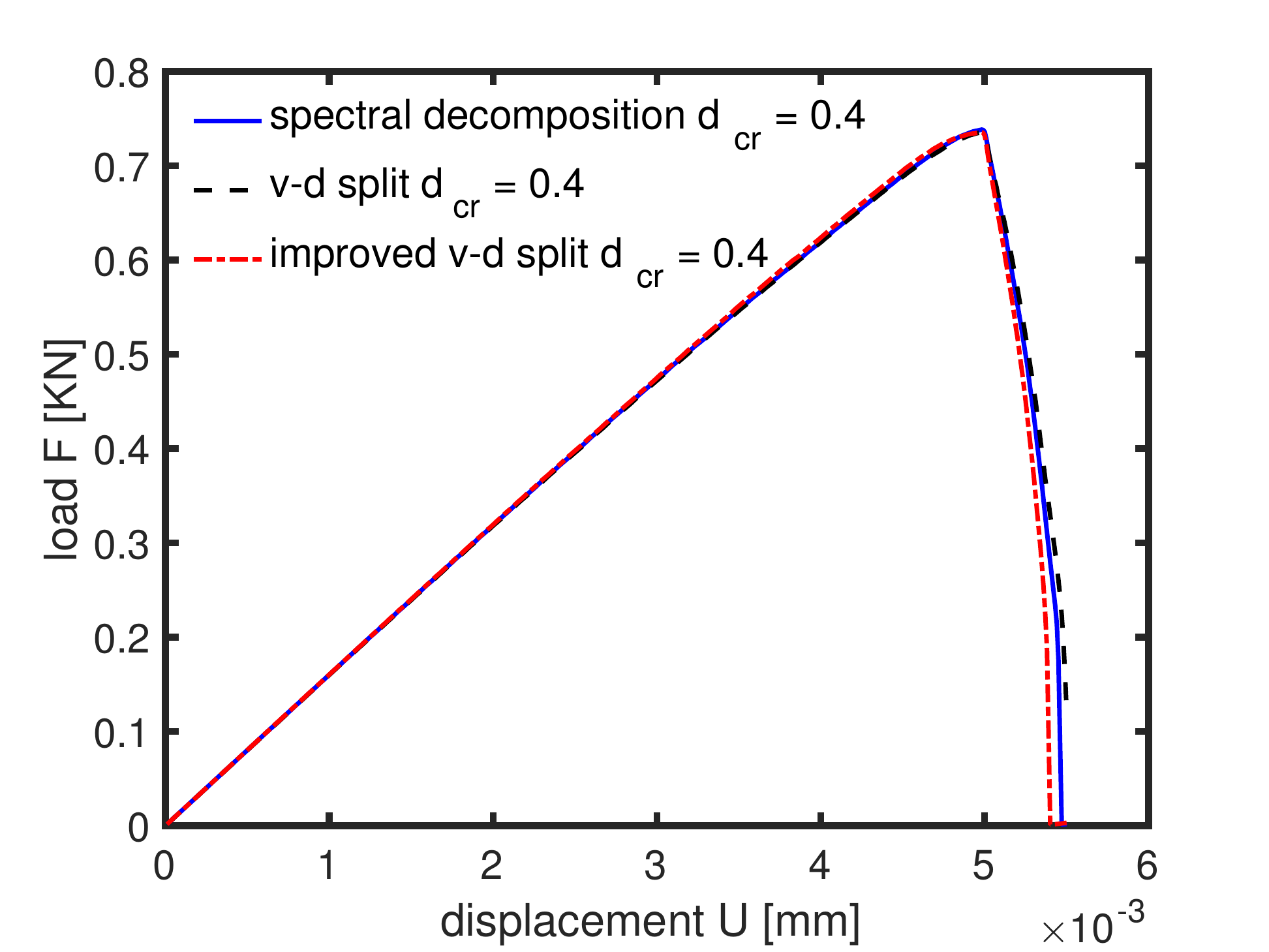}}
\subfigure[$d_{cr} = 0.6$]{\label{fig:subfig:Tld_DiffSplit06}
\includegraphics[width=0.4\linewidth]{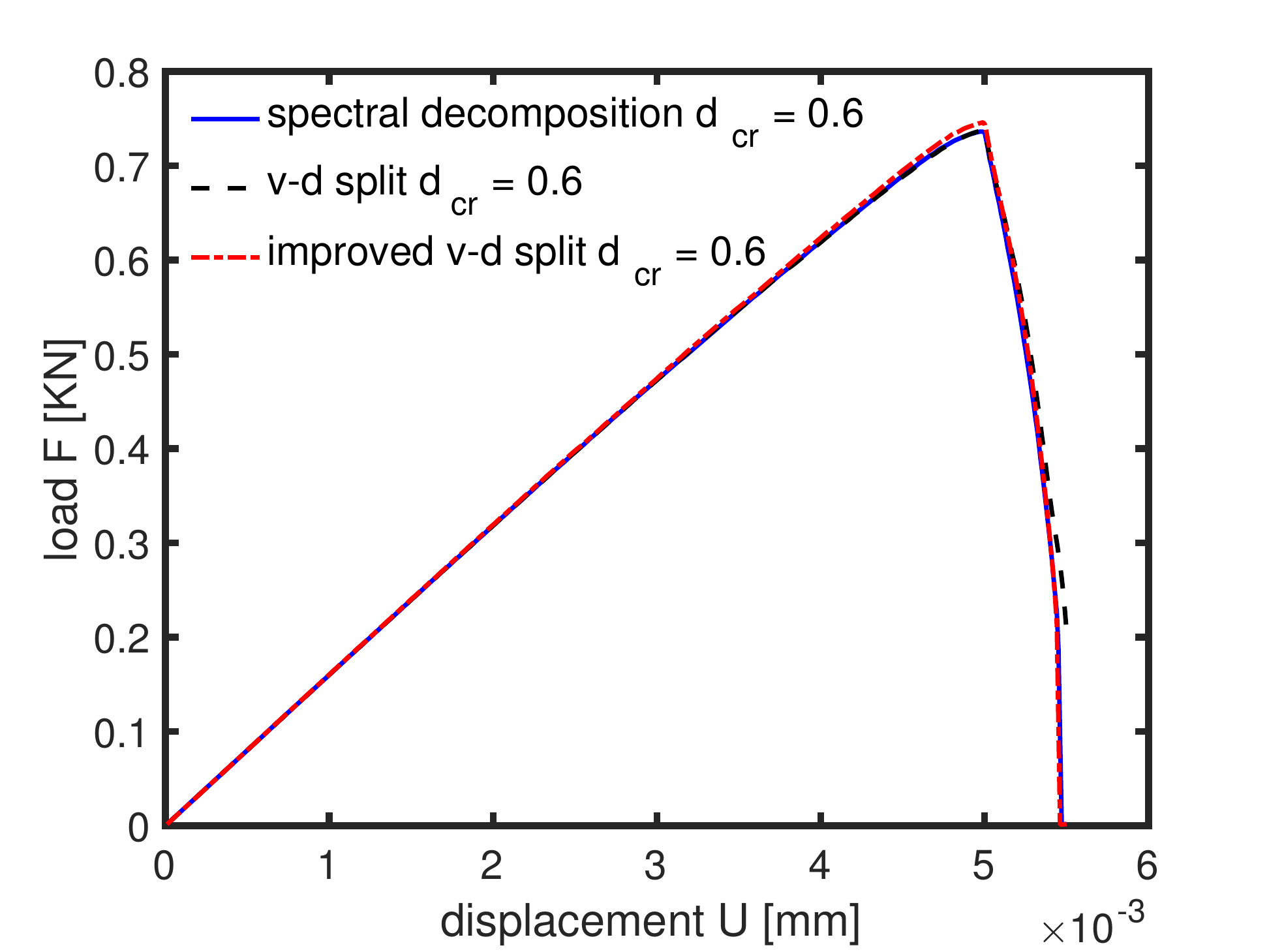}}
\subfigure[$d_{cr} = 0.8$]{\label{fig:subfig:Tld_DiffSplit08}
\includegraphics[width=0.4\linewidth]{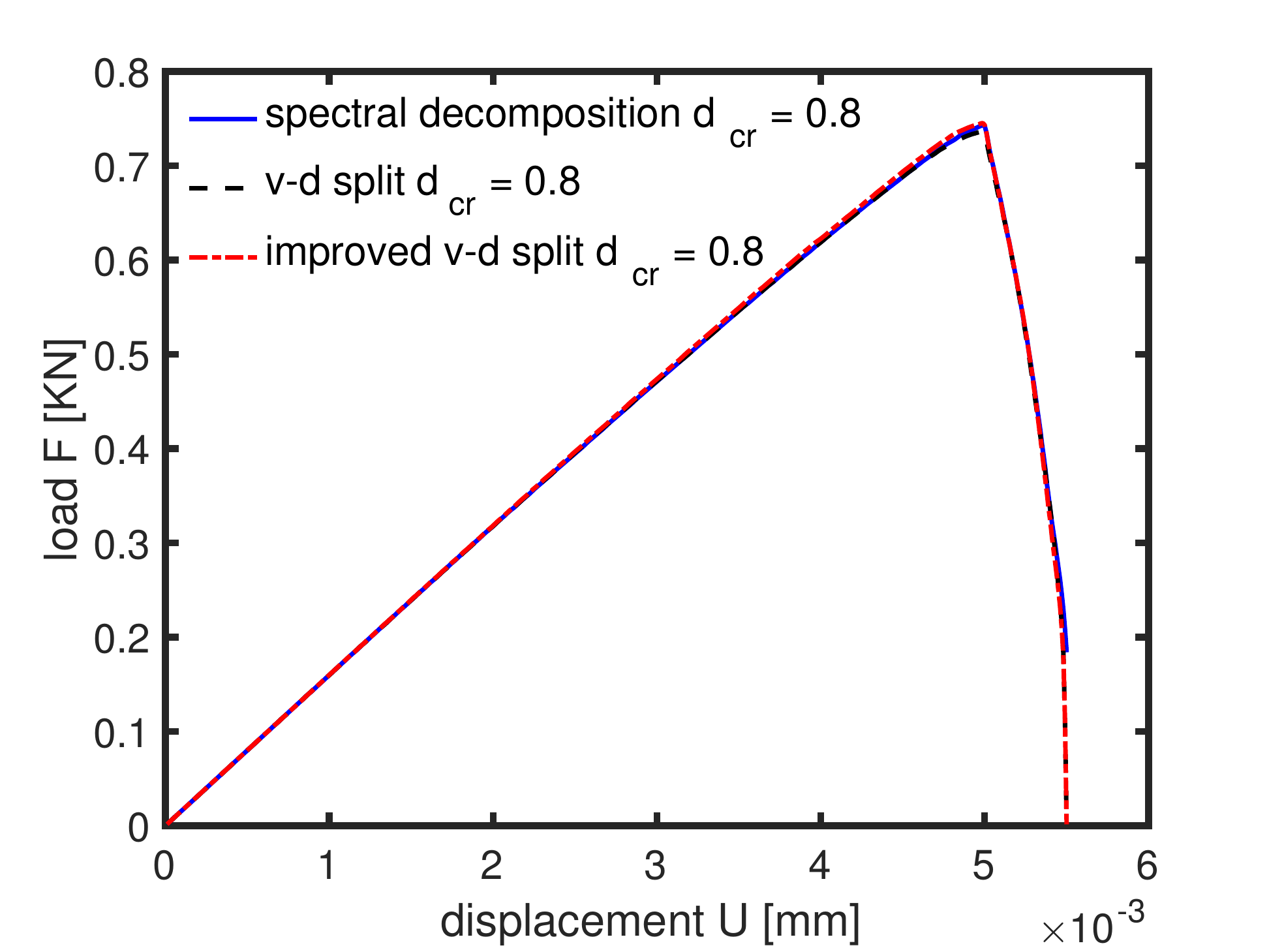}}
\caption{Example 1. The load-deflection curves are compared for three energy decomposition methods with ItCBC.
(a) $d_{cr} = 0.2$;
(b) $d_{cr} = 0.4$;
(c) $d_{cr} = 0.6$;
(d) $d_{cr} = 0.8$.}
\label{fig:tension comparison of the different splits}
\end{figure}

The load-deflection curves for three decomposition models with various values of $d_{cr}$ (0.2, 0.4, 0.6, 0.8)
are shown in Figs. \ref{fig:tension effects of the dcr} and \ref{fig:tension comparison of the different splits}. 
As we can see in Figs.~\ref{fig:subfig:Tld_MieheDiffDcr} and \ref{fig:subfig:Tld_AmorDiffDcr}, there is no significant difference between the original model and the model with ItCBC with various values of $d_{cr}$ for the spectral decomposition and v-d split. However, for the improved v-d split,  smaller $d_{cr}$ values lead to underestimates
of the load after crack starts propagating, see \ref{fig:subfig:Tld_NewDiffDcr}.
The work done by boundary tractions and body forces (external work) on an elastic solid
are stored inside the body in the form of strain energy. According to the Griffith theory,
the energy required to create new crack surfaces is transformed from the stored elastic strain energy.
When $d_{cr}$ is large, more strain energy is degraded and stored energy becomes less.
In this case, more external work is needed to generate new crack,
which means that the peak of reaction force increases with $d_{cr}$
(see Fig.~\ref{fig:tension effects of the dcr}).

Figs. \ref{fig:subfig:Tld_DiffSplit04}, \ref{fig:subfig:Tld_DiffSplit06}, and \ref{fig:subfig:Tld_DiffSplit08} show that
for ItCBC with $d_{cr} \in [0.4, 0.8]$, the load-deflection curves are nearly the same for all three energy decomposition
models. For smaller $d_{cr}$ values ($d_{cr} = 0.2$ in Fig. \ref{fig:subfig:Tld_DiffSplit02}),
the improved v-d split underestimates the load after crack starts propagating.

\subsection{Example 2. A single edge notched shear test}

We now consider a single edge notched shear test. The domain and boundary conditions are shown in
Fig.~\ref{fig:subfig:Shear}. The bottom edge of the domain is fixed and the top edge is fixed along the $y$-direction
while a uniform $x$-displacement $U$ is increased with time to drive the crack propagation.
The material properties are the same as the tension test in Example 1, that is,
$\lambda = 121.15$~kN/mm$^{2}$, $\mu = 80.77$~kN/mm$^{2}$, and $g_c = 2.7 \times 10^{-3}$~kN/mm.
The displacement increment is chosen as $\Delta U = 10^{-5}$~mm for the computation.

\begin{figure} 
\centering 
\subfigure[$U = 2.45 \times 10^{-2}$~mm]{\label{fig:subfig:MD_U2450}
\includegraphics[width=0.25\linewidth]{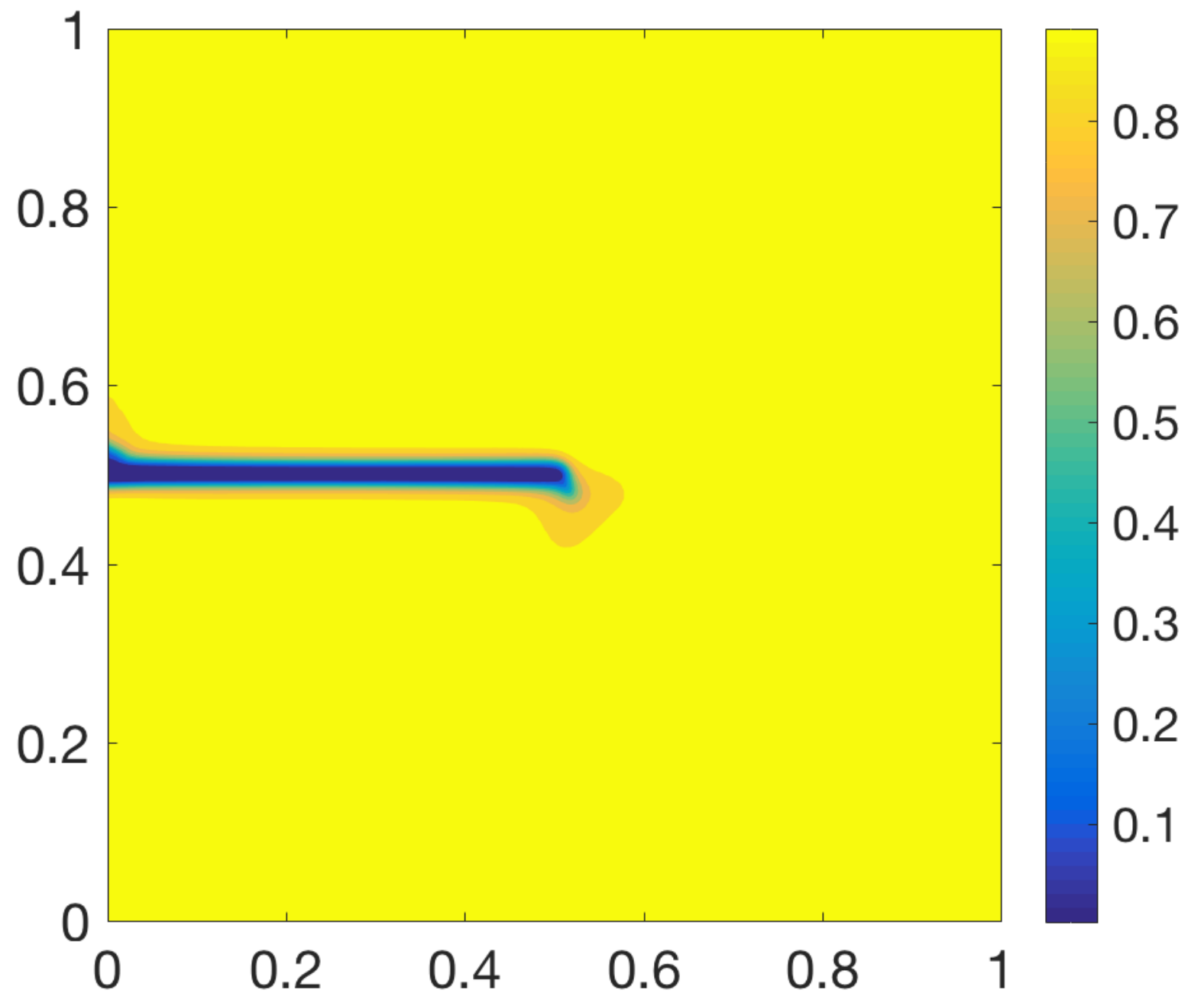}}
\subfigure[$U = 2.55 \times 10^{-2}$~mm]{\label{fig:subfig:MD_U2550}
\includegraphics[width=0.25\linewidth]{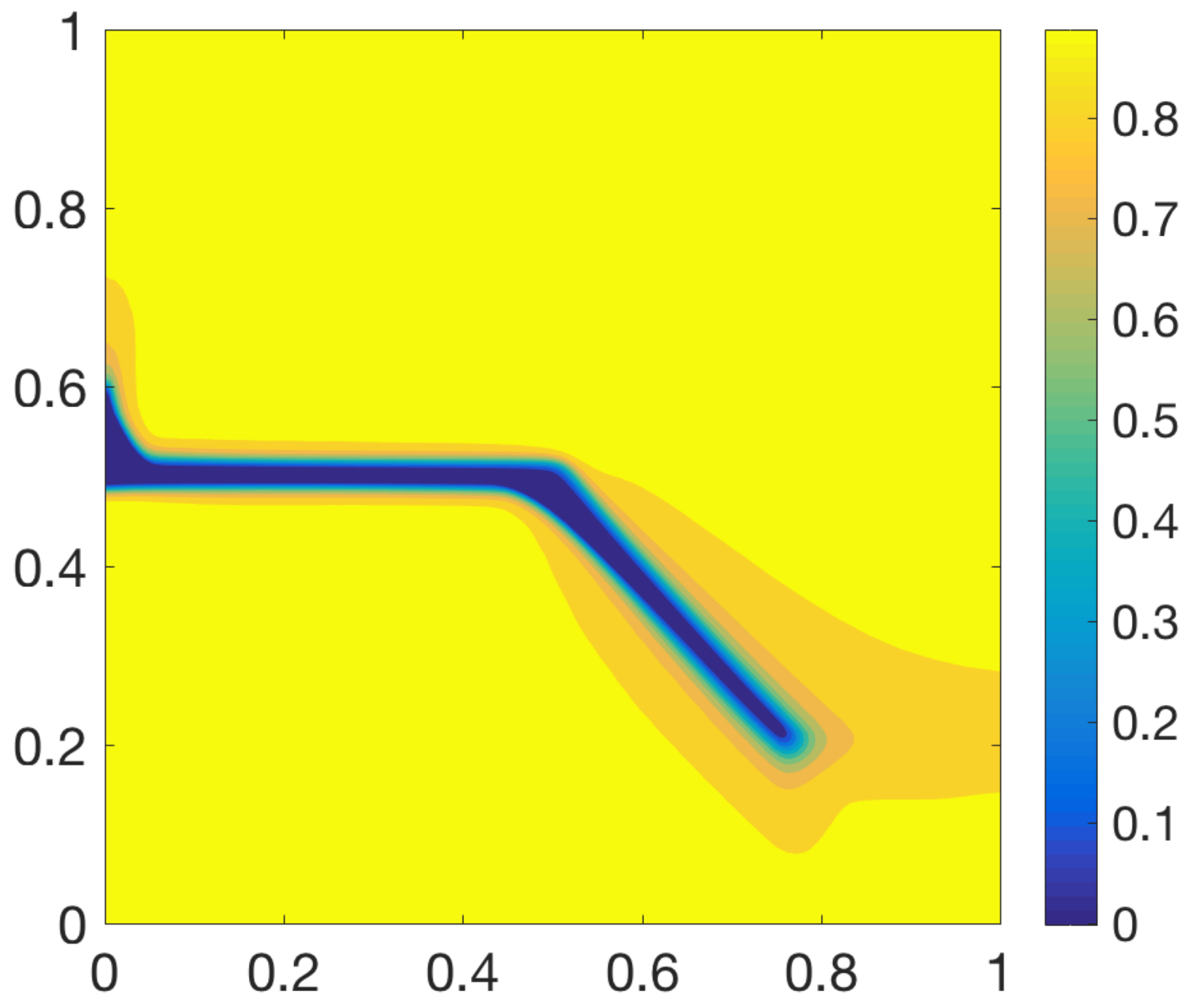}}
\subfigure[$U = 2.6 \times 10^{-2}$~mm]{\label{fig:subfig:MD_U2600}
\includegraphics[width=0.25\linewidth]{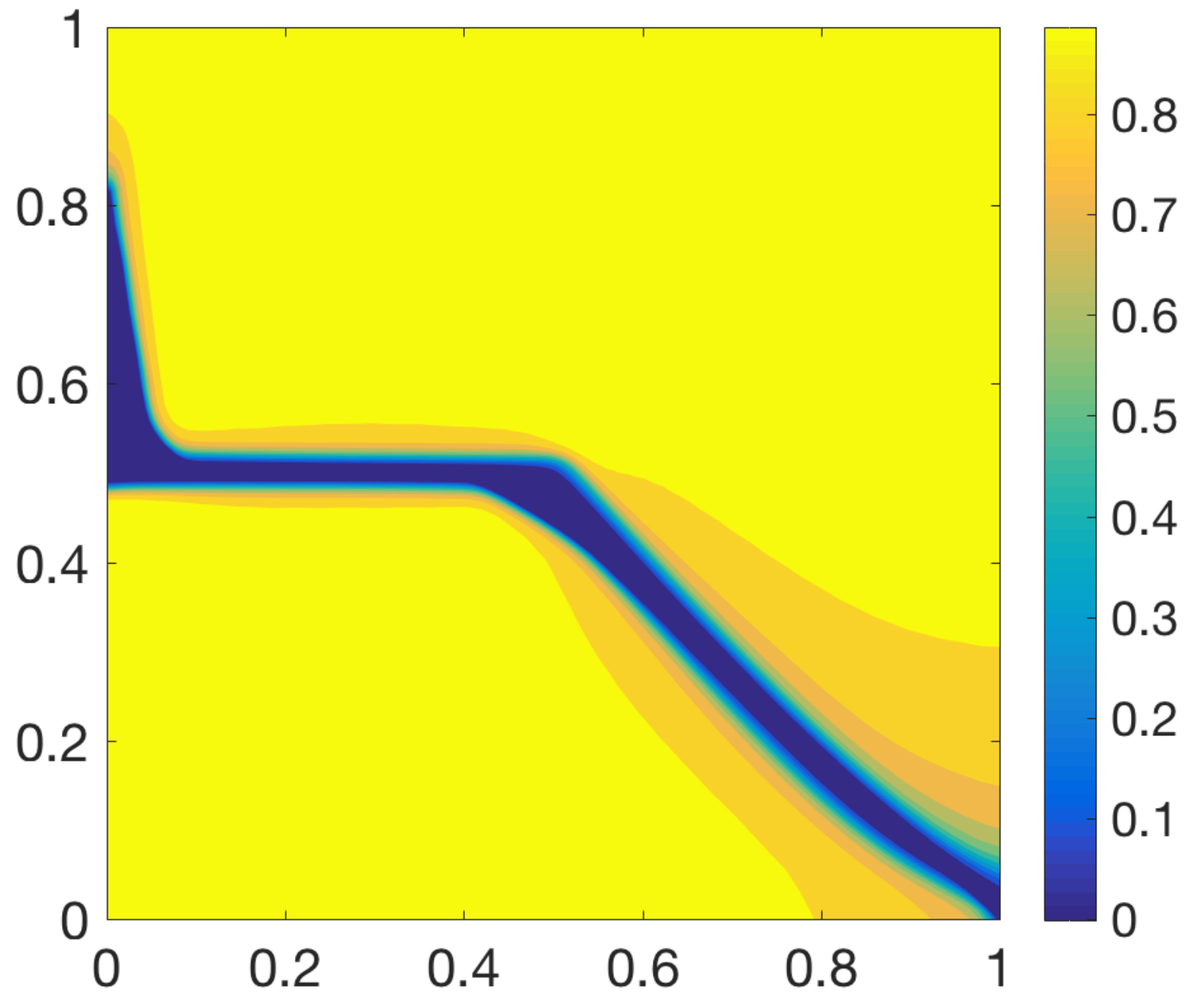}}
\vfill
\subfigure[$U = 2.45 \times 10^{-2}$~mm]{\label{fig:subfig:MM_U2450}
\includegraphics[width=0.25\linewidth]{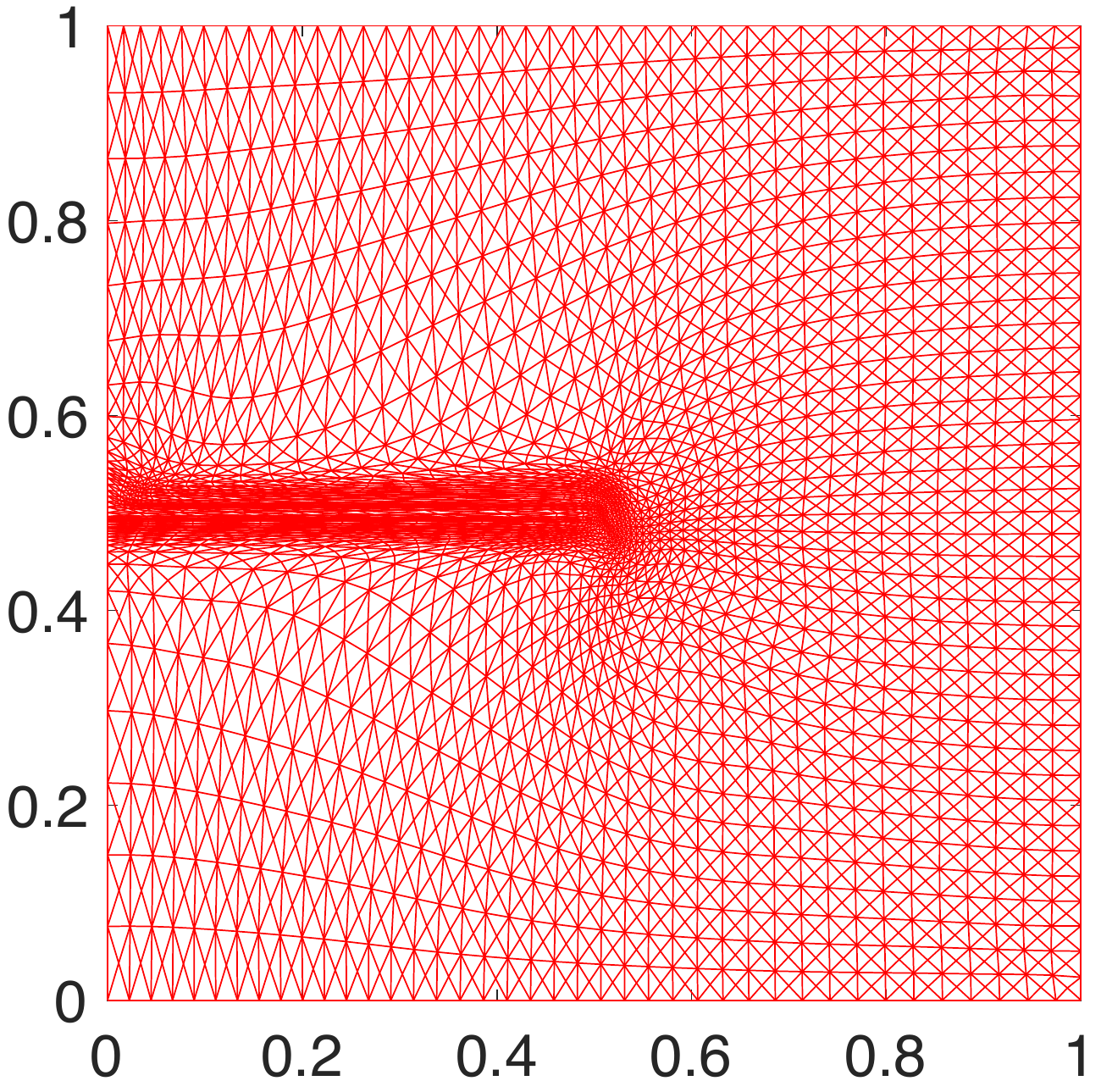}}
\subfigure[$U = 2.55 \times 10^{-2}$~mm]{\label{fig:subfig:MM_U2550}
\includegraphics[width=0.25\linewidth]{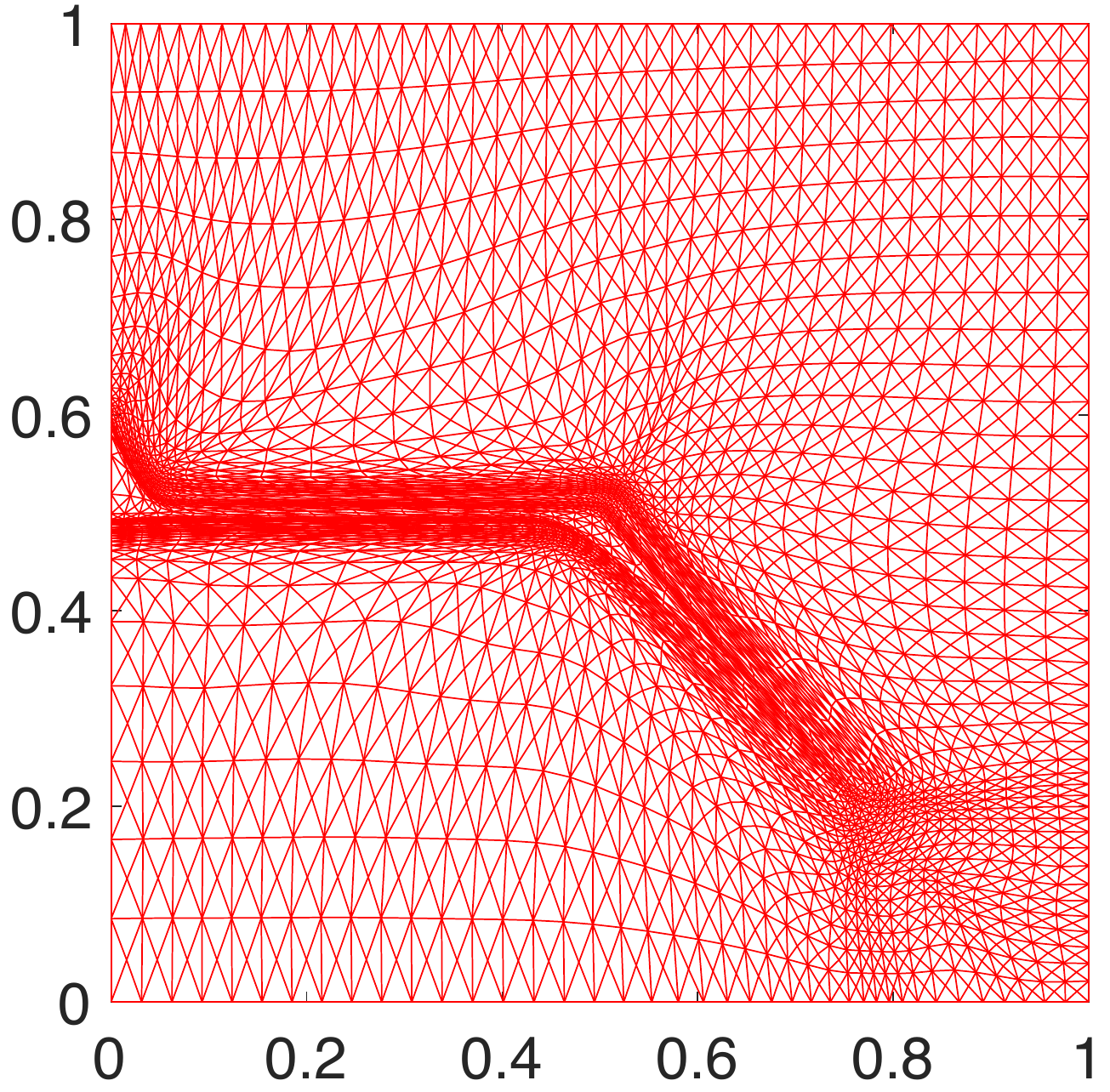}}
\subfigure[$U = 2.6 \times 10^{-2}$~mm]{\label{fig:subfig:MM_U2600}
\includegraphics[width=0.25\linewidth]{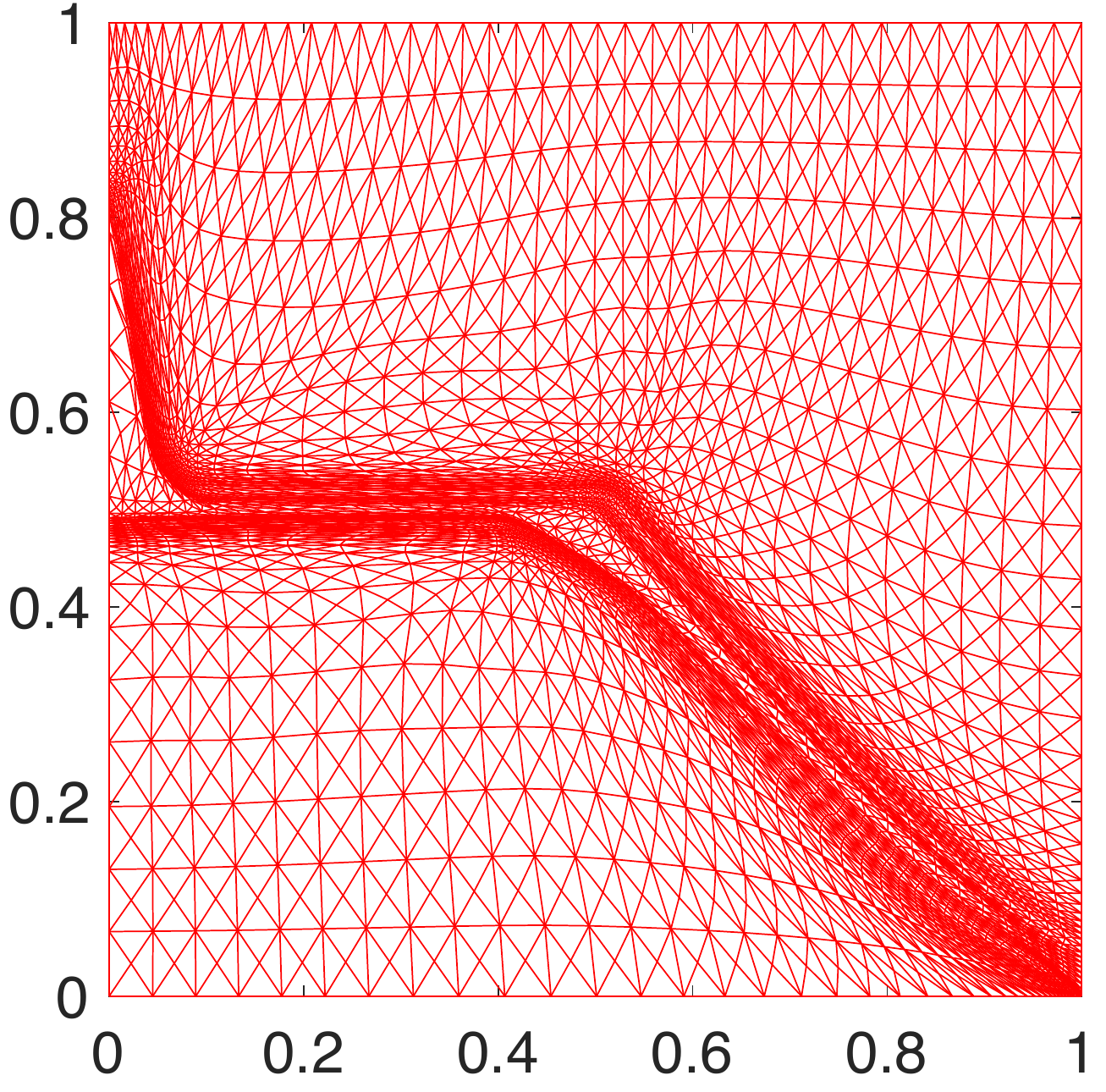}}
\vfill
\subfigure[$U = 2.45 \times 10^{-2}$~mm]{\label{fig:subfig:MS_U2450}
\includegraphics[width=0.25\linewidth]{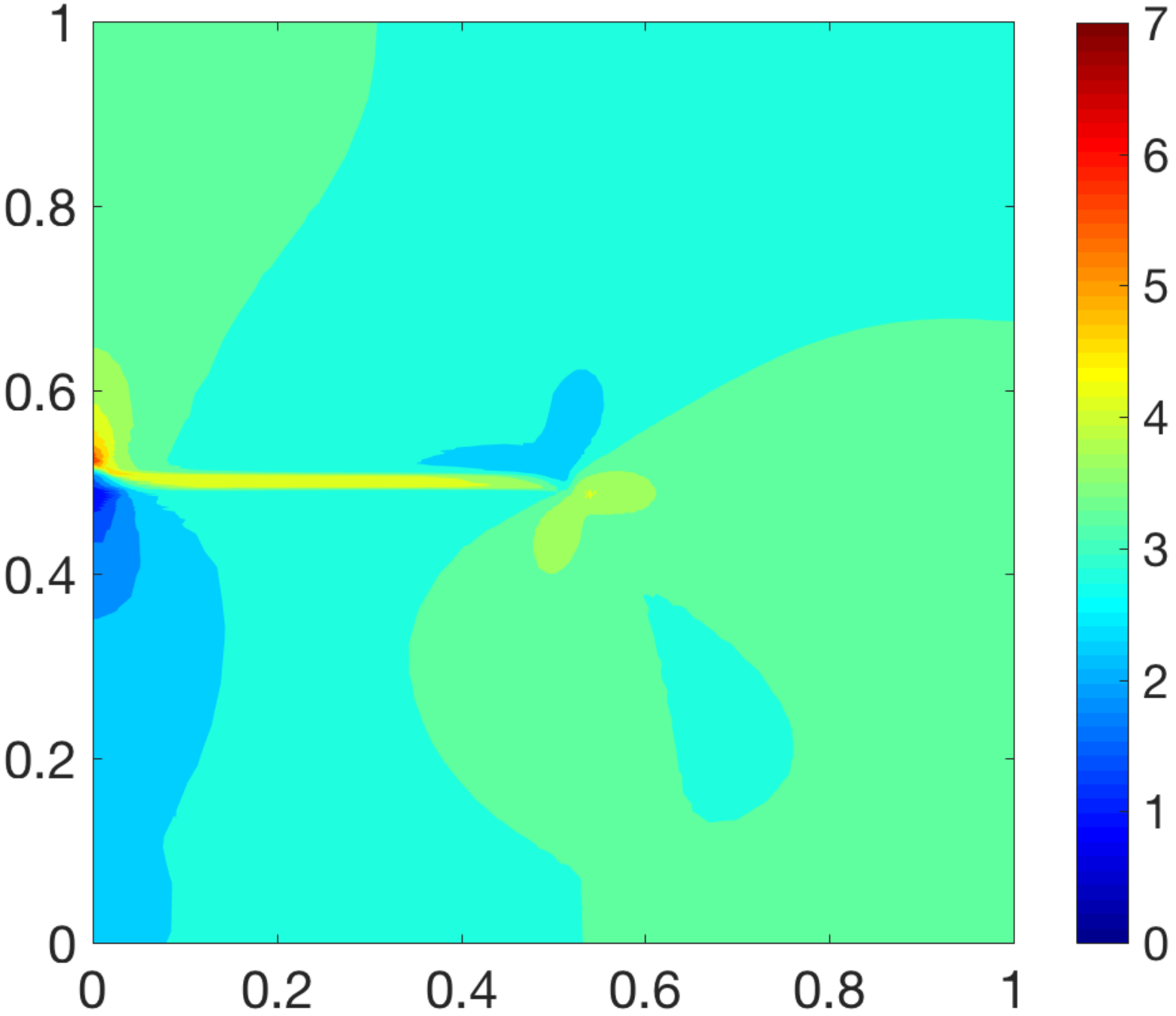}}
\subfigure[$U = 2.55 \times 10^{-2}$~mm]{\label{fig:subfig:MS_U2550}
\includegraphics[width=0.25\linewidth]{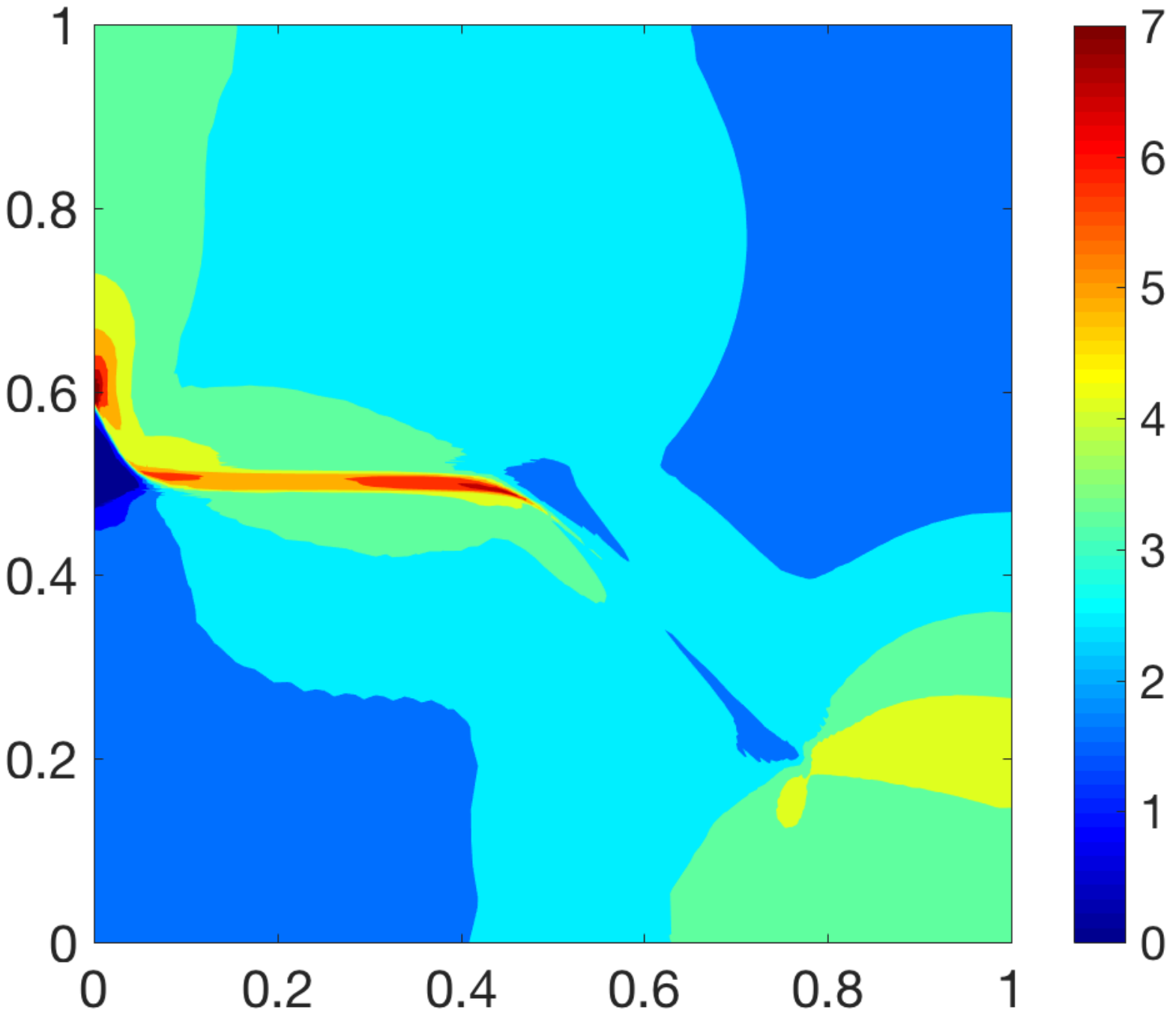}}
\subfigure[$U = 2.6 \times 10^{-2}$~mm]{\label{fig:subfig:MS_U2600}
\includegraphics[width=0.25\linewidth]{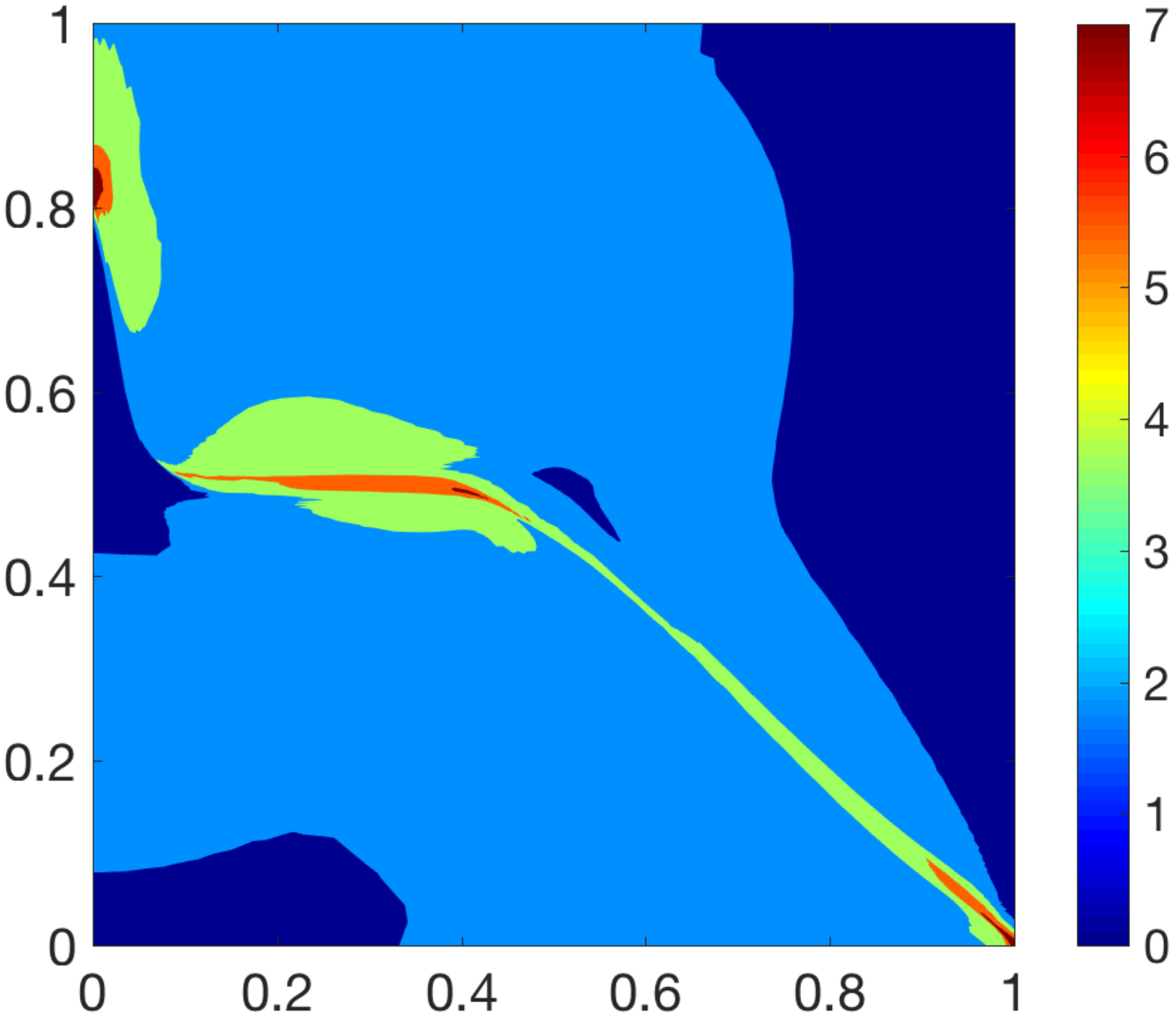}}
\caption{Example 2. Meshes and contours of the phase-field and von Mises stress distribution
are plotted at $U = 2.45\times 10^{-2}$, $2.55\times 10^{-2}$, and $2.60\times 10^{-2}$~mm.
The spectral decomposition model is used.}
\label{fig:shear m's OBC}
\end{figure}

\begin{figure} 
\centering 
\subfigure[$U = 1.1 \times 10^{-2}$~mm]{\label{fig:subfig:MD02_U1100}
\includegraphics[width=0.25\linewidth]{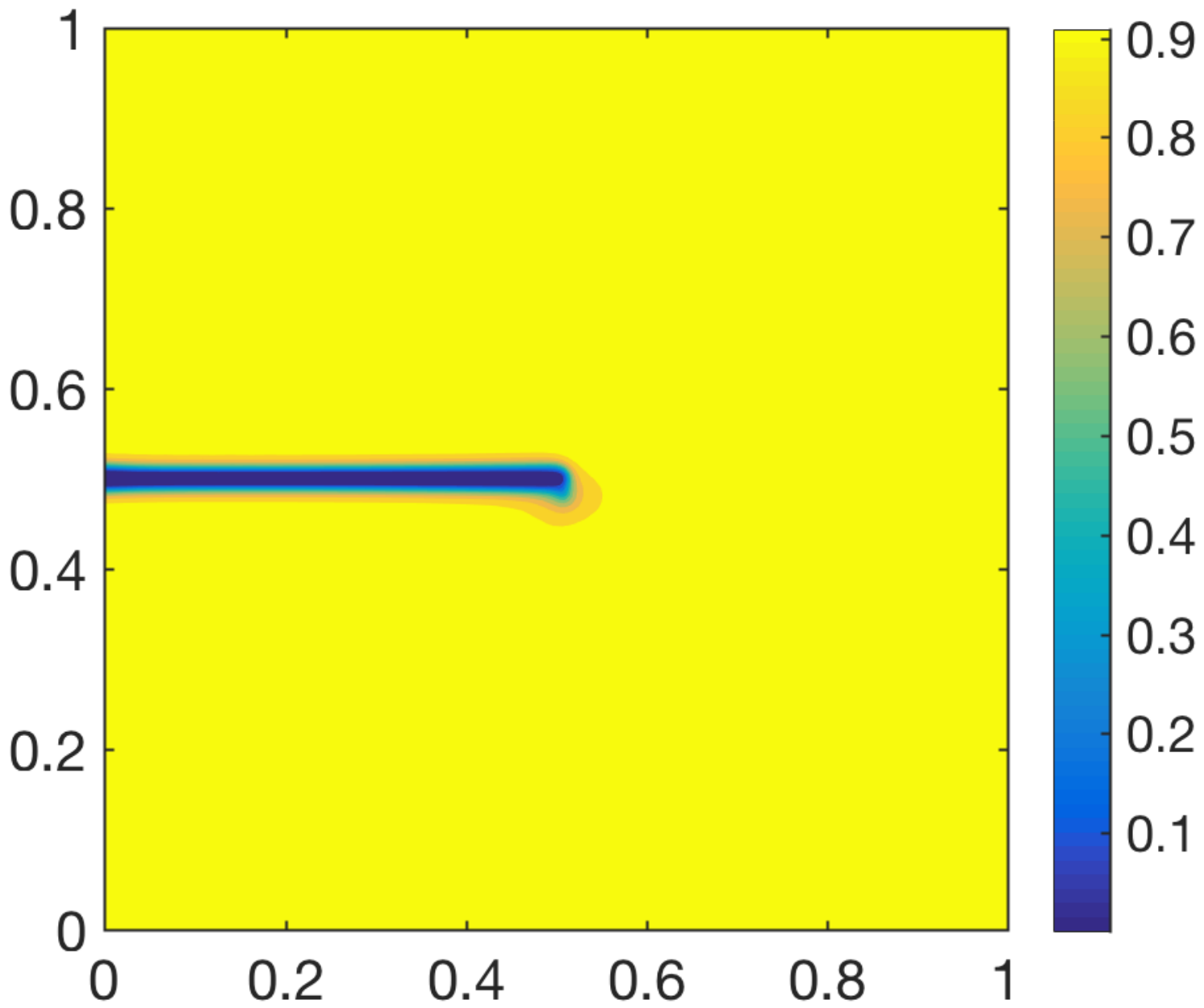}}
\subfigure[$U = 1.3 \times 10^{-2}$~mm]{\label{fig:subfig:MD02_U1300}
\includegraphics[width=0.25\linewidth]{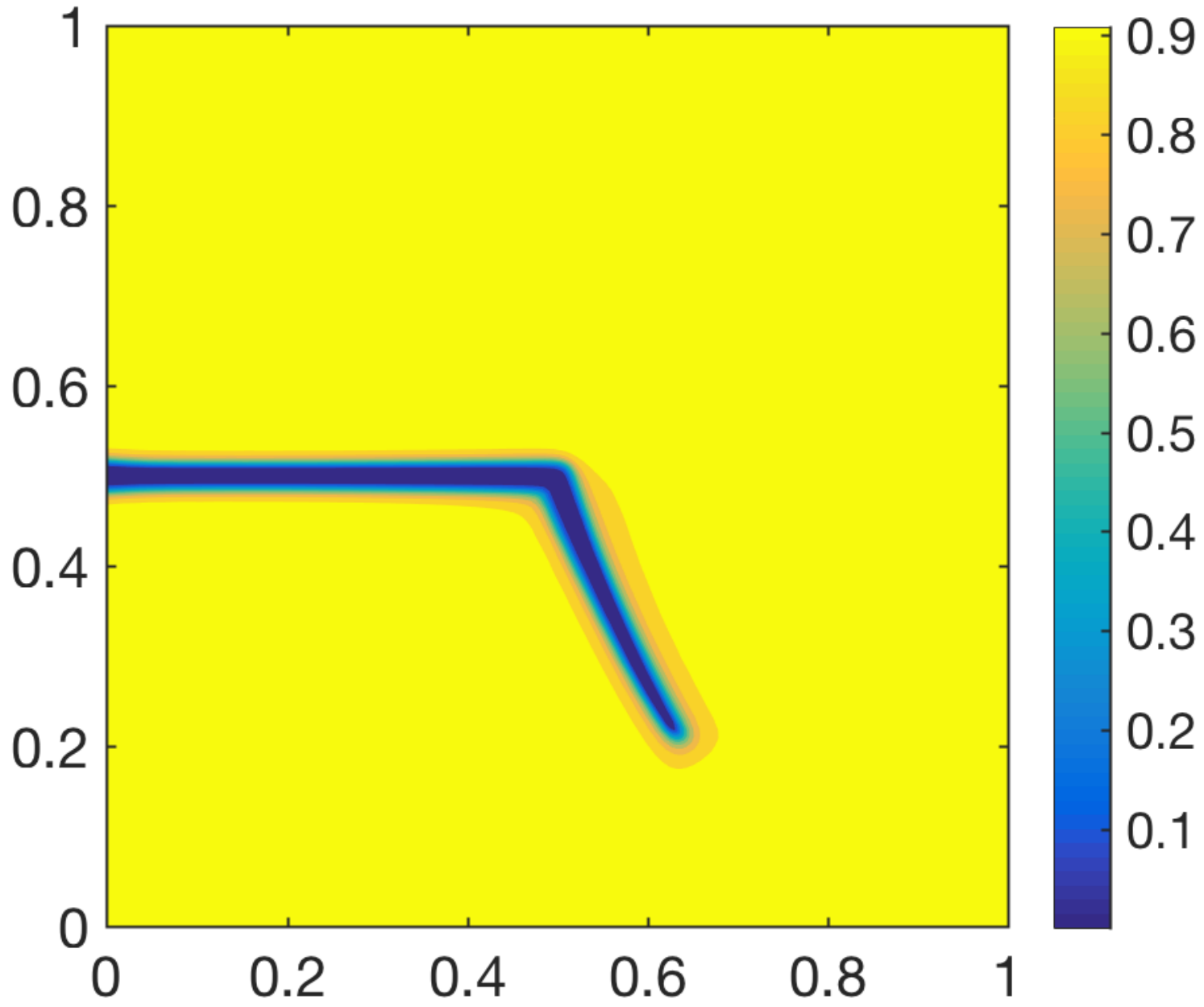}}
\subfigure[$U = 1.45 \times 10^{-2}$~mm]{\label{fig:subfig:MD02_U1450}
\includegraphics[width=0.25\linewidth]{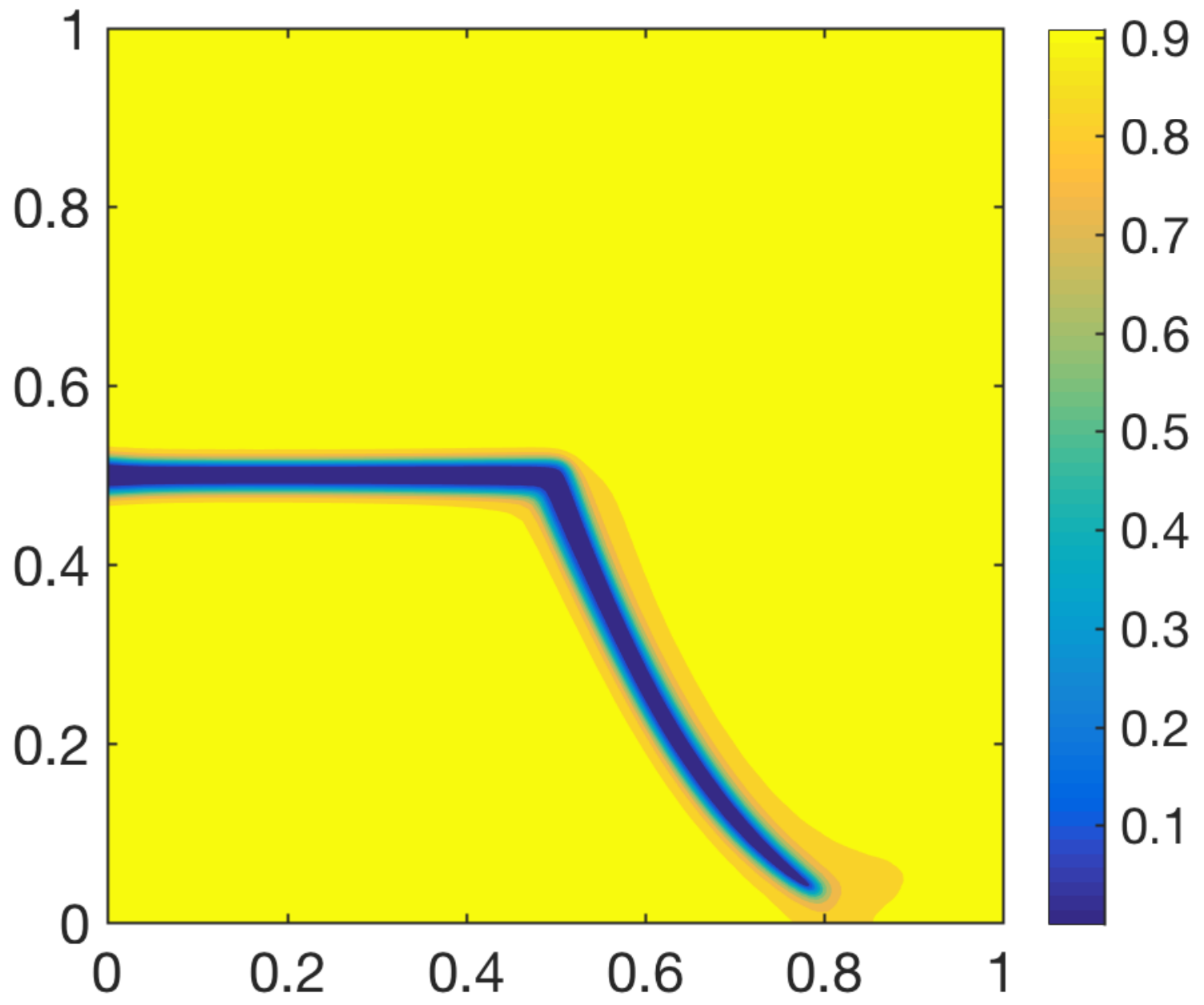}}
\vfill
\subfigure[$U = 1.1 \times 10^{-2}$~mm]{\label{fig:subfig:MM02_U1100}
\includegraphics[width=0.25\linewidth]{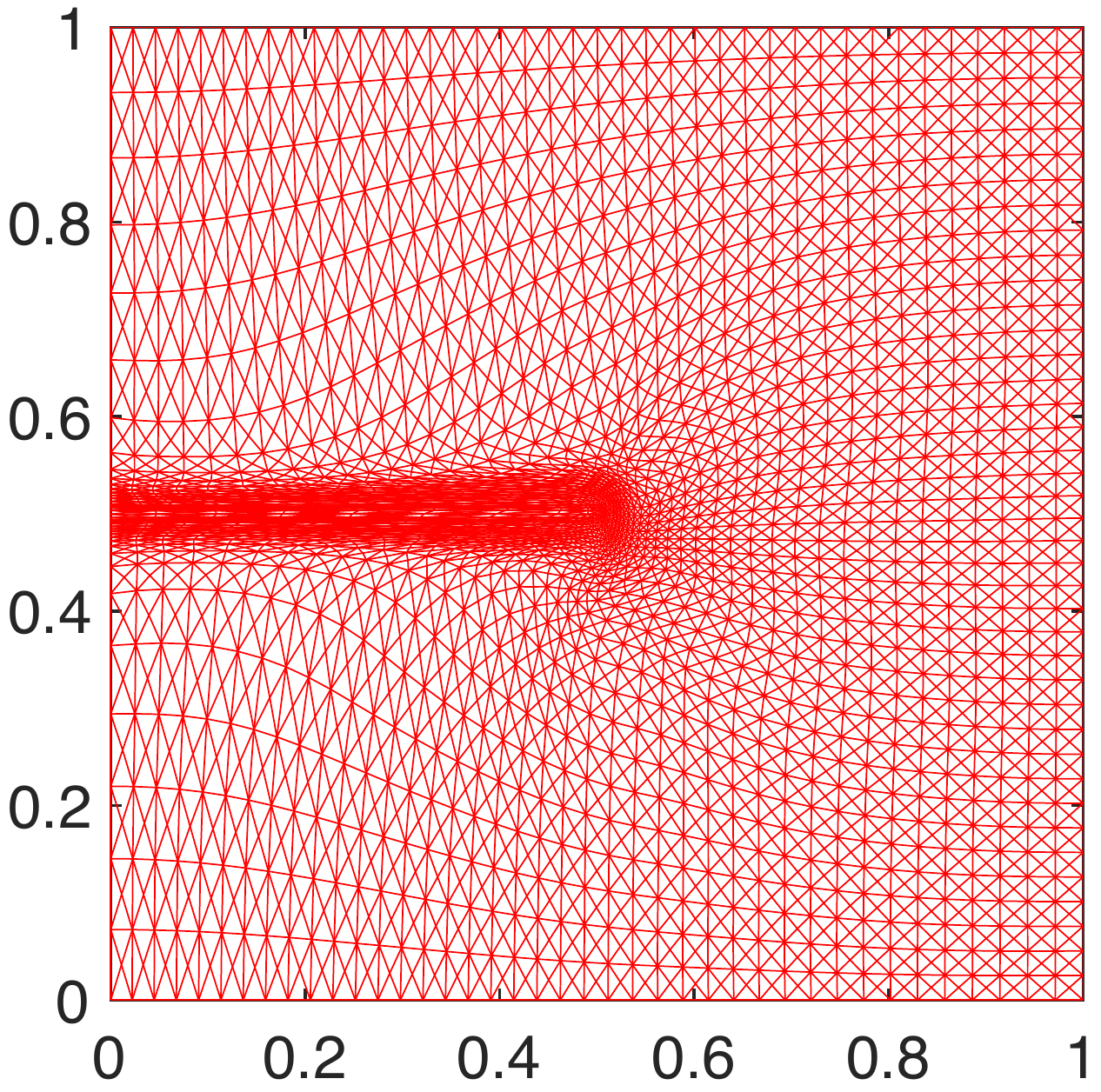}}
\subfigure[$U = 1.3 \times 10^{-2}$~mm]{\label{fig:subfig:MM02_U1300}
\includegraphics[width=0.25\linewidth]{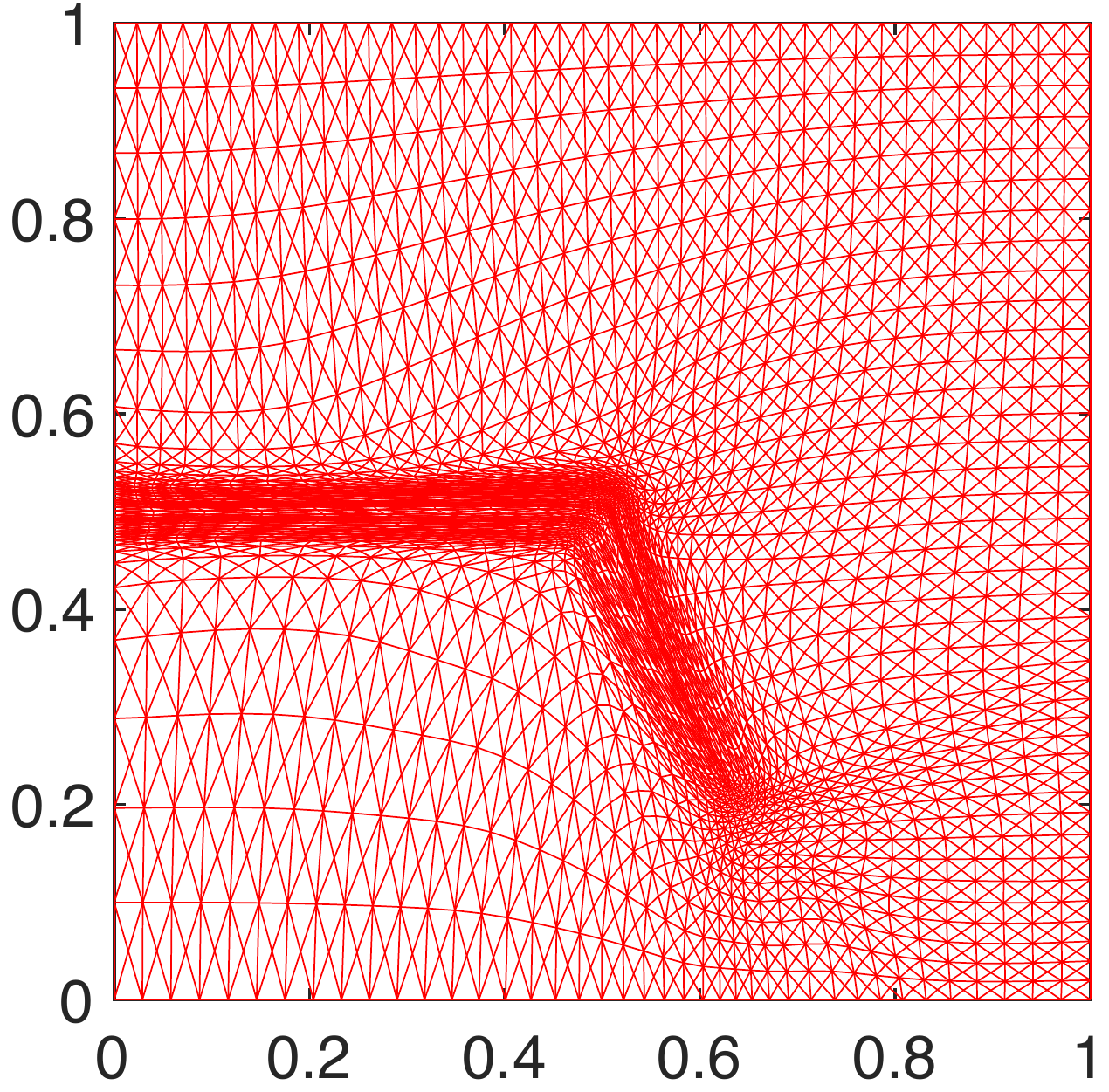}}
\subfigure[$U = 1.45 \times 10^{-2}$~mm]{\label{fig:subfig:MM02_U1450}
\includegraphics[width=0.25\linewidth]{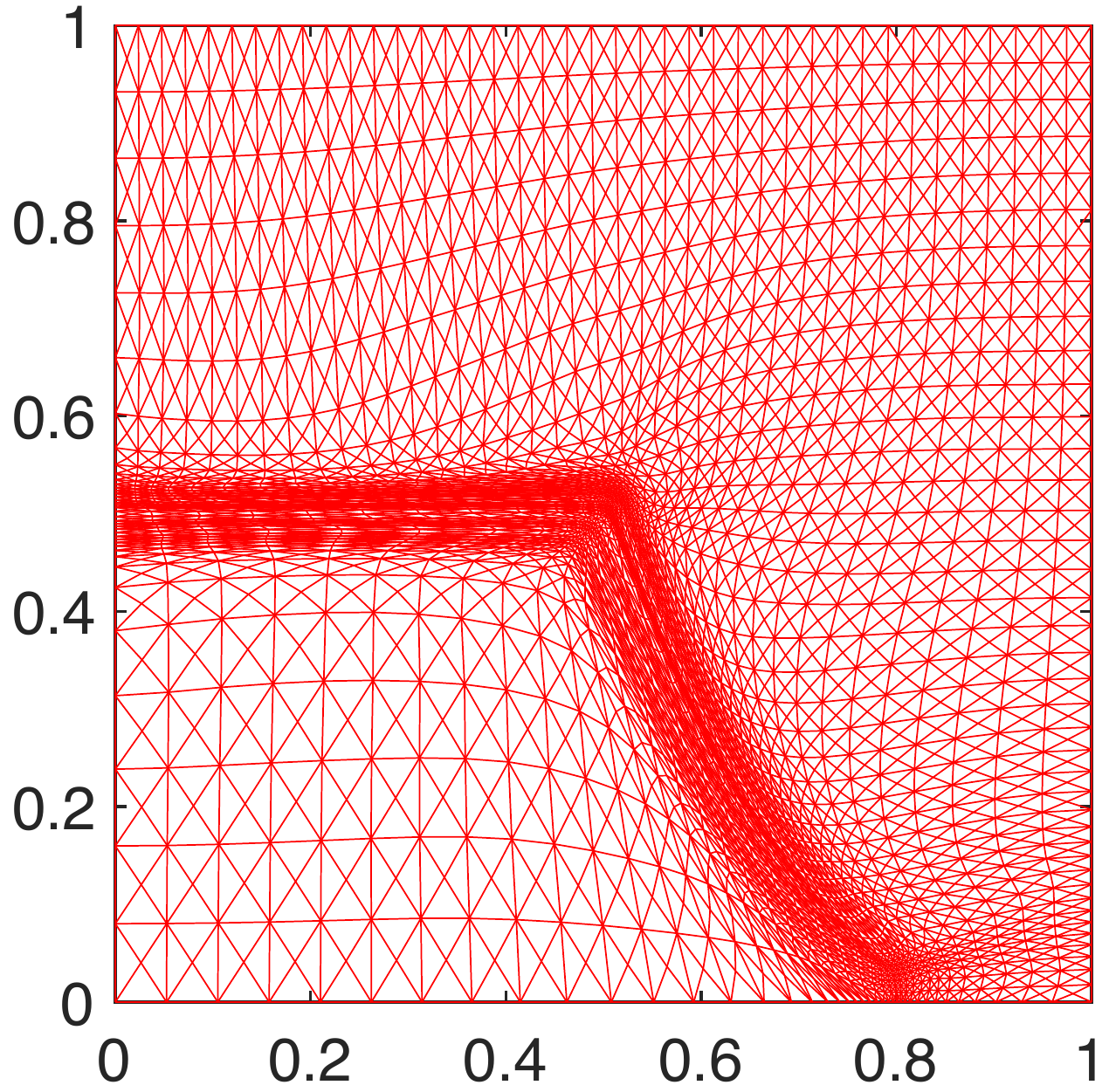}}
\vfill
\subfigure[$U = 1.1 \times 10^{-2}$~mm]{\label{fig:subfig:MS02_U1100}
\includegraphics[width=0.25\linewidth]{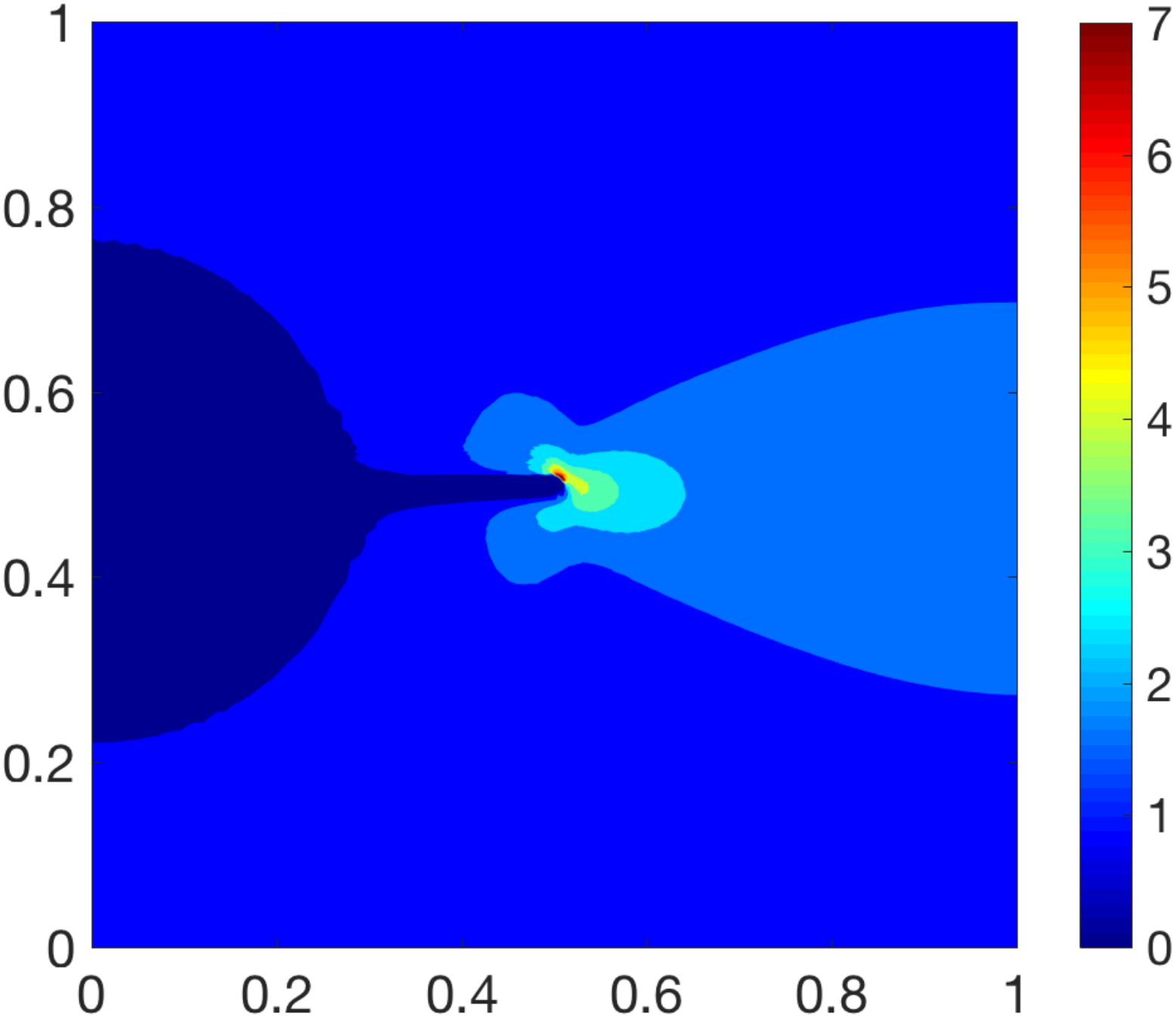}}
\subfigure[$U = 1.3 \times 10^{-2}$~mm]{\label{fig:subfig:MS02_U1300}
\includegraphics[width=0.25\linewidth]{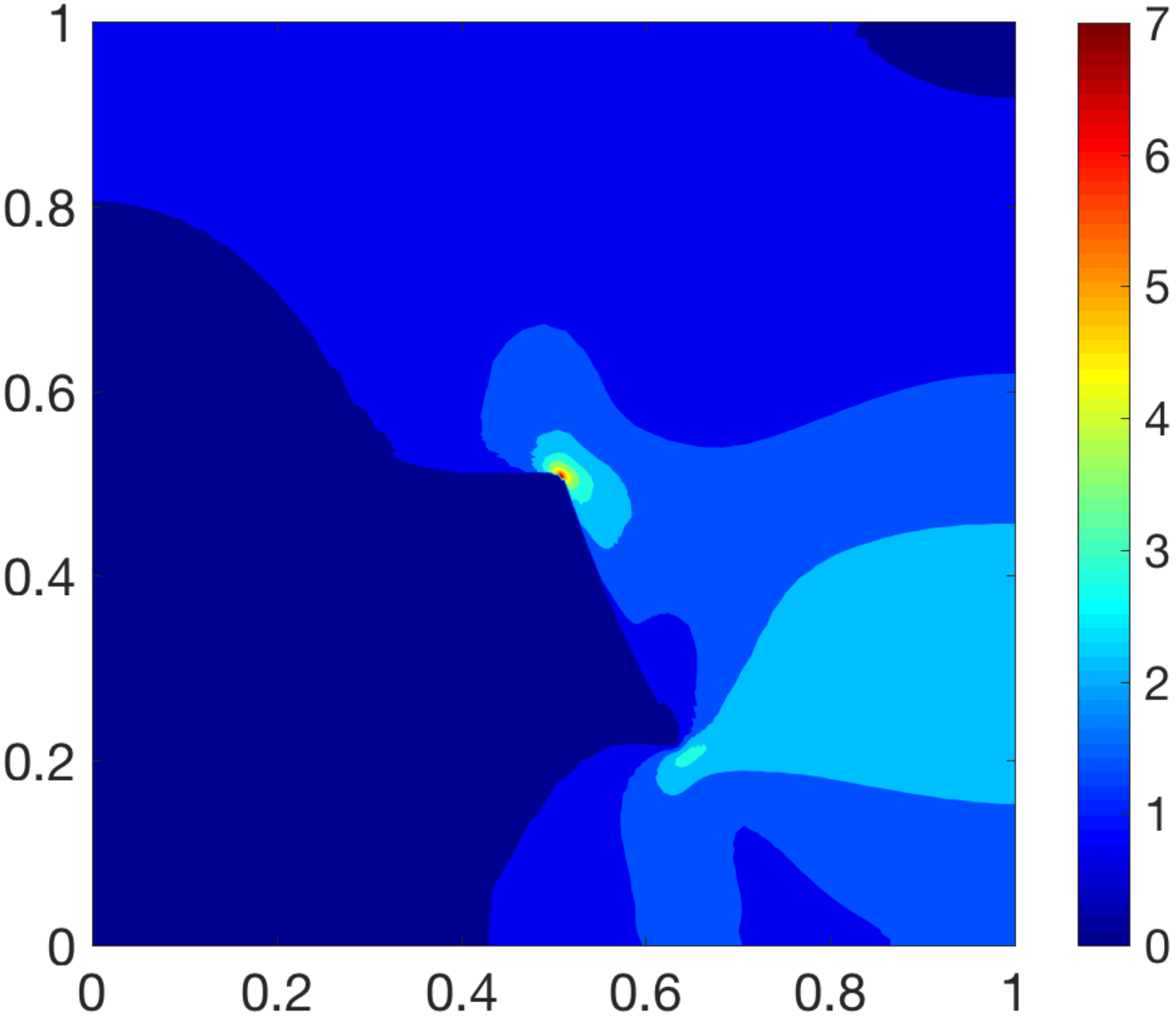}}
\subfigure[$U = 1.45 \times 10^{-2}$~mm]{\label{fig:subfig:MS02_U1450}
\includegraphics[width=0.25\linewidth]{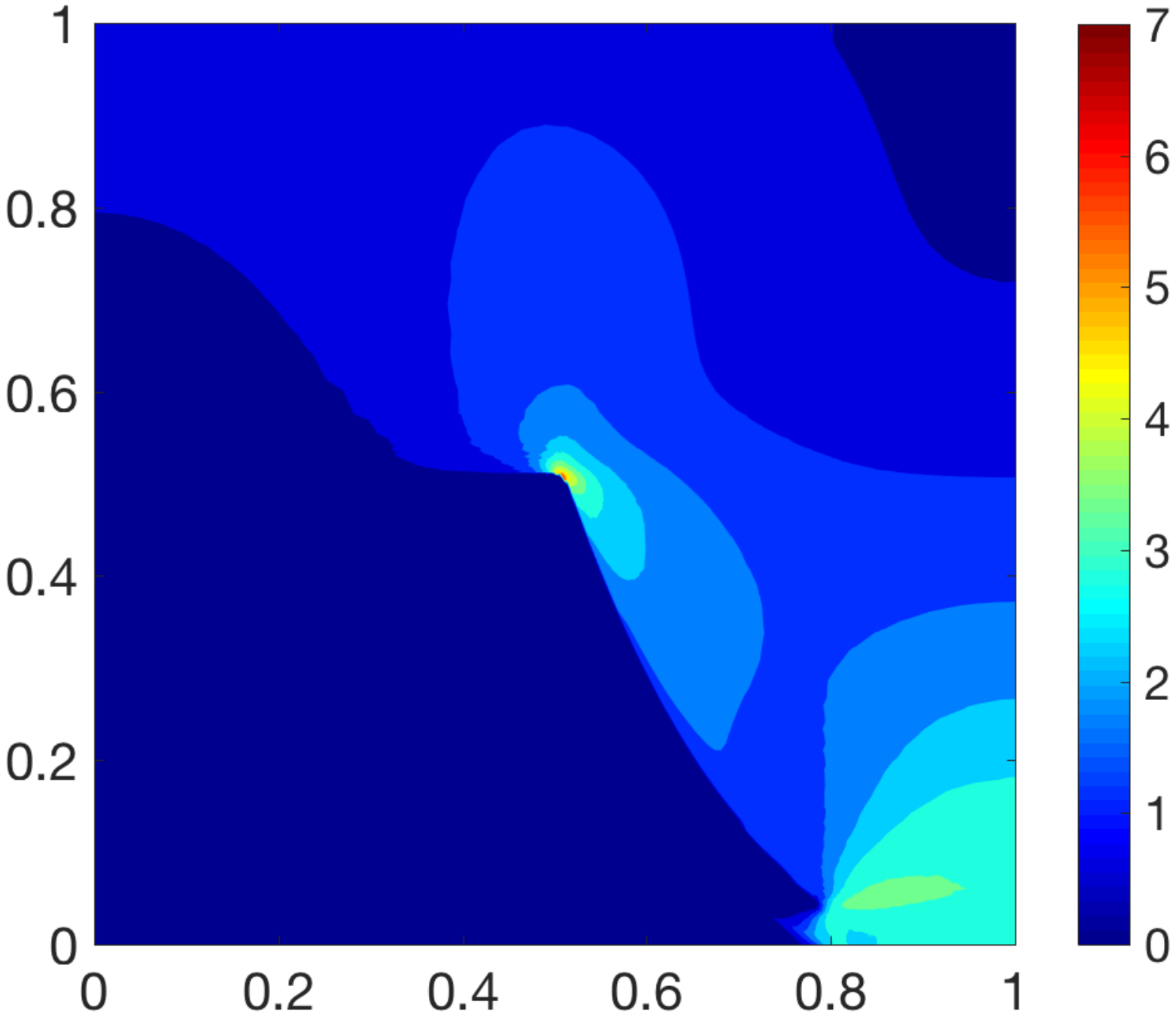}}
\caption{Example 2. Meshes and contours of the phase-field and von Mises stress distribution
are plotted at $U = 1.1\times 10^{-2}$, $1.30\times 10^{-2}$, and $1.45\times 10^{-2}$~mm.
The spectral decomposition model with ItCBC ($d_{cr} = 0.4$) is used.}
\label{fig:shear m's MBC}
\end{figure}

We first investigate the effects of the three decomposition models and ItCBC on the crack propagation
and the distribution of the stress. The spectral decomposition model without ItCBC
leads to the development of a secondary crack along the left edge of the domain and the concentration
of the stress in the damaged region; see Fig.~\ref{fig:shear m's OBC}.
This (unphysical) phenomenon has also been observed by May et al. \cite{MVB15}
where the initial notch is modeled with $d = 0$. 
The results with ItCBC modification ($d_{cr} = 0.4$) are shown in Fig.~\ref{fig:shear m's MBC}.
It can be seen that a secondary crack does not occur and the stress concentrates at the turning point and
crack tip only. The crack path agrees well with that in \cite{MHW10} where the initial crack is modeled
as a discrete crack in the geometry.

\begin{figure} 
\centering 
\subfigure[$U = 1.2 \times 10^{-2}$~mm]{\label{fig:subfig:AD_U1200}
\includegraphics[width=0.25\linewidth]{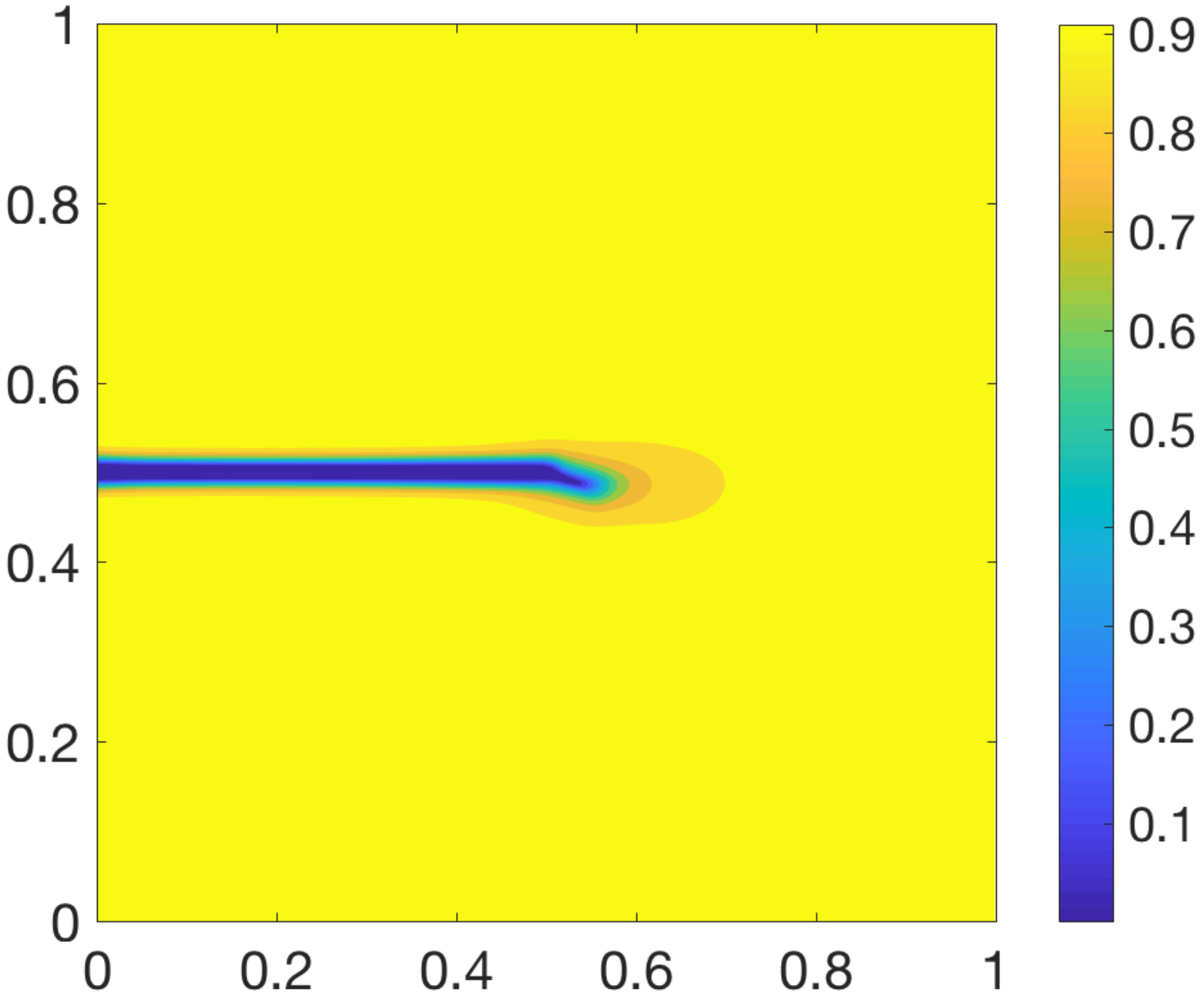}}
\subfigure[$U = 1.25 \times 10^{-2}$~mm]{\label{fig:subfig:AD_U1250}
\includegraphics[width=0.25\linewidth]{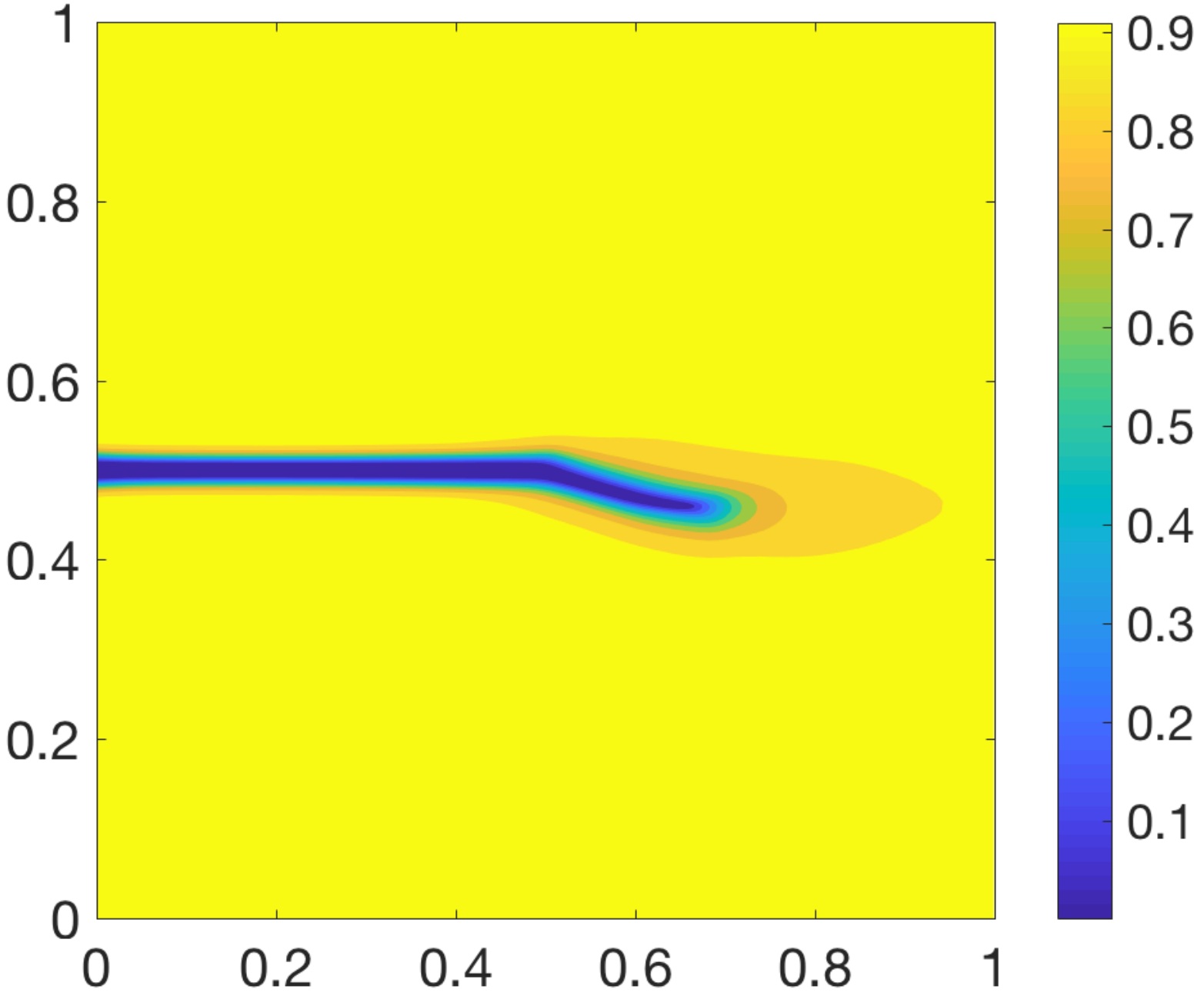}}
\subfigure[$U = 1.3 \times 10^{-2}$~mm]{\label{fig:subfig:AD_U1300}
\includegraphics[width=0.25\linewidth]{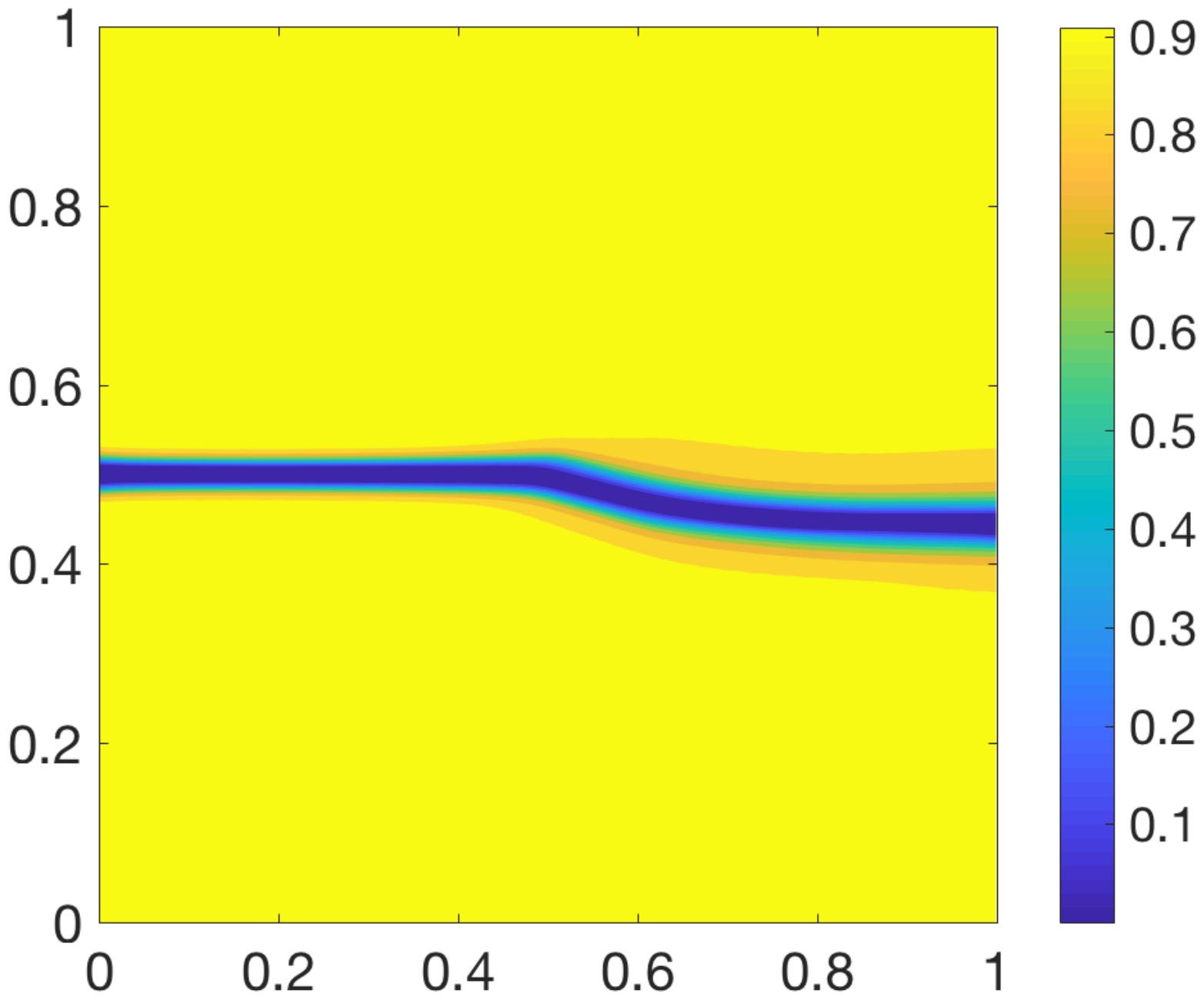}}
\vfill
\subfigure[$U = 1.2 \times 10^{-2}$~mm]{\label{fig:subfig:AM_U1200}
\includegraphics[width=0.25\linewidth]{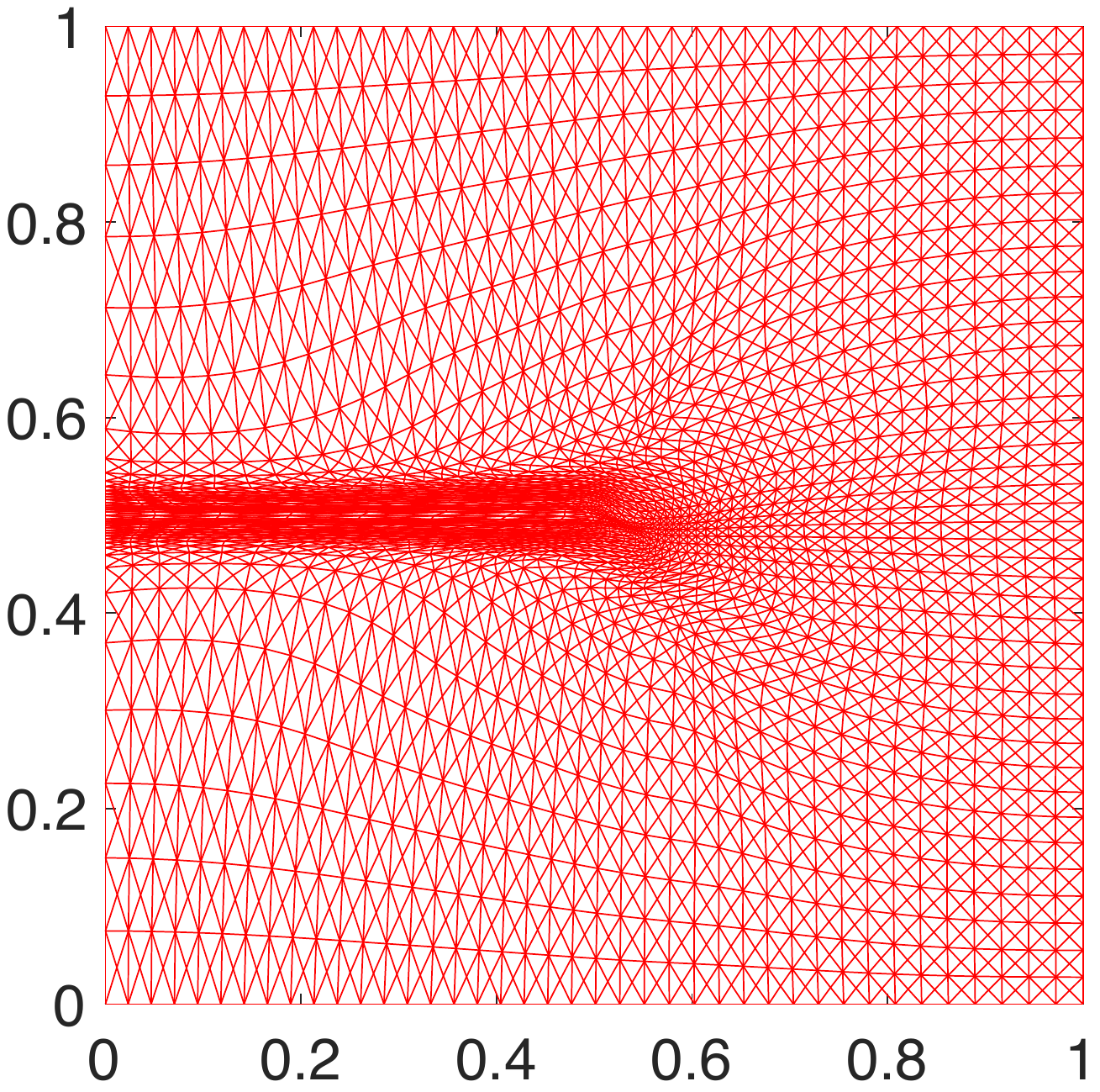}}
\subfigure[$U = 1.25 \times 10^{-2}$~mm]{\label{fig:subfig:AM_U1250}
\includegraphics[width=0.25\linewidth]{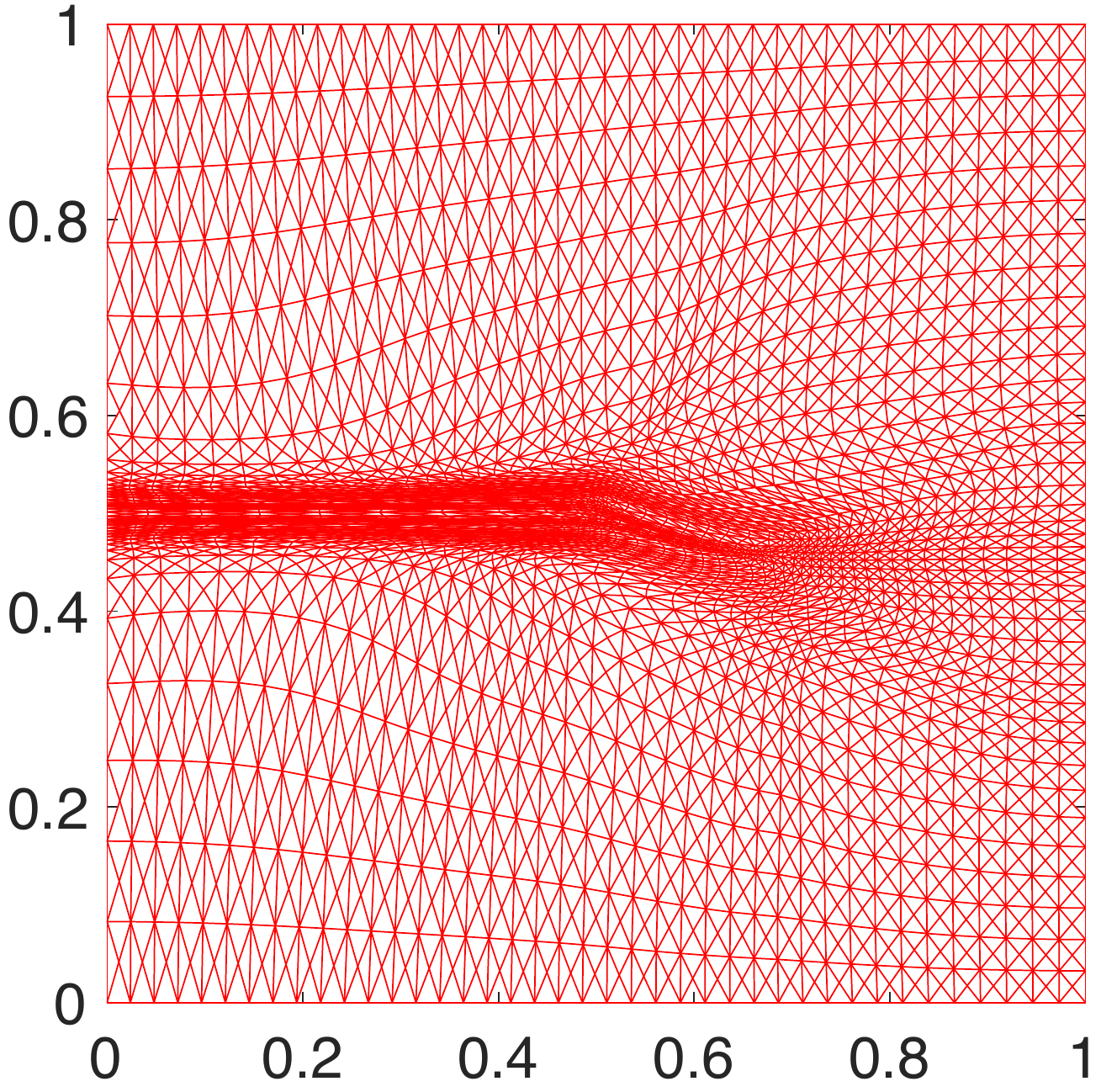}}
\subfigure[$U = 1.3 \times 10^{-2}$~mm]{\label{fig:subfig:AM_U1300}
\includegraphics[width=0.25\linewidth]{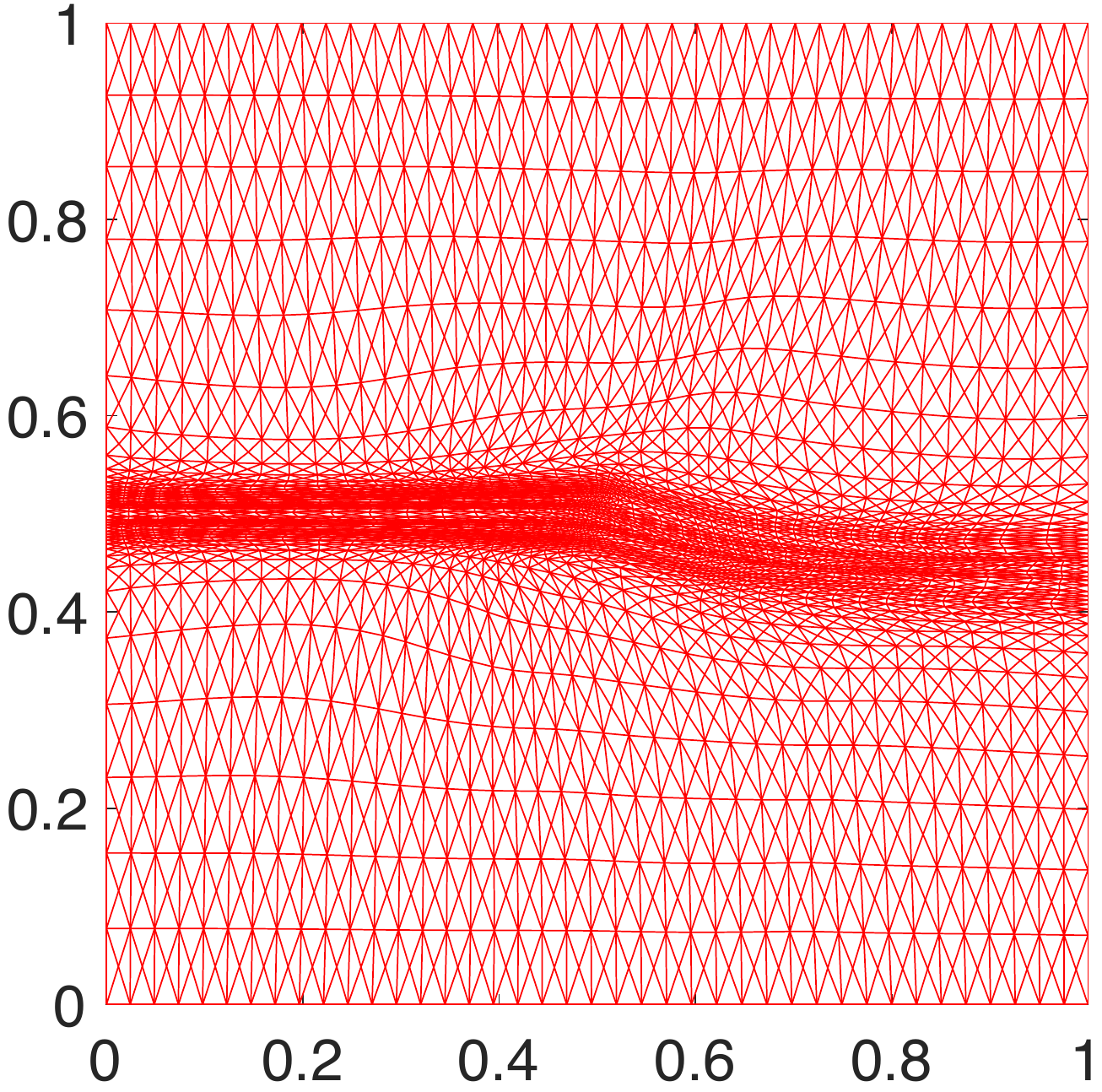}}
\vfill
\subfigure[$U = 1.2 \times 10^{-2}$~mm]{\label{fig:subfig:AS_U1200}
\includegraphics[width=0.25\linewidth]{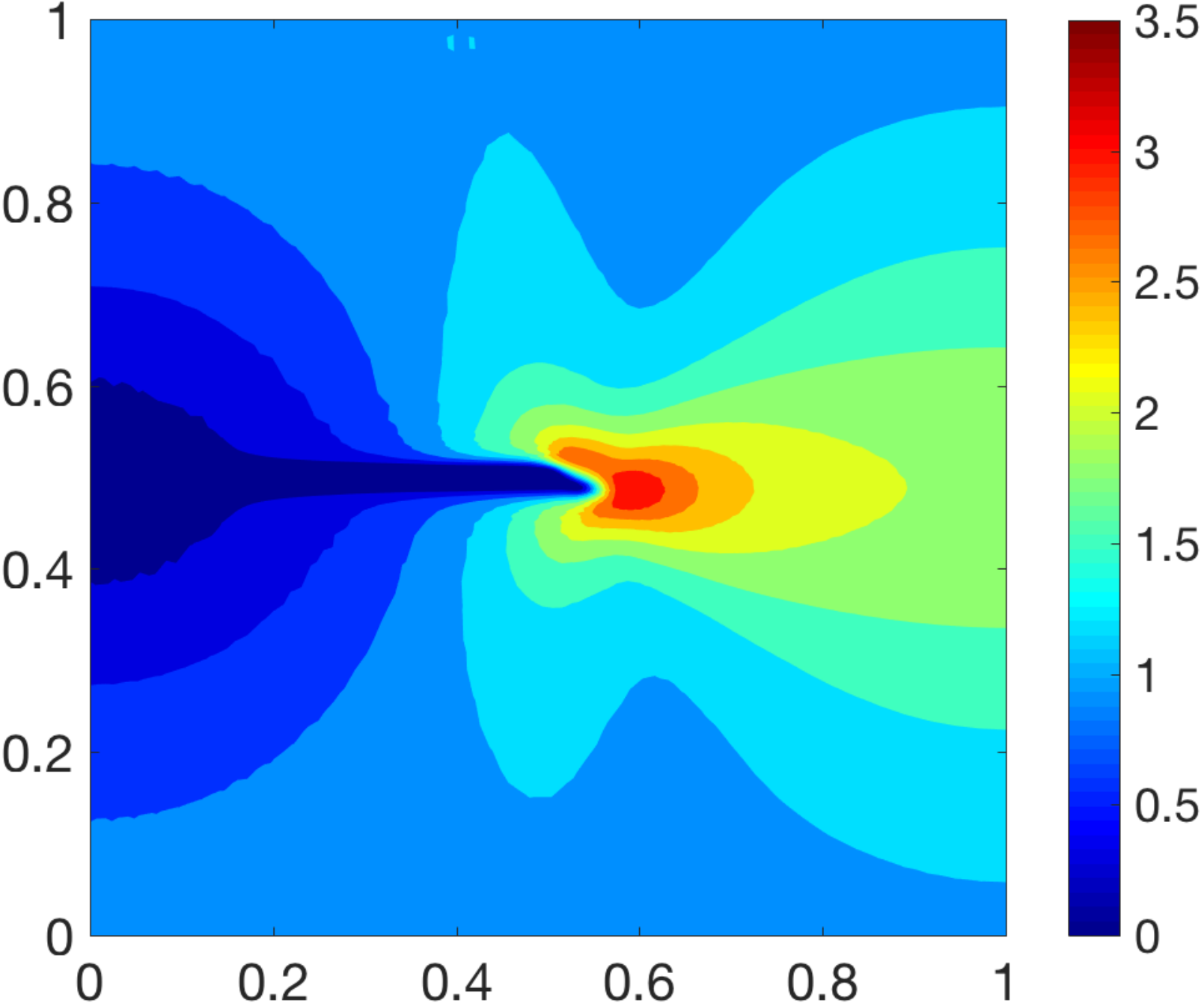}}
\subfigure[$U = 1.25 \times 10^{-2}$~mm]{\label{fig:subfig:AS_U1250}
\includegraphics[width=0.25\linewidth]{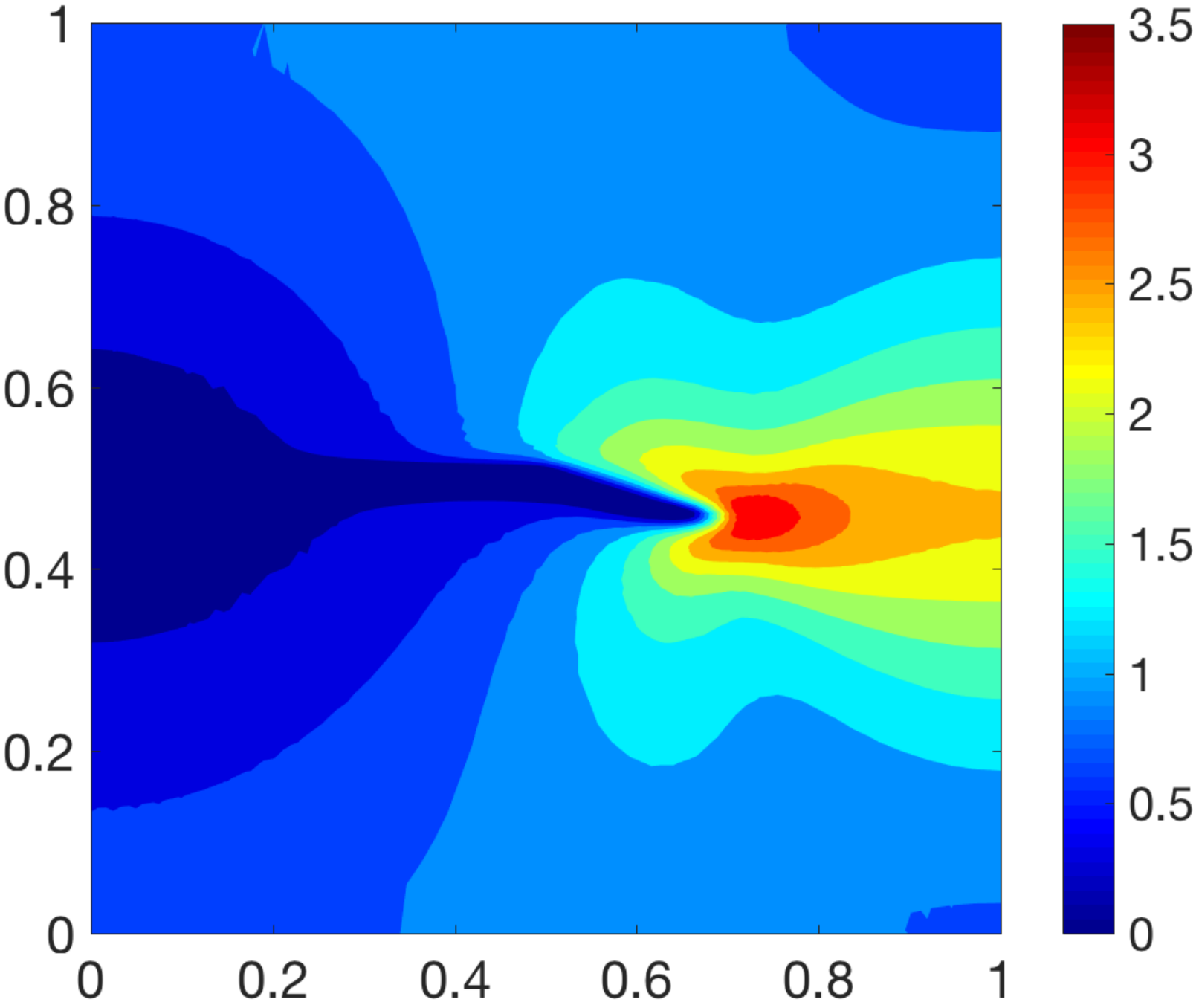}}
\subfigure[$U = 1.3 \times 10^{-2}$~mm]{\label{fig:subfig:AS_U1300}
\includegraphics[width=0.25\linewidth]{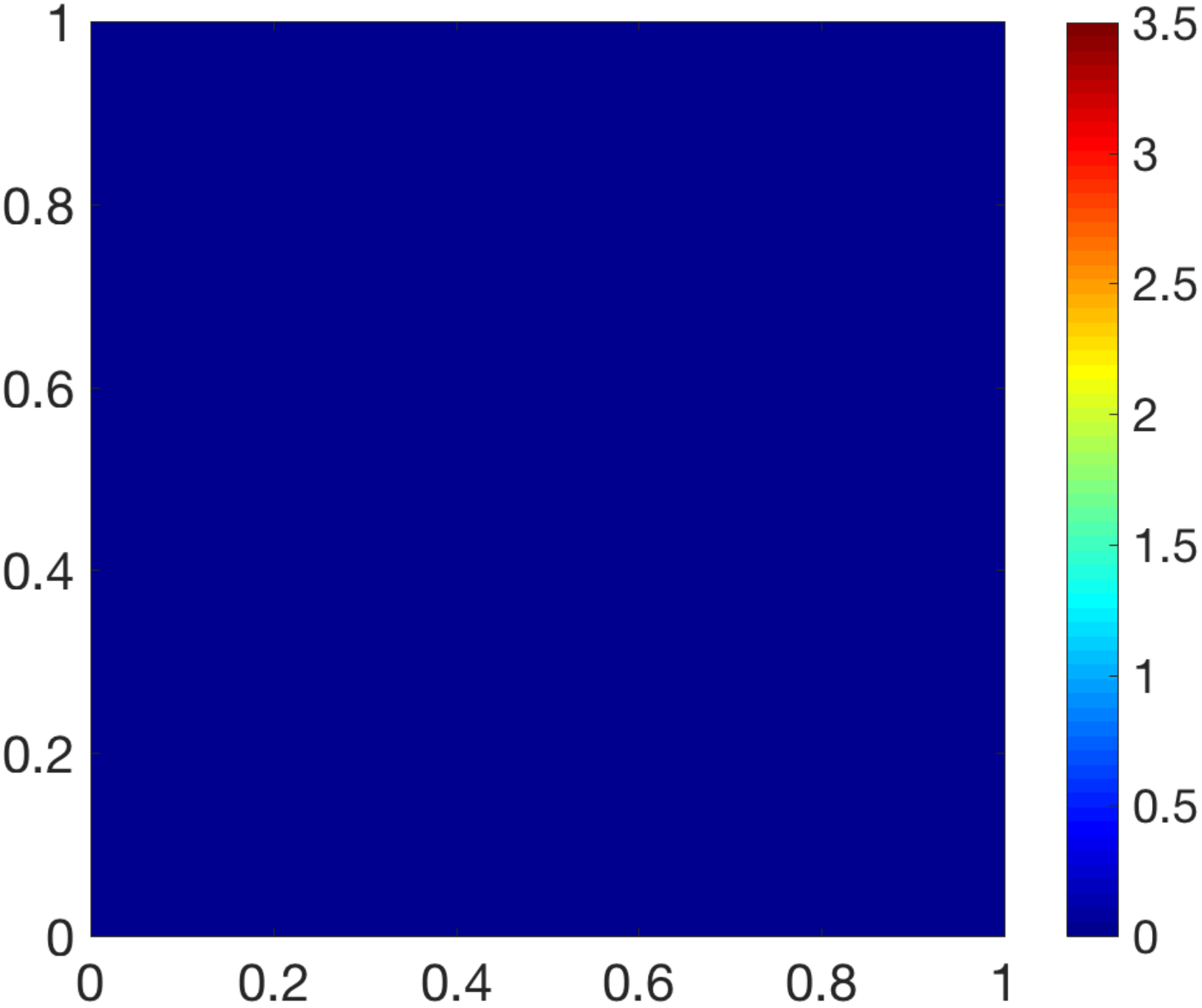}}
\caption{Example 2. Meshes and contours of the phase-field and von Mises stress distribution
are plotted at $U = 1.2\times 10^{-2}$, $1.25\times 10^{-2}$, and $1.30\times 10^{-2}$~mm.
The v-d split model is used.}
\label{fig:shear a's OBC}
\end{figure}

\begin{figure} 
\centering 
\subfigure[$U = 1.2 \times 10^{-2}$~mm]{\label{fig:subfig:AD02_U1200}
\includegraphics[width=0.25\linewidth]{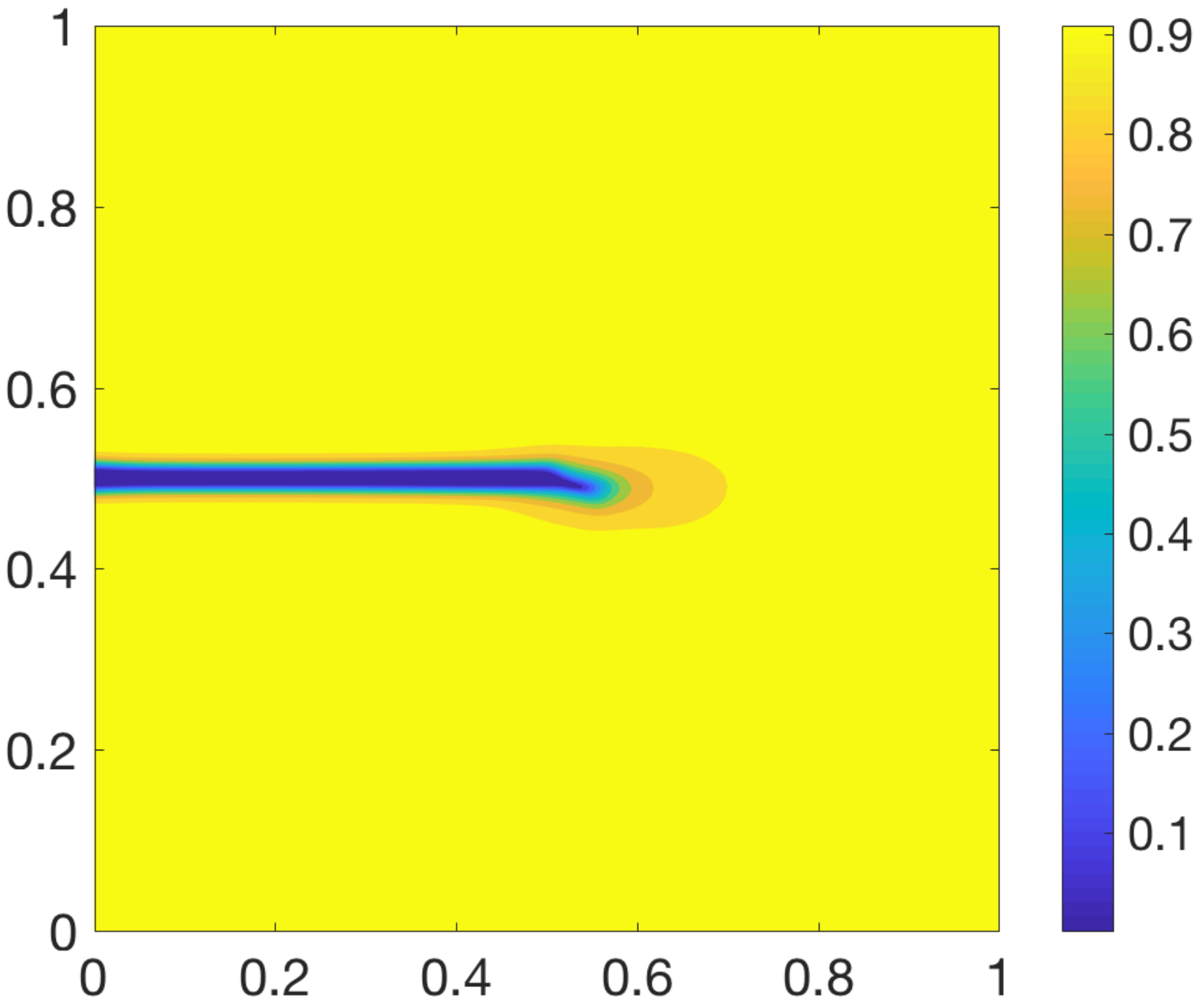}}
\subfigure[$U = 1.25 \times 10^{-2}$~mm]{\label{fig:subfig:AD02_U1250}
\includegraphics[width=0.25\linewidth]{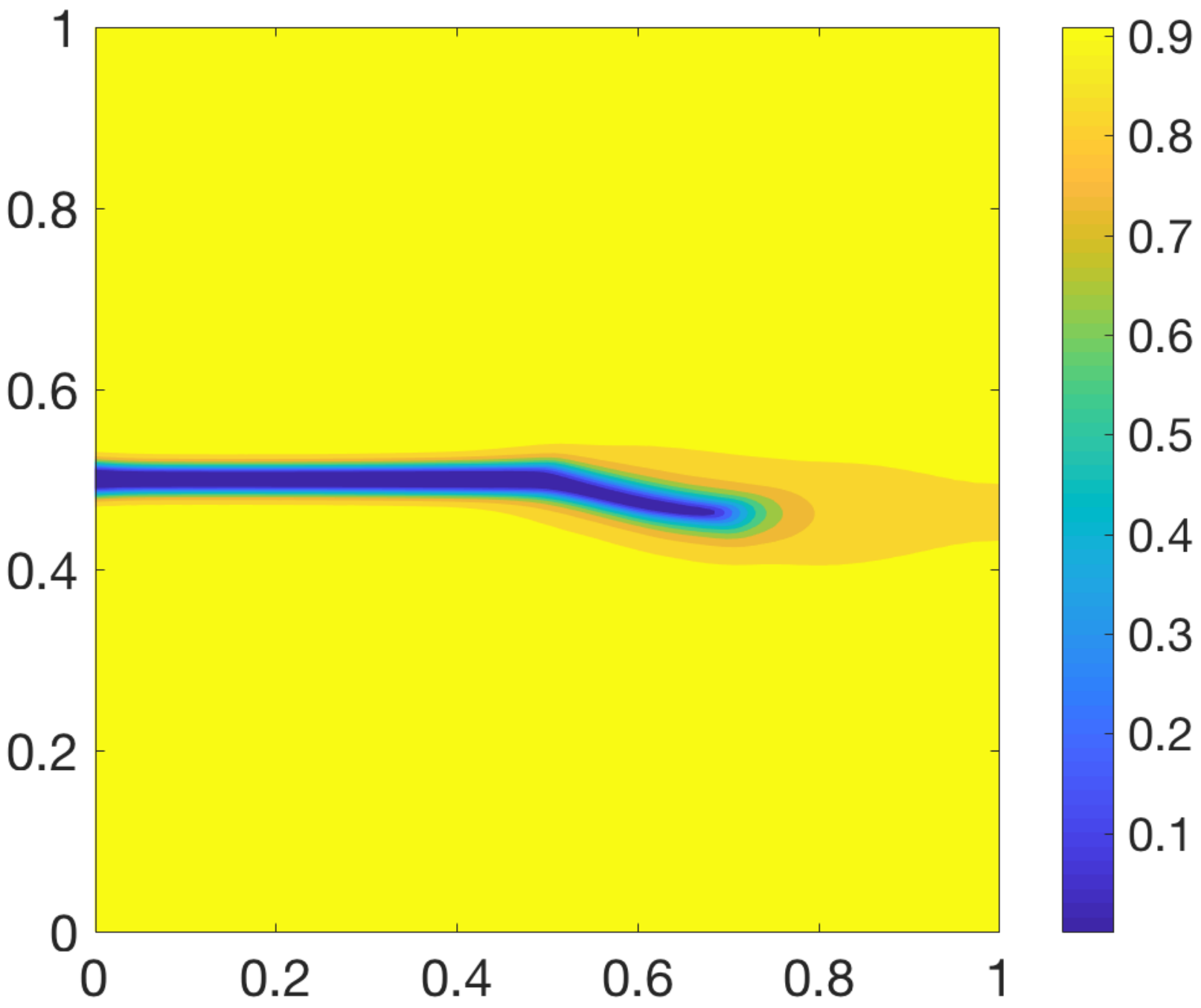}}
\subfigure[$U = 1.3 \times 10^{-2}$~mm]{\label{fig:subfig:AD02_U1300}
\includegraphics[width=0.25\linewidth]{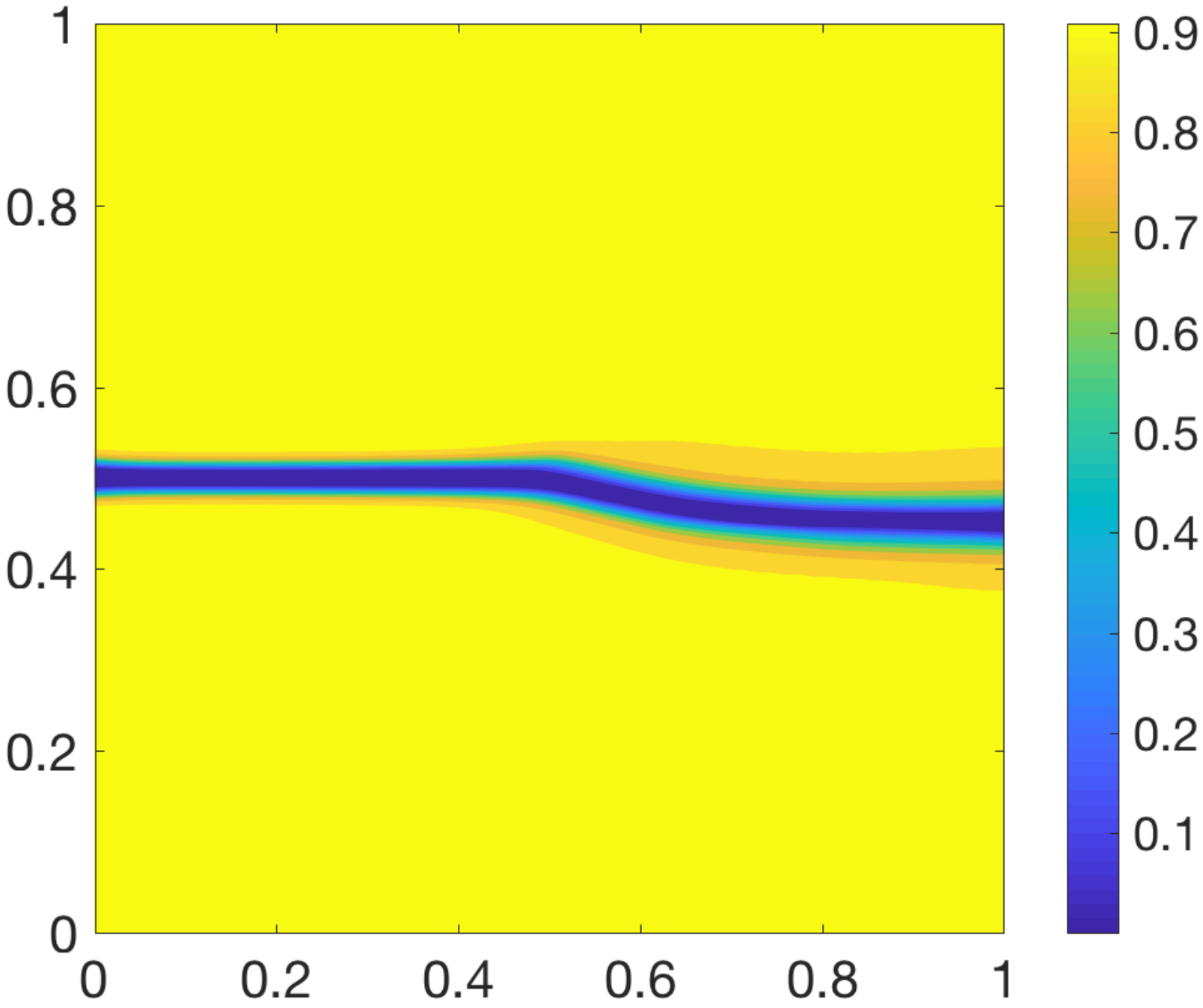}}
\vfill
\subfigure[$U = 1.2 \times 10^{-2}$~mm]{\label{fig:subfig:AM02_U1200}
\includegraphics[width=0.25\linewidth]{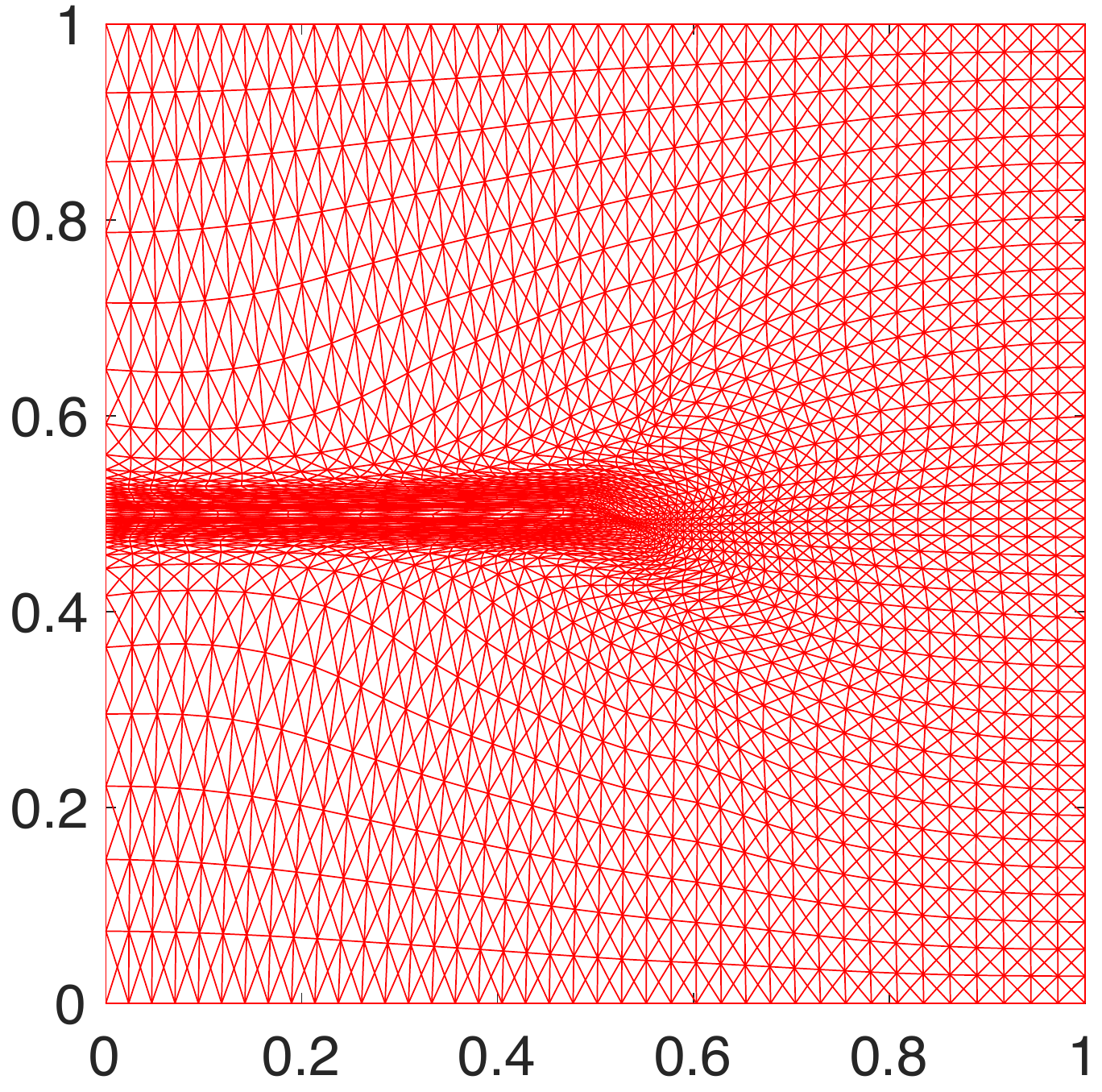}}
\subfigure[$U = 1.25 \times 10^{-2}$~mm]{\label{fig:subfig:AM02_U1250}
\includegraphics[width=0.25\linewidth]{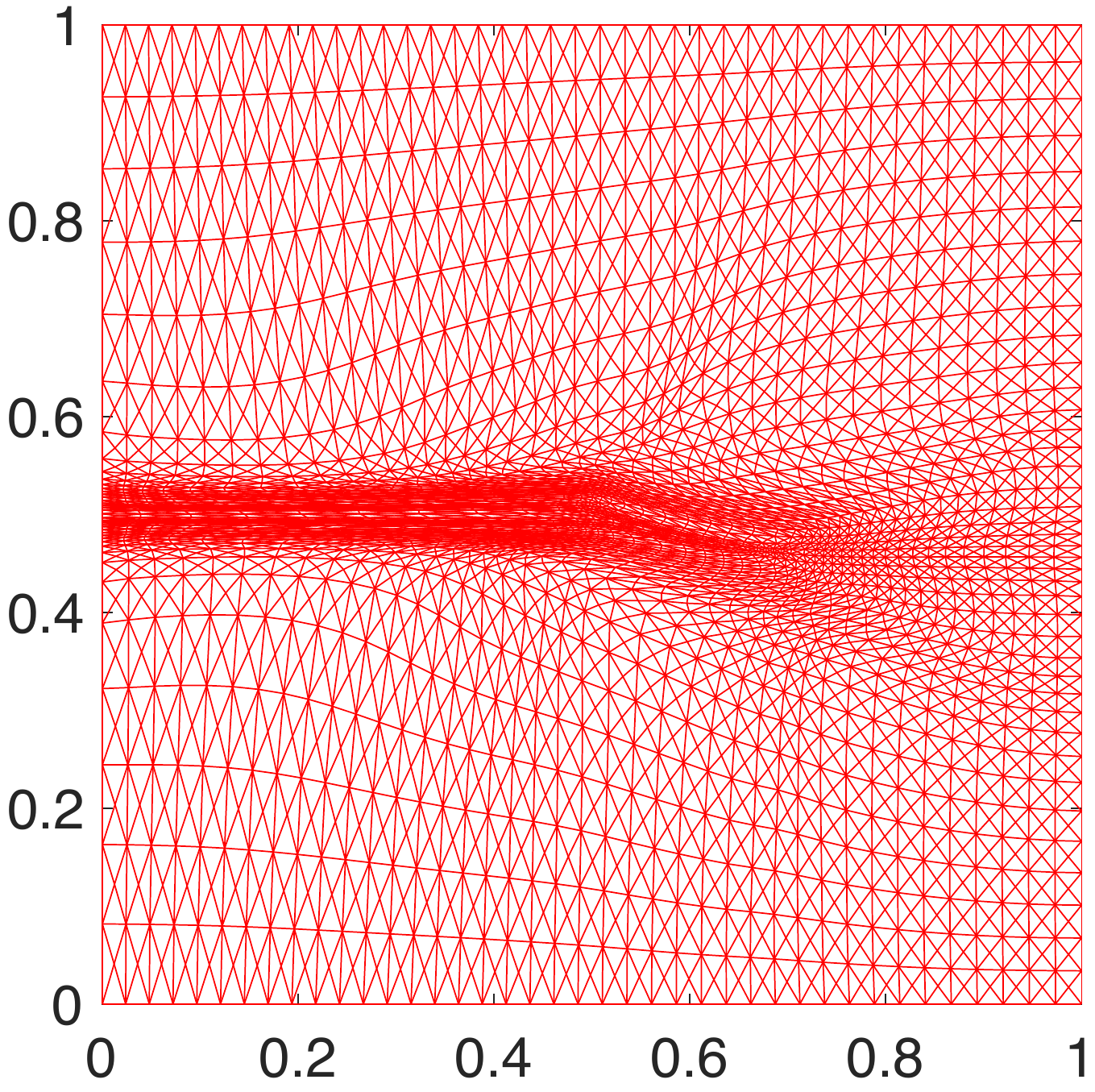}}
\subfigure[$U = 1.3 \times 10^{-2}$~mm]{\label{fig:subfig:AM02_U1300}
\includegraphics[width=0.25\linewidth]{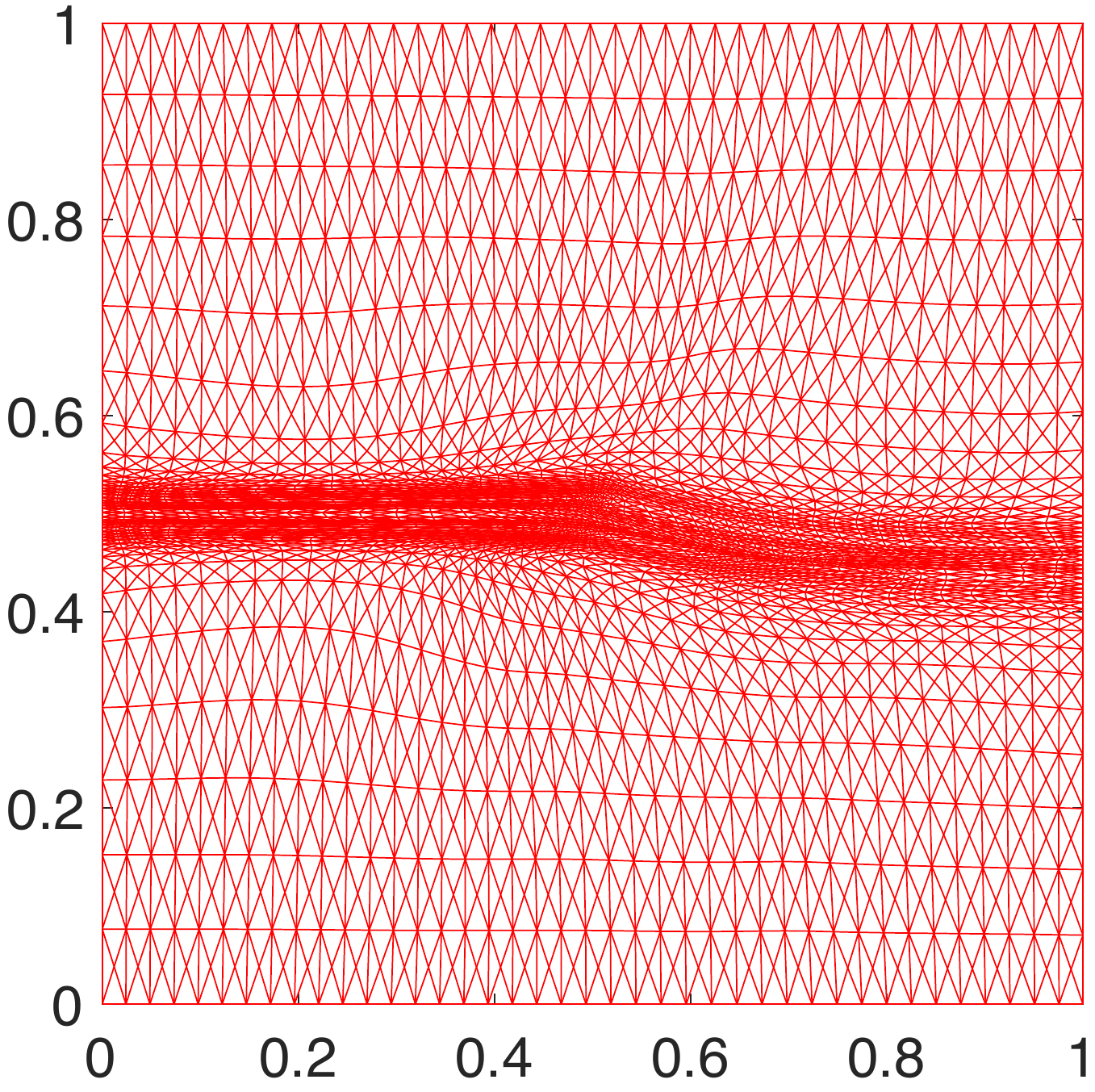}}
\vfill
\subfigure[$U = 1.2 \times 10^{-2}$~mm]{\label{fig:subfig:AS02_U1200}
\includegraphics[width=0.25\linewidth]{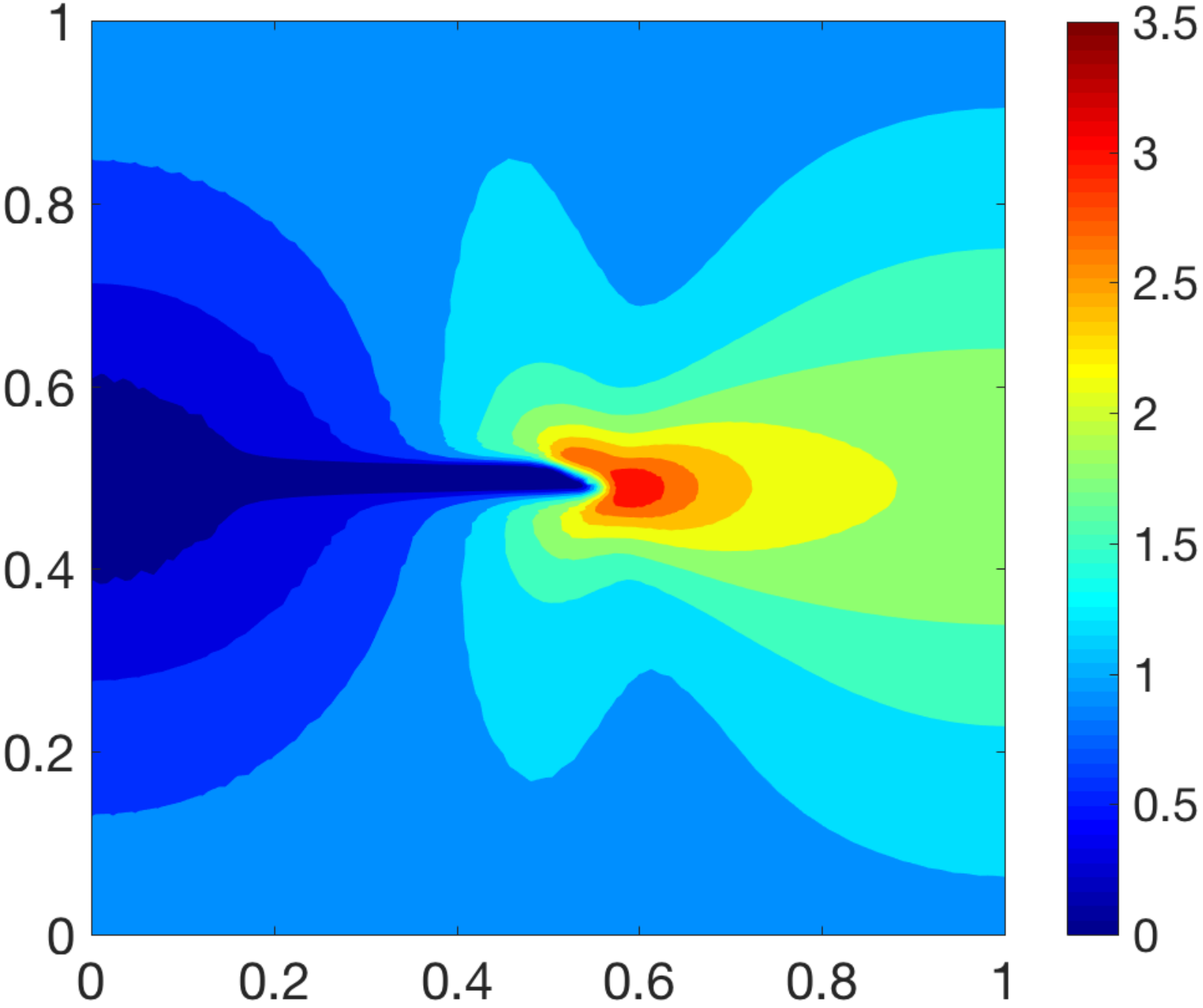}}
\subfigure[$U = 1.25 \times 10^{-2}$~mm]{\label{fig:subfig:AS02_U1250}
\includegraphics[width=0.25\linewidth]{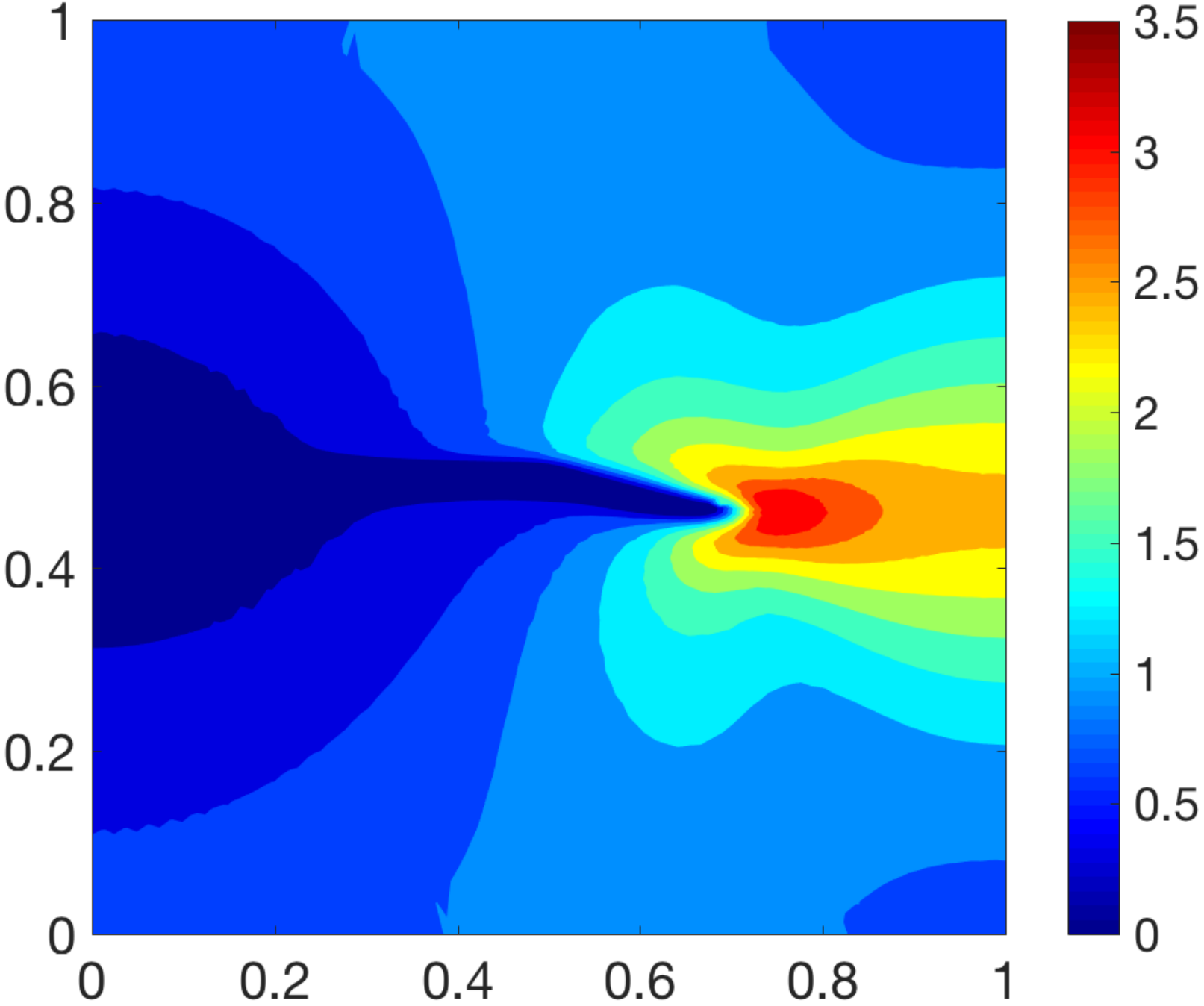}}
\subfigure[$U = 1.3 \times 10^{-2}$~mm]{\label{fig:subfig:AS02_U1300}
\includegraphics[width=0.25\linewidth]{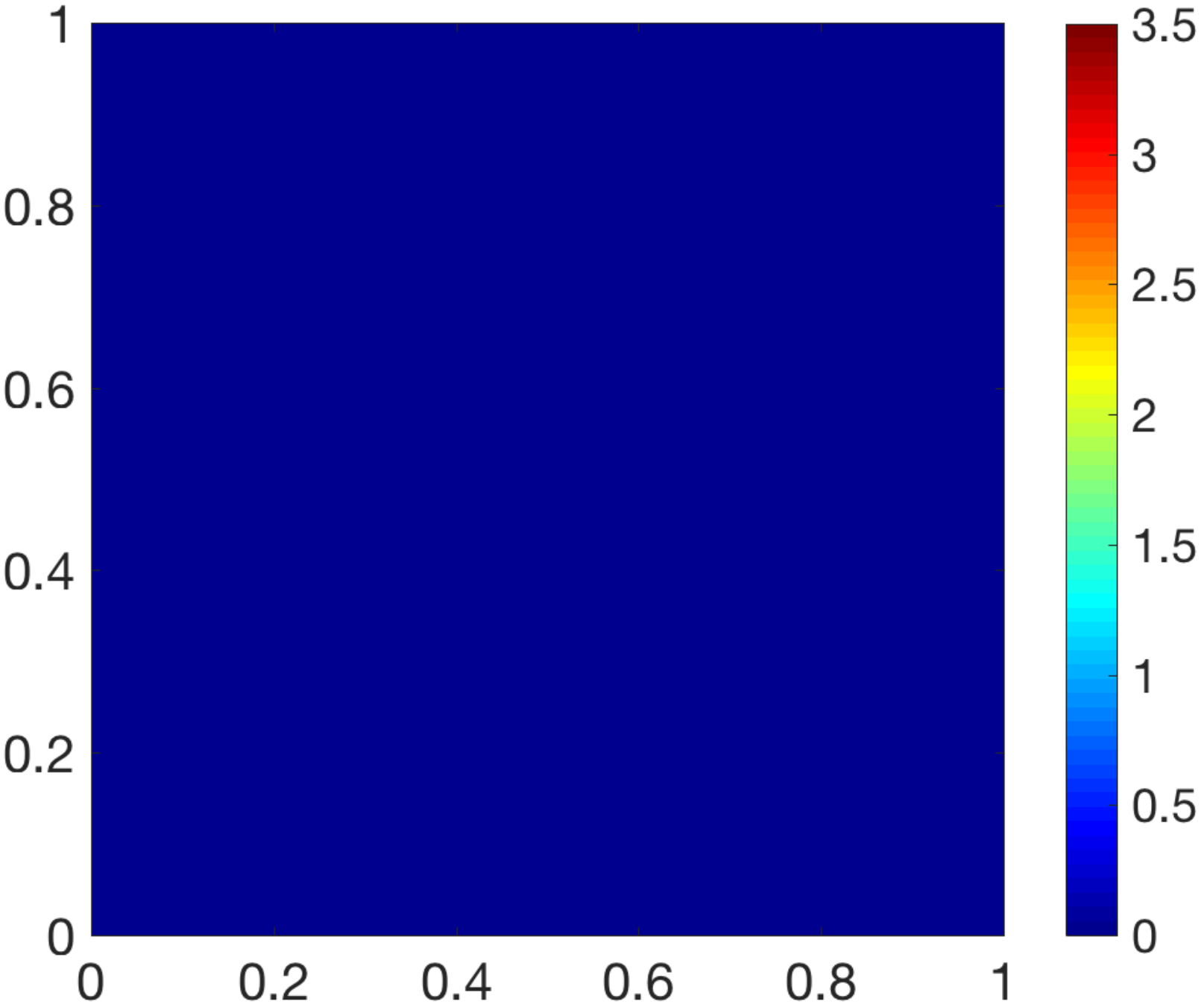}}
\caption{Example 2. Meshes and contours of the phase-field and von Mises stress distribution
are plotted at $U = 1.20\times 10^{-2}$, $1.25\times 10^{-2}$, and $1.30\times 10^{-2}$~mm.
The v-d split model with ItCBC ($d_{cr} = 0.4$) is used.}
\label{fig:shear a's MBC}
\end{figure}

\begin{figure} 
\centering 
\subfigure[$U = 3.0 \times 10^{-2}$~mm]{\label{fig:subfig:ID_U3000}
\includegraphics[width=0.25\linewidth]{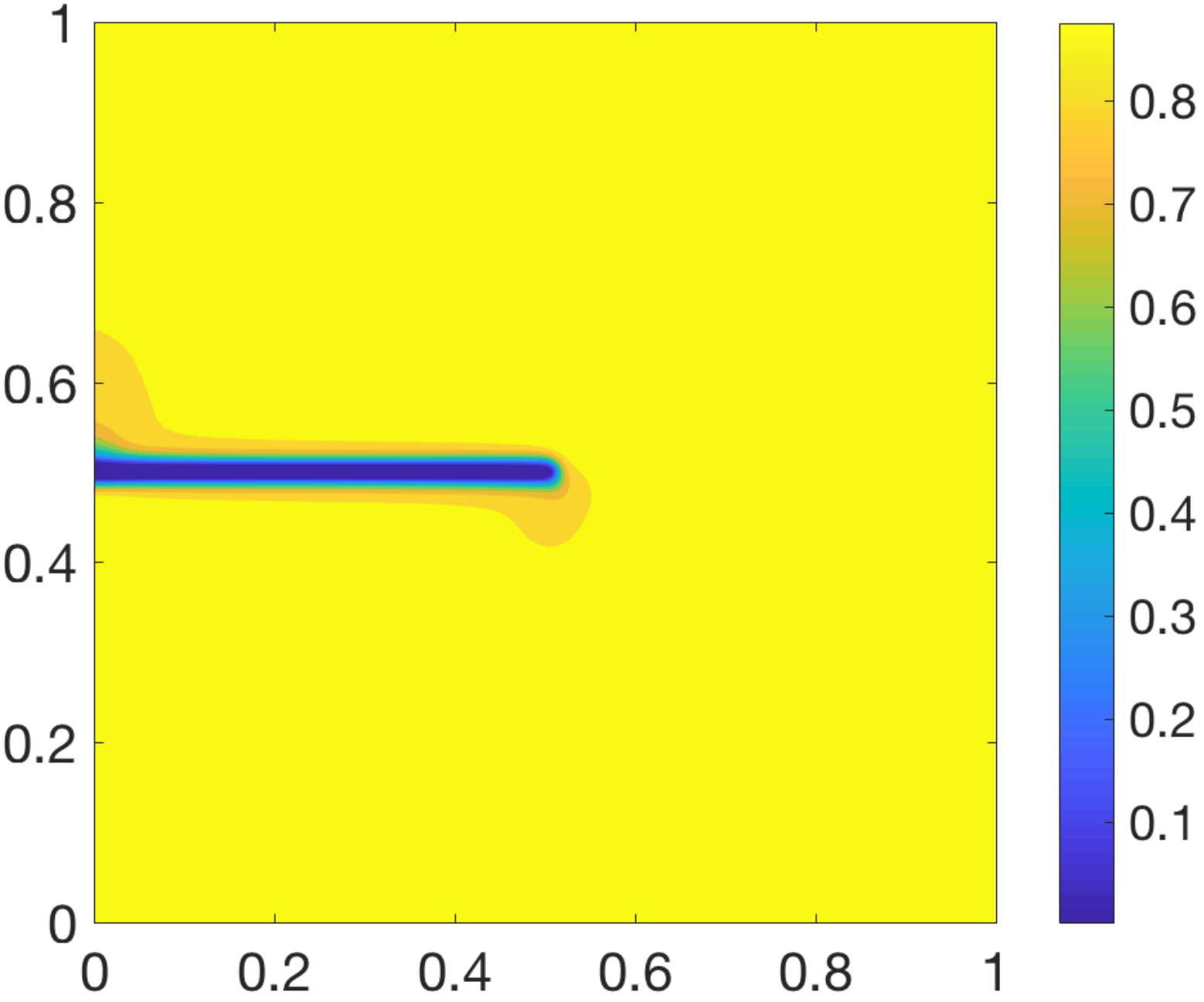}}
\subfigure[$U = 3.2 \times 10^{-2}$~mm]{\label{fig:subfig:ID_U3200}
\includegraphics[width=0.25\linewidth]{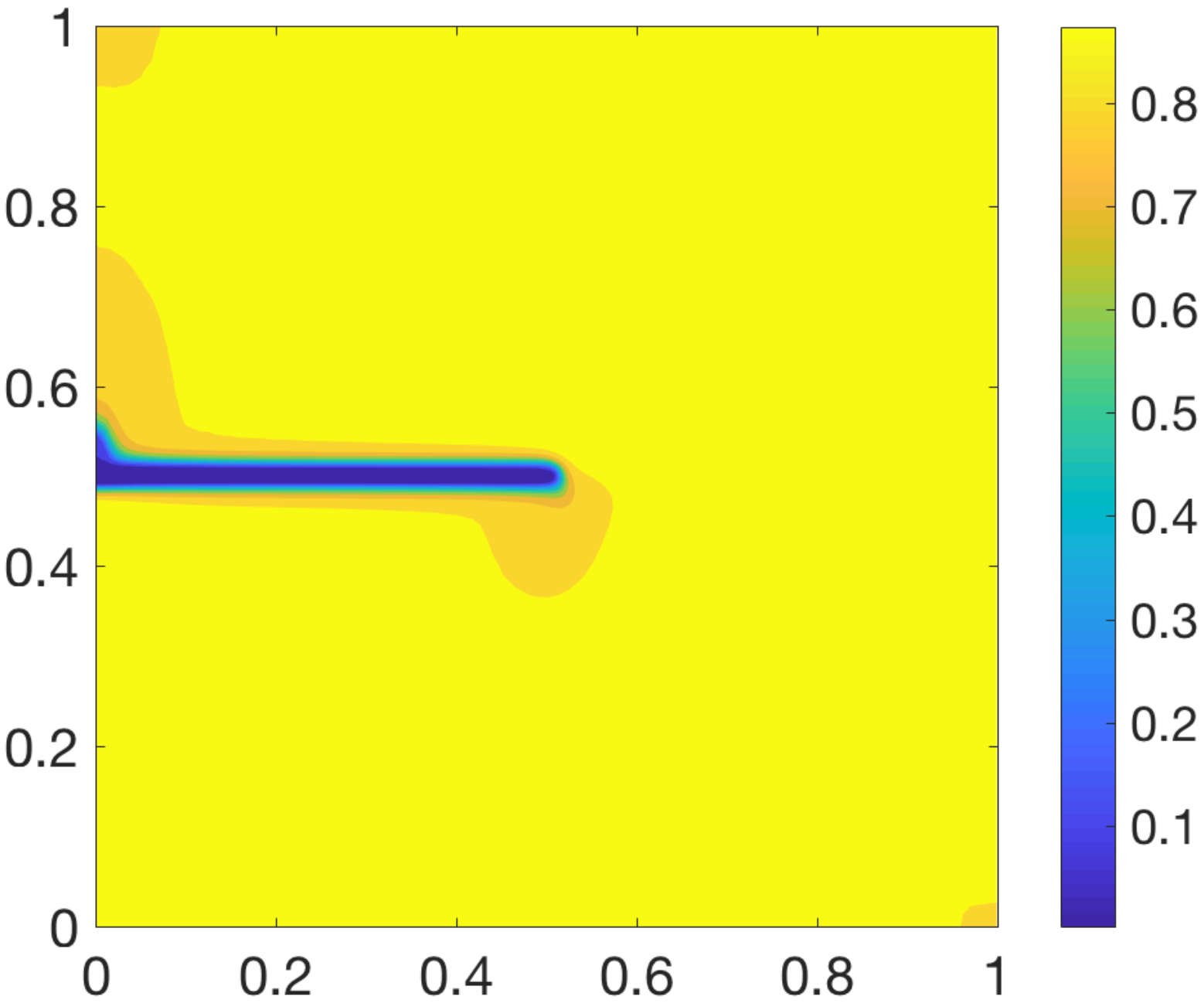}}
\subfigure[$U = 3.25 \times 10^{-2}$~mm]{\label{fig:subfig:ID_U3250}
\includegraphics[width=0.25\linewidth]{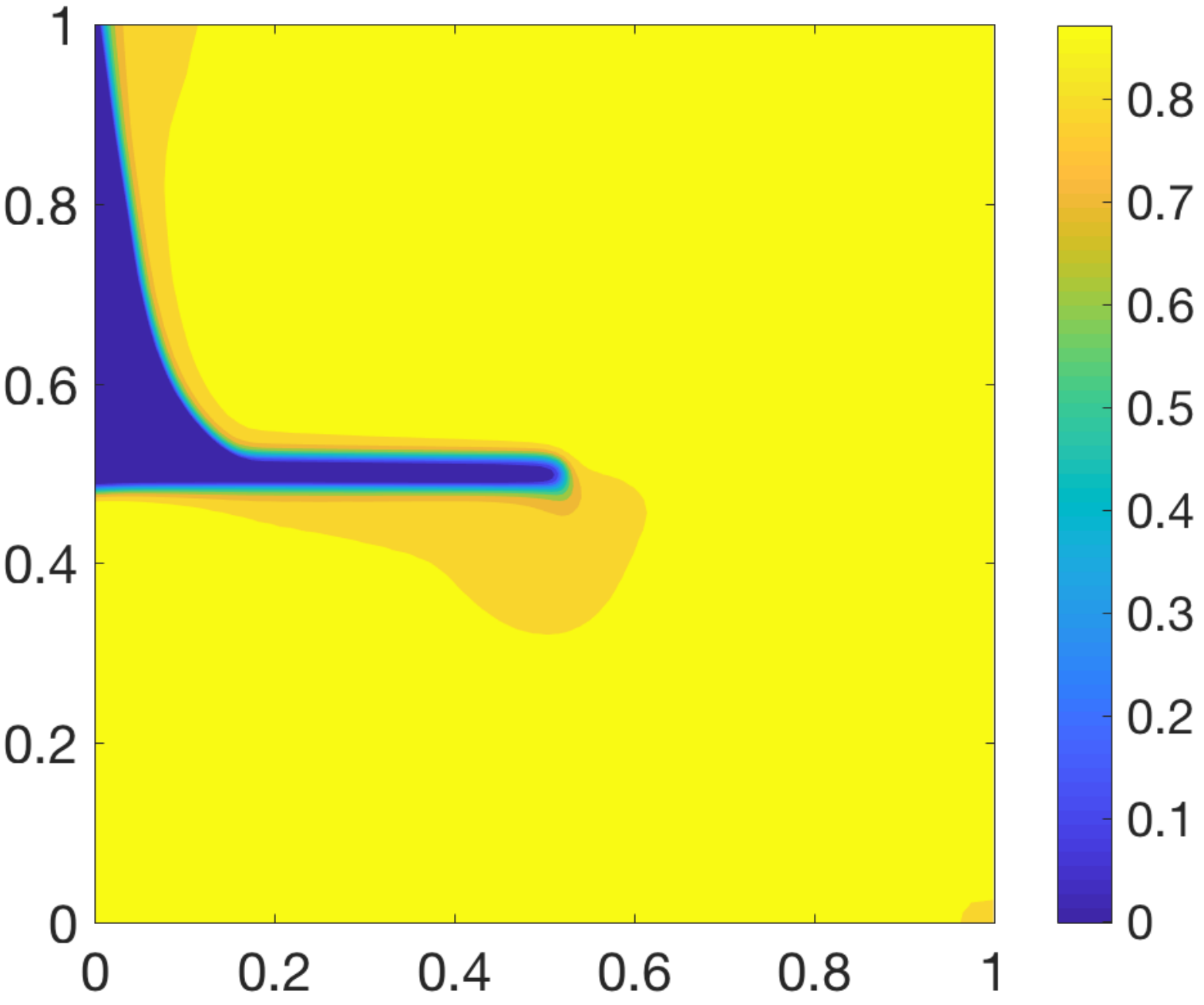}}
\vfill
\subfigure[$U = 3.0 \times 10^{-2}$~mm]{\label{fig:subfig:ID_U3000}
\includegraphics[width=0.25\linewidth]{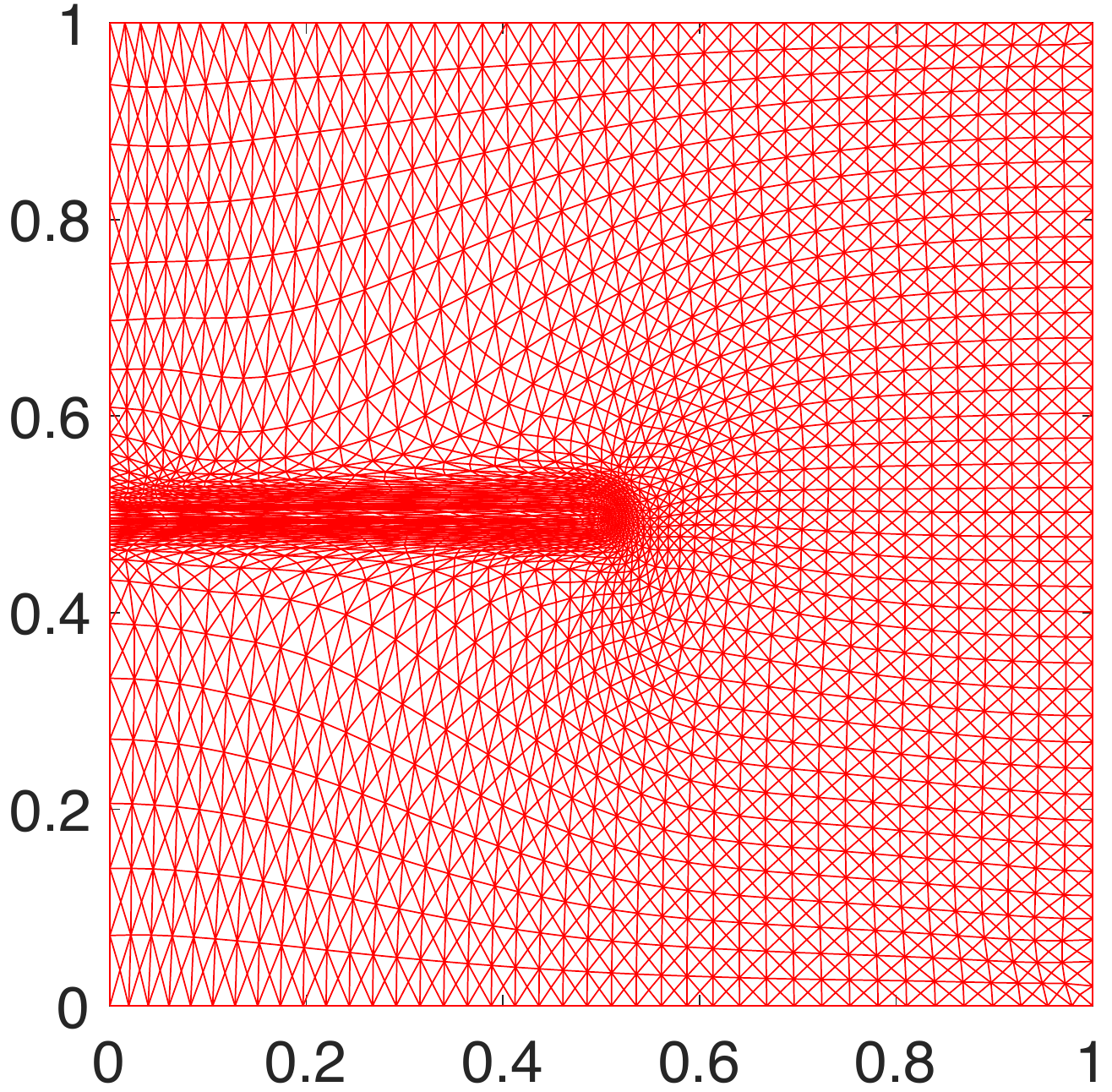}}
\subfigure[$U = 3.2 \times 10^{-2}$~mm]{\label{fig:subfig:ID_U3200}
\includegraphics[width=0.25\linewidth]{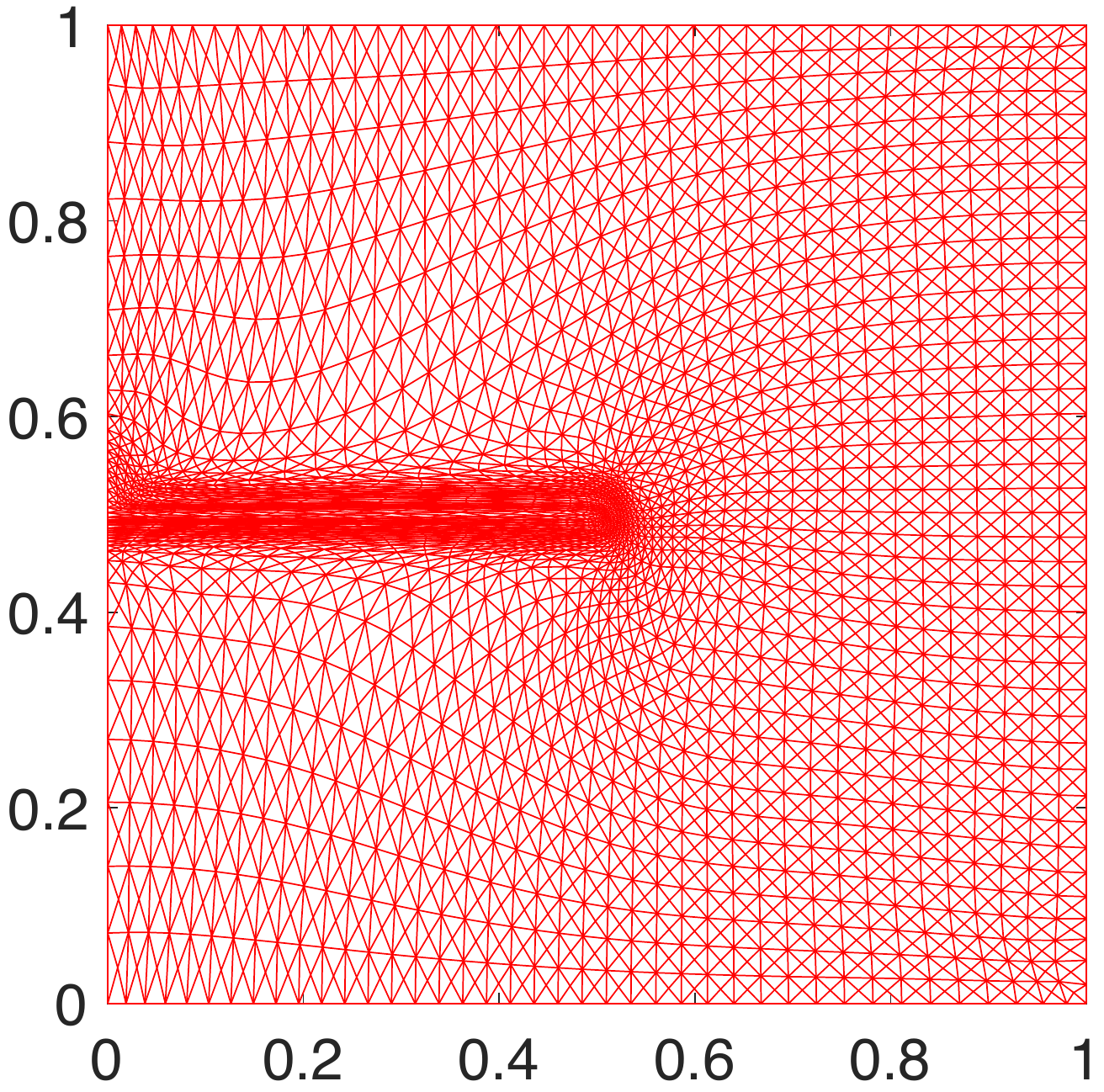}}
\subfigure[$U = 3.25 \times 10^{-2}$~mm]{\label{fig:subfig:ID_U3250}
\includegraphics[width=0.25\linewidth]{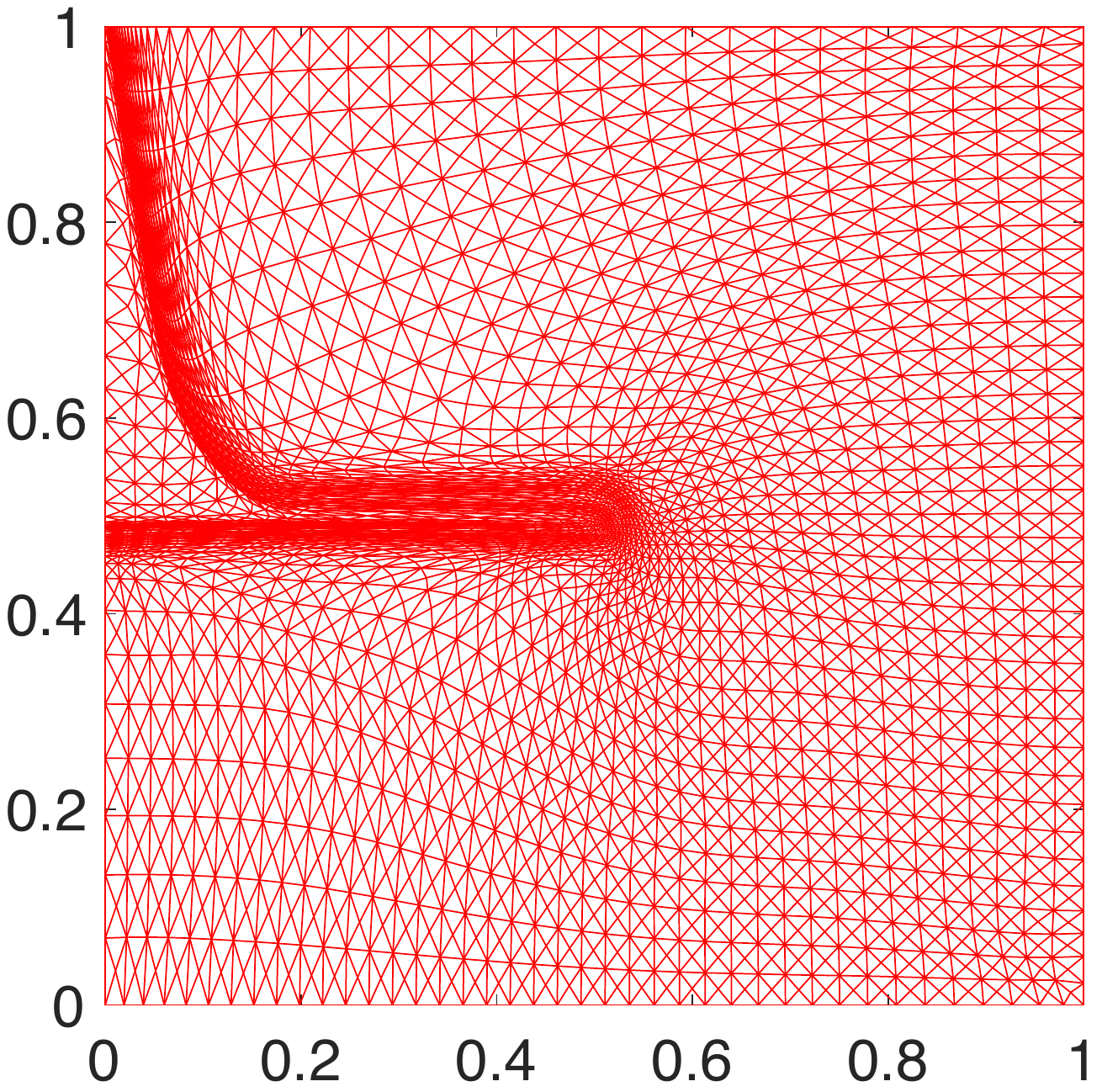}}
\vfill
\subfigure[$U = 3.0 \times 10^{-2}$~mm]{\label{fig:subfig:ID_U3000}
\includegraphics[width=0.25\linewidth]{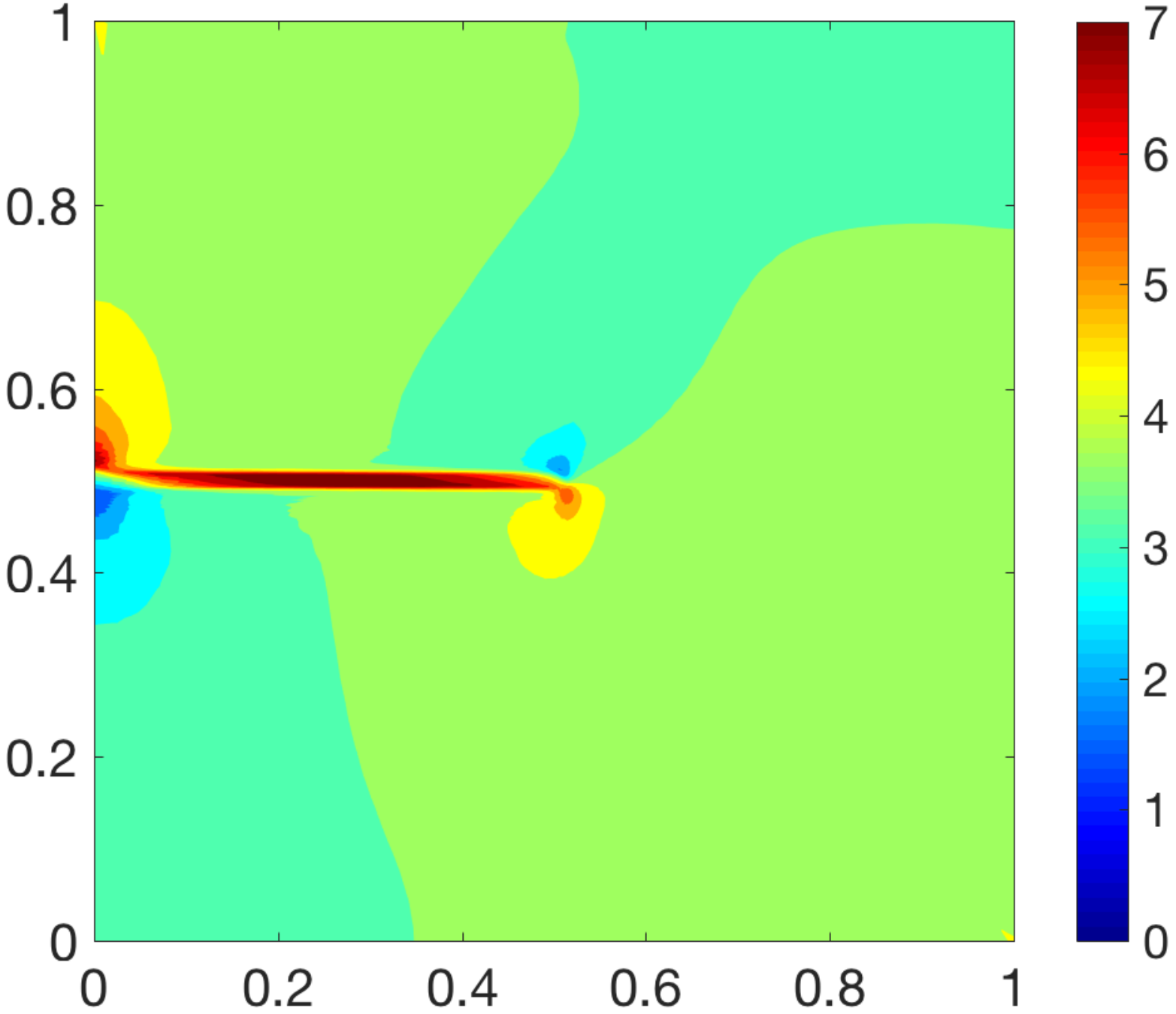}}
\subfigure[$U = 3.2 \times 10^{-2}$~mm]{\label{fig:subfig:ID_U3200}
\includegraphics[width=0.25\linewidth]{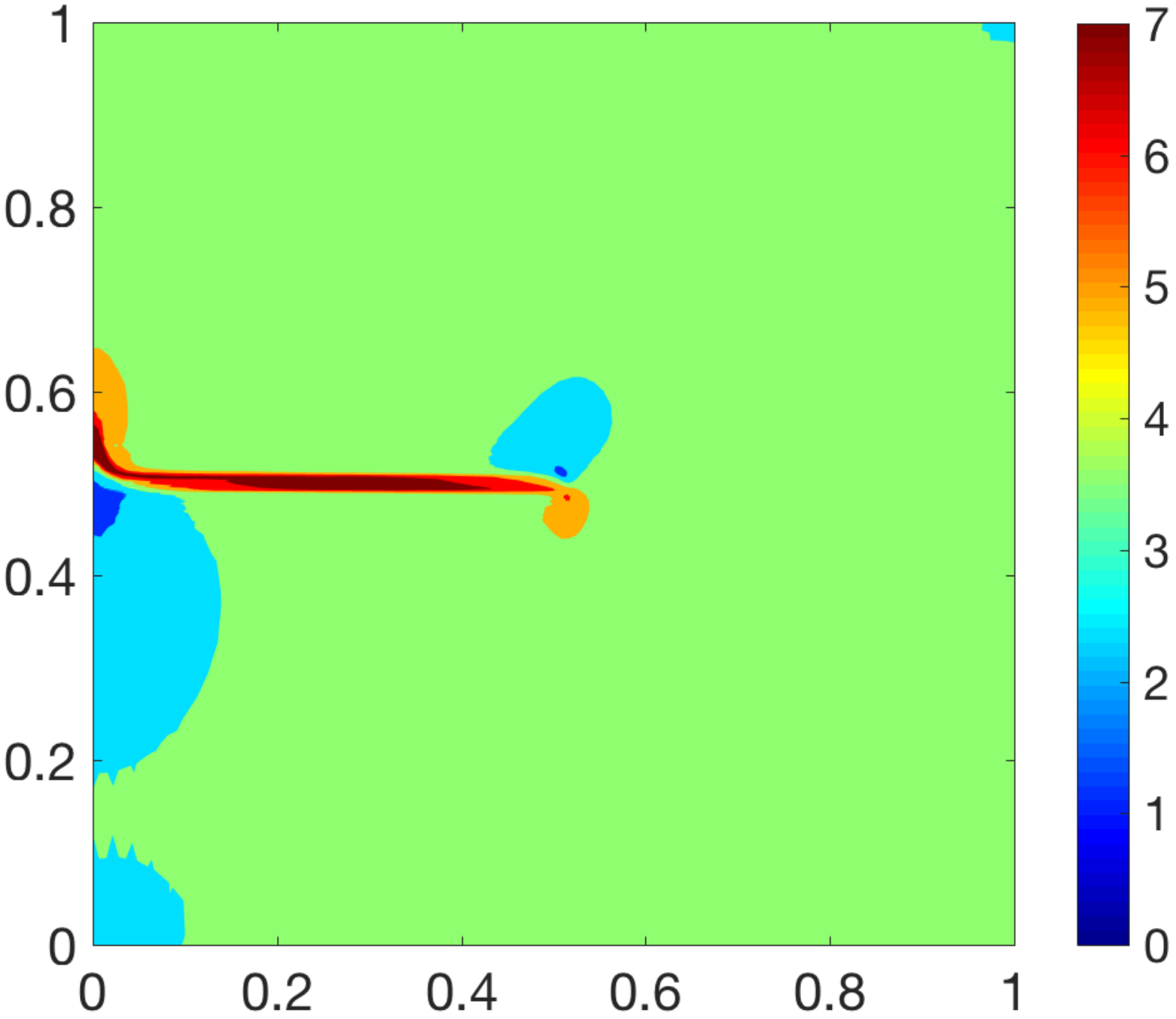}}
\subfigure[$U = 3.25 \times 10^{-2}$~mm]{\label{fig:subfig:ID_U3250}
\includegraphics[width=0.25\linewidth]{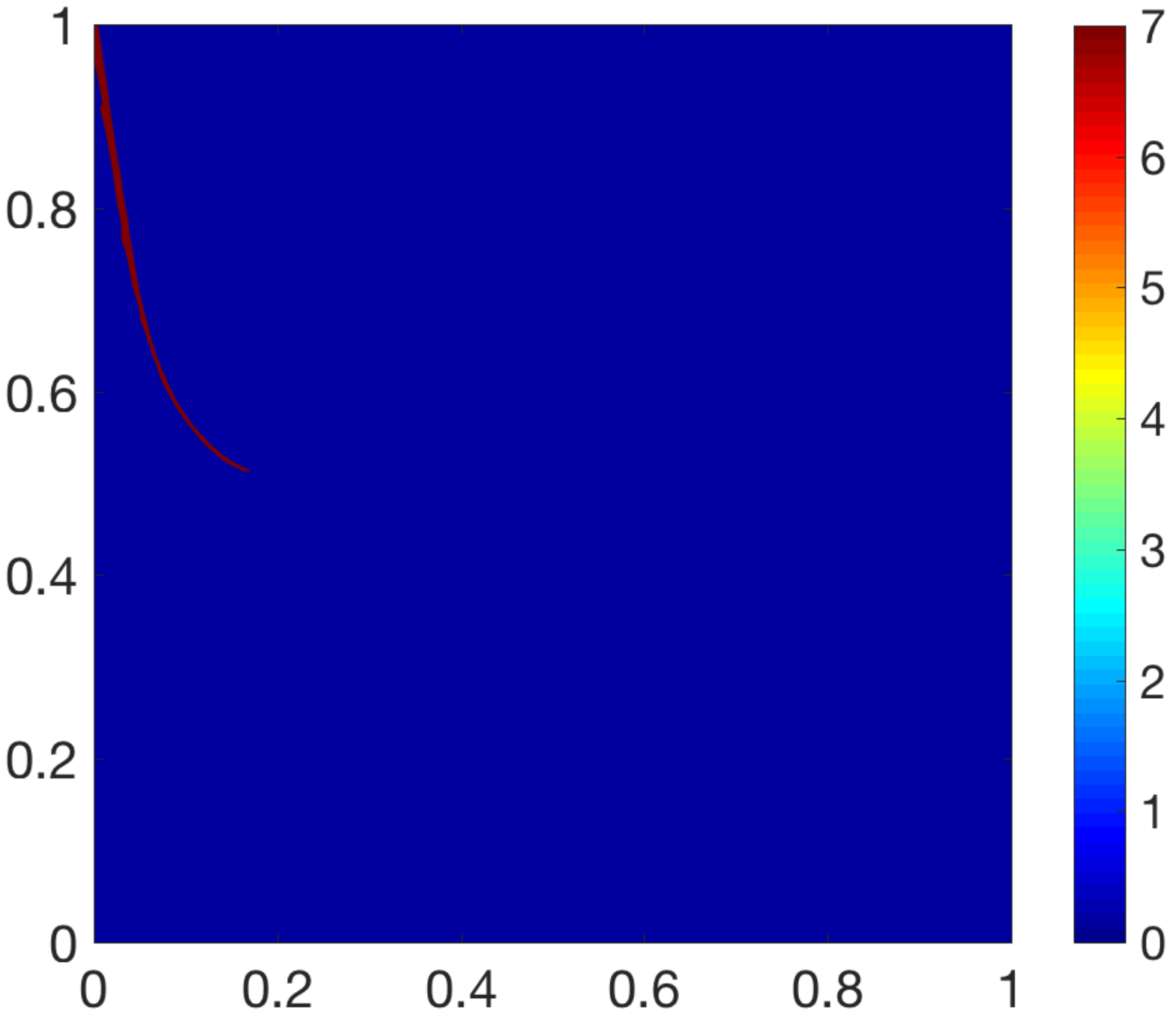}}
\caption{Example 2. Meshes and contours of the phase-field and von Mises stress distribution
are plotted at $U = 3.0\times 10^{-2}$, $3.20\times 10^{-2}$, and $3.25\times 10^{-2}$~mm.
The improved v-d split model is used.}
\label{fig:shear n's OBC}
\end{figure}

\begin{figure} 
\centering 
\subfigure[$U = 1.1 \times 10^{-2}$~mm]{\label{fig:subfig:MD02_U1100}
\includegraphics[width=0.25\linewidth]{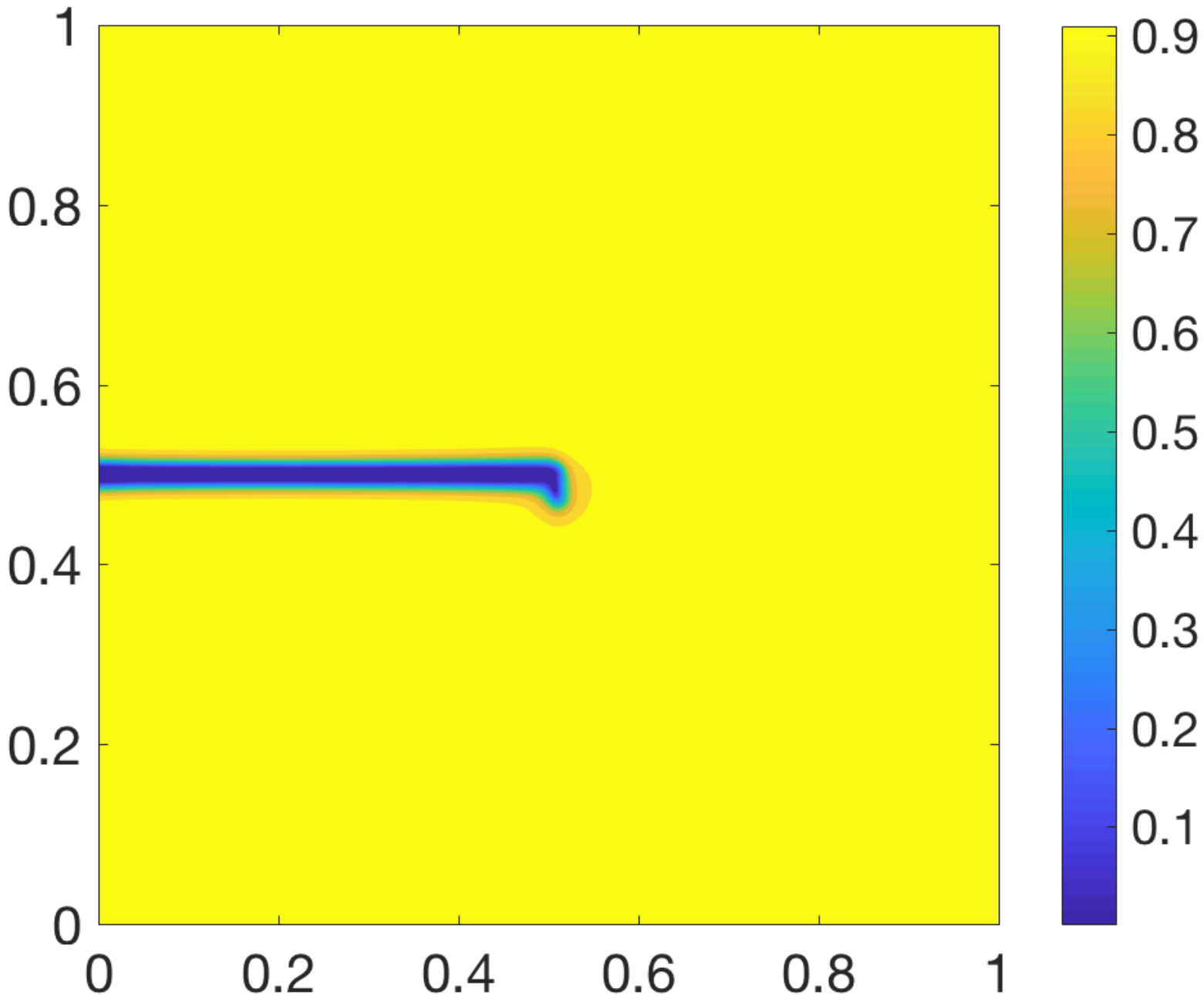}}
\subfigure[$U = 1.3 \times 10^{-2}$~mm]{\label{fig:subfig:MD02_U1300}
\includegraphics[width=0.25\linewidth]{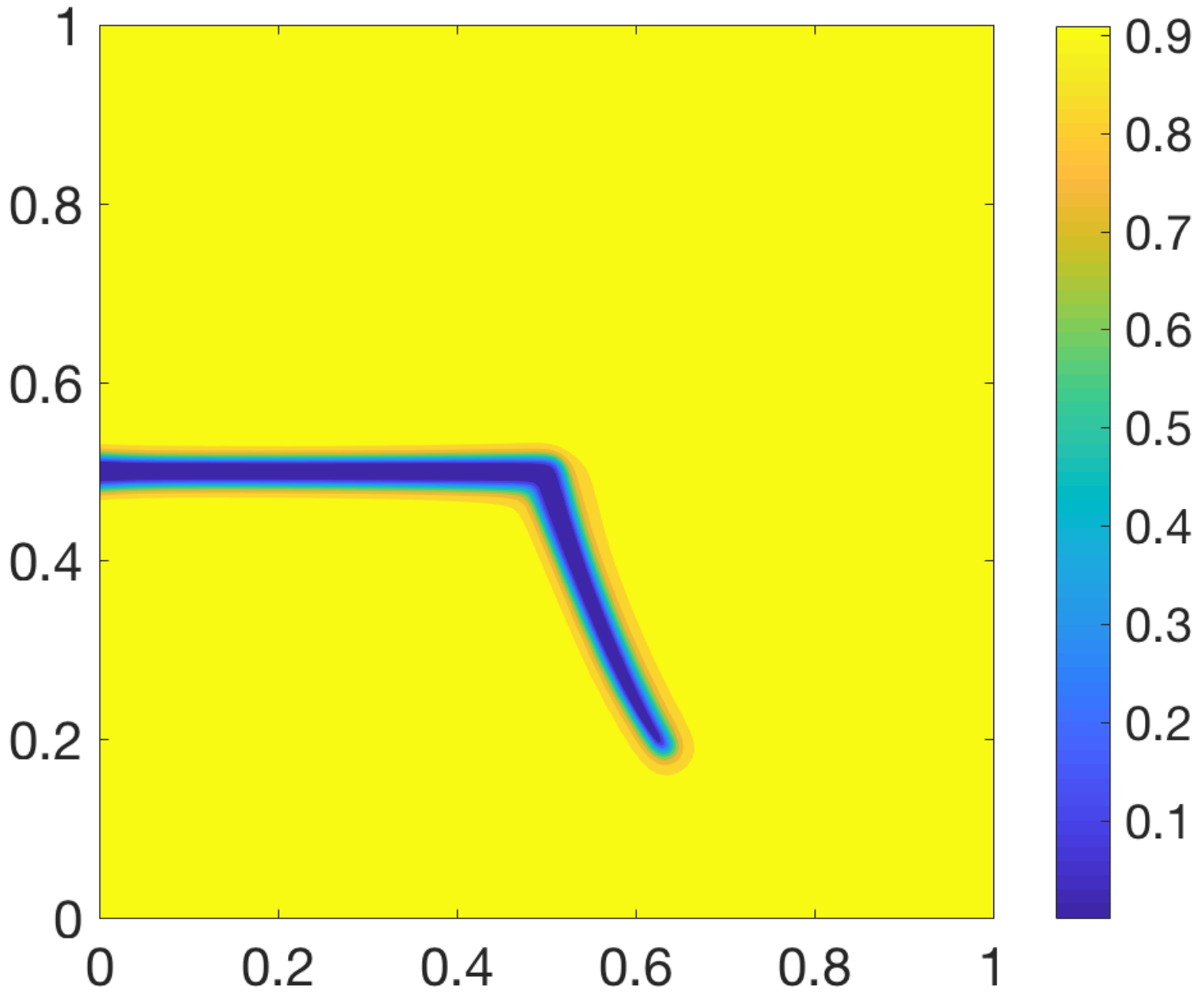}}
\subfigure[$U = 1.45 \times 10^{-2}$~mm]{\label{fig:subfig:MD02_U1450}
\includegraphics[width=0.25\linewidth]{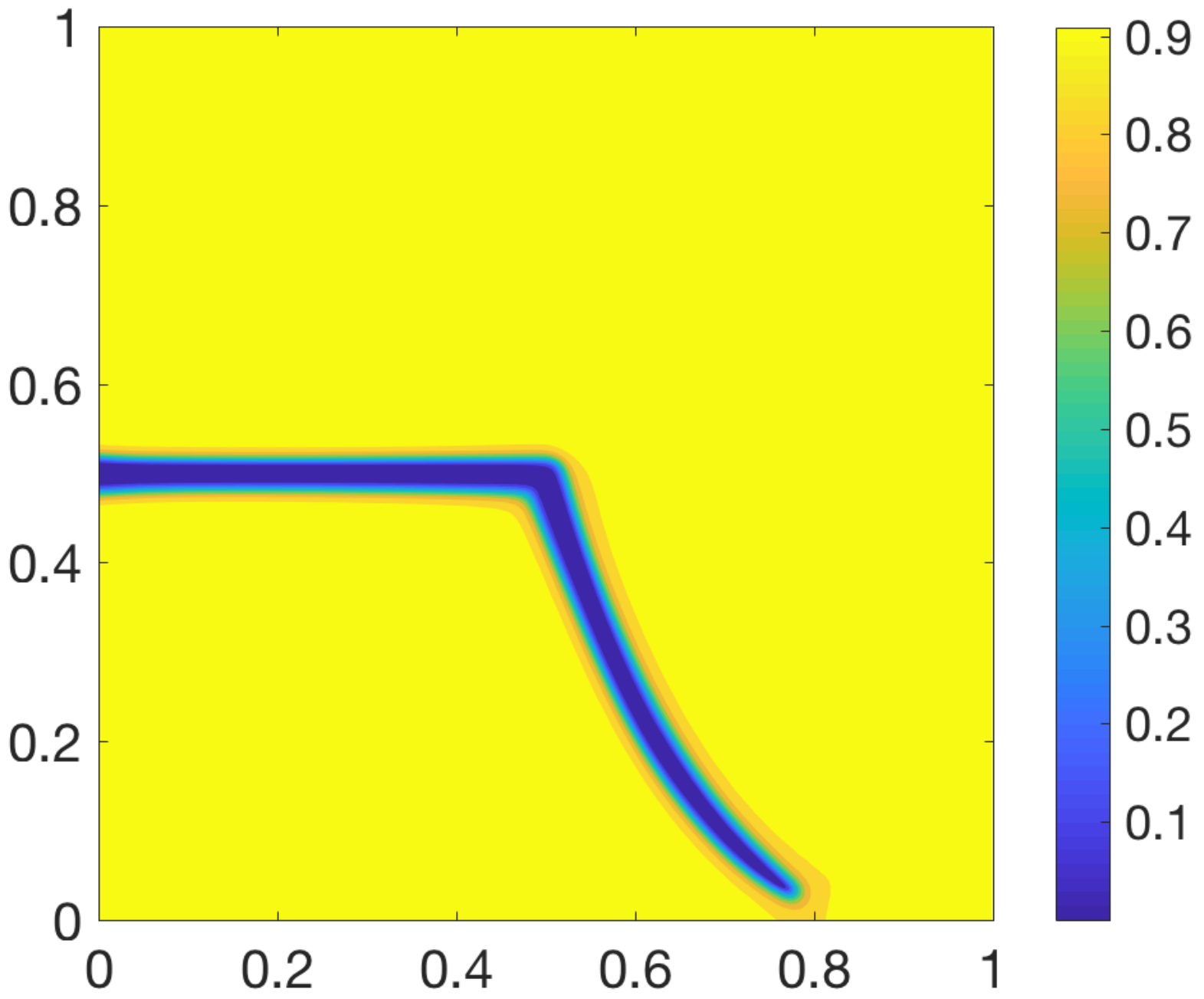}}
\vfill
\subfigure[$U = 1.1 \times 10^{-2}$~mm]{\label{fig:subfig:MM02_U1100}
\includegraphics[width=0.25\linewidth]{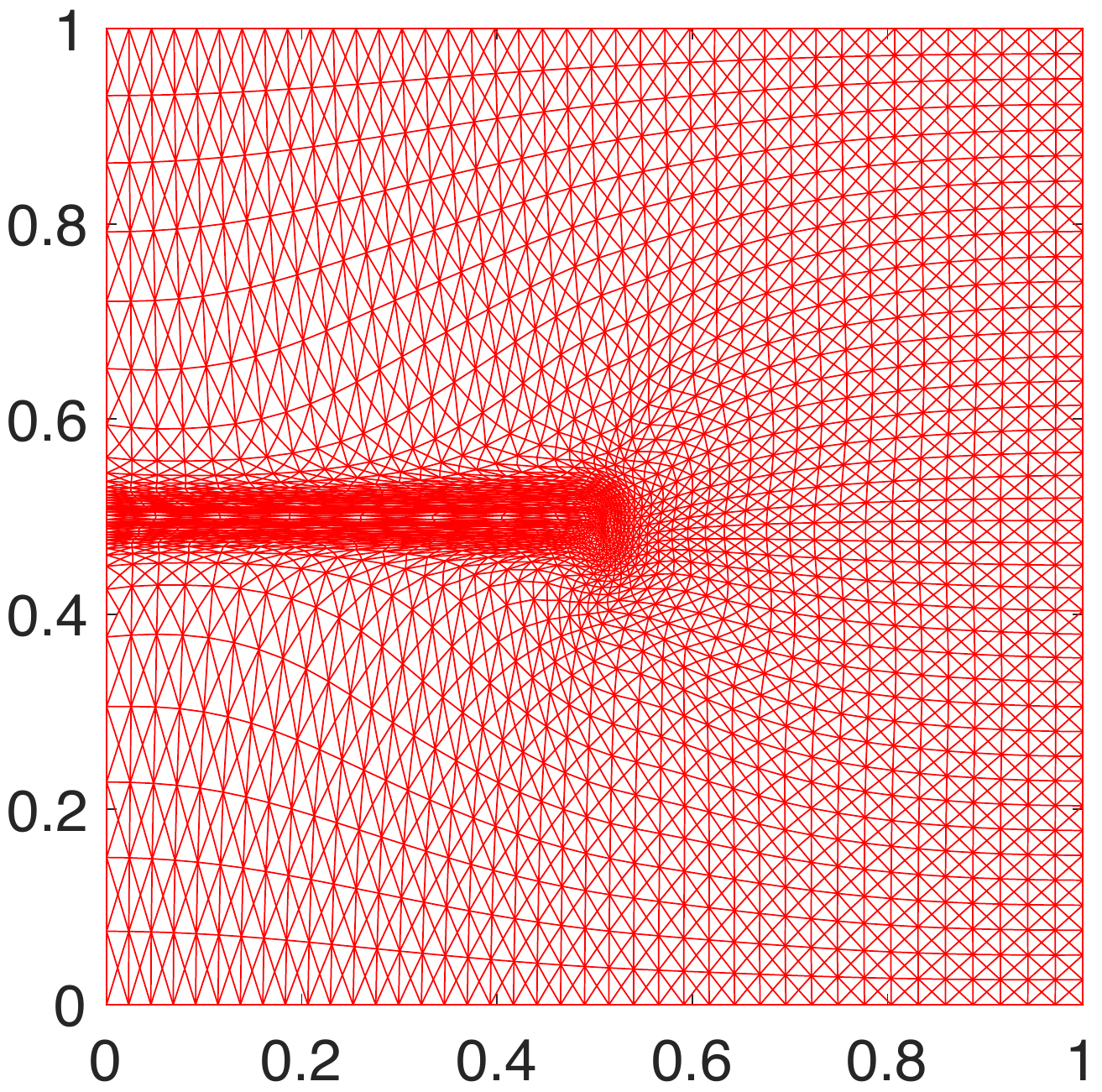}}
\subfigure[$U = 1.3 \times 10^{-2}$~mm]{\label{fig:subfig:MM02_U1300}
\includegraphics[width=0.25\linewidth]{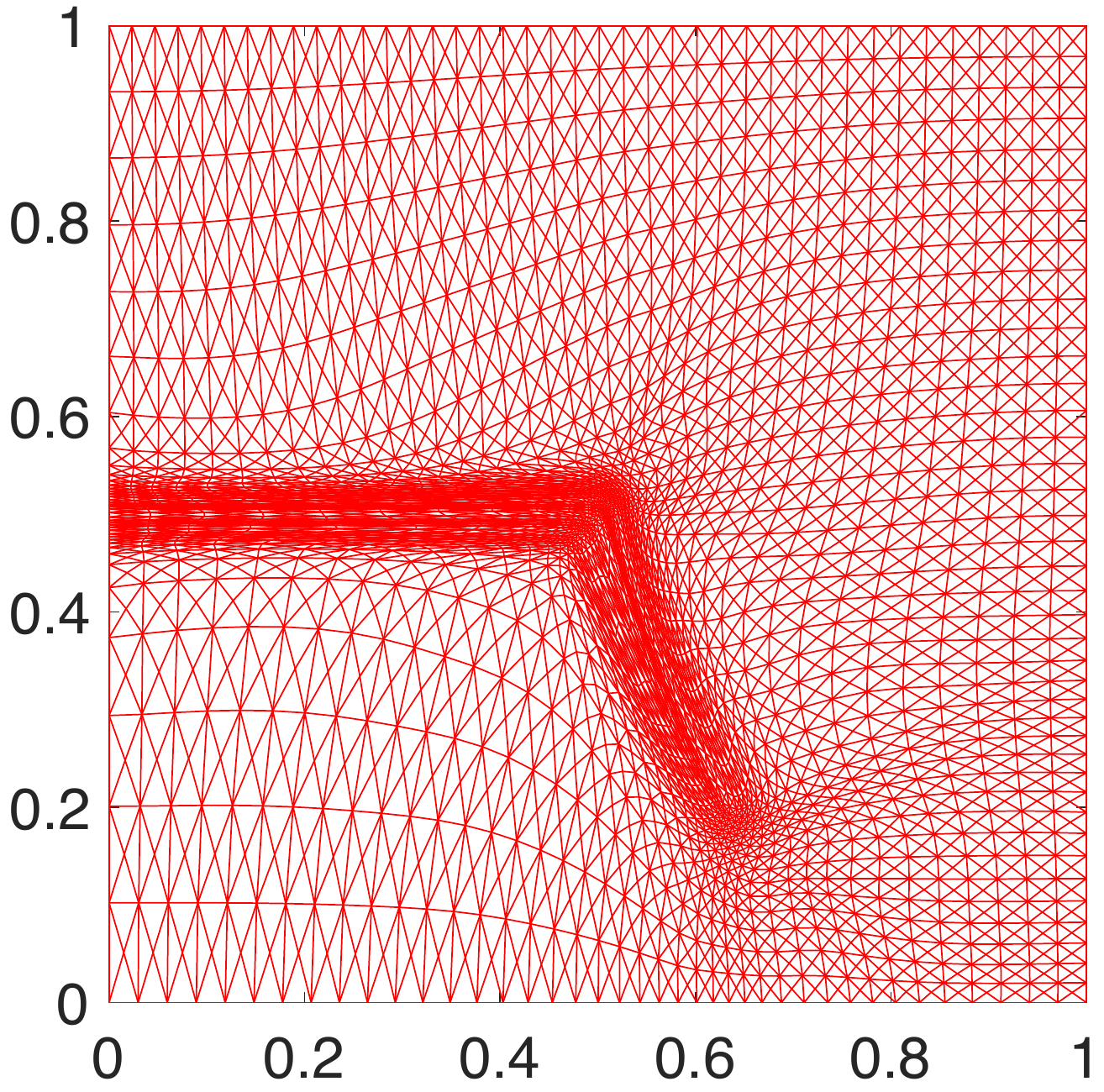}}
\subfigure[$U = 1.45 \times 10^{-2}$~mm]{\label{fig:subfig:MM02_U1450}
\includegraphics[width=0.25\linewidth]{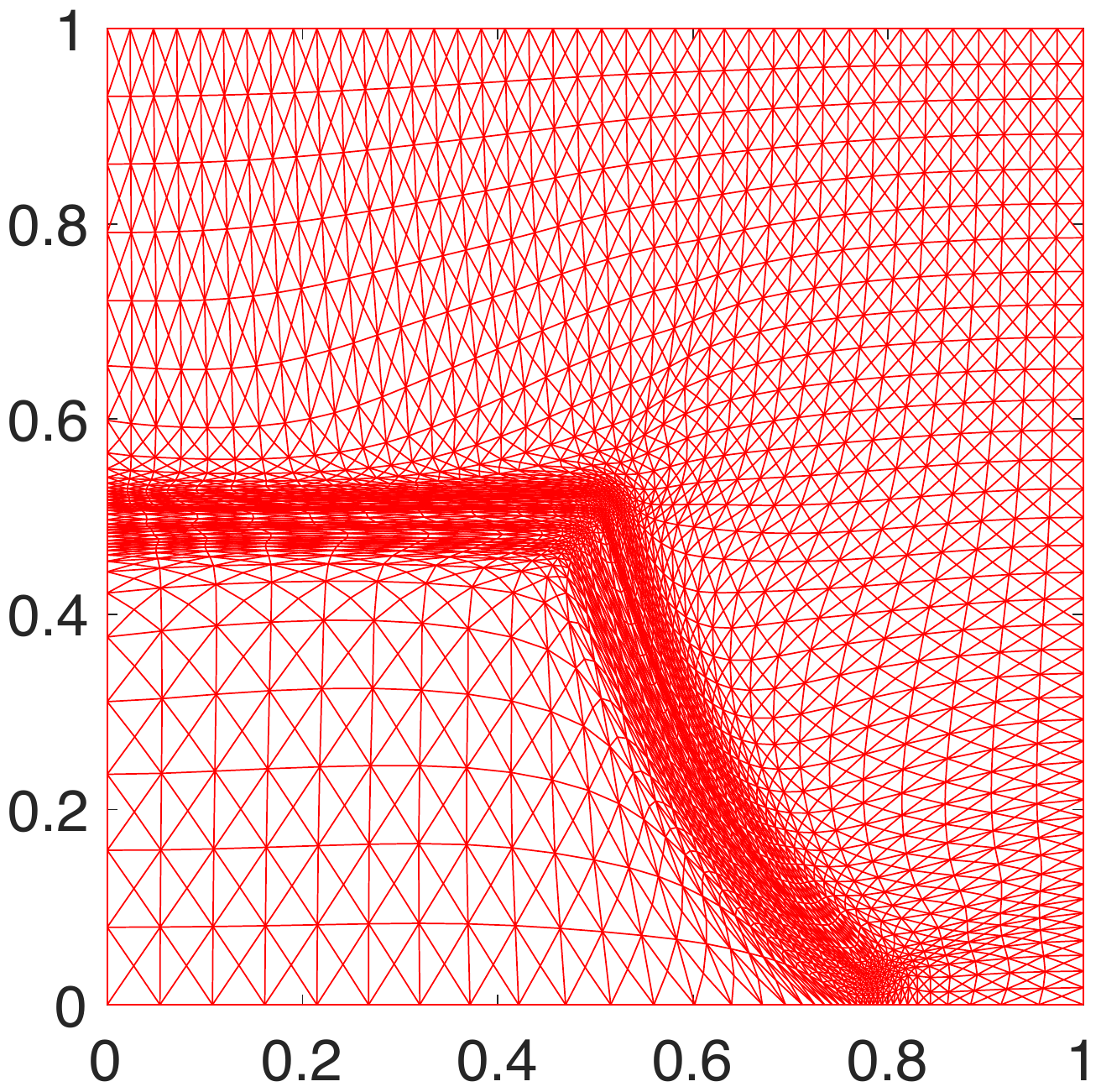}}
\vfill
\subfigure[$U = 1.1 \times 10^{-2}$~mm]{\label{fig:subfig:MS02_U1100}
\includegraphics[width=0.25\linewidth]{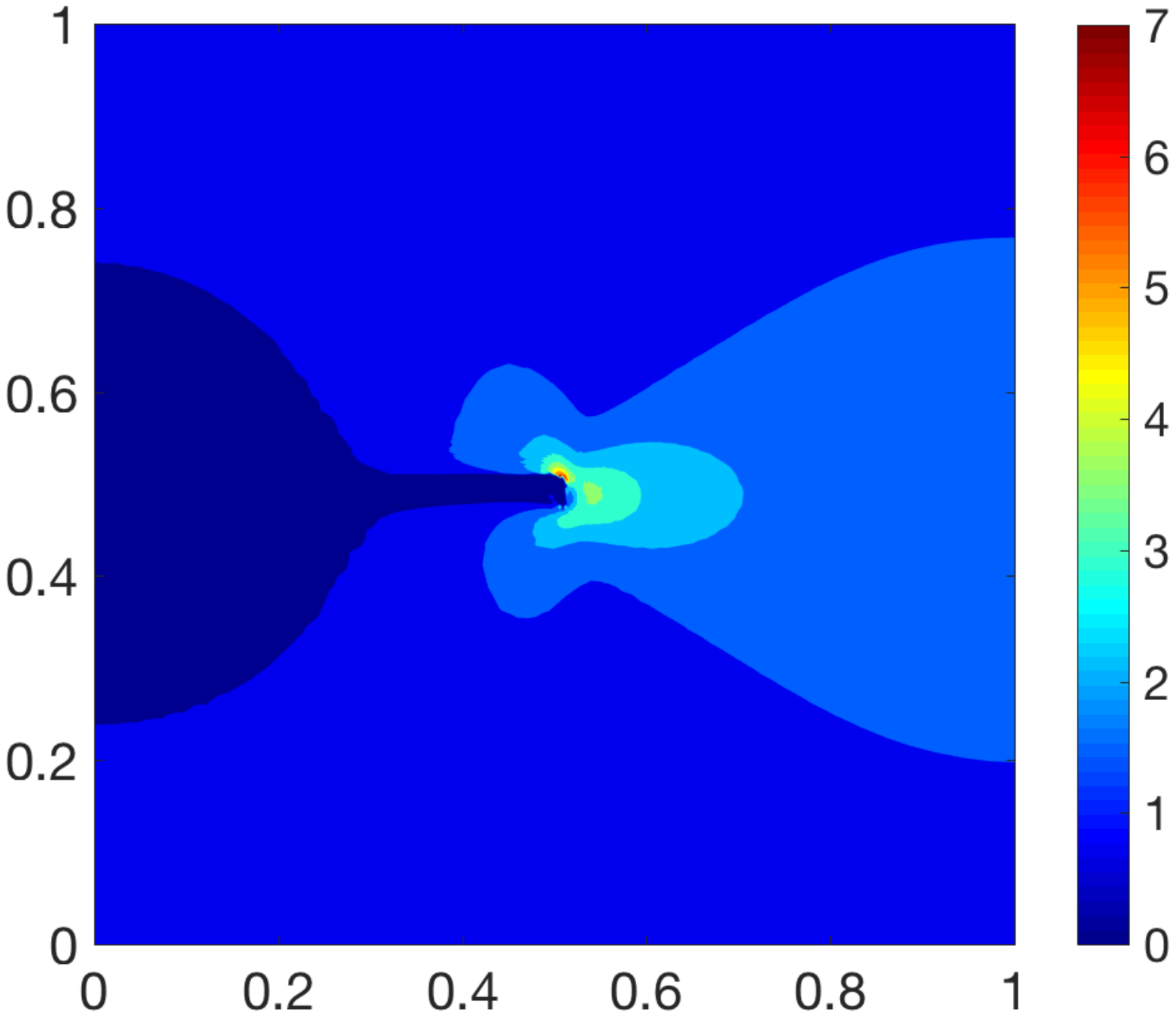}}
\subfigure[$U = 1.3 \times 10^{-2}$~mm]{\label{fig:subfig:MS02_U1300}
\includegraphics[width=0.25\linewidth]{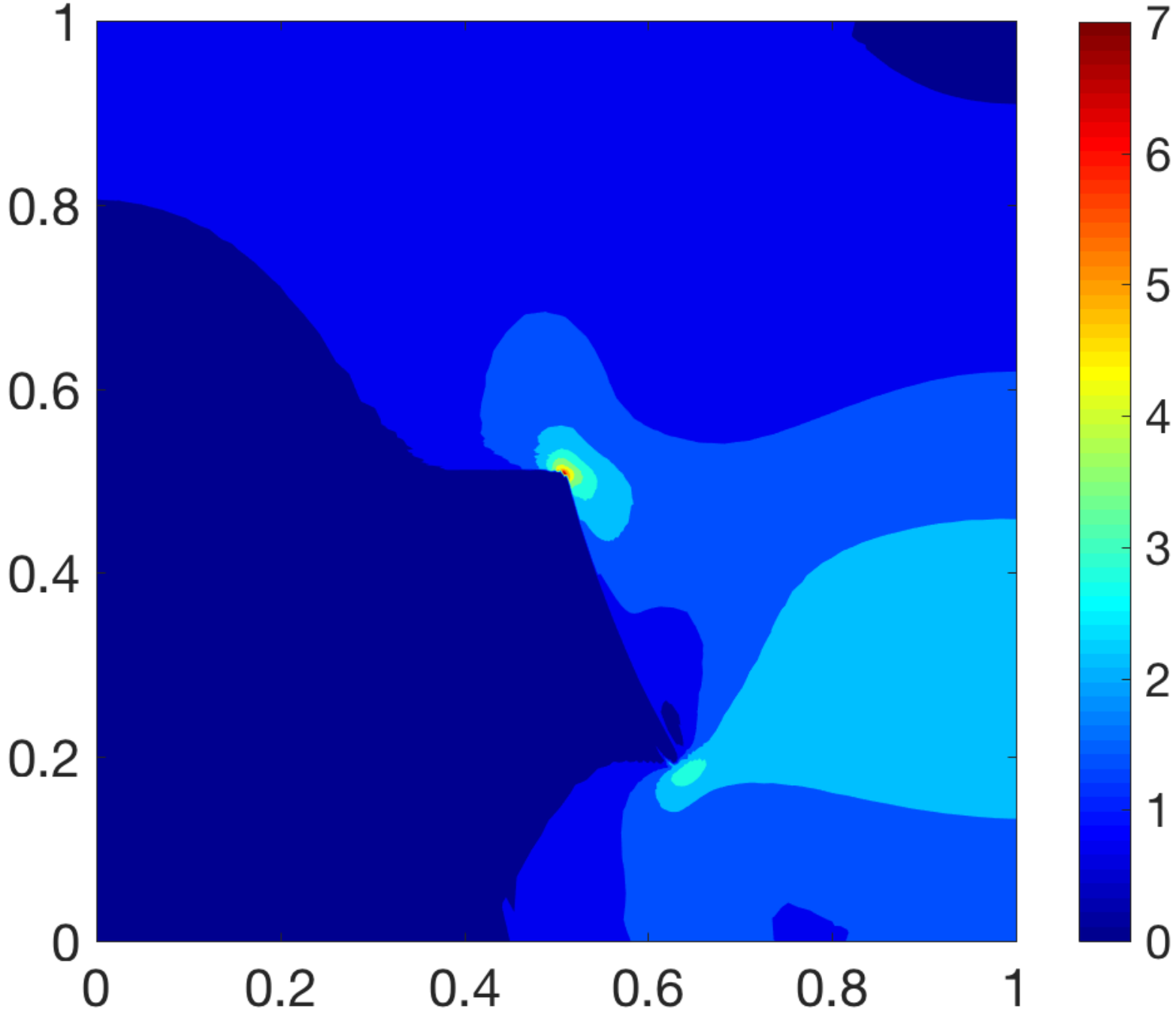}}
\subfigure[$U = 1.45 \times 10^{-2}$~mm]{\label{fig:subfig:MS02_U1450}
\includegraphics[width=0.25\linewidth]{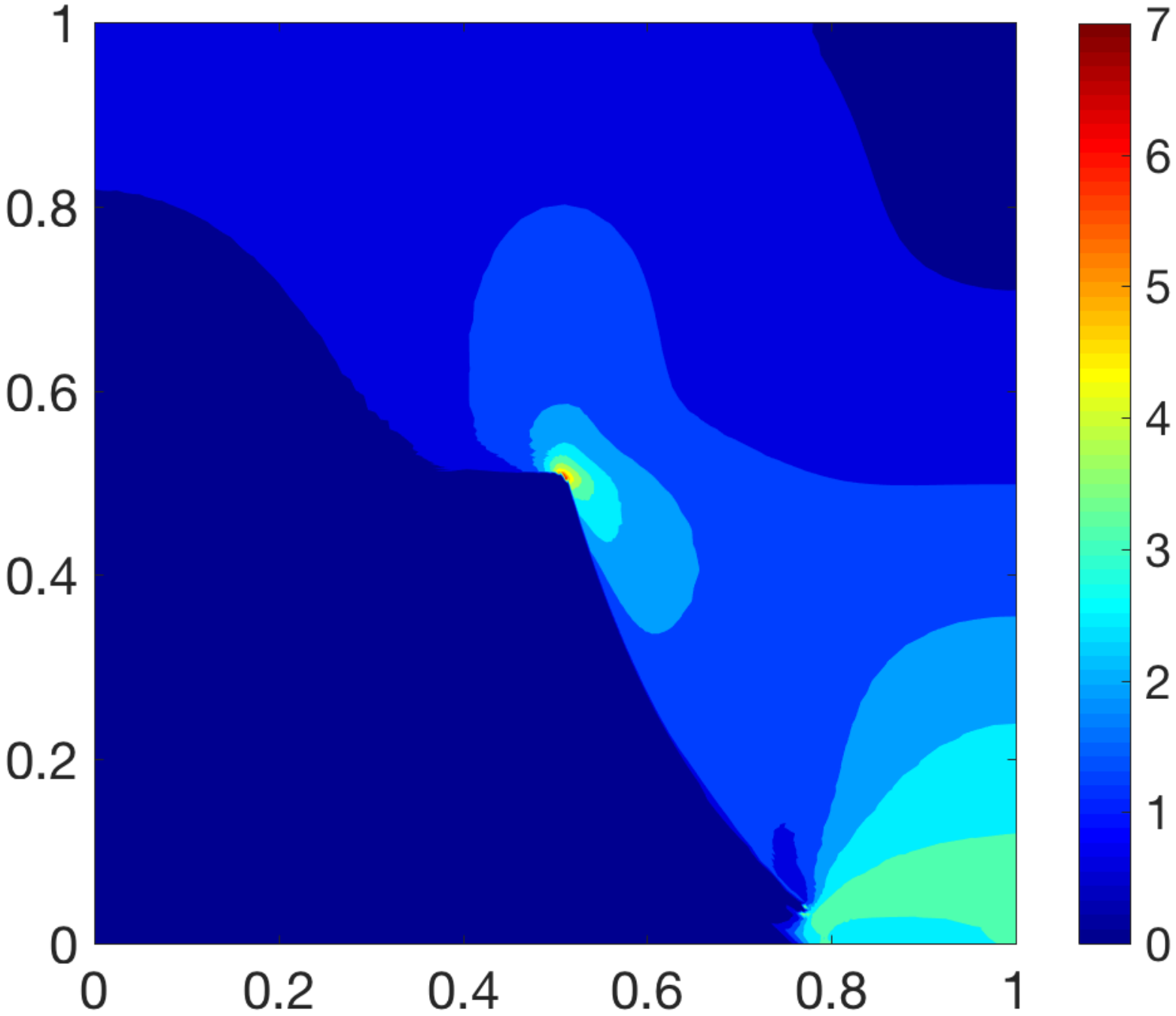}}
\caption{Example 2. Meshes and contours of the phase-field and von Mises stress distribution
are plotted at $U = 1.1\times 10^{-2}$, $1.30\times 10^{-2}$, and $1.45\times 10^{-2}$~mm.
The improved v-d split model with ItCBC ($d_{cr} = 0.4$) is used.}
\label{fig:shear n's MBC}
\end{figure}

Figs.~\ref{fig:shear a's OBC} and \ref{fig:shear a's MBC} show the results using the v-d split model without and with
ItCBC ($d_{cr} = 0.4$), respectively. The crack evolution paths are similar to each other. However, they are
very different from that with the spectral decomposition model with ItCBC.
The reason is that the v-d split model accounts the compressive deviatoric strain energy
for contributing to the damage process and this leads to unphysical crack evolution. 
On the other hand, the improved v-d split model takes this energy component off the damage process
and, together with ItCBC ($d_{cr} = 0.4$), leads to a correct crack propagation path;
see Fig.~\ref{fig:shear n's MBC}.  The results are comparable to those in Fig.~\ref{fig:shear m's MBC}
for spectral decomposition with ItCBC. Fig.~\ref{fig:shear n's OBC} shows the results for the improved v-d split
model without ItCBC. Once again, the stress concentrates in the damaged region. Moreover, the crack
hardly propagates at the tip. Instead, a secondary crack develops along the left edge of the domain.

\begin{figure} 
\centering 
\subfigure[spectral decomposition]{\label{fig:subfig:Sld_MieheModified}
\includegraphics[width=0.4\linewidth]{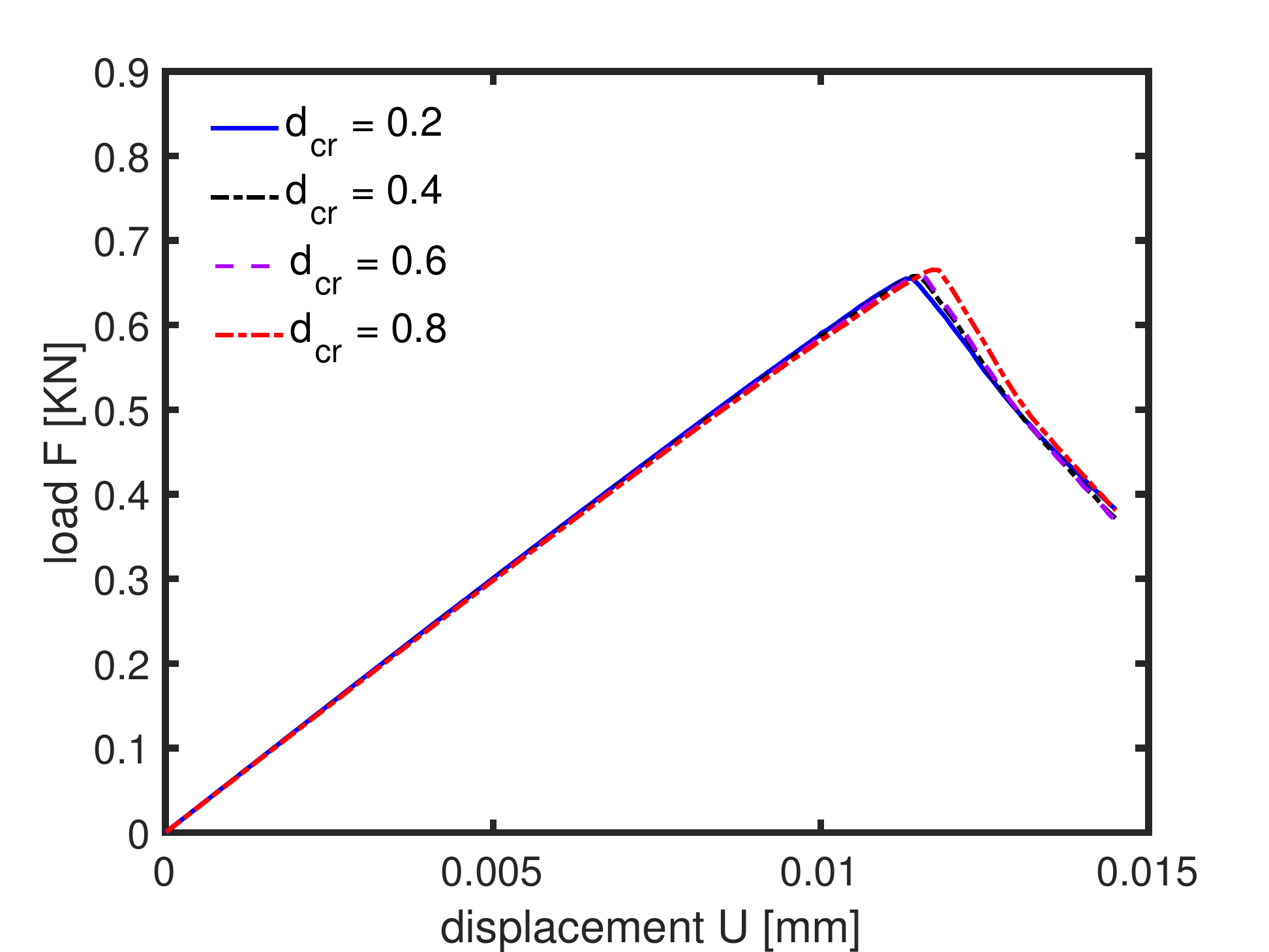}}
\subfigure[improved v-d split]{\label{fig:subfig:Sld_NewModified}
\includegraphics[width=0.4\linewidth]{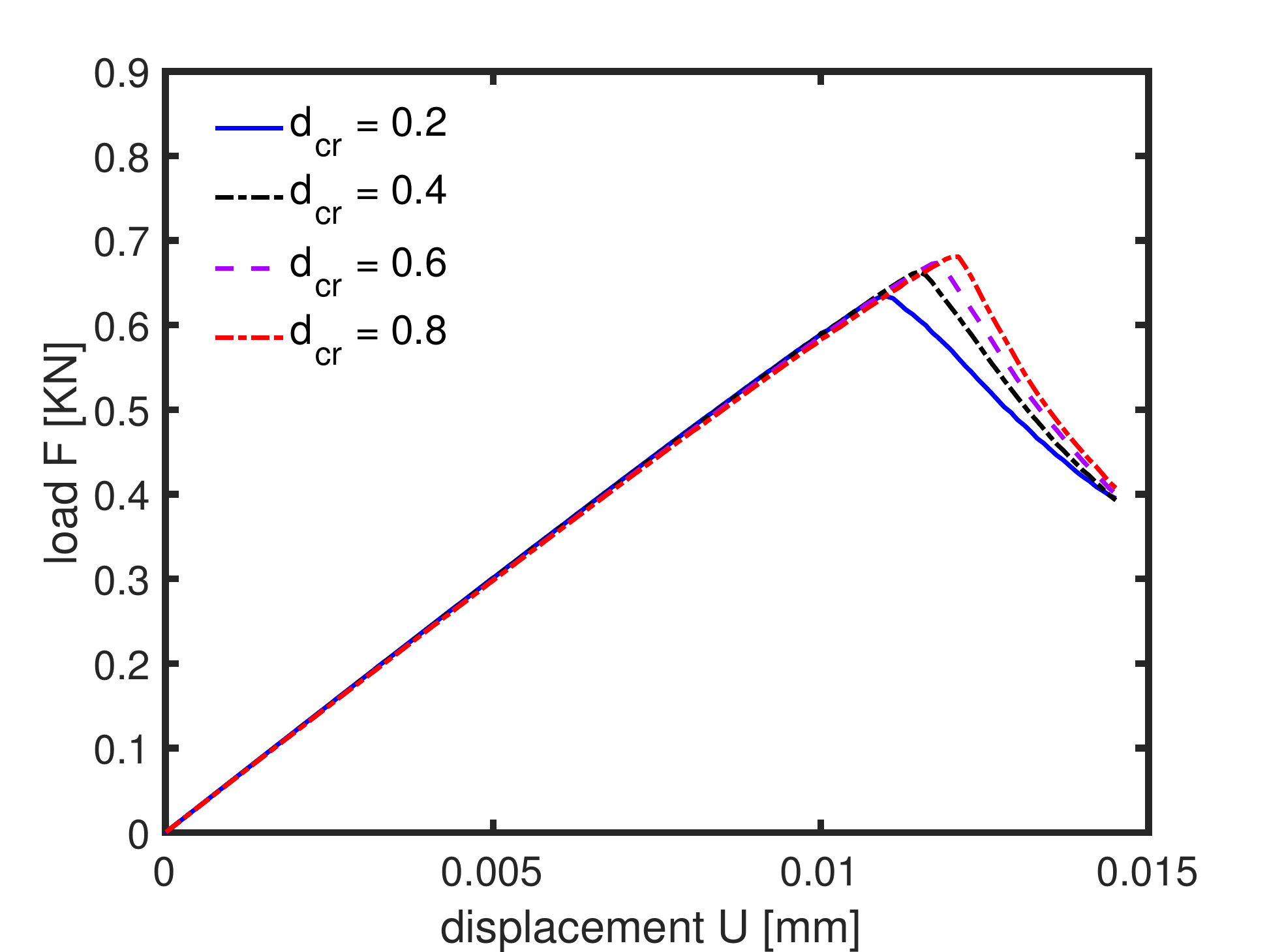}}
\caption{Example 2. The load-deflection curves are plotted for two decomposition models.
(a) The spectral decomposition model with ItCBC (with various $d_{cr}$);
(b) The improved v-d split model with ItCBC (with various $d_{cr}$).}
\label{fig:shear dcr-1}
\end{figure}

\begin{figure} 
\centering 
\subfigure[$d_{cr} = 0.2$]{\label{fig:subfig:Sld_DiffSplit02}
\includegraphics[width=0.4\linewidth]{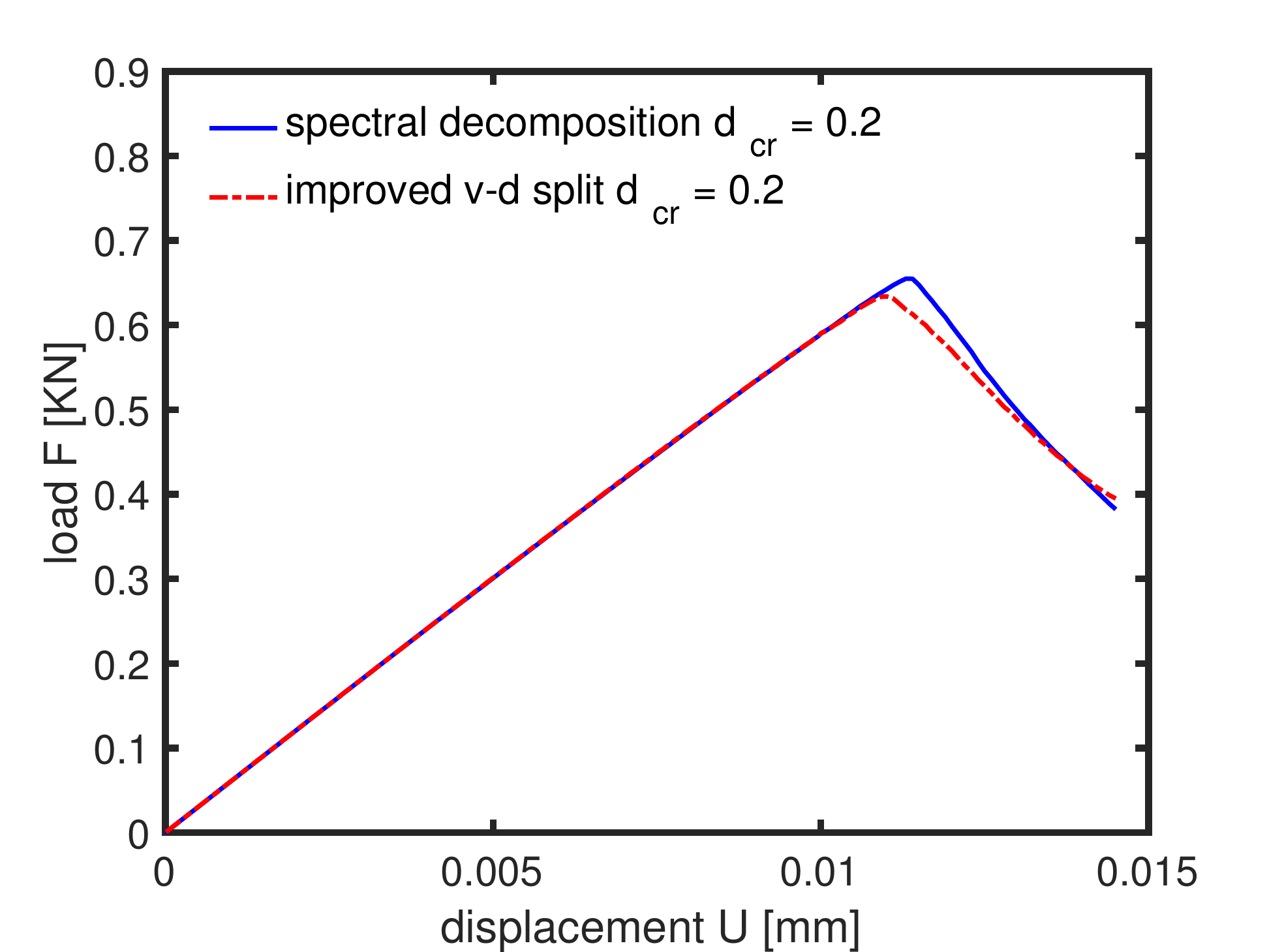}}
\subfigure[$d_{cr} = 0.4$]{\label{fig:subfig:Sld_DiffSplit04}
\includegraphics[width=0.4\linewidth]{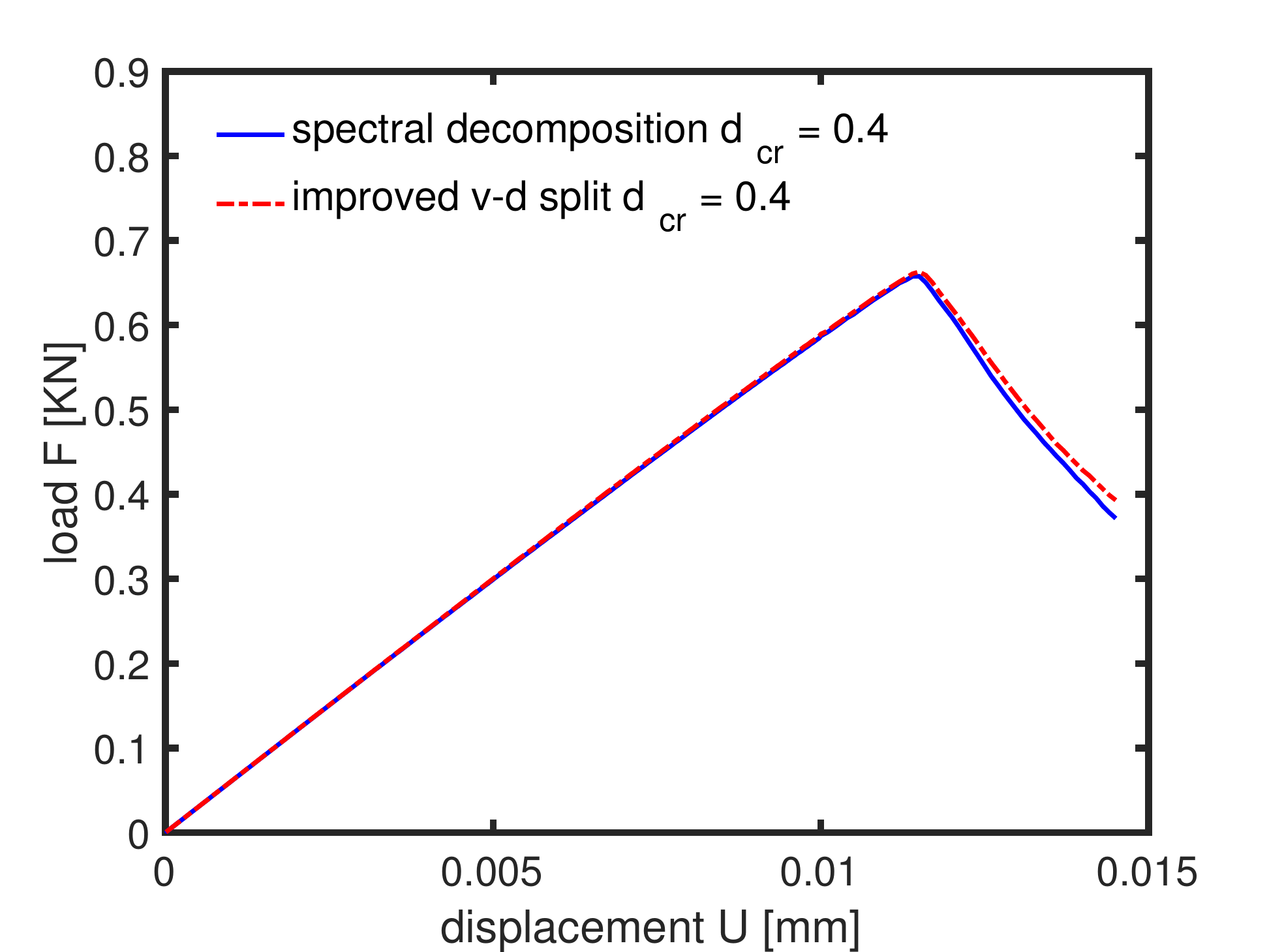}}
\subfigure[$d_{cr} = 0.6$]{\label{fig:subfig:Sld_DiffSplit06}
\includegraphics[width=0.4\linewidth]{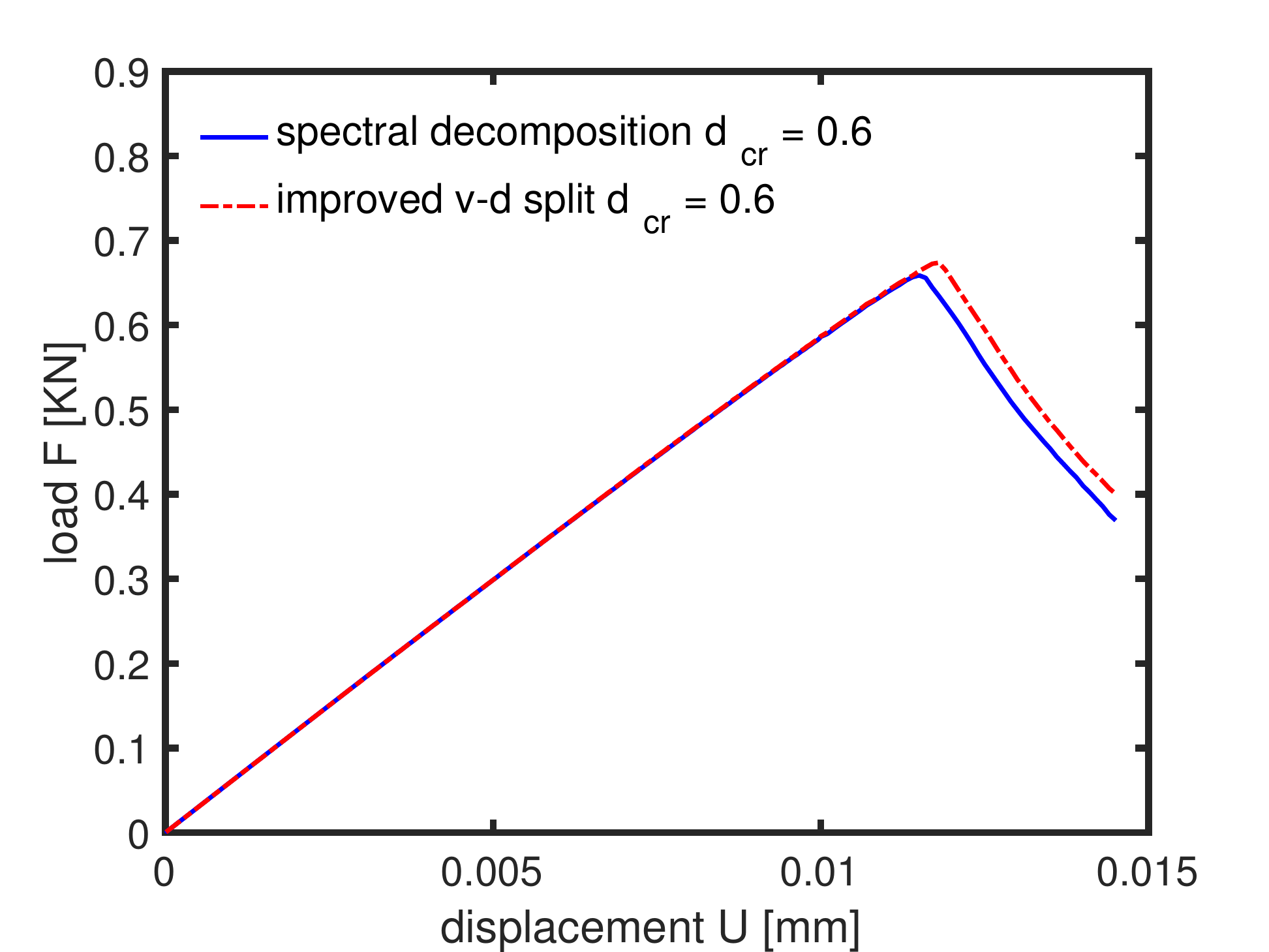}}
\subfigure[$d_{cr} = 0.8$]{\label{fig:subfig:Sld_DiffSplit08}
\includegraphics[width=0.4\linewidth]{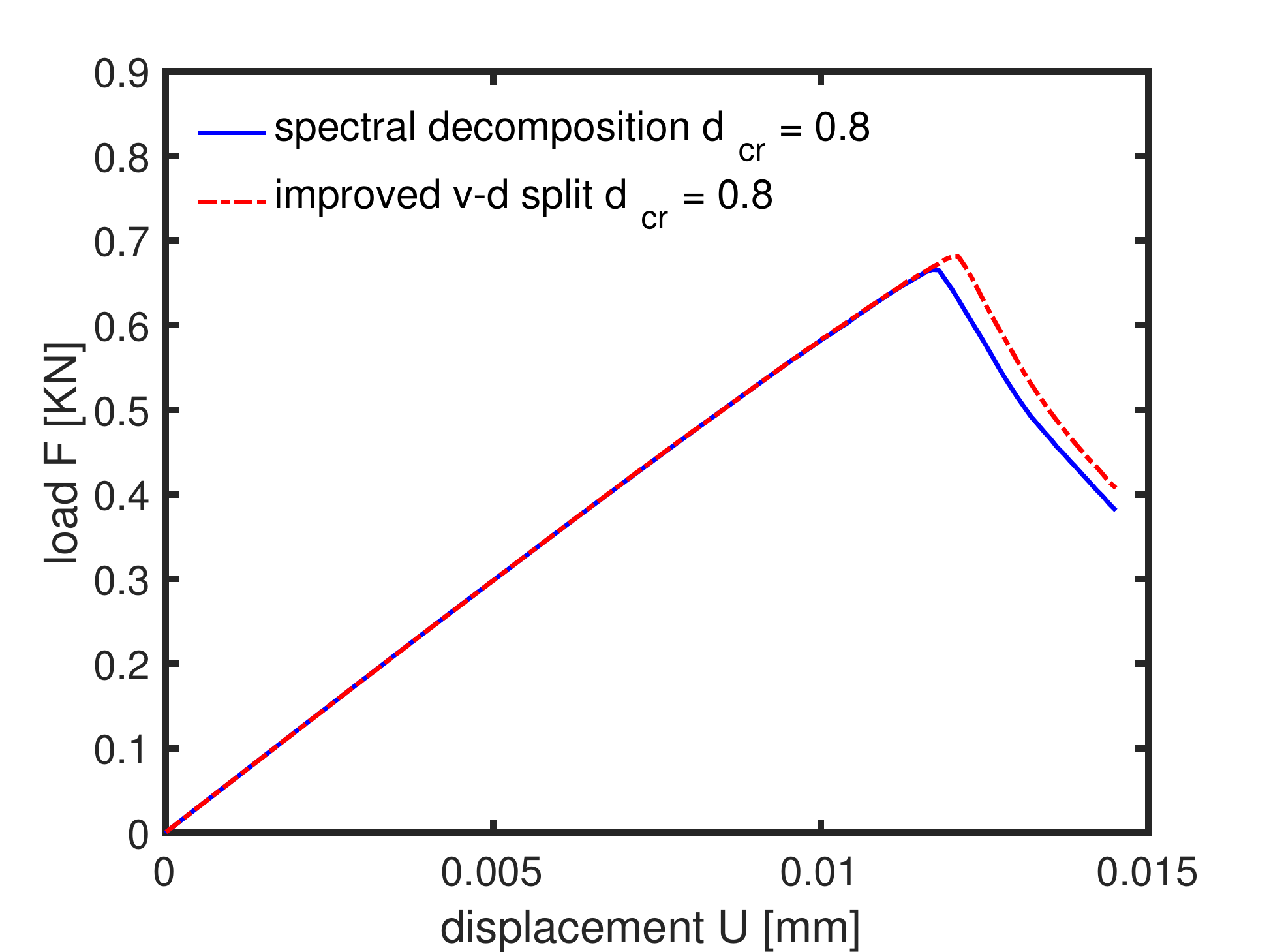}}
\caption{Example 2. The load-deflection curves are compared for two decomposition models
with ItCBC.
(a) $d_{cr} = 0.2$;
(b) $d_{cr} = 0.4$;
(c) $d_{cr} = 0.6$;
(d) $d_{cr} = 0.8$.}
\label{fig:shear dcr-2}
\end{figure}

Next, we investigate the effects of the choice of $d_{cr}$.
The load-deflection curves obtained using the spectral decomposition and improved v-d split model with
various values of $d_{cr}$ are plotted in Figs.~\ref{fig:shear dcr-1} and \ref{fig:shear dcr-2}.

As can be seen in Fig. \ref{fig:subfig:Sld_MieheModified}, the effects of $d_{cr}$ on the load-deflection
are small for the spectral decomposition model when $d_{cr}$ is small ($d_{cr} \le 0.6$).
The load is overestimated after crack starts propagating for large values of $d_{cr}$ (e.g., $d_{cr} = 0.8$).
For the improved v-d split method, the peak of the load-deflection increases with $d_{cr}$; see Fig. \ref{fig:subfig:Sld_NewModified}.
It is also interesting to observe that the curves for the two decomposition methods with $d_{cr} = 0.4$ are nearly identical, as can be seen in Fig. \ref{fig:subfig:Sld_DiffSplit04}.

\subsection{Example 3. A test with multiple cracks}

\begin{figure} 
\centering 
\subfigure[junction between two cracks]{\label{fig:subfig:TwoCrack}
\includegraphics[width=0.4\linewidth]{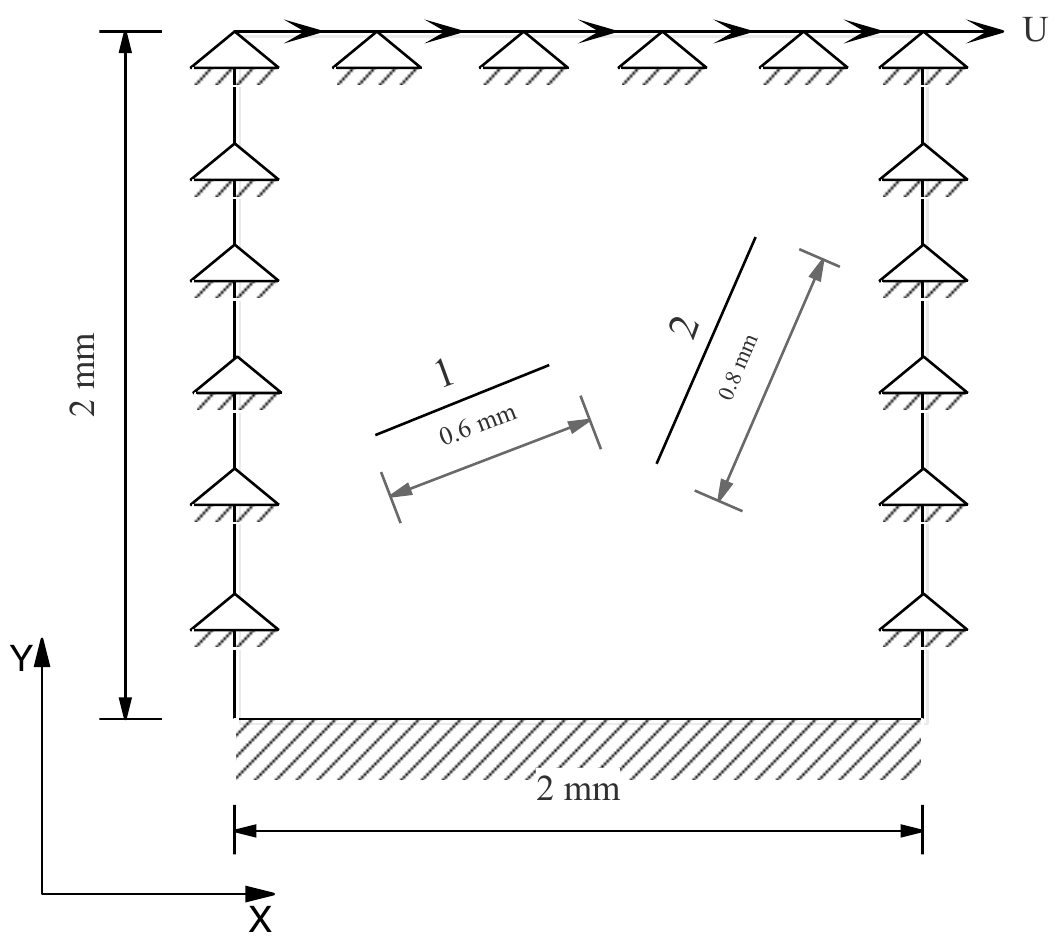}}
\subfigure[random distribution of cracks]{\label{fig:subfig:FiveCrack}
\includegraphics[width=0.4\linewidth]{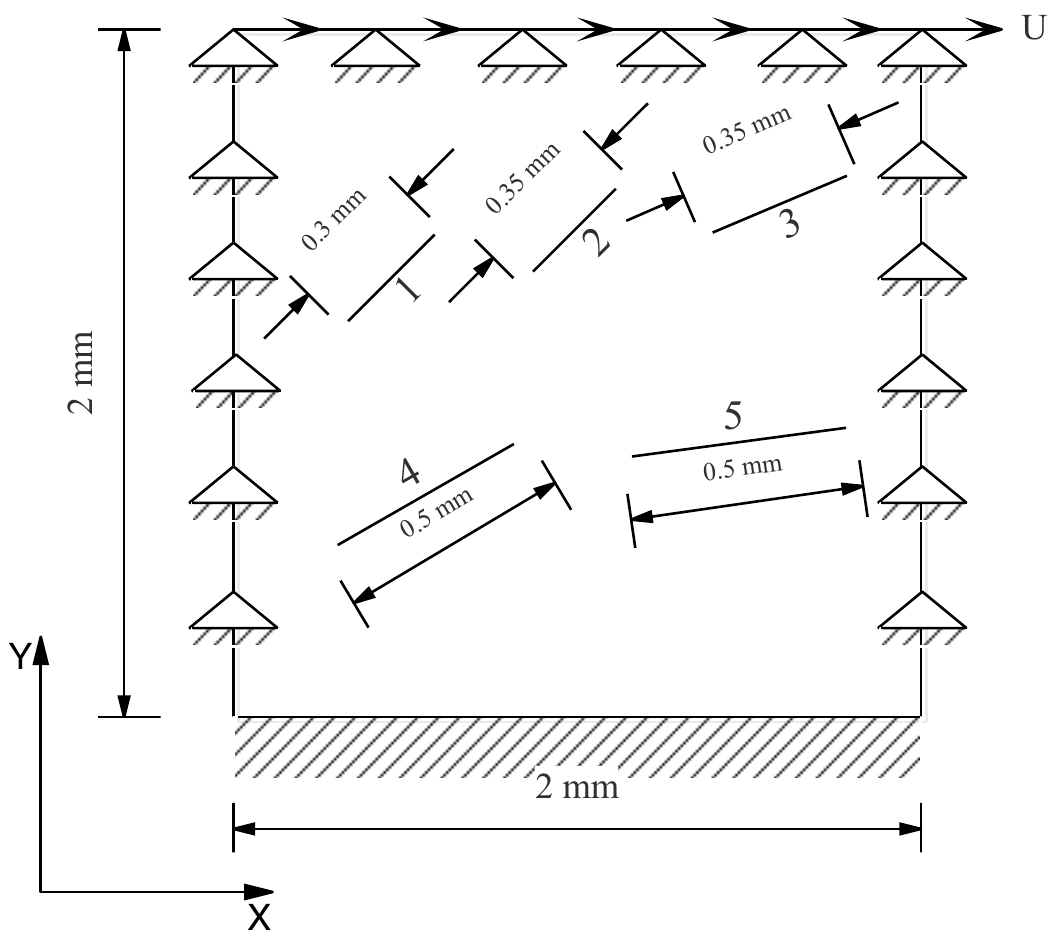}}
\caption{Example 3. Domain and boundary conditions for the shear test with multiple cracks,
(a) two cracks, (b) five cracks.}
\label{fig:sketch_MultipleCracks}
\end{figure}

In this example, we test the modeling of complex crack systems. We consider a square plate of width $2$~mm with two or five cracks. The domain and boundary conditions are shown in Figs. \ref{fig:subfig:TwoCrack} and \ref{fig:subfig:FiveCrack}, respectively. 
For both problems, the bottom edge of the domain is fixed. The top edge is fixed along $y$-direction while a uniform $x$-displacement $U$ is increased with time. The material parameter are the same as in Example 1 except $g_c = 2.7 \times 10^{-4}$~kN/mm for the five-crack problem. 
For the two-crack problem, the lengths of Crack 1 and Crack 2 are $0.6$~mm and $0.8$~mm, with polar angle $9^{\circ}$ and $64.8^{\circ}$, respectively.

\begin{figure} 
\centering 
\subfigure[$U = 6.0 \times 10^{-2}$~mm]{\label{fig:subfig:TOMD_U600}
\includegraphics[width=0.225\linewidth]{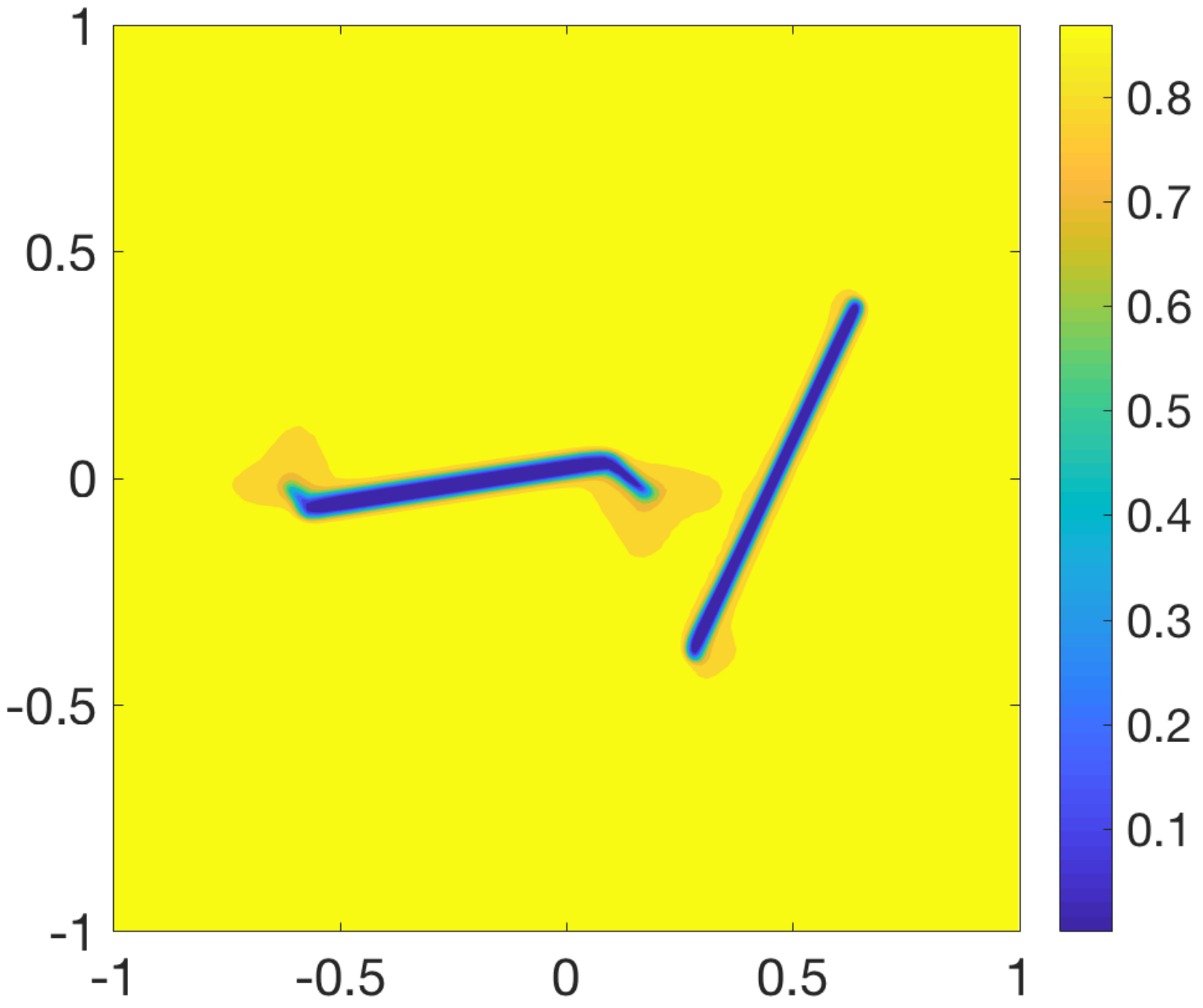}}
\subfigure[$U = 6.2 \times 10^{-2}$~mm]{\label{fig:subfig:TOMD_U620}
\includegraphics[width=0.225\linewidth]{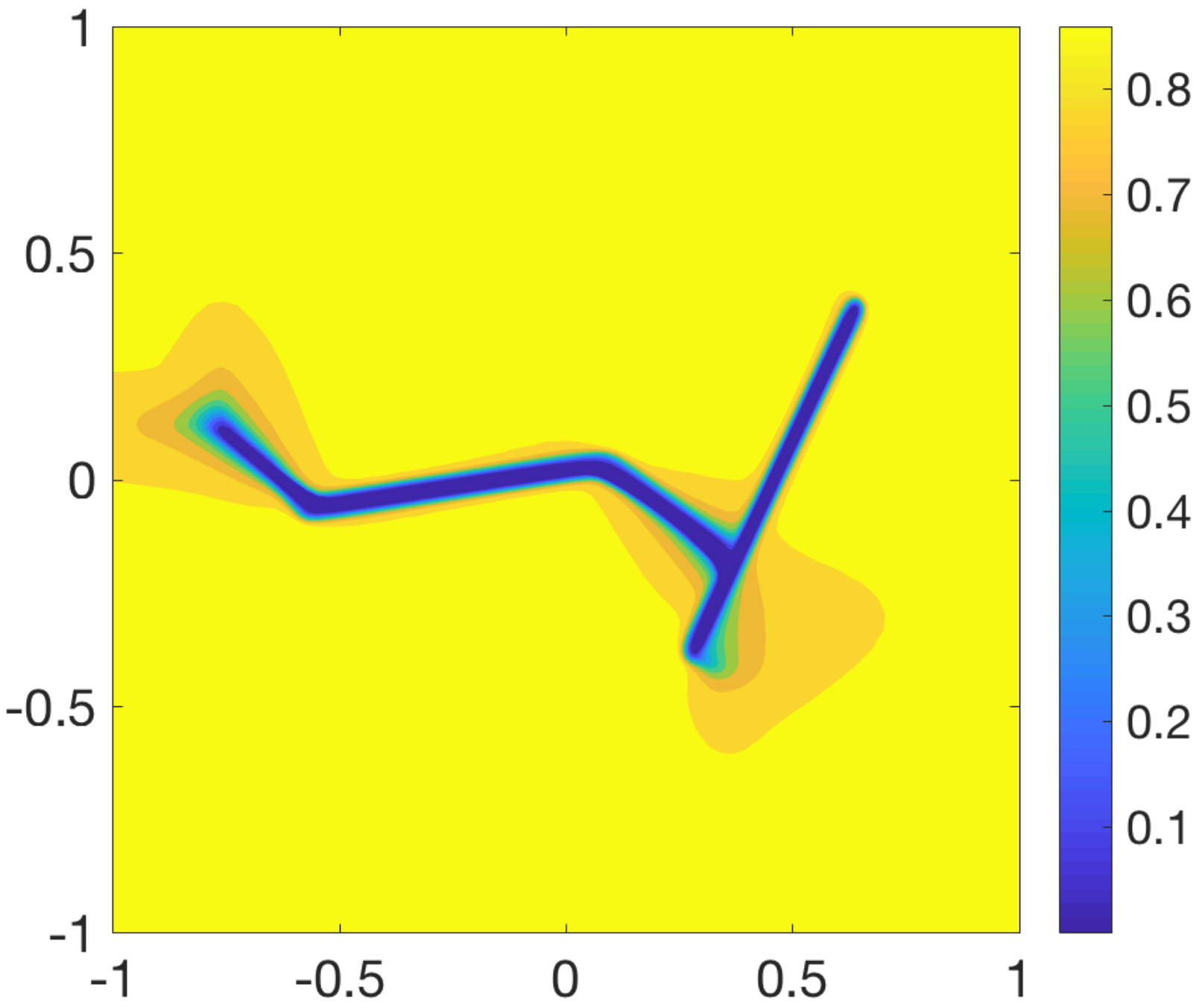}}
\subfigure[$U = 6.4 \times 10^{-2}$~mm]{\label{fig:subfig:TOMD_U640}
\includegraphics[width=0.225\linewidth]{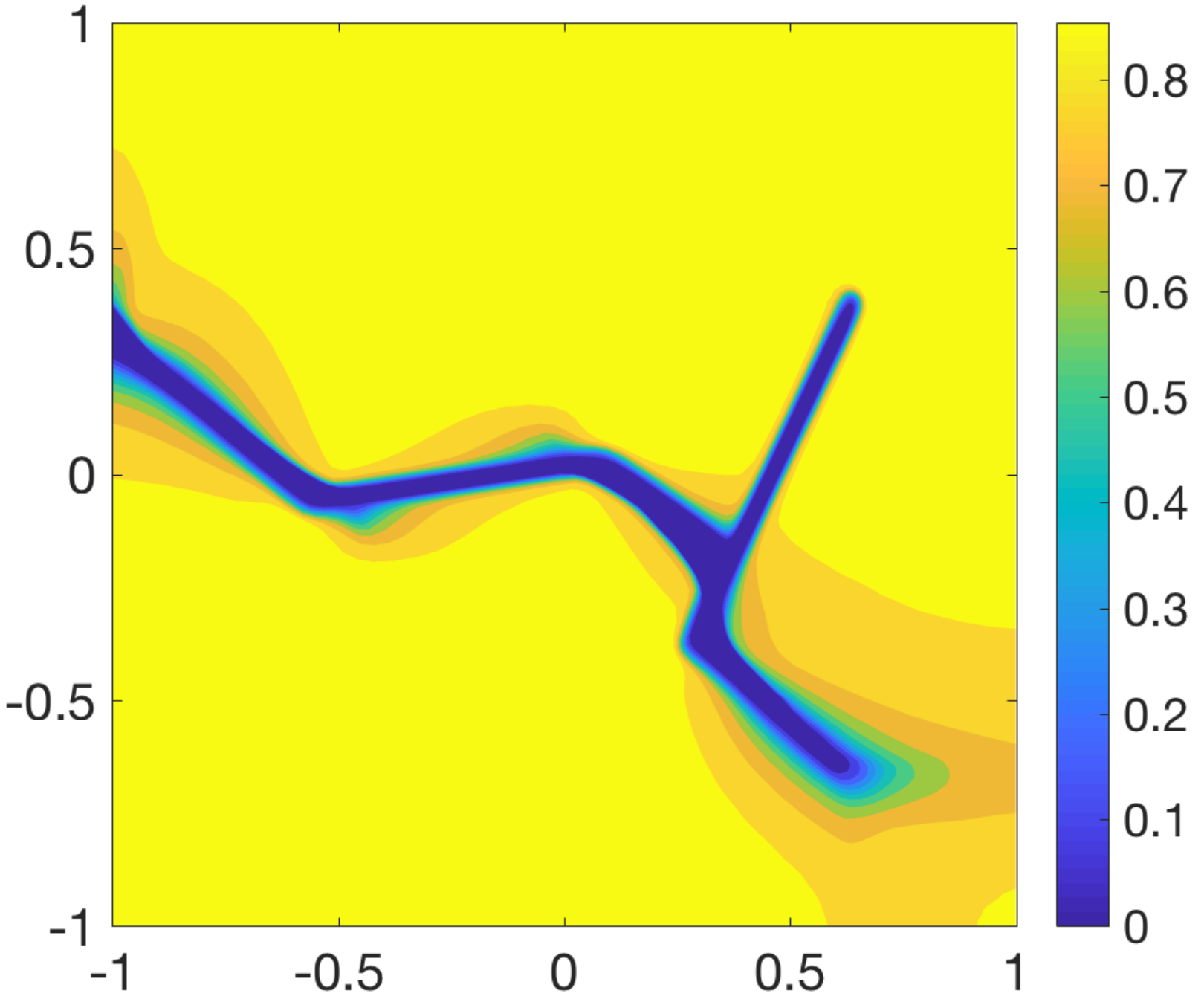}}
\subfigure[$U = 6.6 \times 10^{-2}$~mm]{\label{fig:subfig:TOMD_U660}
\includegraphics[width=0.225\linewidth]{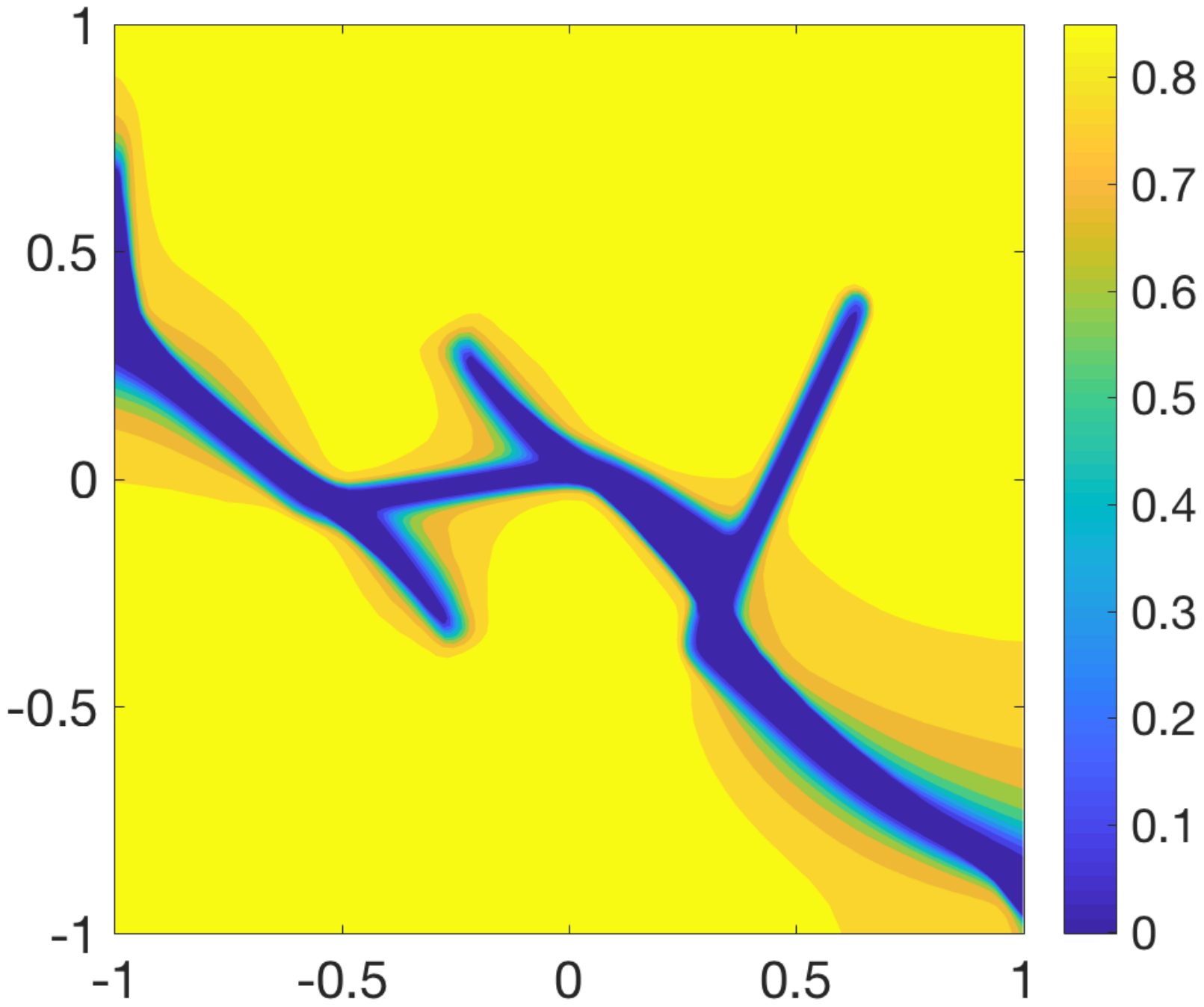}}
\vfill
\subfigure[$U = 6.0 \times 10^{-2}$~mm]{\label{fig:subfig:TOMM_U600}
\includegraphics[width=0.225\linewidth]{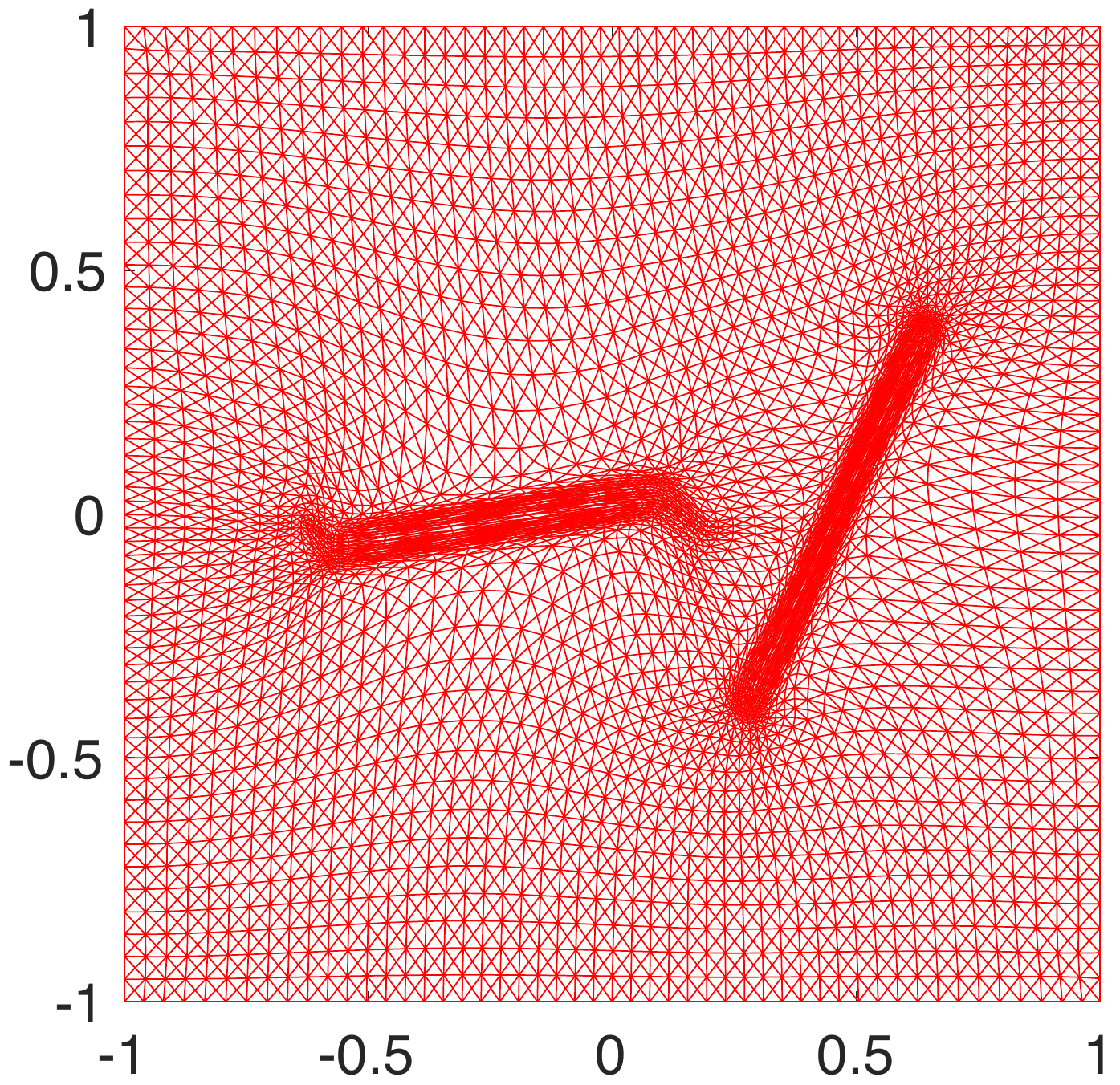}}
\subfigure[$U = 6.2 \times 10^{-2}$~mm]{\label{fig:subfig:TOMM_U620}
\includegraphics[width=0.225\linewidth]{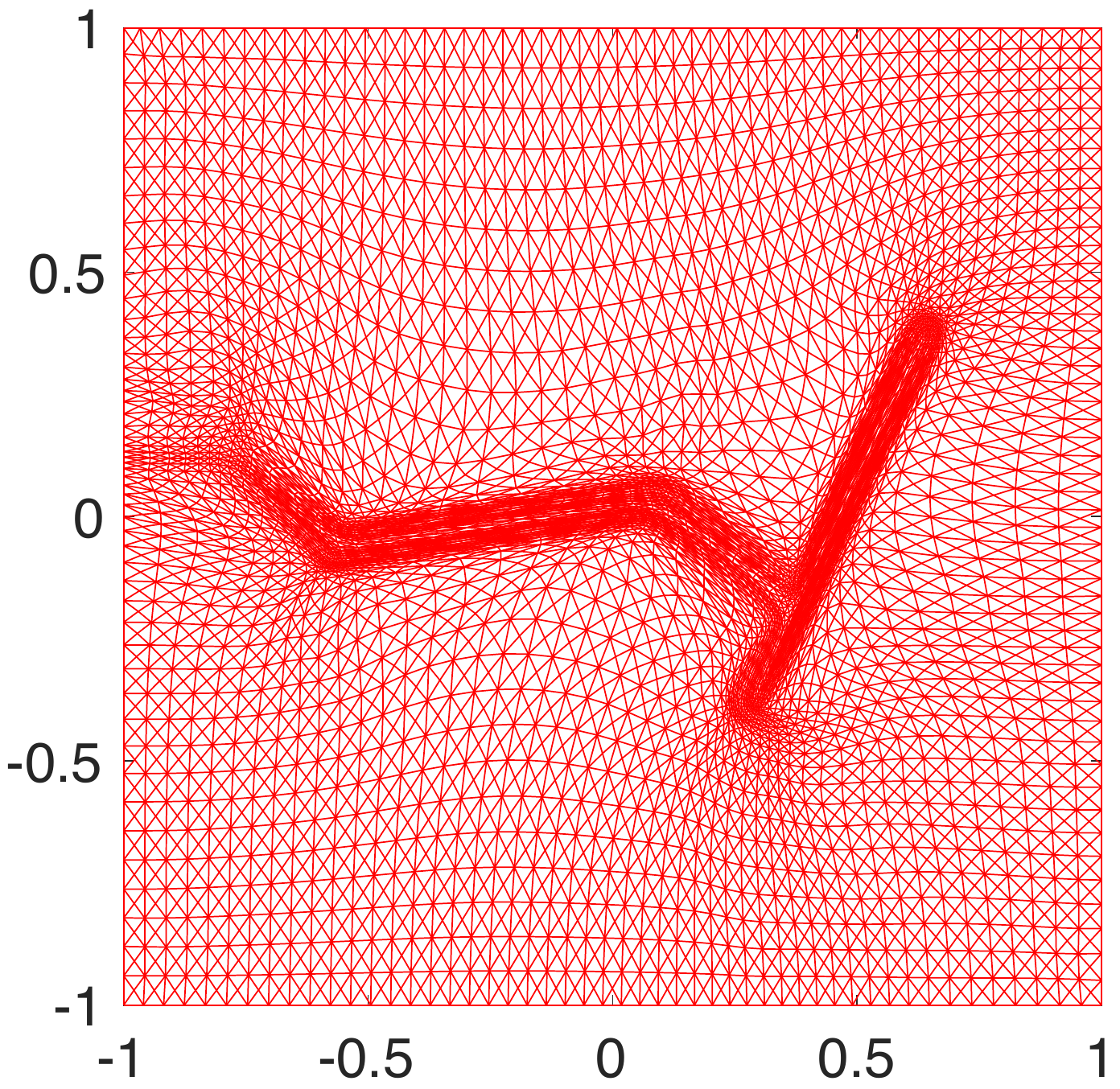}}
\subfigure[$U = 6.4 \times 10^{-2}$~mm]{\label{fig:subfig:TOMM_U640}
\includegraphics[width=0.225\linewidth]{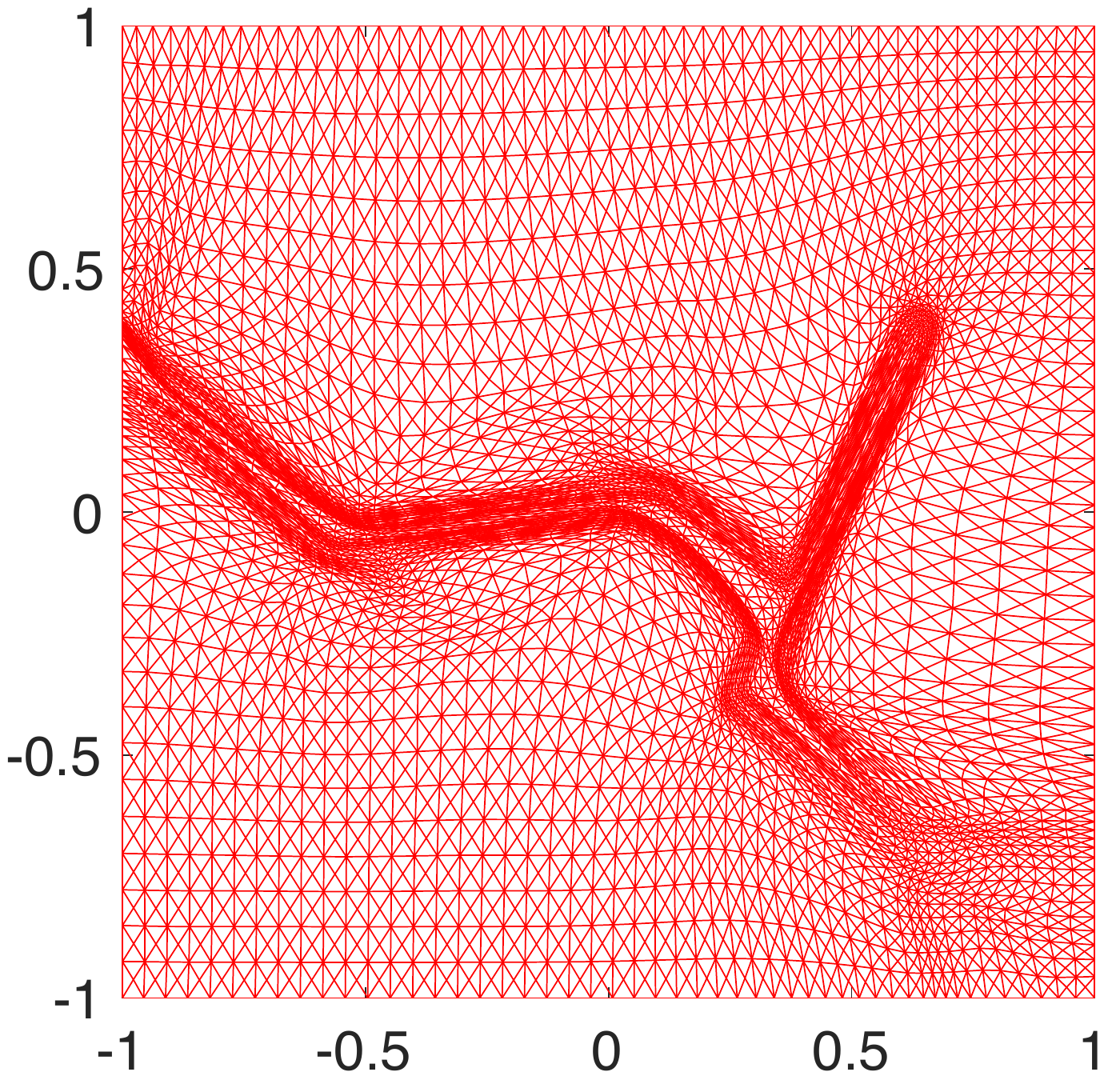}}
\subfigure[$U = 6.6 \times 10^{-2}$~mm]{\label{fig:subfig:TOMM_U660}
\includegraphics[width=0.225\linewidth]{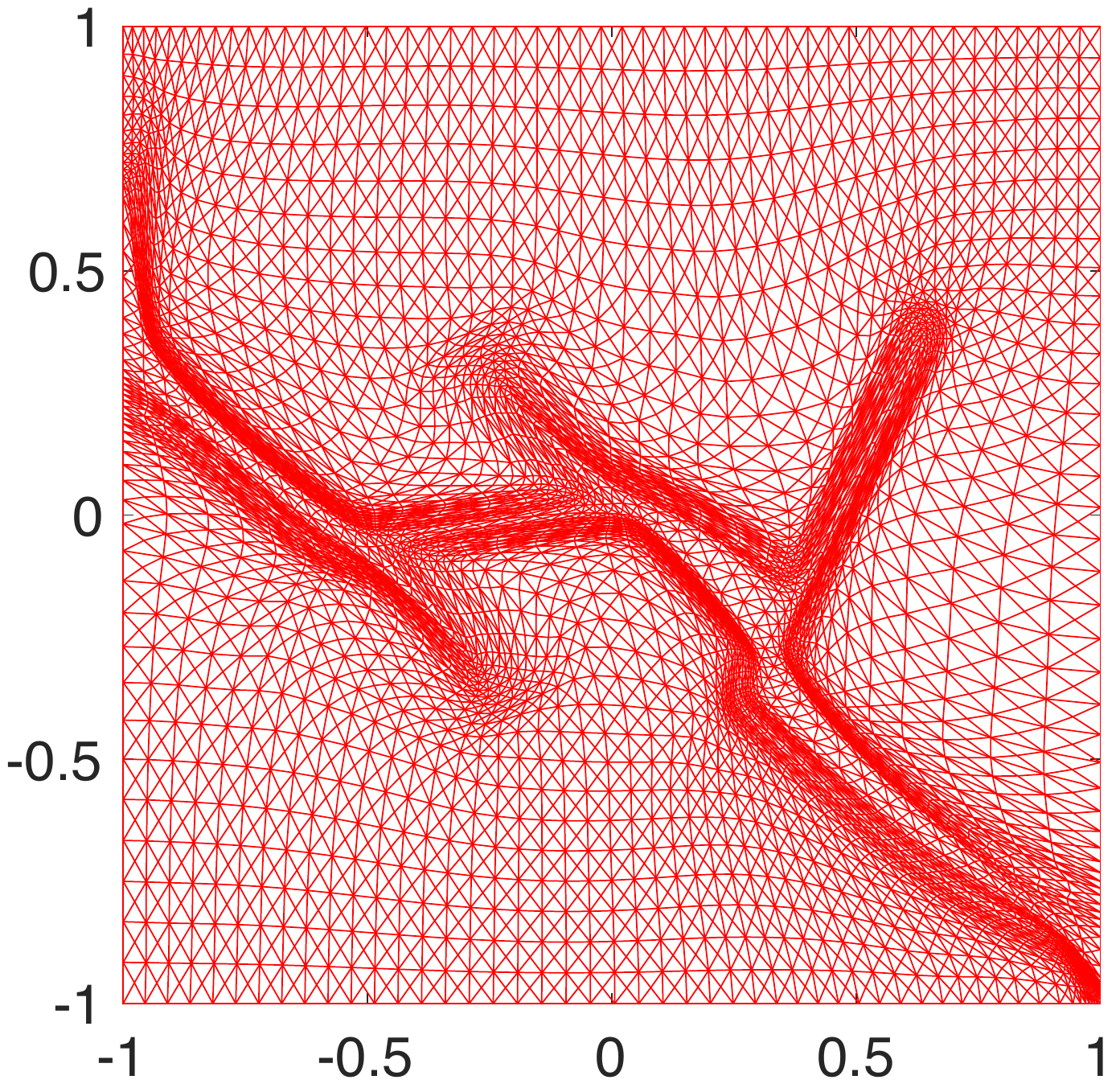}}
\vfill
\subfigure[$U = 6.0 \times 10^{-2}$~mm]{\label{fig:subfig:TOMS_U600}
\includegraphics[width=0.225\linewidth]{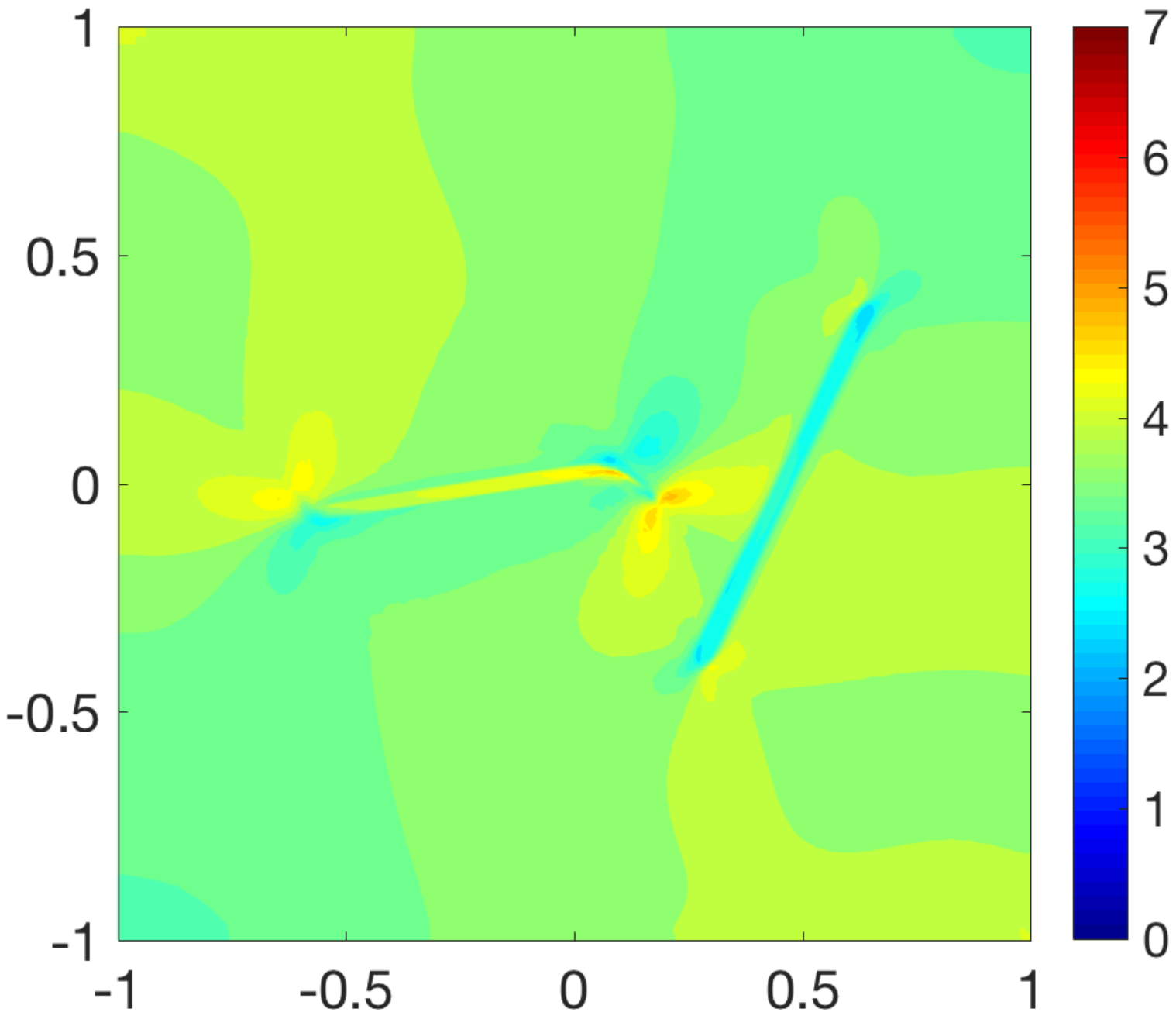}}
\subfigure[$U = 6.2 \times 10^{-2}$~mm]{\label{fig:subfig:TOMS_U620}
\includegraphics[width=0.225\linewidth]{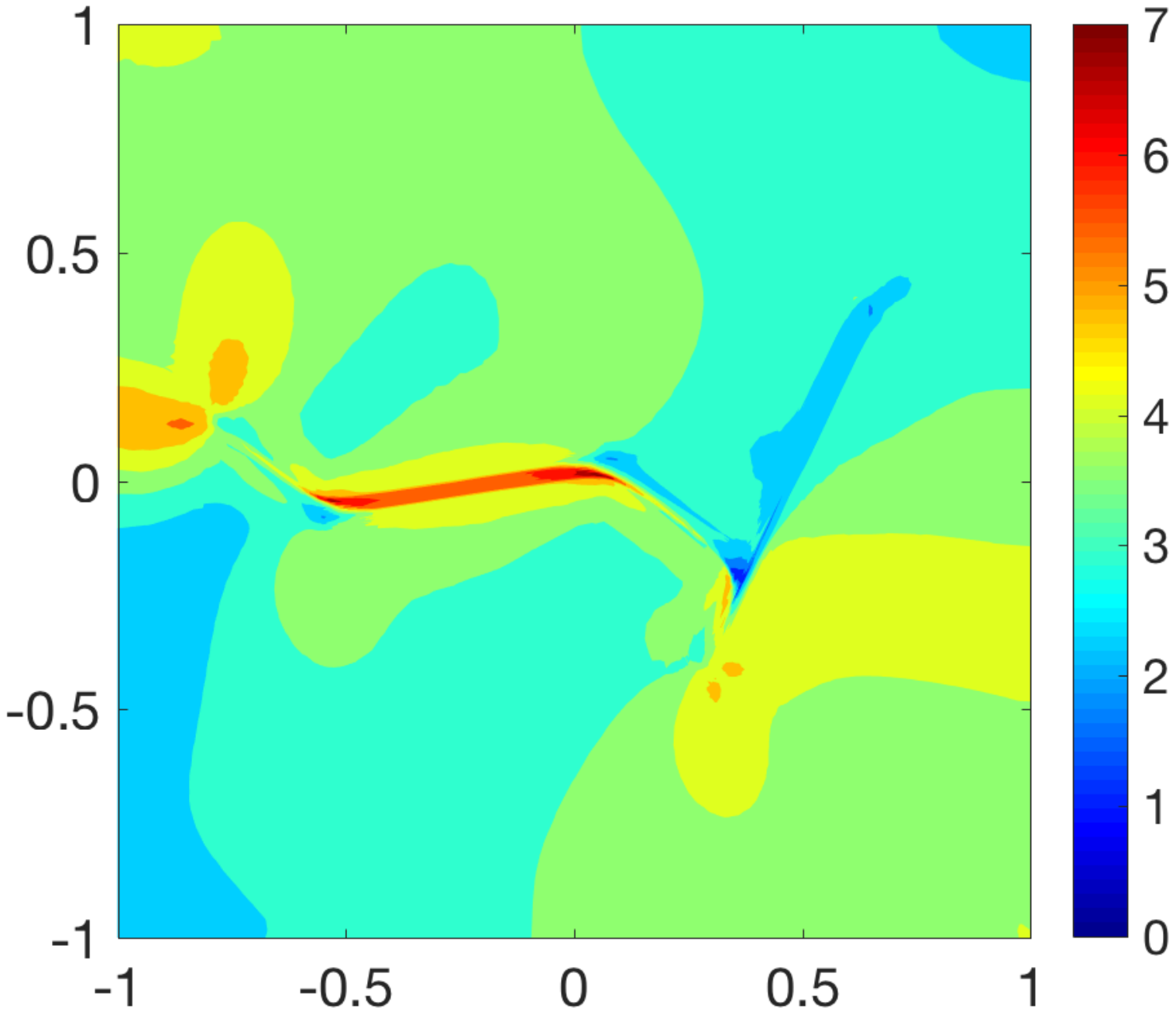}}
\subfigure[$U = 6.4 \times 10^{-2}$~mm]{\label{fig:subfig:TOMS_U640}
\includegraphics[width=0.225\linewidth]{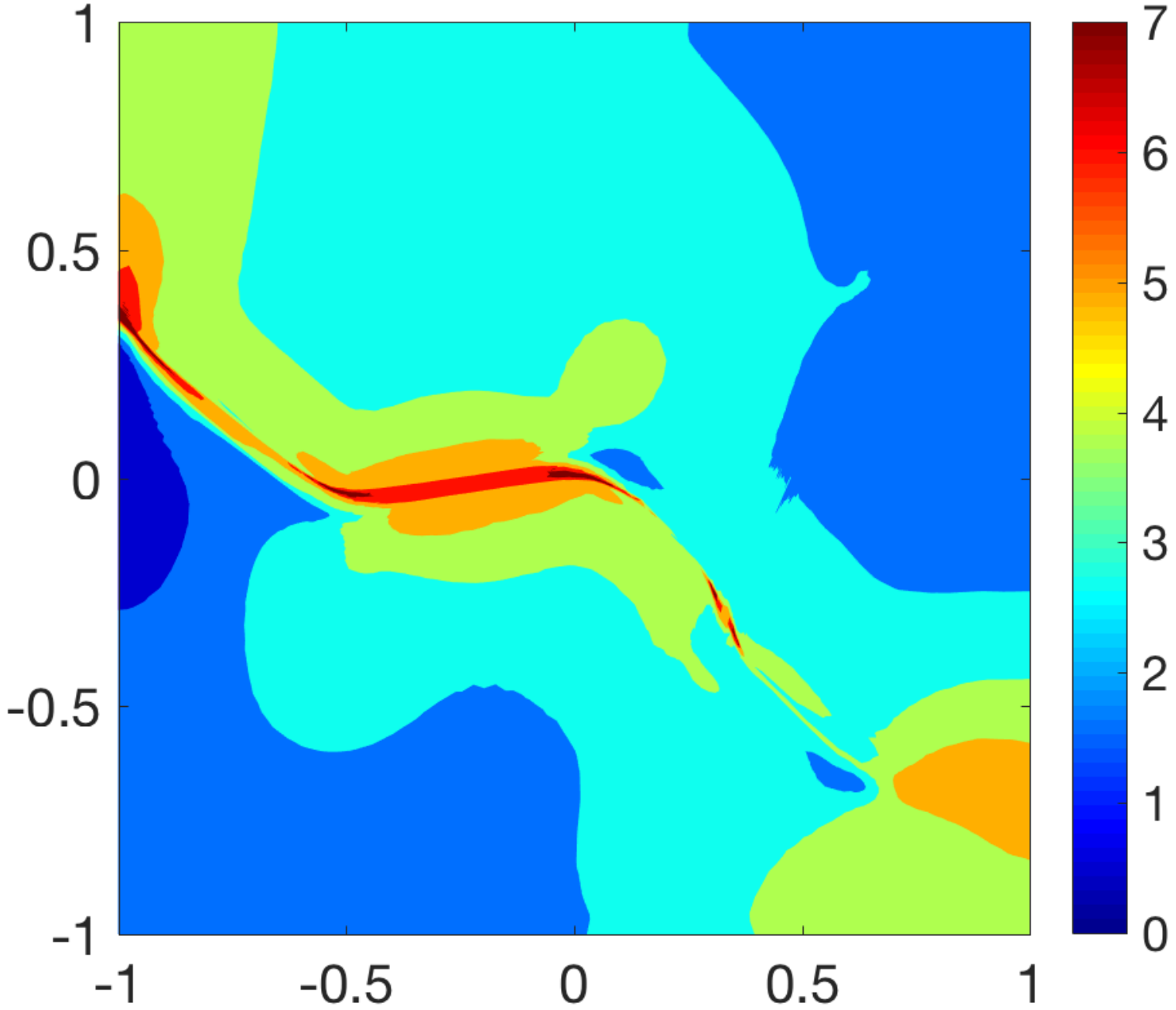}}
\subfigure[$U = 6.6 \times 10^{-2}$~mm]{\label{fig:subfig:TOMS_U660}
\includegraphics[width=0.225\linewidth]{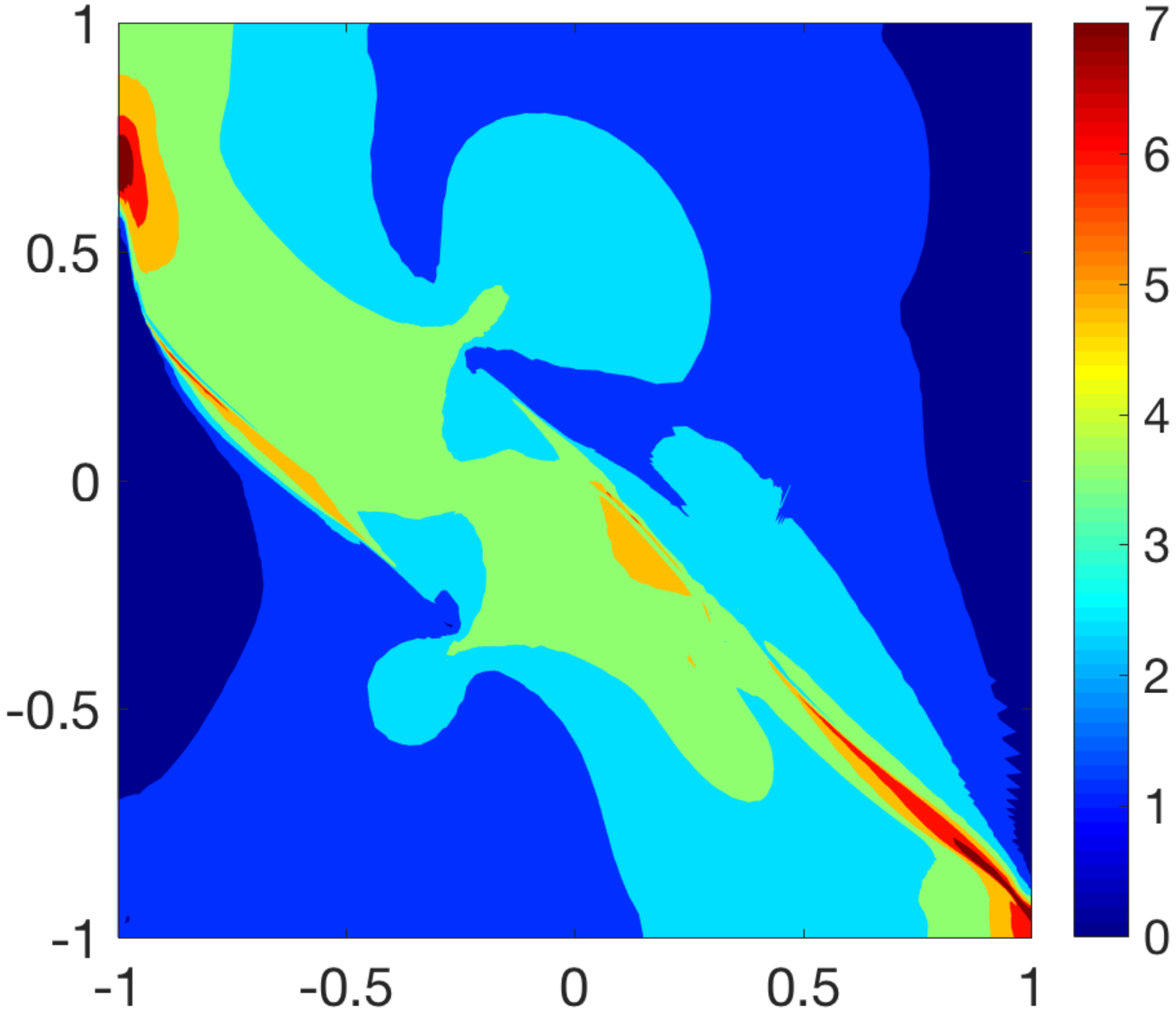}}
\caption{Example 3. The mesh and contours of the phase-field and von Mises stress distribution during
crack evolution for the two-crack shear test with $l = 0.00375$~mm, $N = 10,000\;  (51\times51)$.
(spectral decomposition with original crack boundary conditions)}
\label{fig:SSOBC}
\end{figure}

\begin{figure} 
\centering 
\subfigure[$U = 2.8 \times 10^{-2}$~mm]{\label{fig:subfig:TMMD_U280}
\includegraphics[width=0.225\linewidth]{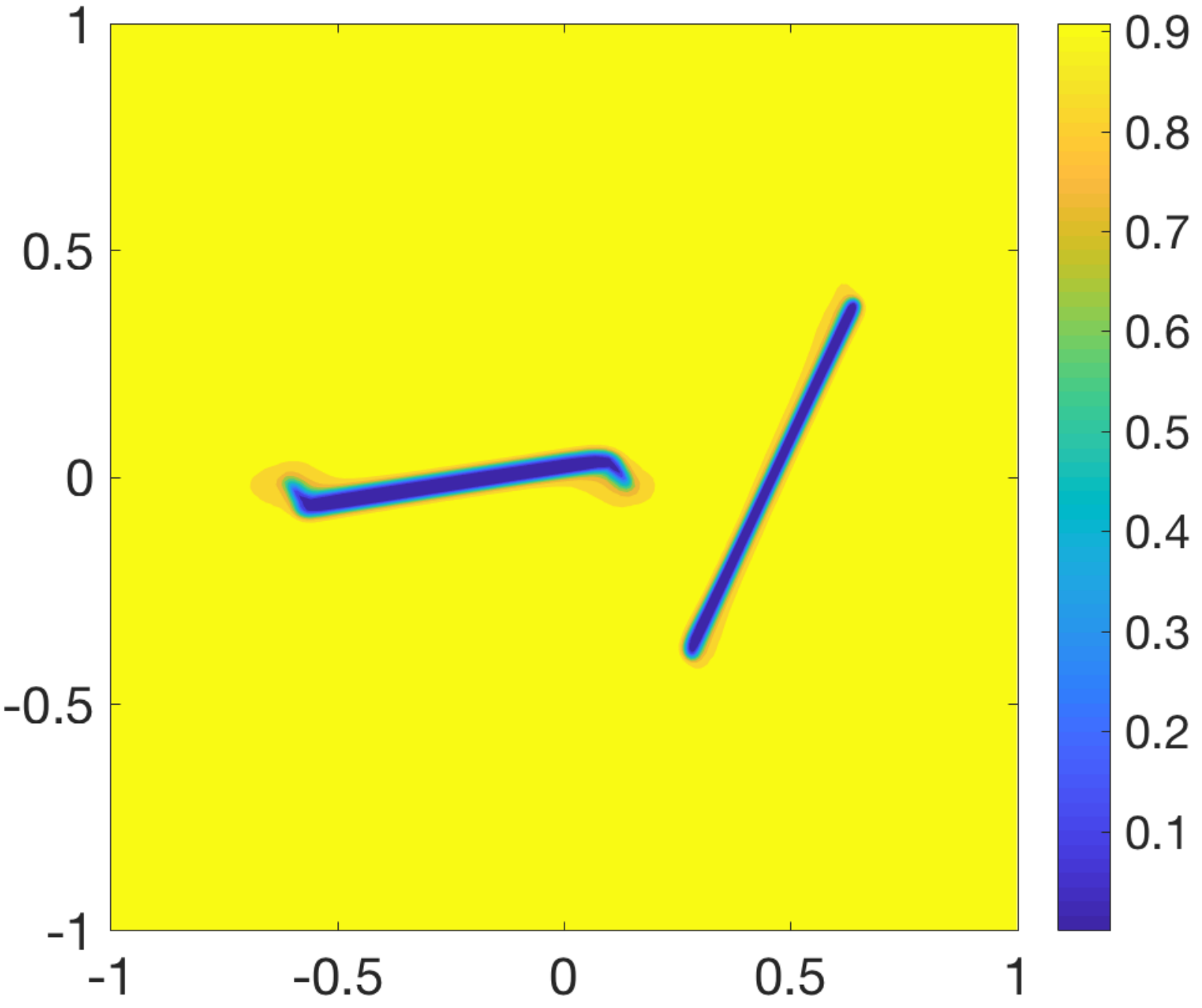}}
\subfigure[$U = 3.0 \times 10^{-2}$~mm]{\label{fig:subfig:TMMD_U300}
\includegraphics[width=0.225\linewidth]{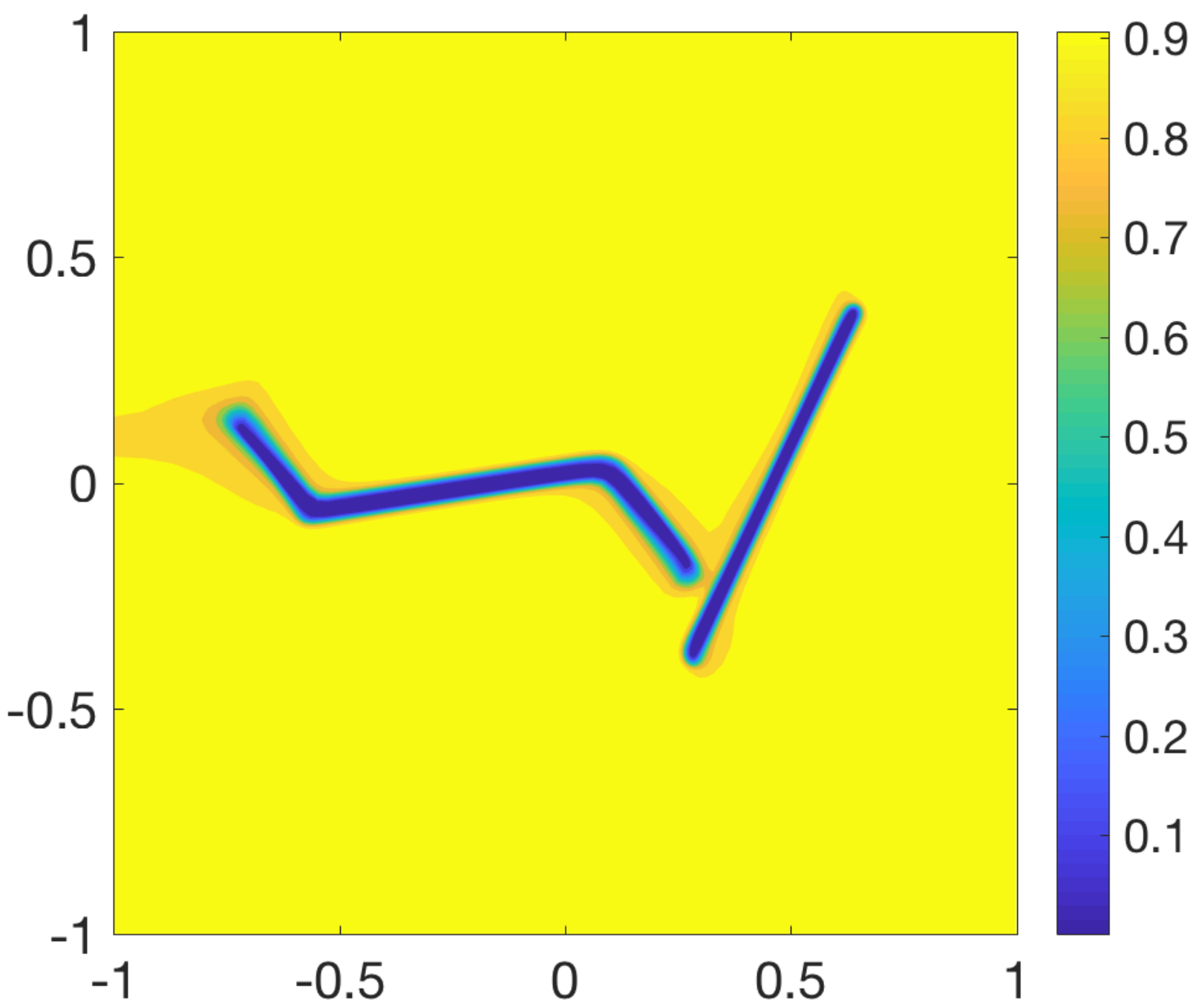}}
\subfigure[$U = 3.4 \times 10^{-2}$~mm]{\label{fig:subfig:TMMD_U340}
\includegraphics[width=0.225\linewidth]{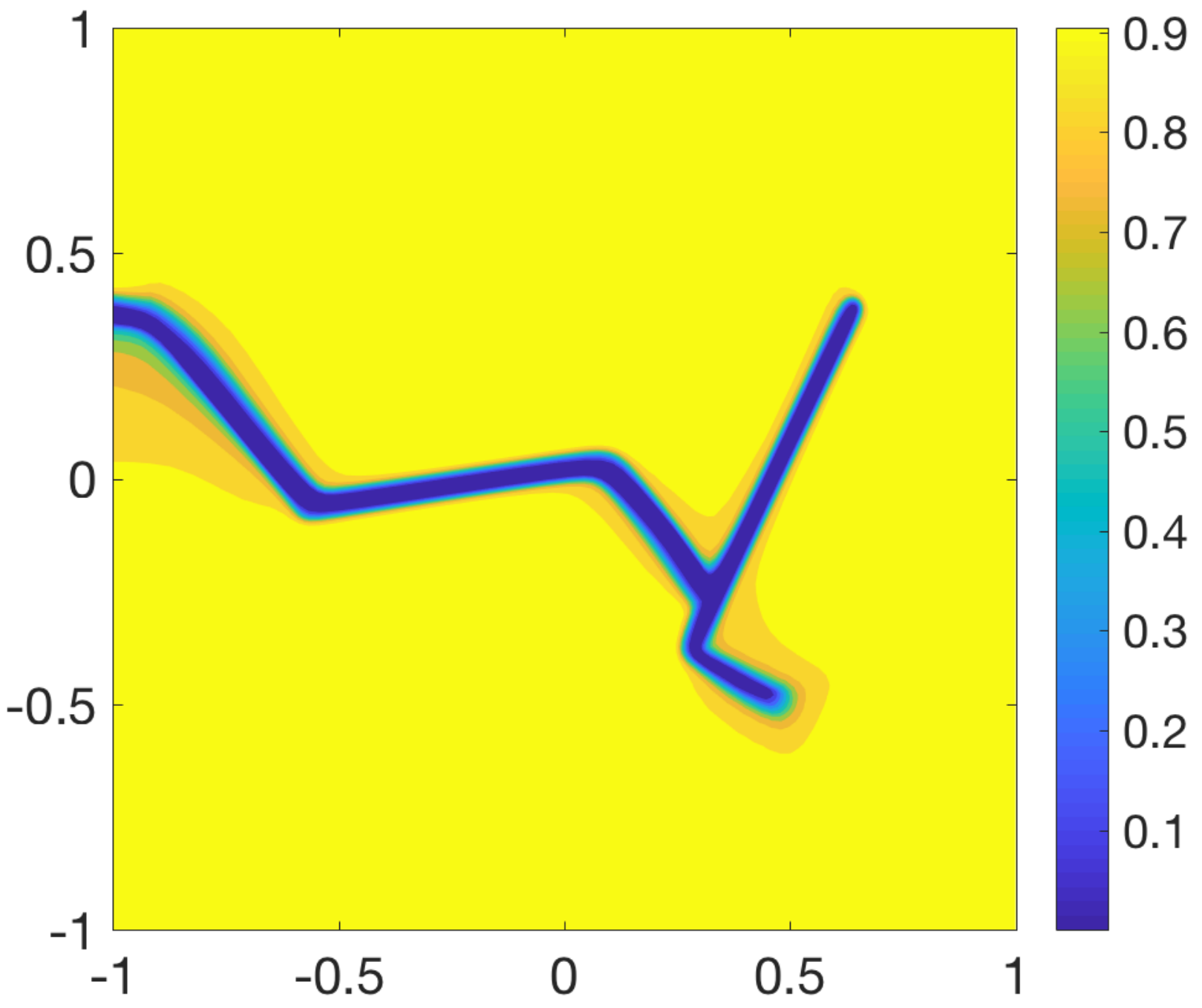}}
\subfigure[$U = 4.2 \times 10^{-2}$~mm]{\label{fig:subfig:TMMD_U420}
\includegraphics[width=0.225\linewidth]{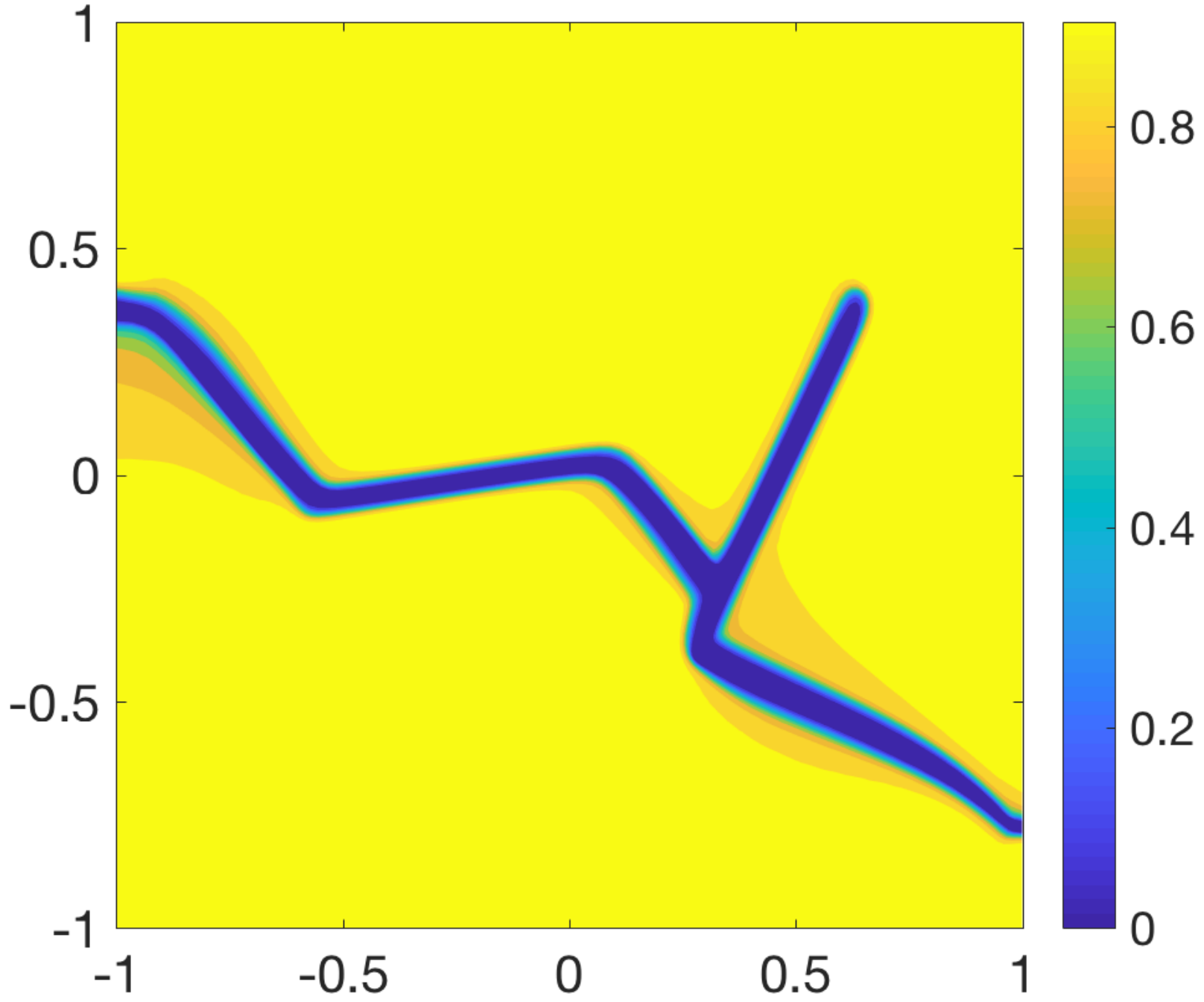}}
\vfill
\subfigure[$U = 2.8 \times 10^{-2}$~mm]{\label{fig:subfig:TMMM_U280}
\includegraphics[width=0.225\linewidth]{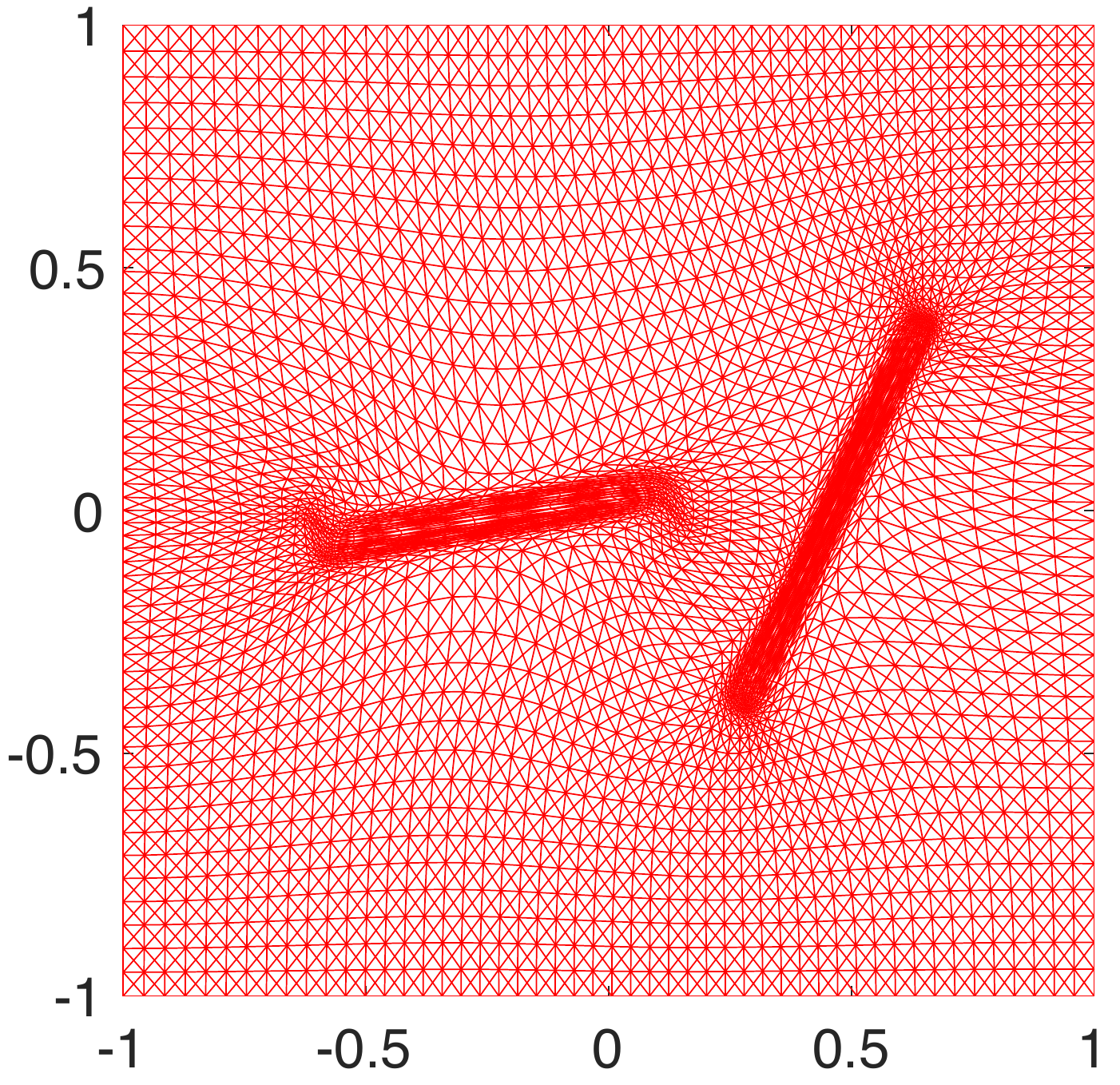}}
\subfigure[$U = 3.0 \times 10^{-2}$~mm]{\label{fig:subfig:TMMM_U300}
\includegraphics[width=0.225\linewidth]{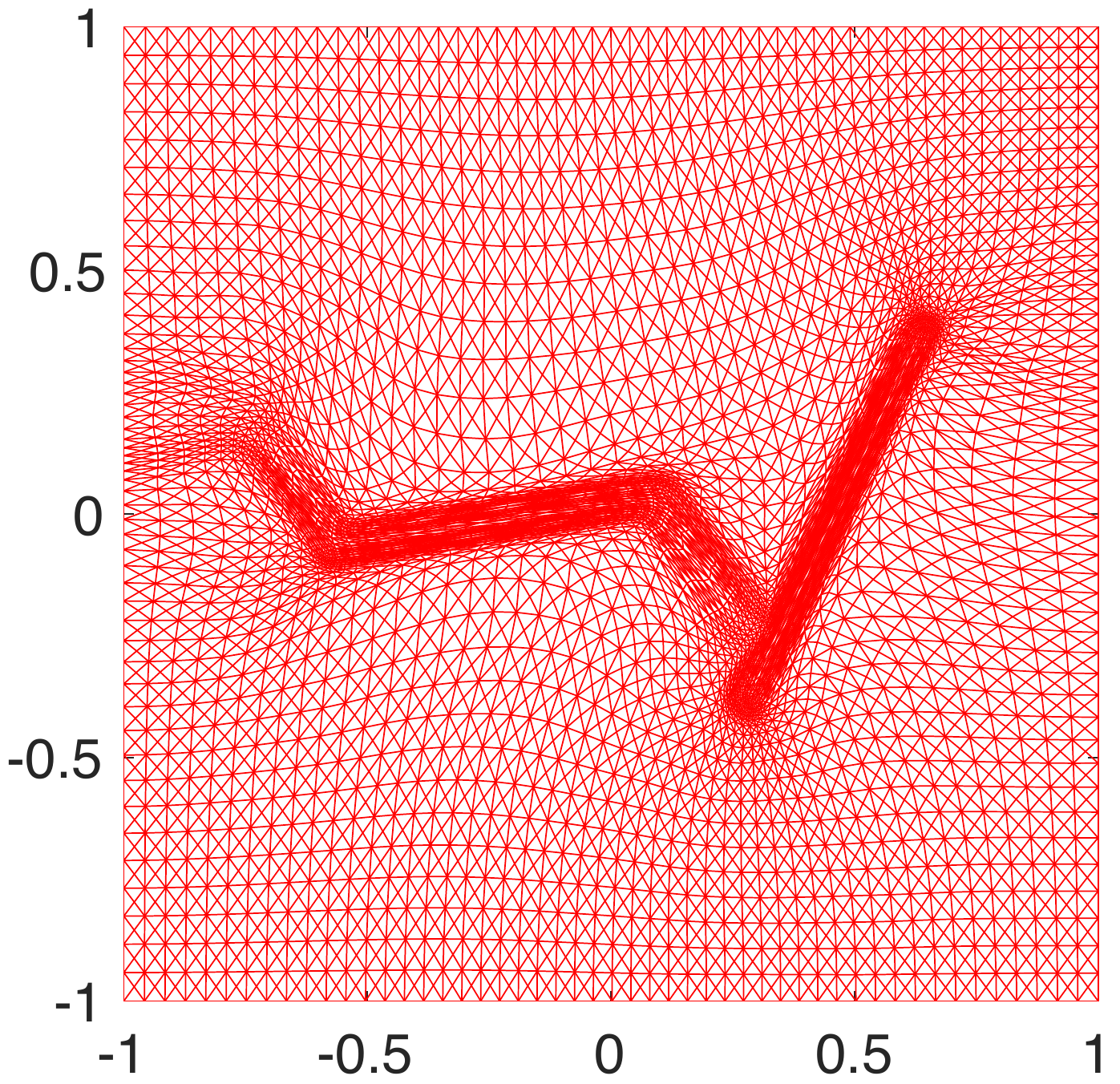}}
\subfigure[$U = 3.4 \times 10^{-2}$~mm]{\label{fig:subfig:TMMM_U340}
\includegraphics[width=0.225\linewidth]{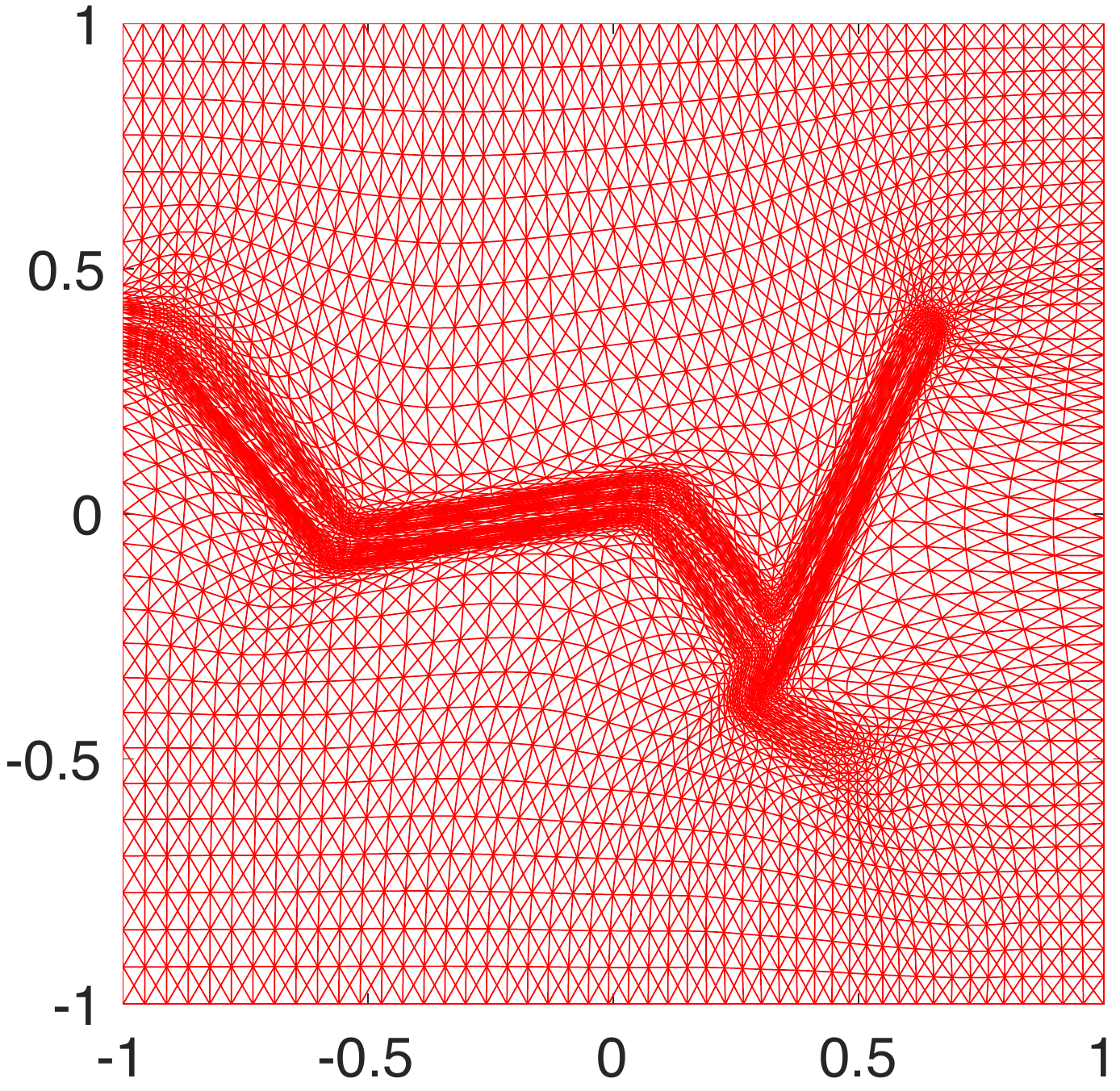}}
\subfigure[$U = 4.2 \times 10^{-2}$~mm]{\label{fig:subfig:TMMM_U420}
\includegraphics[width=0.225\linewidth]{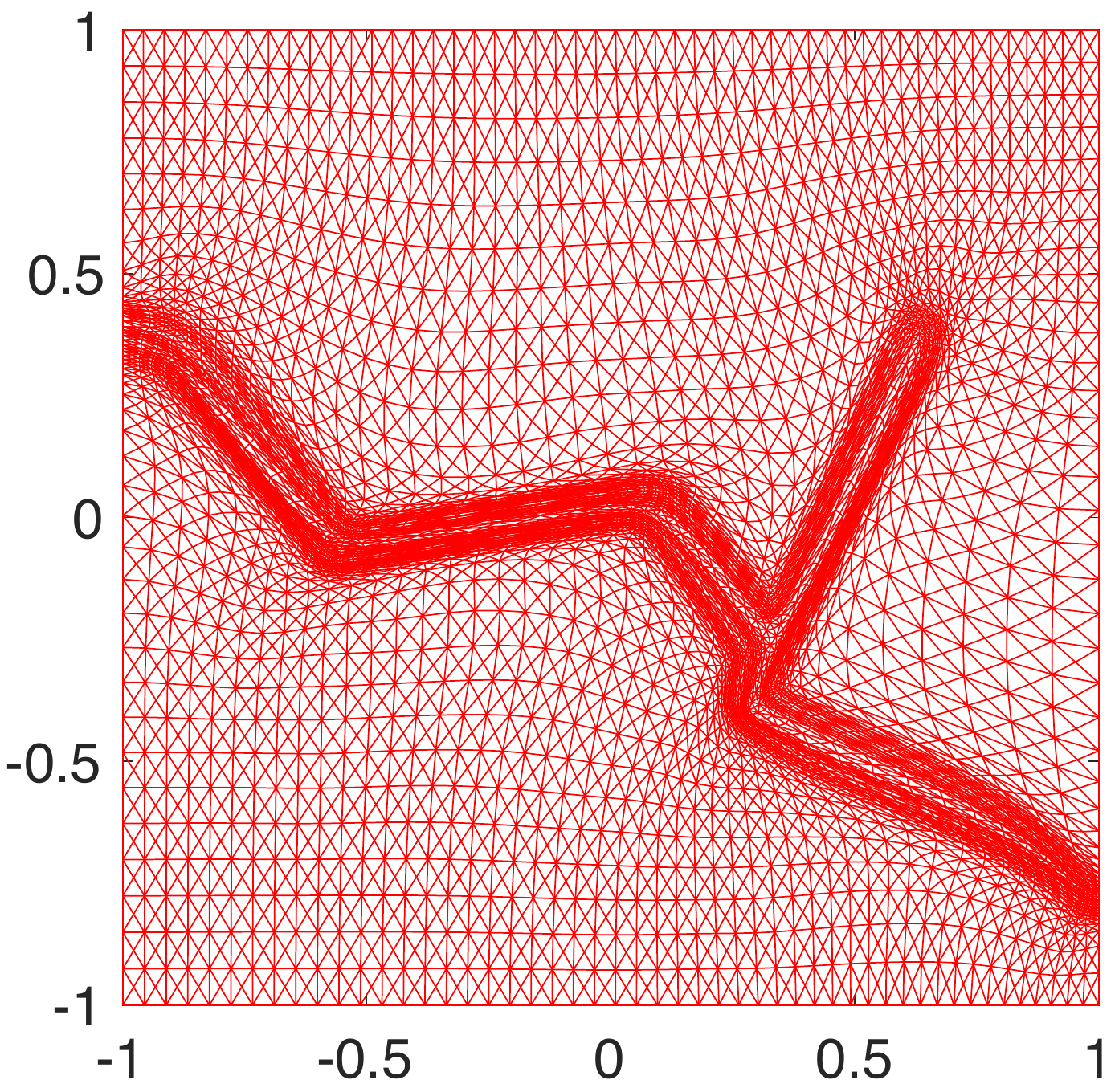}}
\vfill
\subfigure[$U = 2.8 \times 10^{-2}$~mm]{\label{fig:subfig:TMMS_U280}
\includegraphics[width=0.225\linewidth]{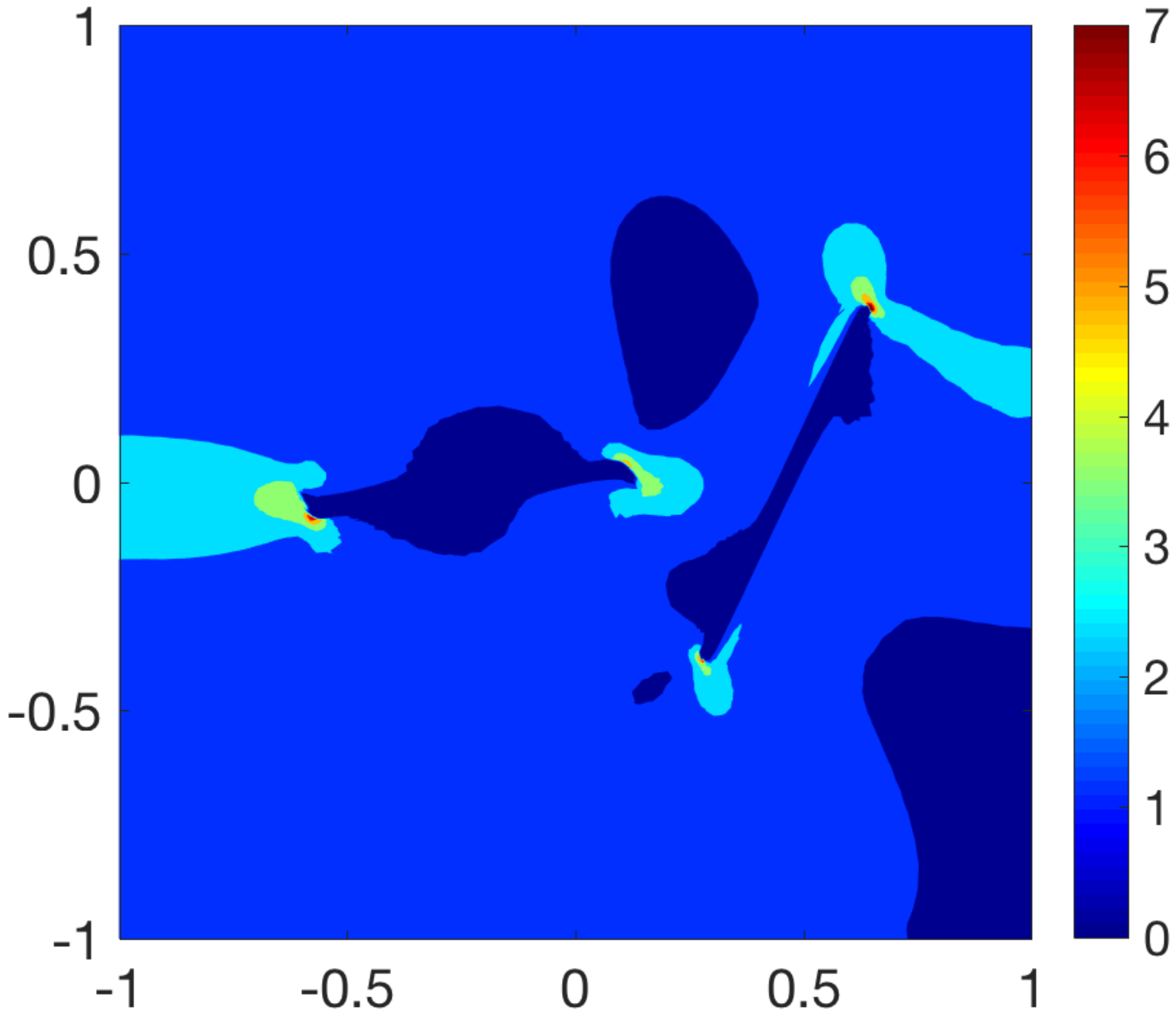}}
\subfigure[$U = 3.0 \times 10^{-2}$~mm]{\label{fig:subfig:TMMS_U300}
\includegraphics[width=0.225\linewidth]{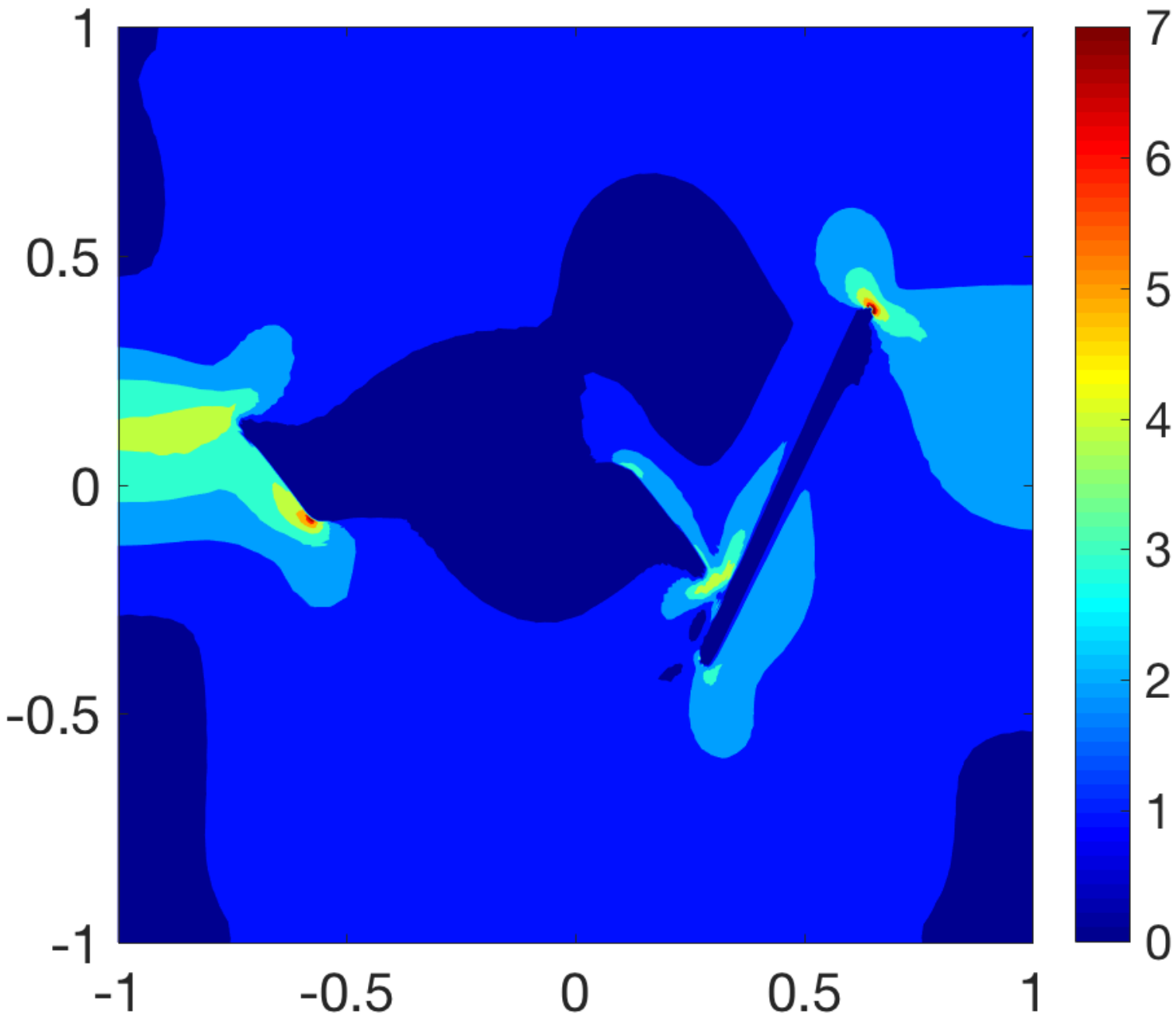}}
\subfigure[$U = 3.4 \times 10^{-2}$~mm]{\label{fig:subfig:TMMS_U340}
\includegraphics[width=0.225\linewidth]{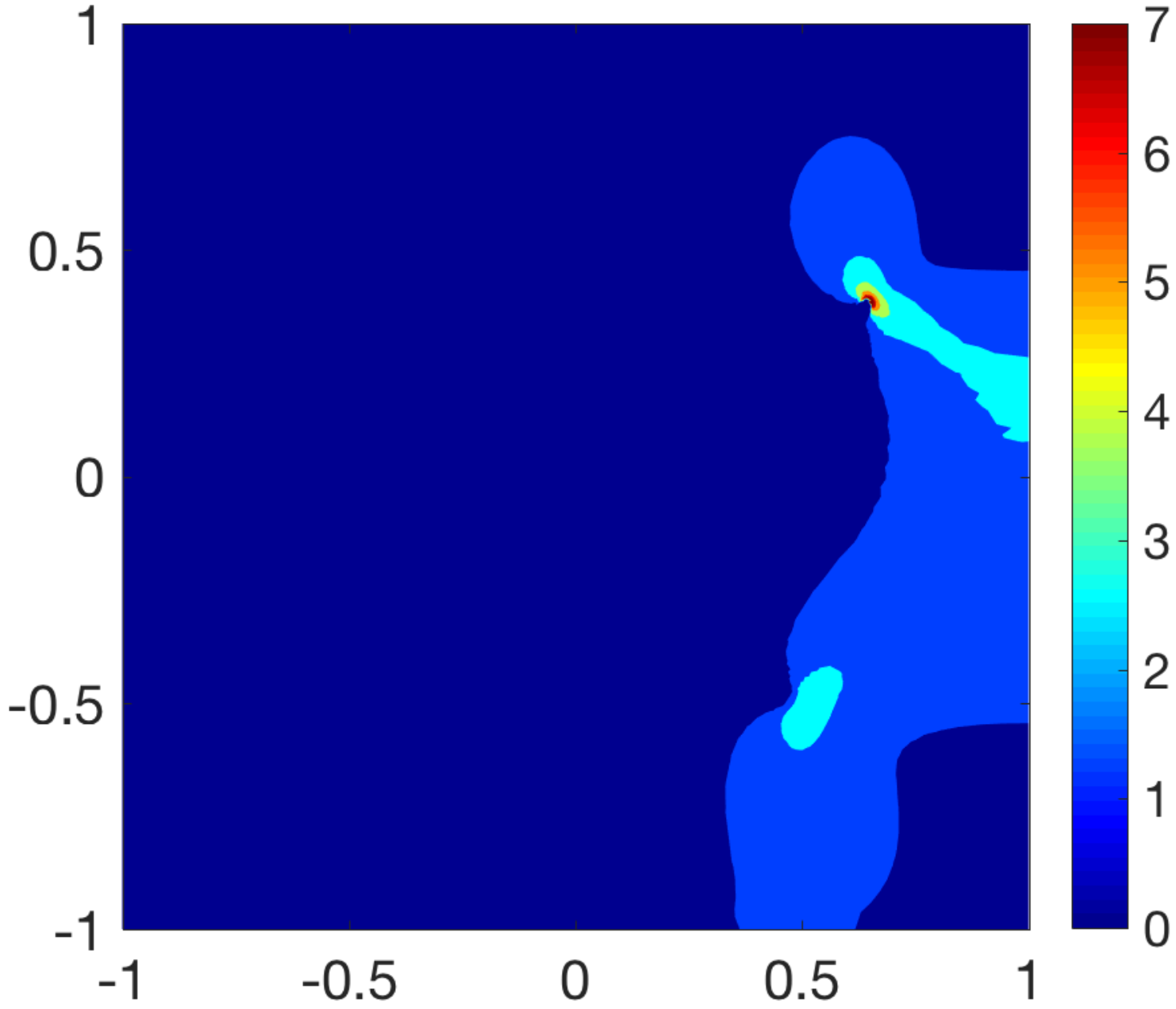}}
\subfigure[$U = 4.2 \times 10^{-2}$~mm]{\label{fig:subfig:TMMS_U420}
\includegraphics[width=0.225\linewidth]{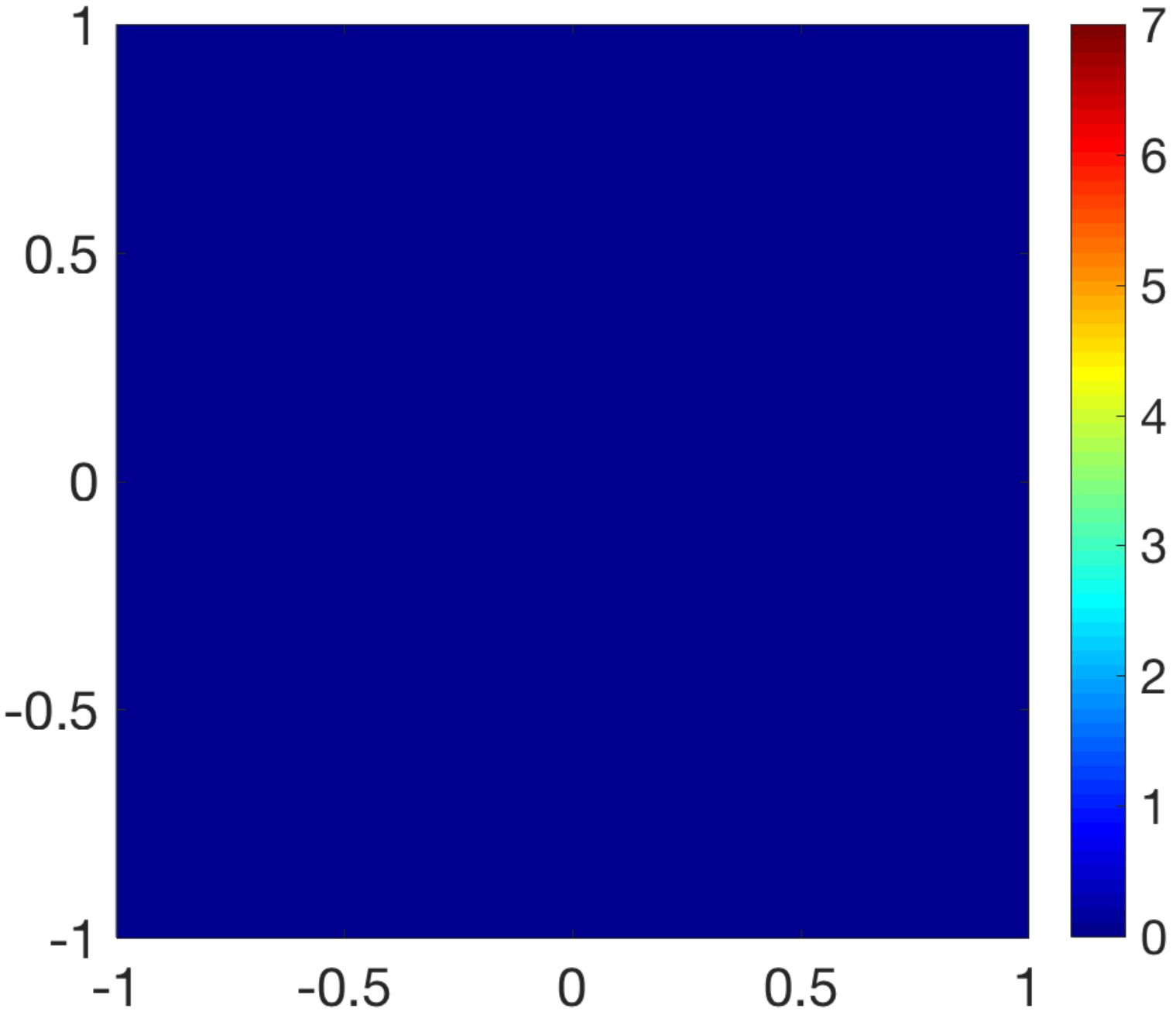}}
\caption{Example 3. The mesh and contours of the phase-field and von Mises stress distribution during
crack evolution for the two-crack shear test with $l = 0.00375$~mm, $N = 10,000\; (51\times51)$.
(spectral decomposition with ItCBC)}
\label{fig:SSMBC}
\end{figure}

\begin{figure} 
\centering 
\subfigure[$U = 2.4 \times 10^{-2}$~mm]{\label{fig:subfig:SMAD_U240}
\includegraphics[width=0.225\linewidth]{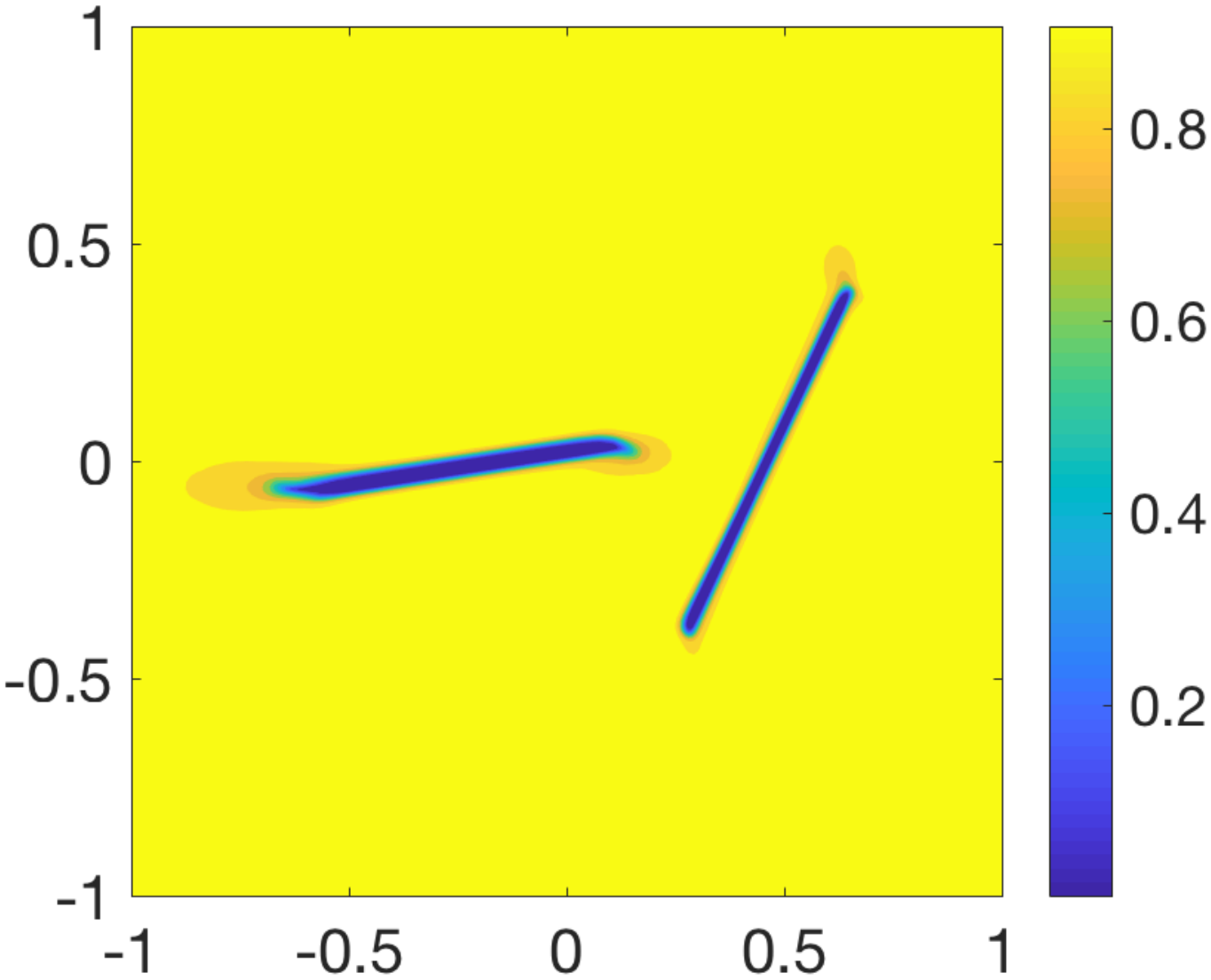}}
\subfigure[$U = 2.6 \times 10^{-2}$~mm]{\label{fig:subfig:SMAD_U260}
\includegraphics[width=0.225\linewidth]{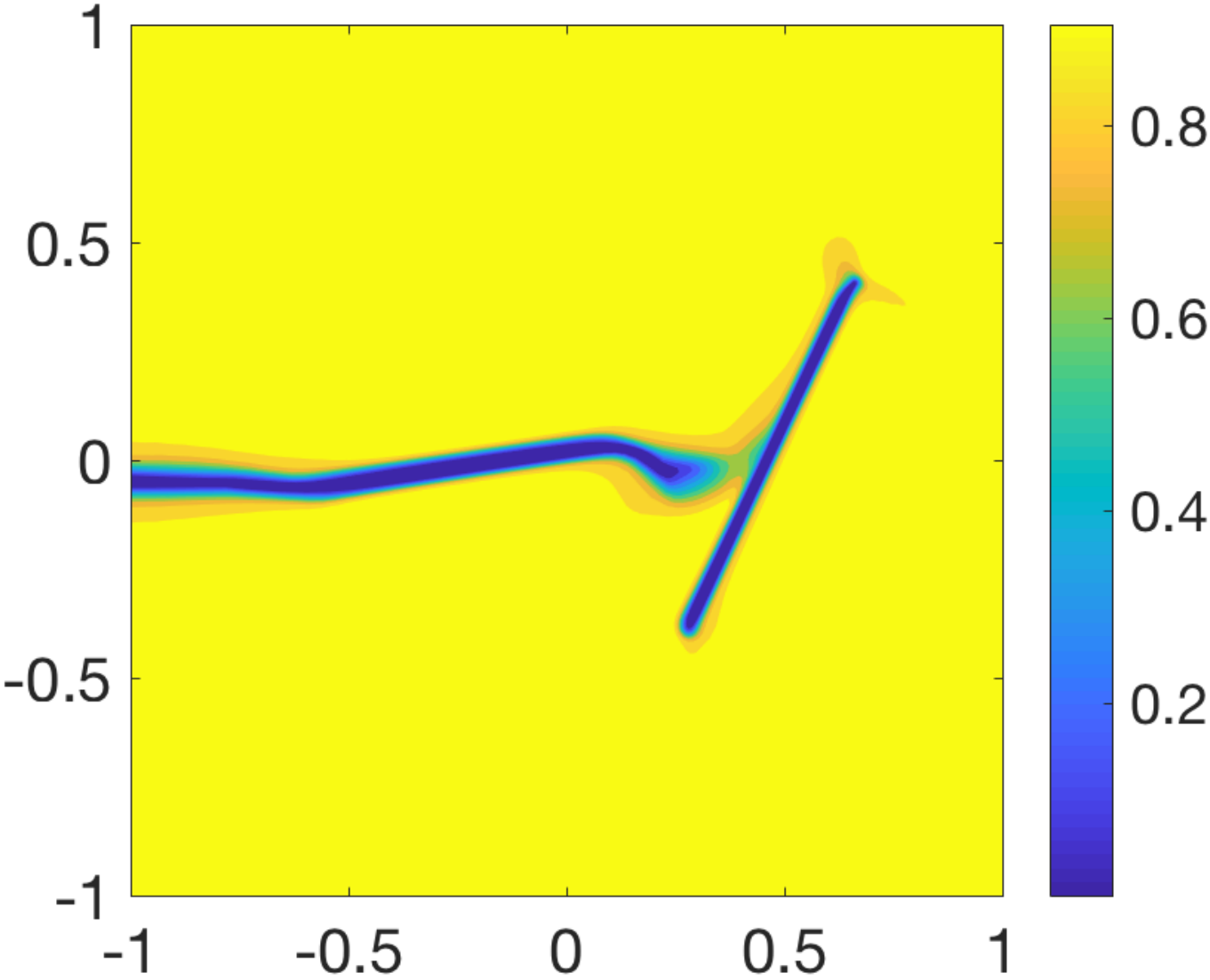}}
\subfigure[$U = 2.8 \times 10^{-2}$~mm]{\label{fig:subfig:SMAD_U280}
\includegraphics[width=0.225\linewidth]{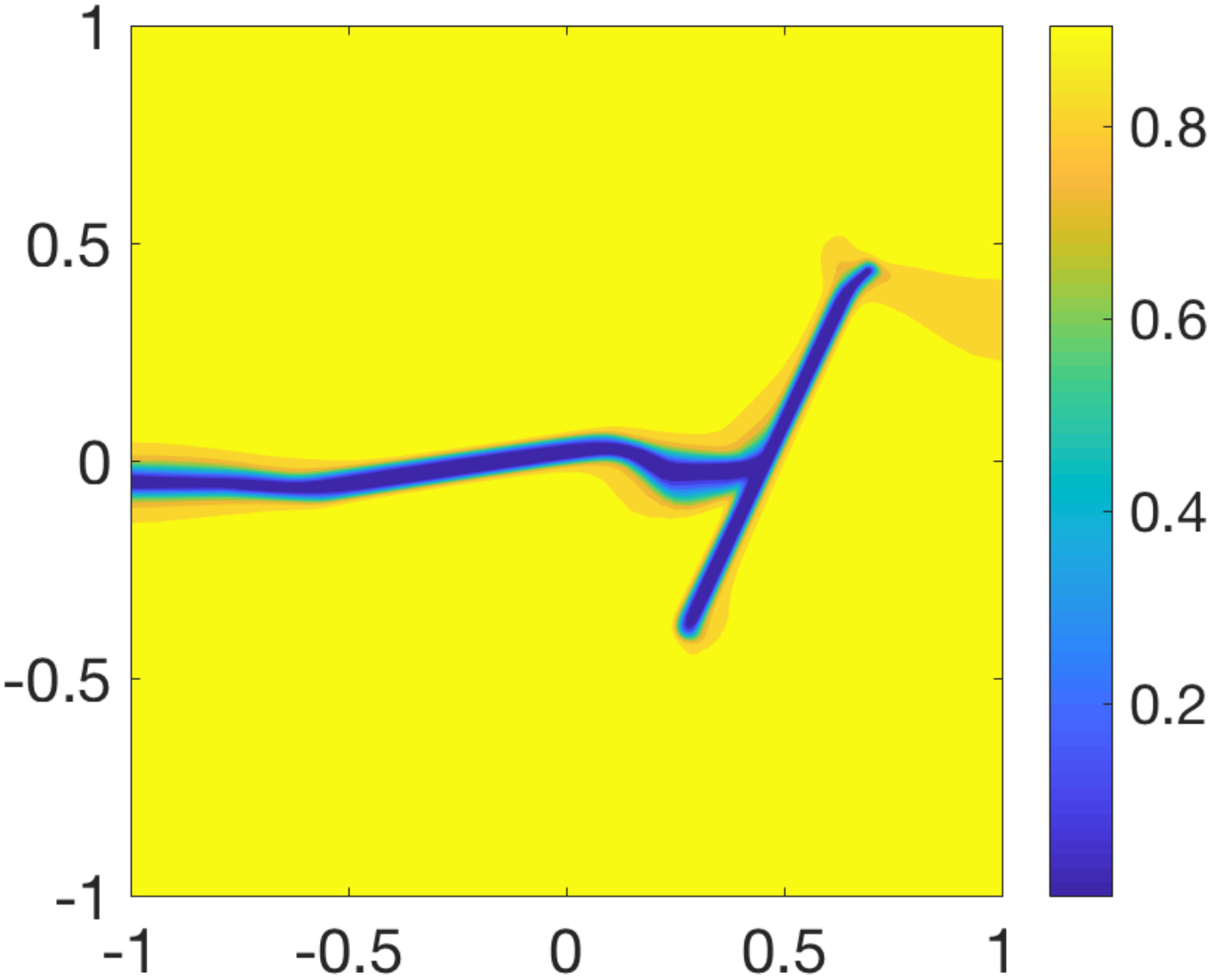}}
\subfigure[$U = 3.2 \times 10^{-2}$~mm]{\label{fig:subfig:SMAD_U320}
\includegraphics[width=0.225\linewidth]{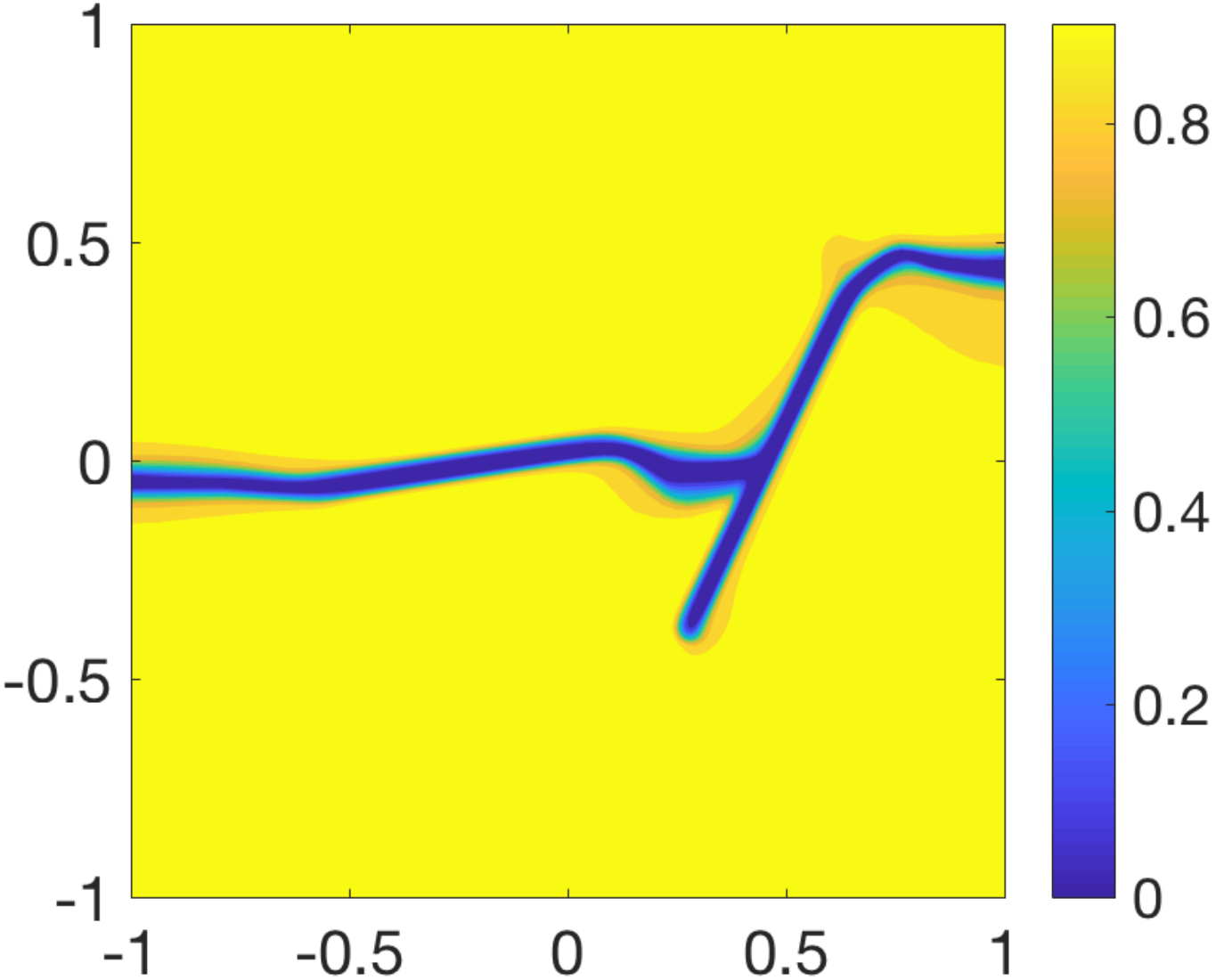}}
\vfill
\subfigure[$U = 2.4 \times 10^{-2}$~mm]{\label{fig:subfig:SMAM_U240}
\includegraphics[width=0.225\linewidth]{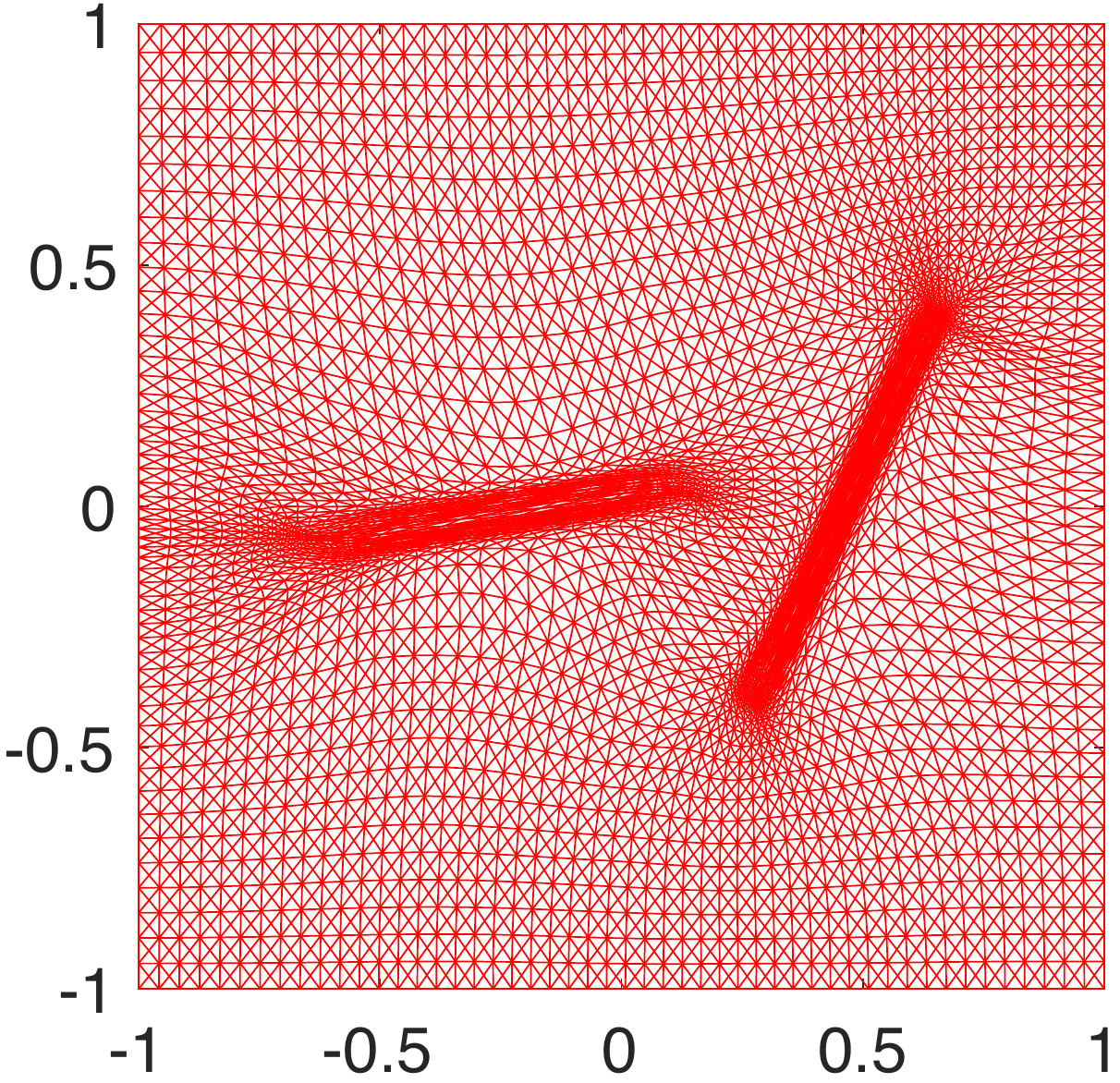}}
\subfigure[$U = 3.6 \times 10^{-2}$~mm]{\label{fig:subfig:SMAM_U260}
\includegraphics[width=0.225\linewidth]{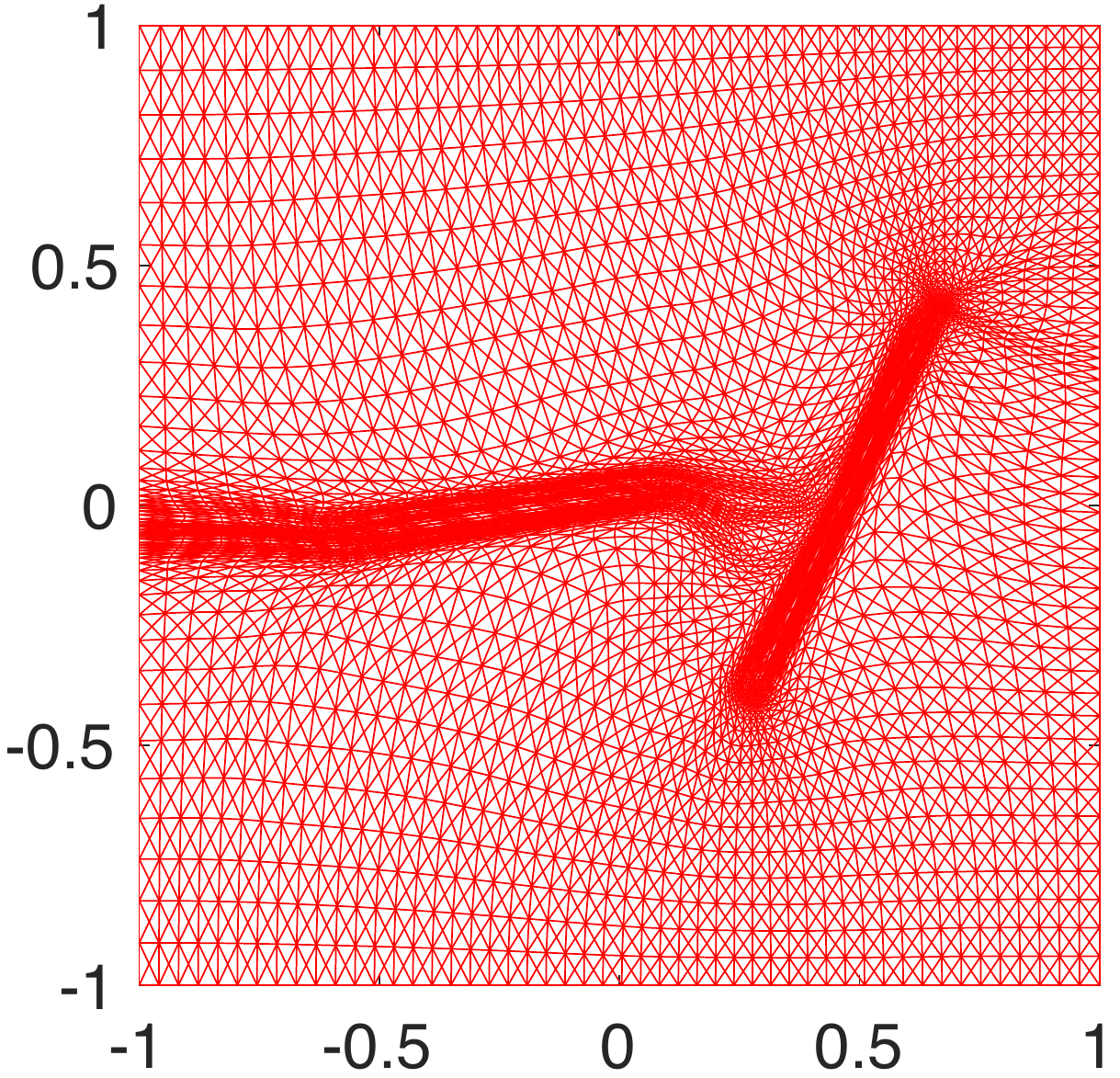}}
\subfigure[$U = 2.8 \times 10^{-2}$~mm]{\label{fig:subfig:SMAM_U280}
\includegraphics[width=0.225\linewidth]{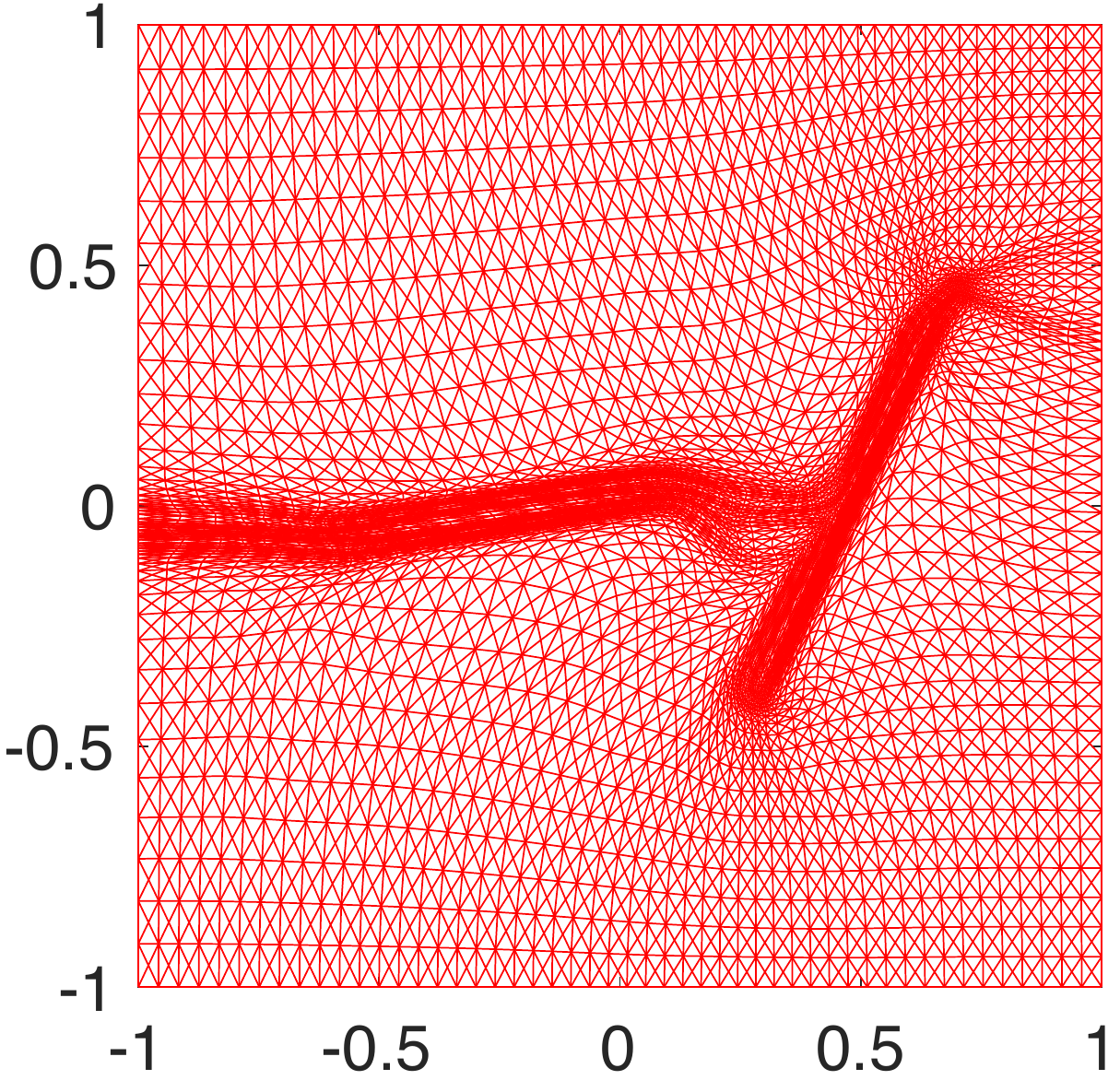}}
\subfigure[$U = 3.2 \times 10^{-2}$~mm]{\label{fig:subfig:SMAM_U320}
\includegraphics[width=0.225\linewidth]{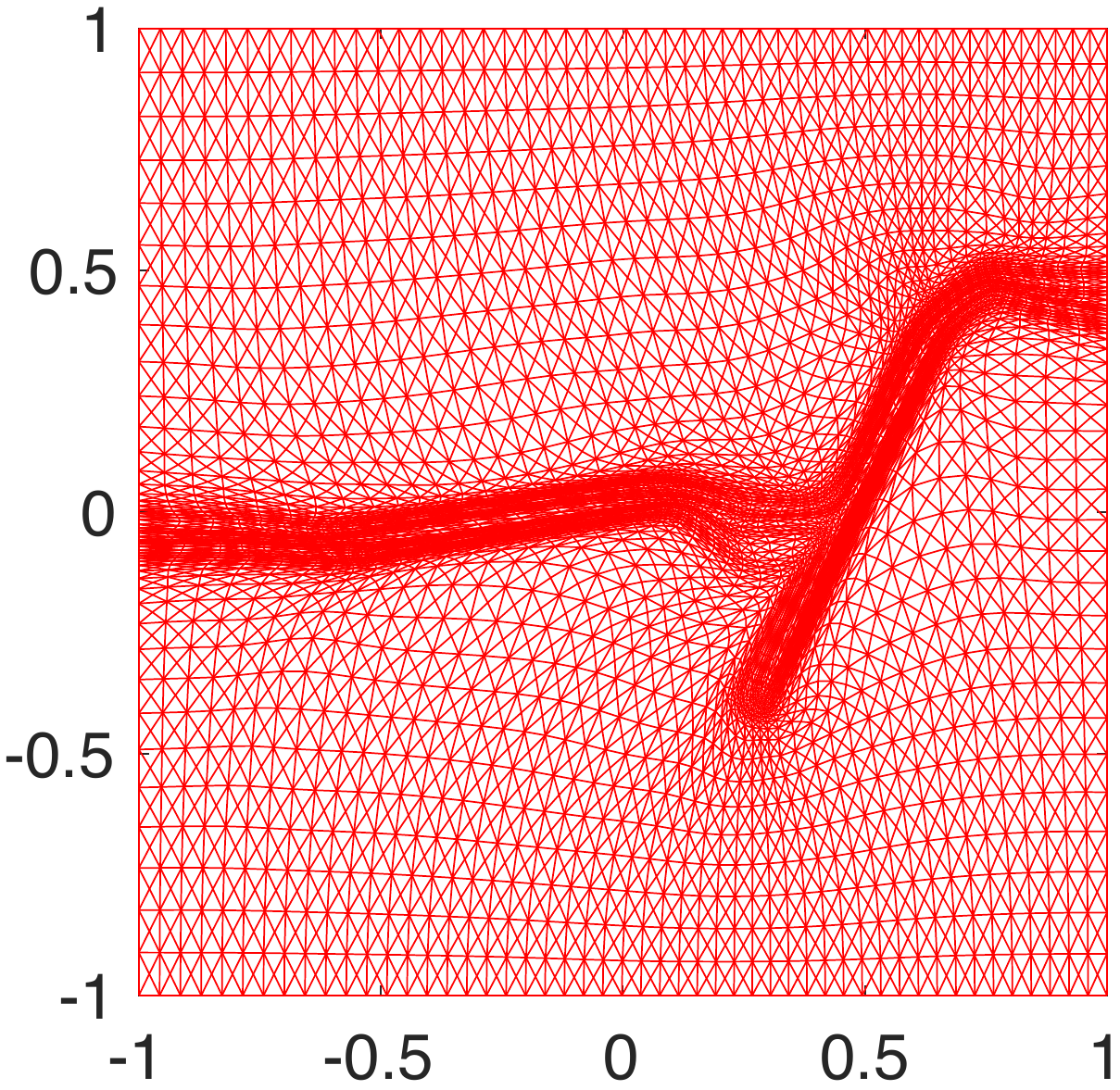}}
\vfill
\subfigure[$U = 2.4 \times 10^{-2}$~mm]{\label{fig:subfig:SMAS_U240}
\includegraphics[width=0.225\linewidth]{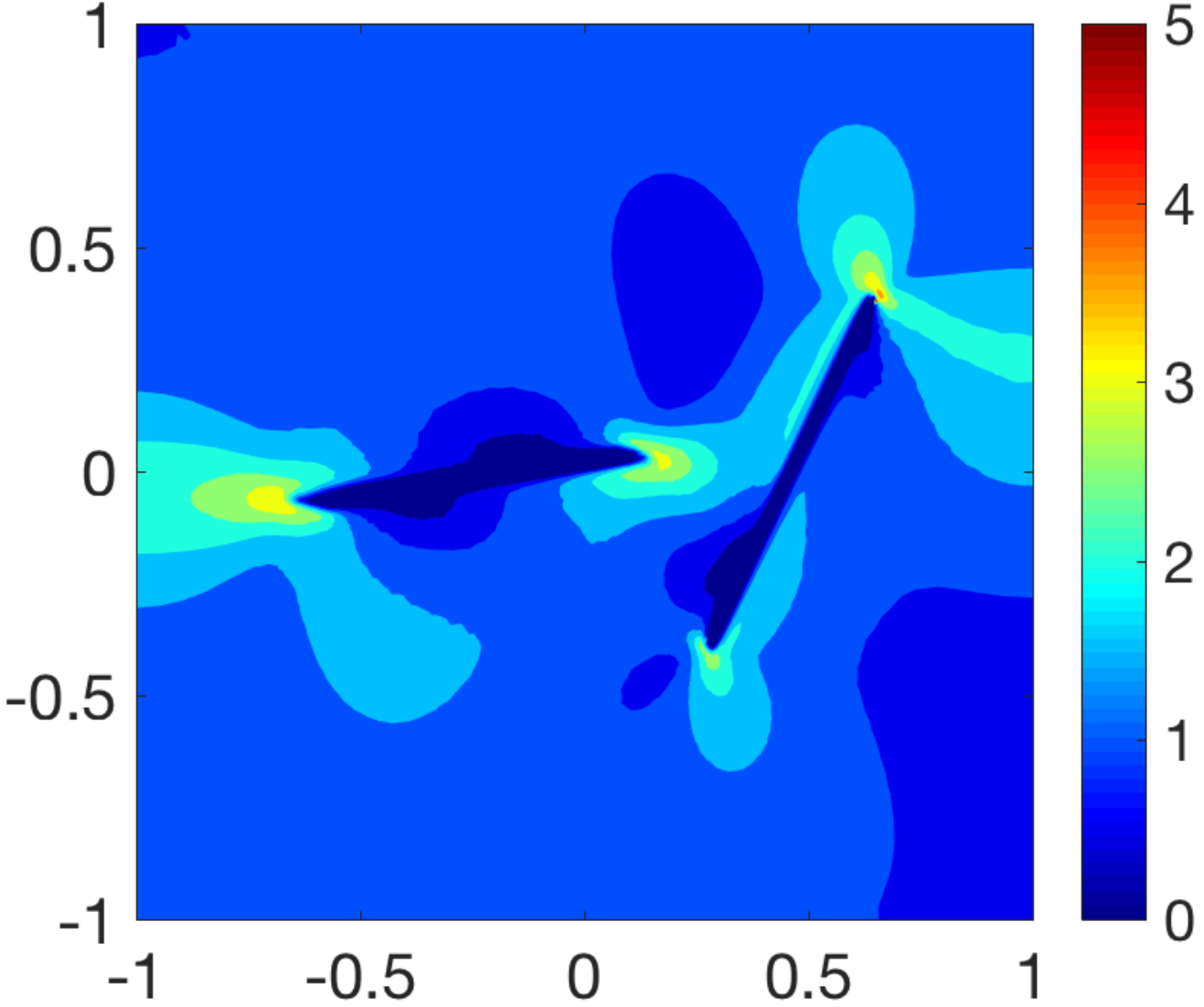}}
\subfigure[$U = 2.6 \times 10^{-2}$~mm]{\label{fig:subfig:SMAS_U260}
\includegraphics[width=0.225\linewidth]{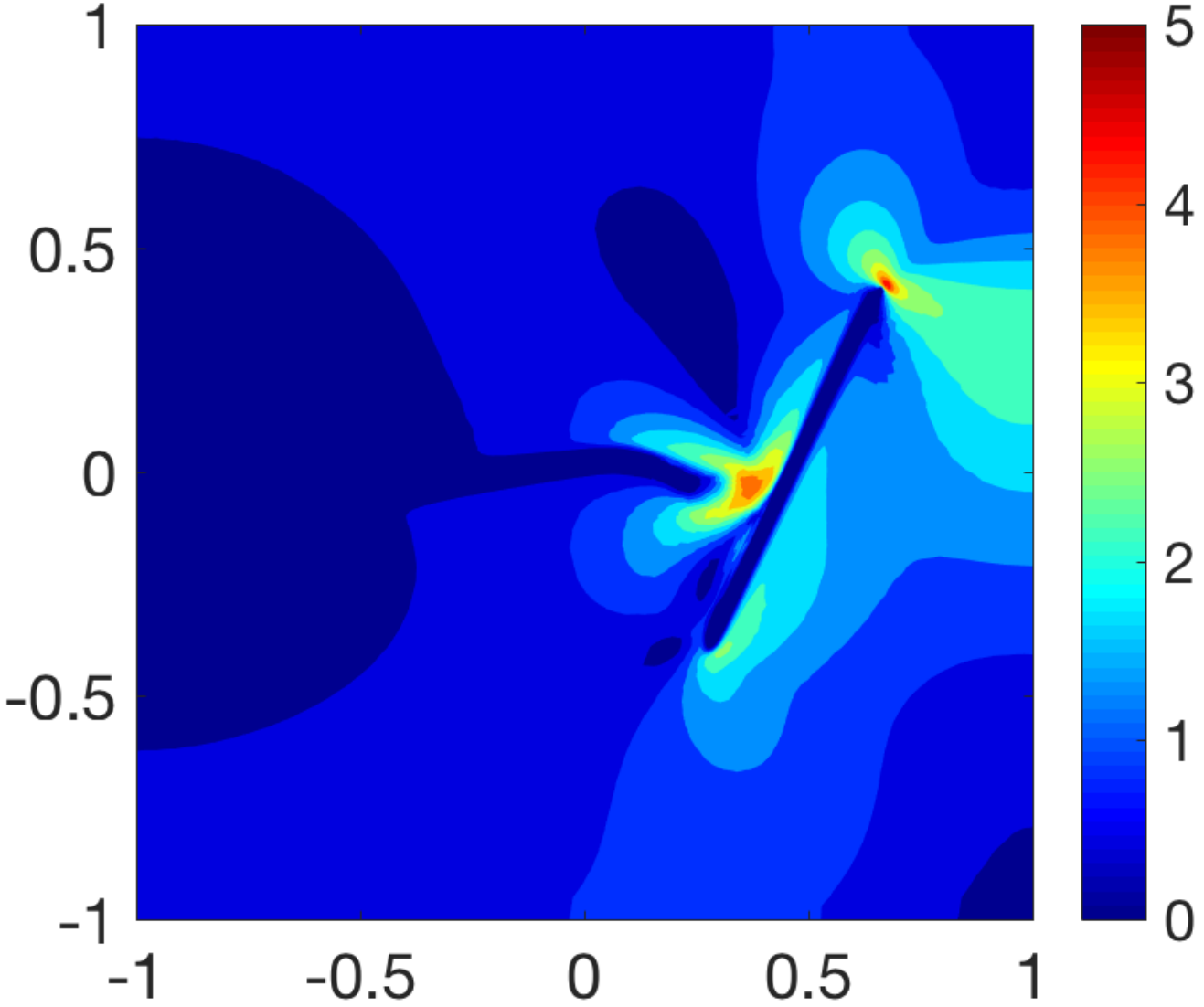}}
\subfigure[$U = 2.8 \times 10^{-2}$~mm]{\label{fig:subfig:SMAS_U280}
\includegraphics[width=0.225\linewidth]{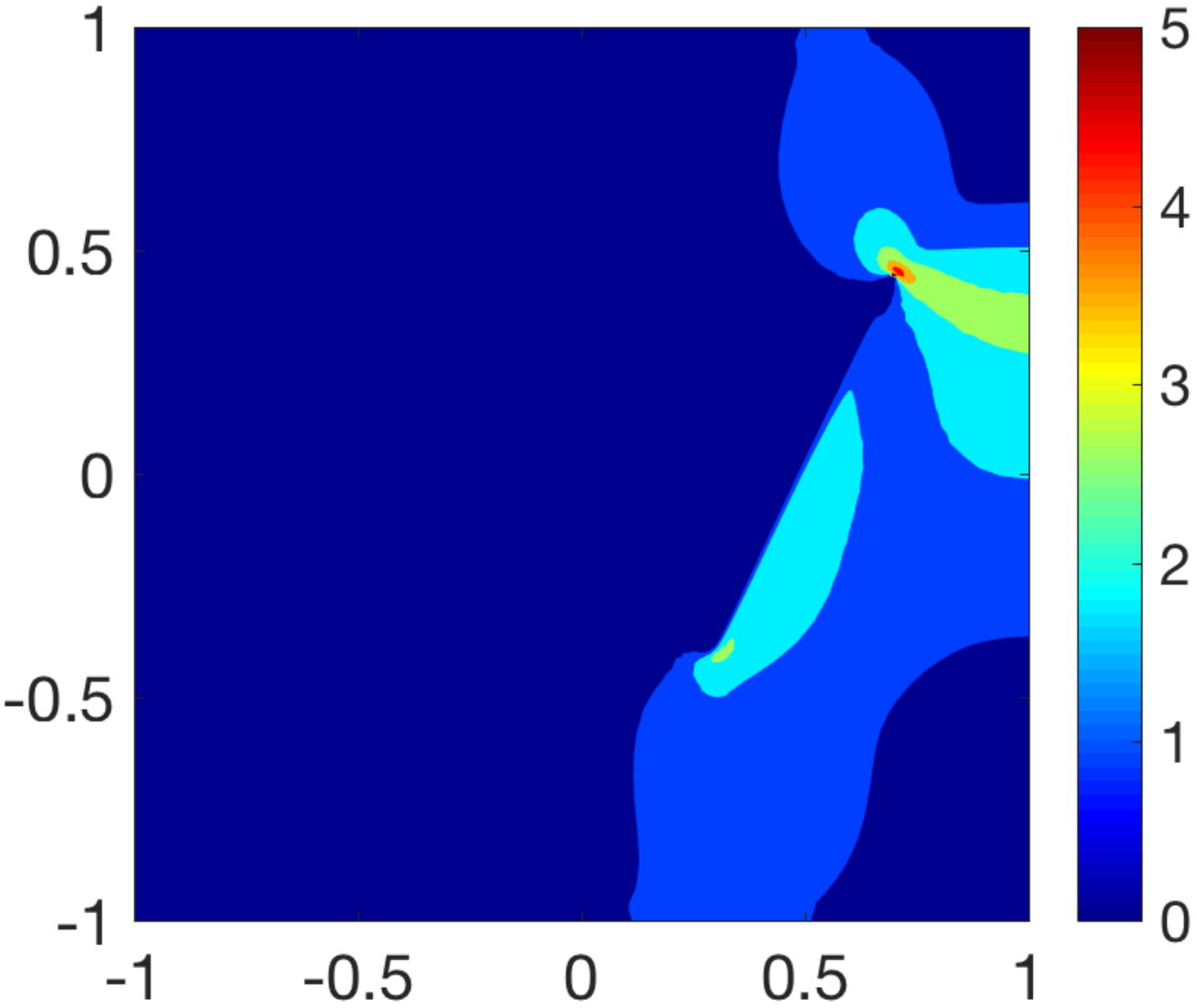}}
\subfigure[$U = 3.2 \times 10^{-2}$~mm]{\label{fig:subfig:SMAS_U320}
\includegraphics[width=0.225\linewidth]{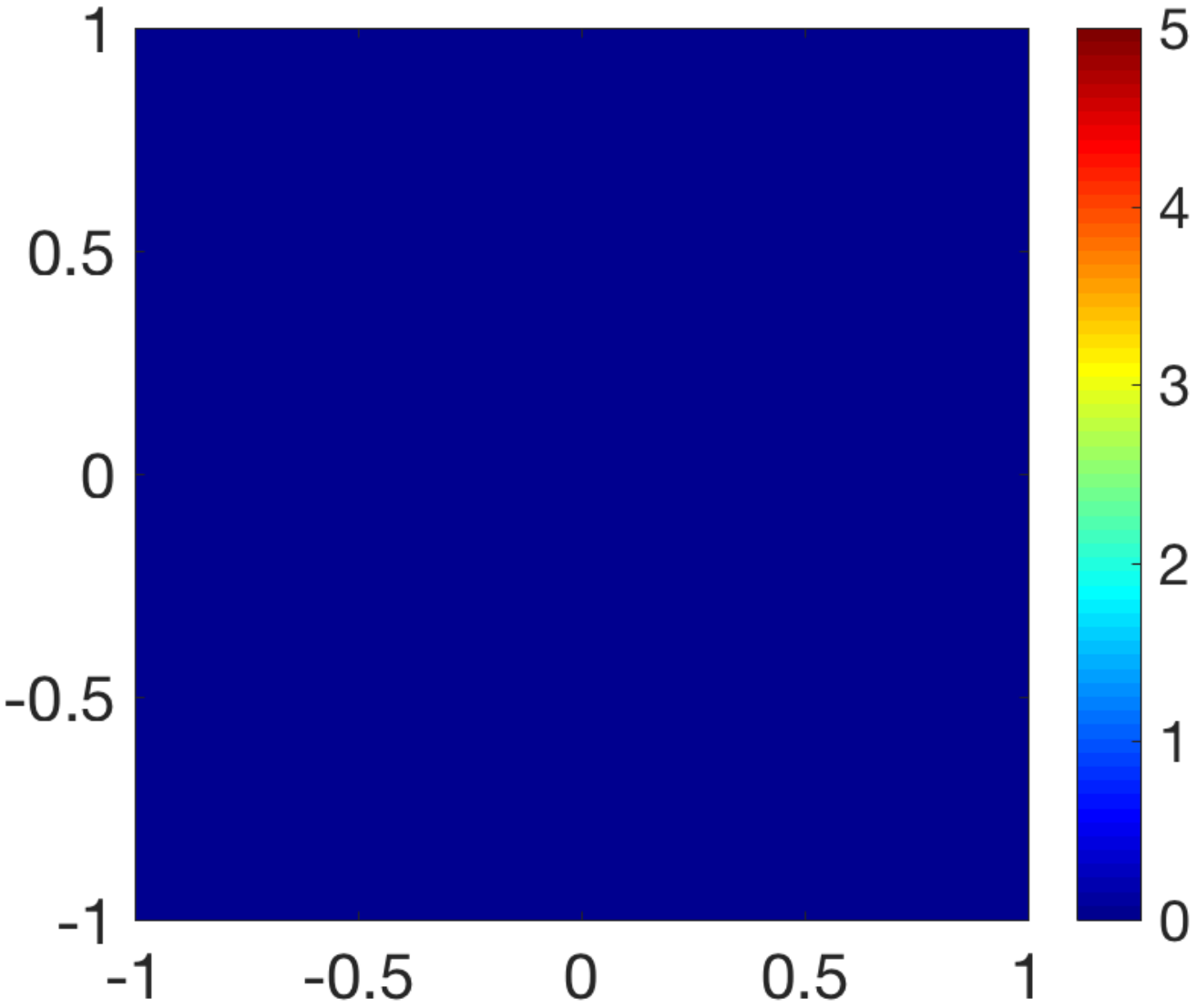}}
\caption{Example 3. The mesh and contours of the phase-field and von Mises stress distribution during
crack evolution for the two-crack shear test with $l = 0.00375$~mm, $N = 10,000\; (51\times51)$.
(v-d split with ItCBC)}
\label{fig:SAMBC}
\end{figure}

\begin{figure} 
\centering 
\subfigure[$U = 2.8 \times 10^{-2}$~mm]{\label{fig:subfig:TMND_U280}
\includegraphics[width=0.225\linewidth]{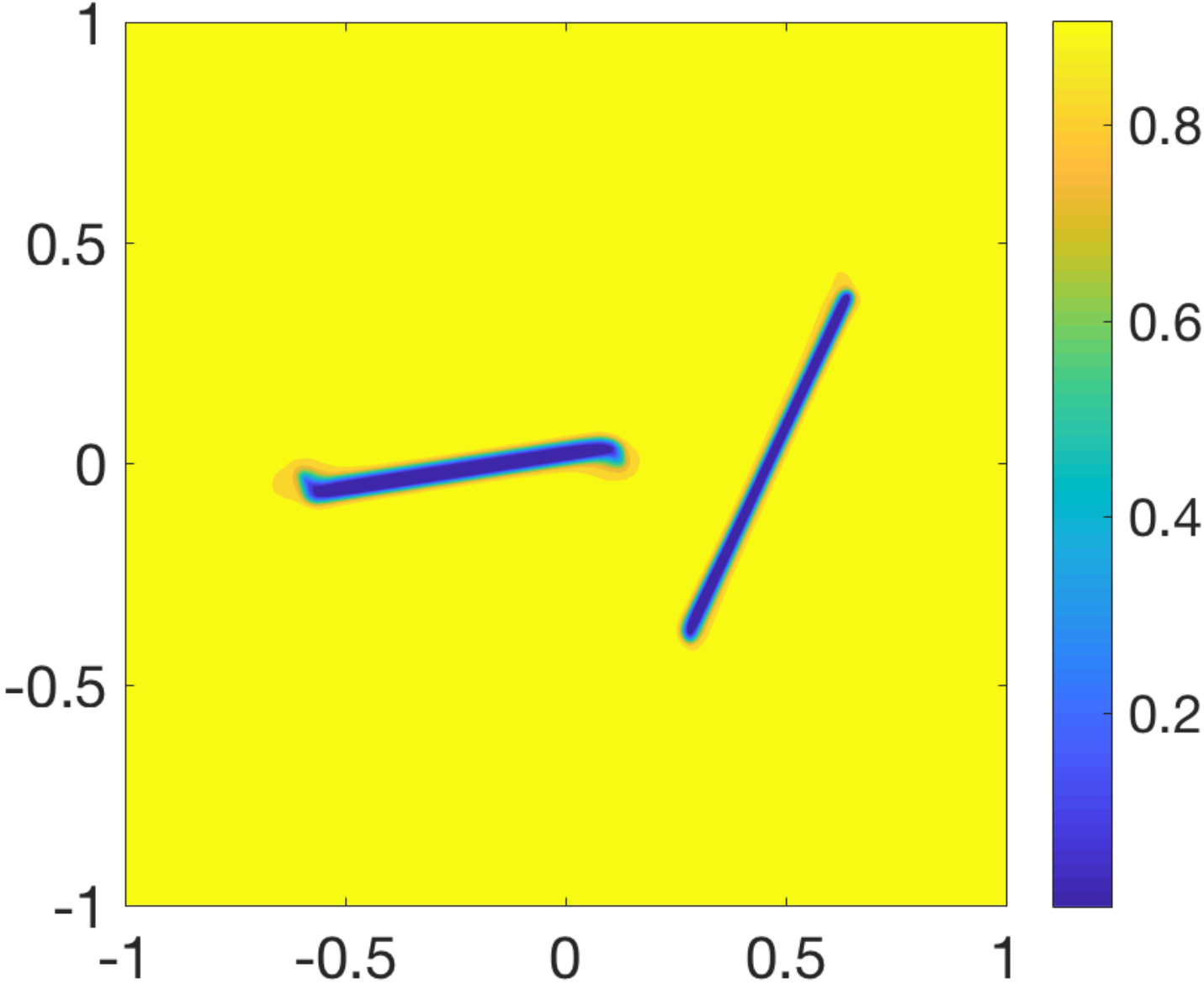}}
\subfigure[$U = 3.0 \times 10^{-2}$~mm]{\label{fig:subfig:TMND_U300}
\includegraphics[width=0.225\linewidth]{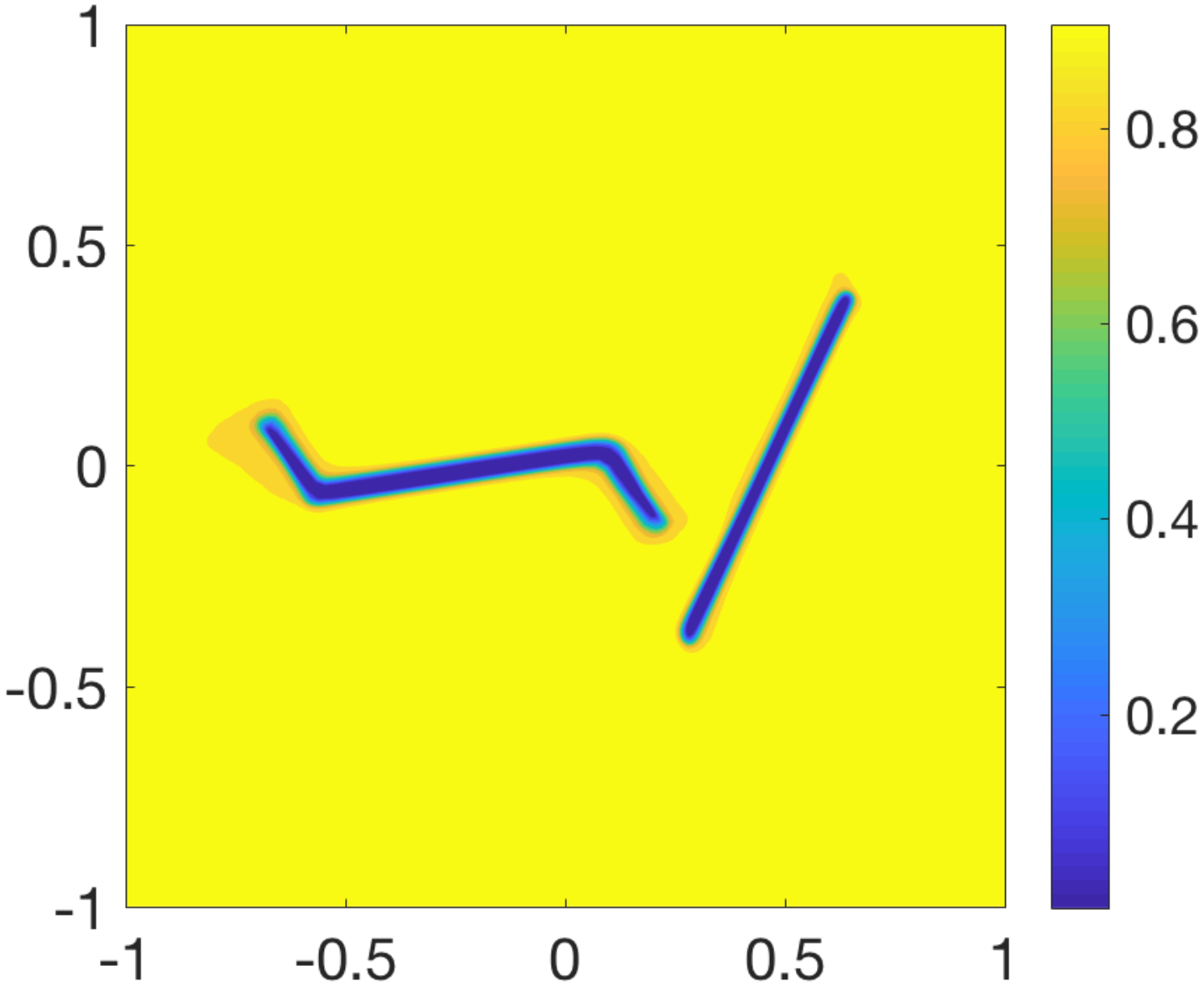}}
\subfigure[$U = 3.4 \times 10^{-2}$~mm]{\label{fig:subfig:TMND_U340}
\includegraphics[width=0.225\linewidth]{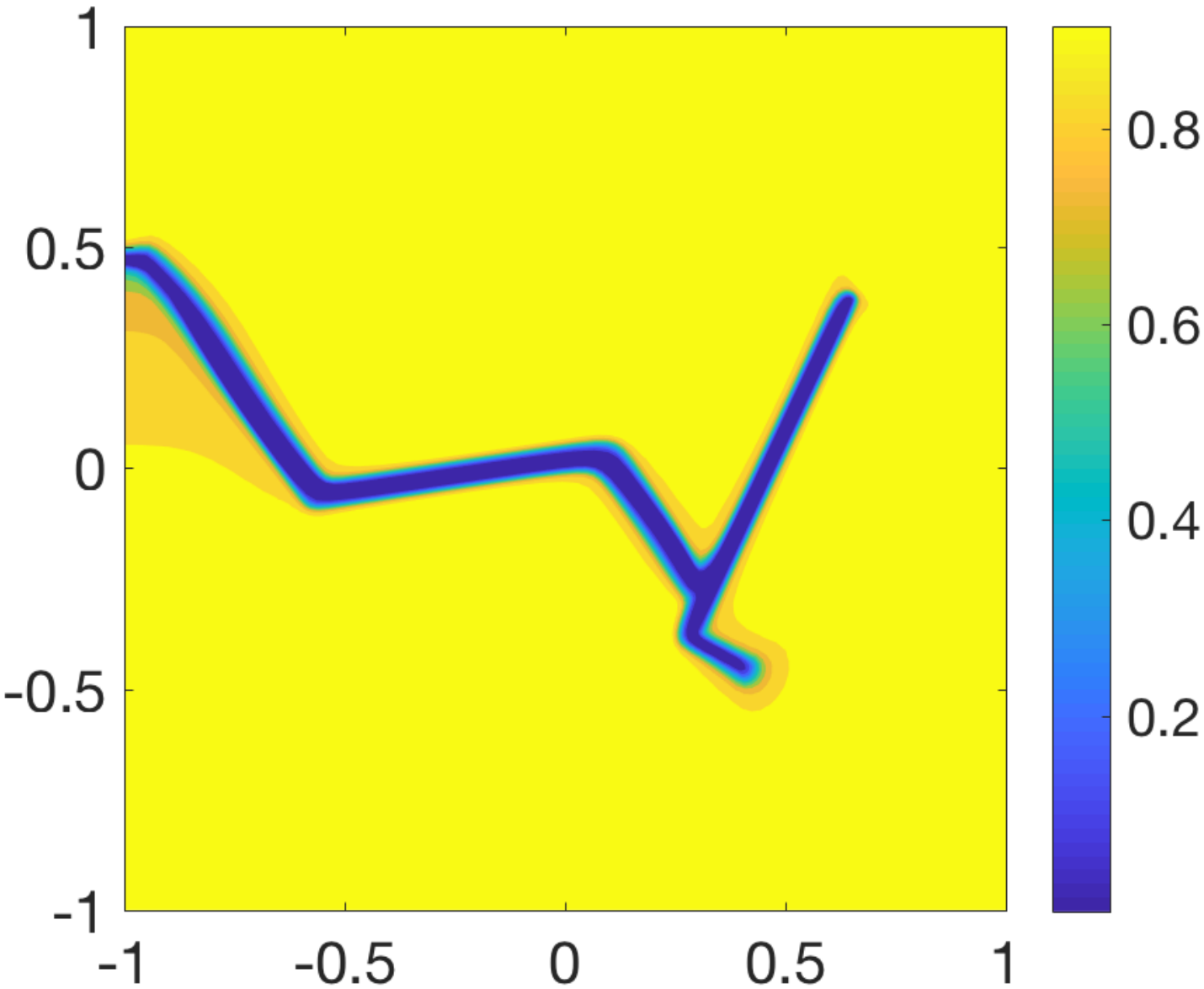}}
\subfigure[$U = 4.2 \times 10^{-2}$~mm]{\label{fig:subfig:TMND_U420}
\includegraphics[width=0.225\linewidth]{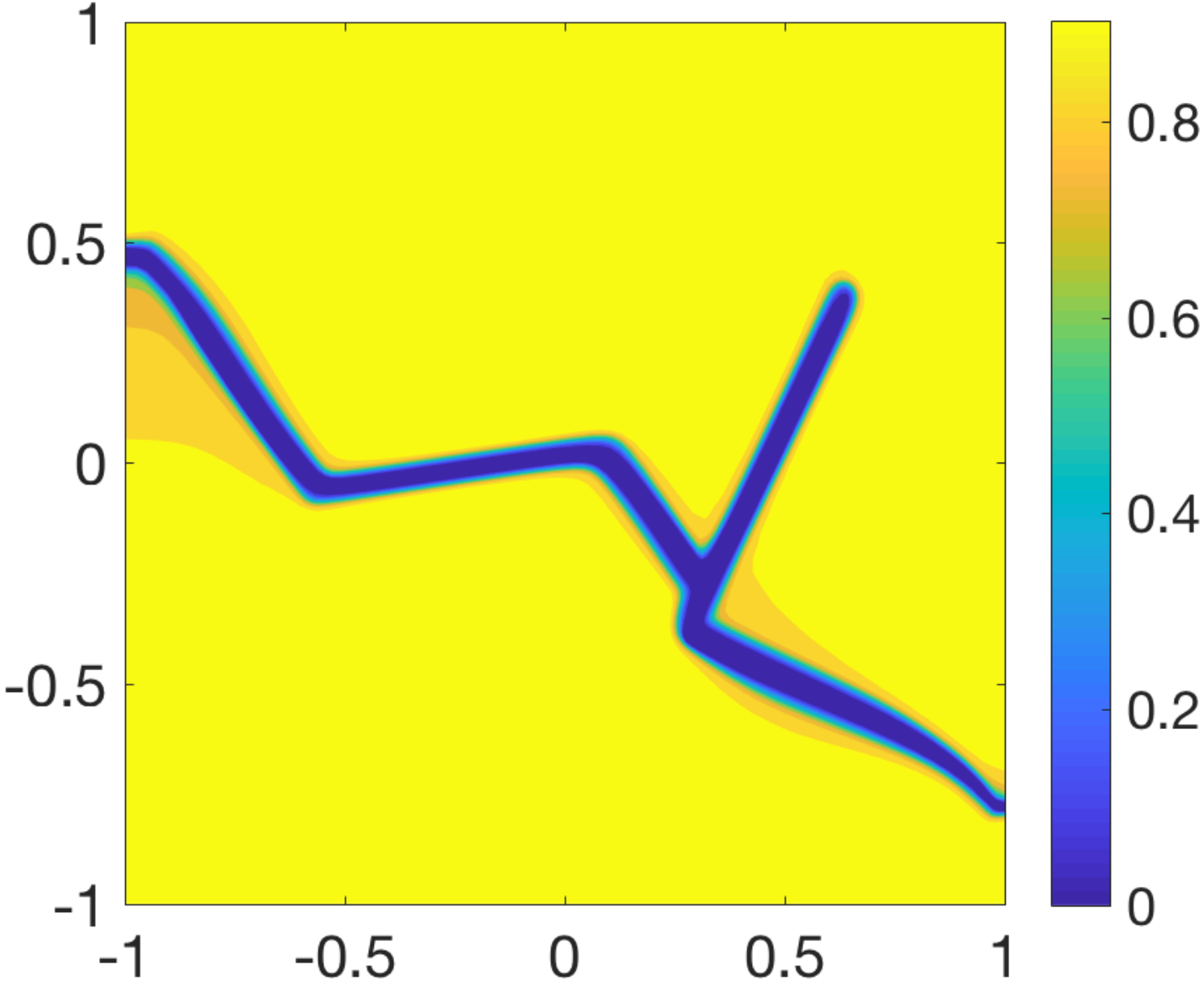}}
\vfill
\subfigure[$U = 2.8 \times 10^{-2}$~mm]{\label{fig:subfig:SMNM_U280}
\includegraphics[width=0.225\linewidth]{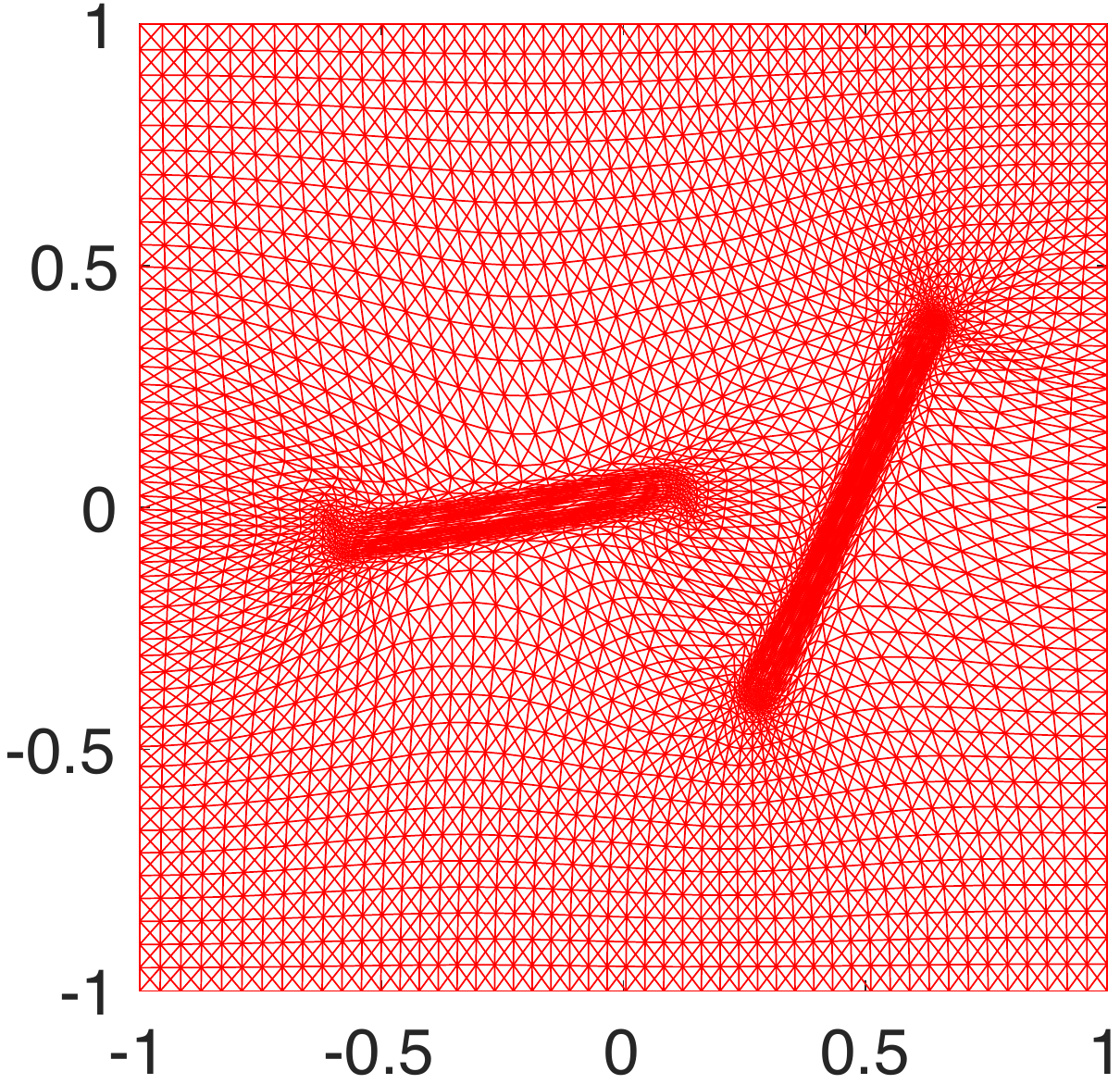}}
\subfigure[$U = 3.0 \times 10^{-2}$~mm]{\label{fig:subfig:SMNM_U300}
\includegraphics[width=0.225\linewidth]{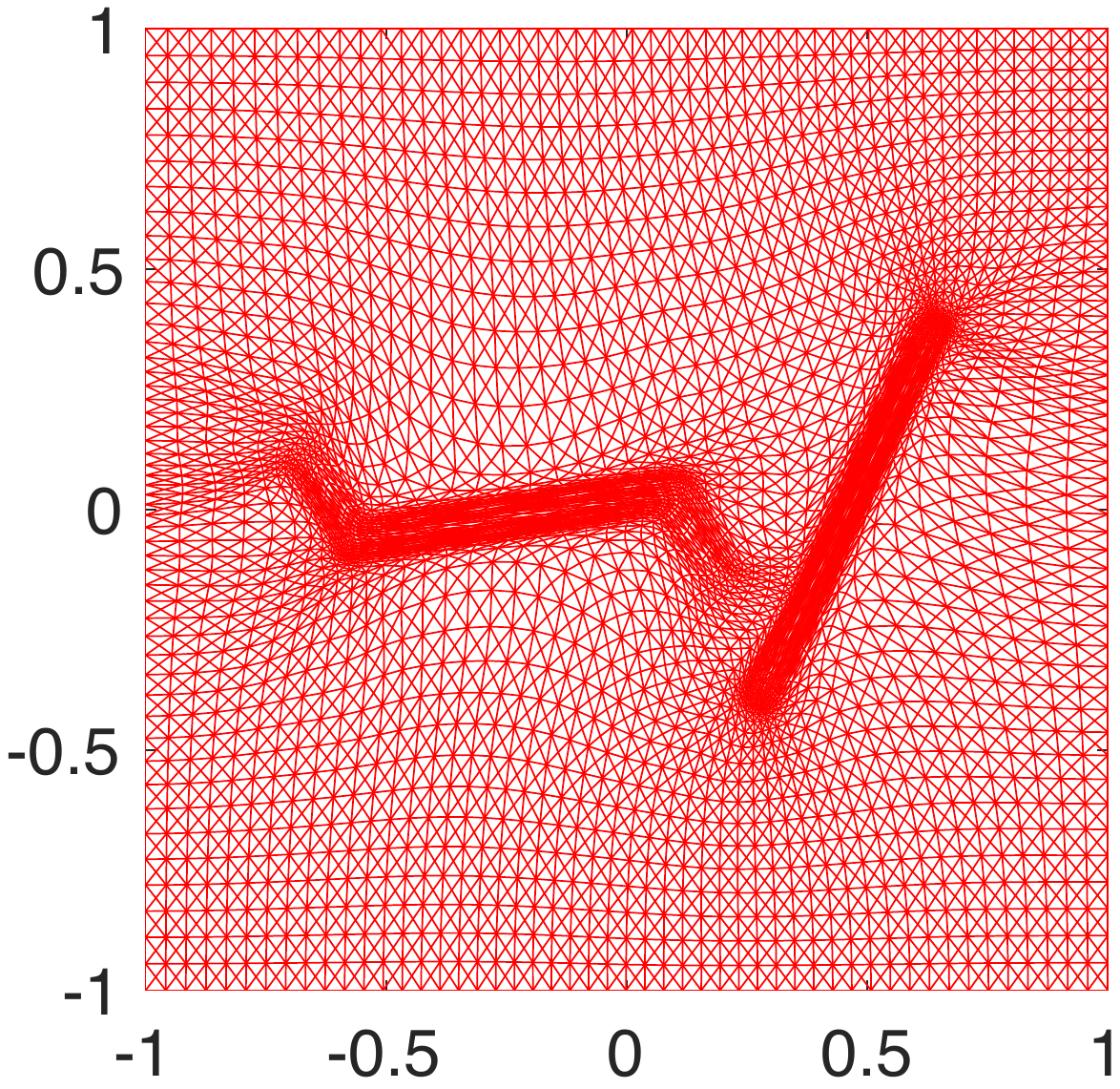}}
\subfigure[$U = 3.4 \times 10^{-2}$~mm]{\label{fig:subfig:SMNM_U340}
\includegraphics[width=0.225\linewidth]{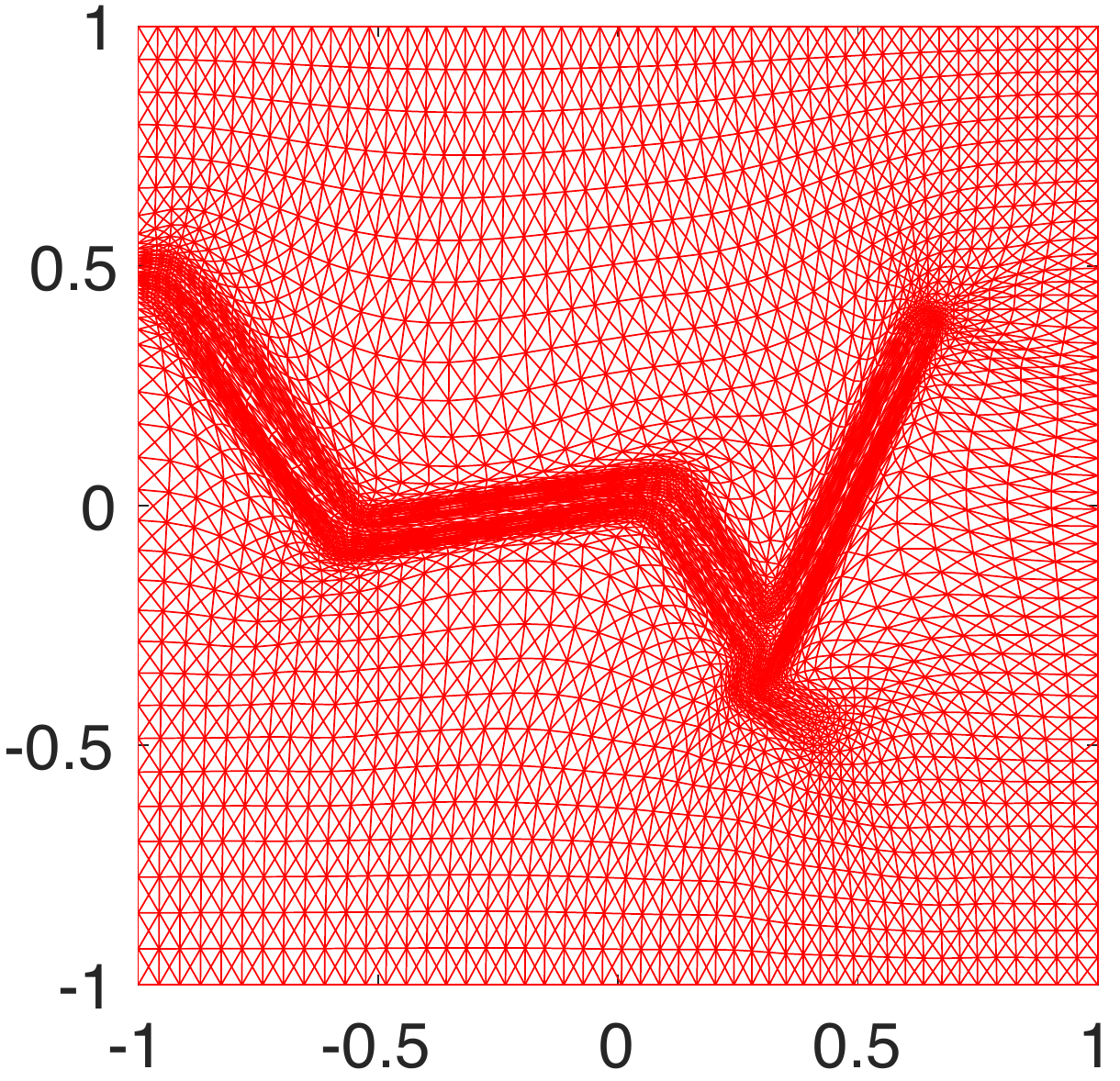}}
\subfigure[$U = 4.2 \times 10^{-2}$~mm]{\label{fig:subfig:SMNM_U420}
\includegraphics[width=0.225\linewidth]{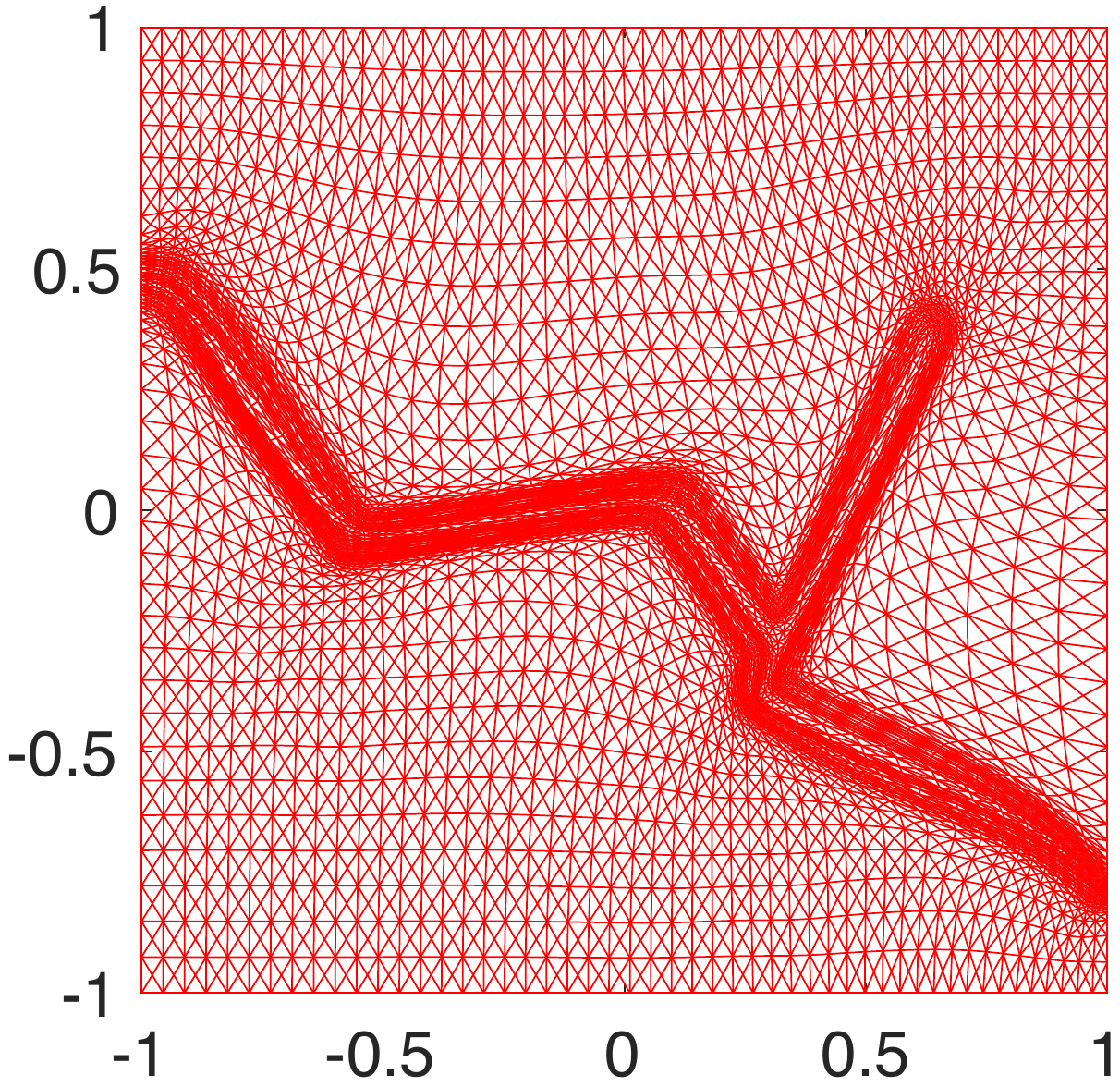}}
\vfill
\subfigure[$U = 2.8 \times 10^{-2}$~mm]{\label{fig:subfig:TMNS_U280}
\includegraphics[width=0.225\linewidth]{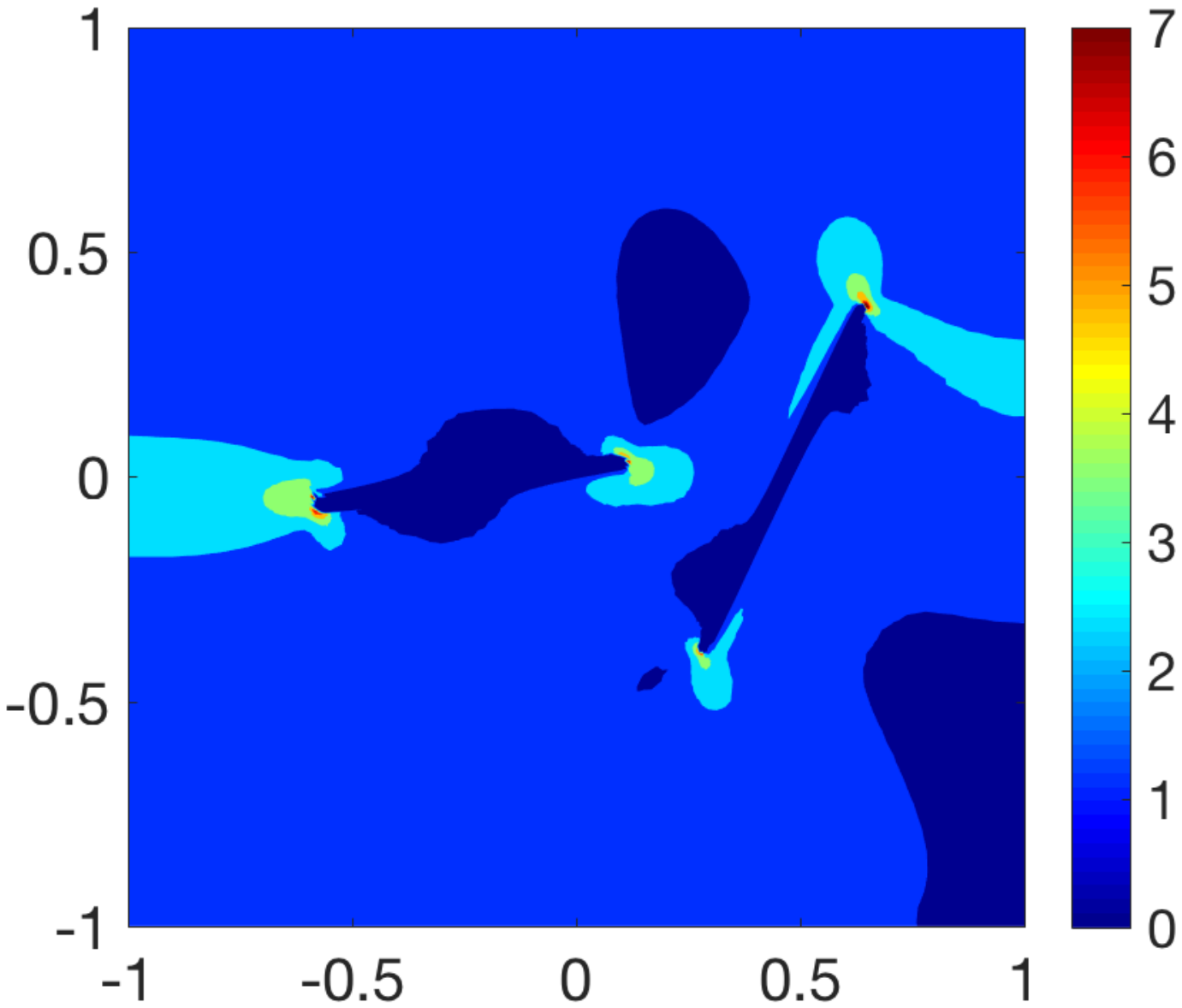}}
\subfigure[$U = 3.0 \times 10^{-2}$~mm]{\label{fig:subfig:TMNS_U300}
\includegraphics[width=0.225\linewidth]{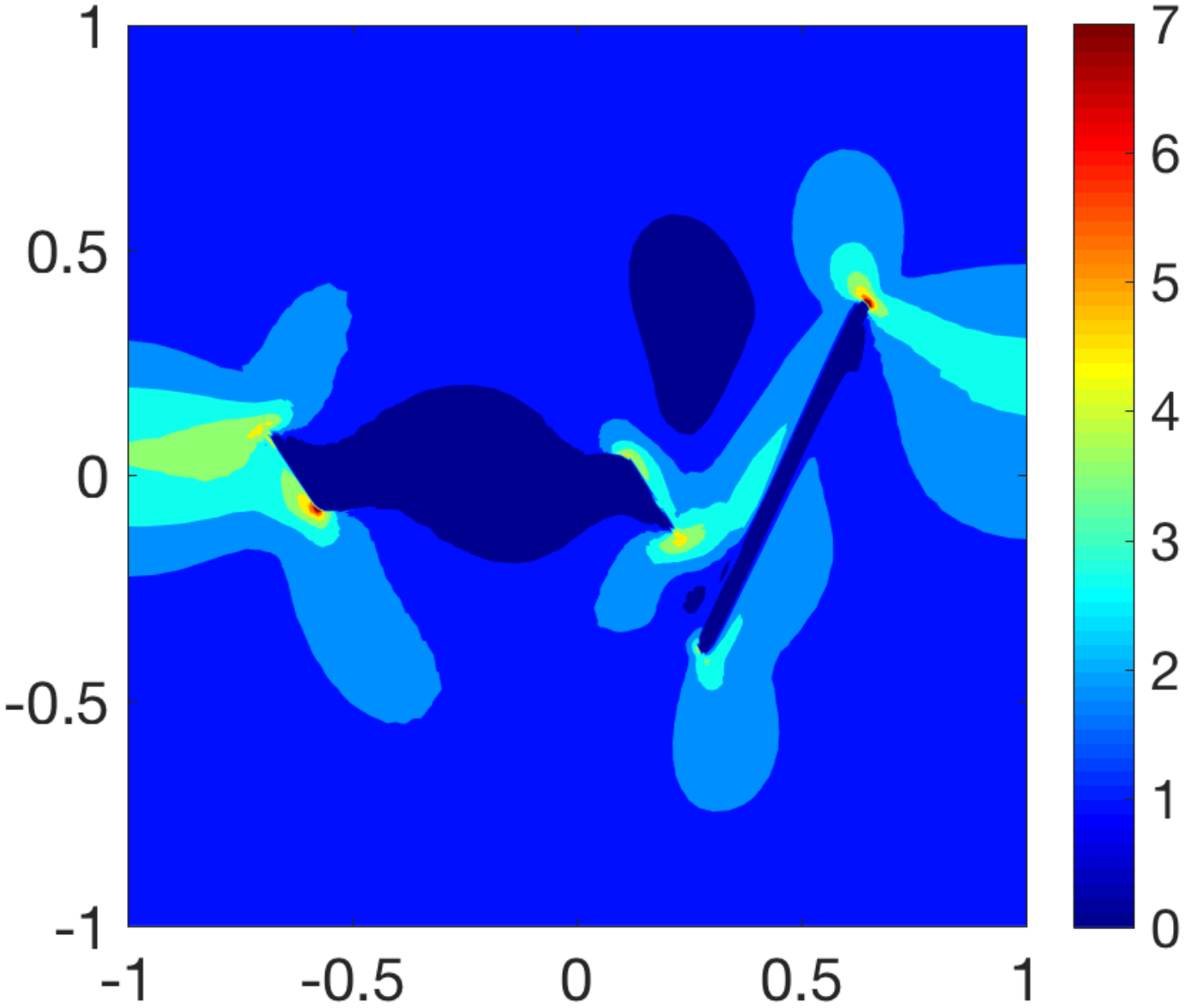}}
\subfigure[$U = 3.4 \times 10^{-2}$~mm]{\label{fig:subfig:TMNS_U340}
\includegraphics[width=0.225\linewidth]{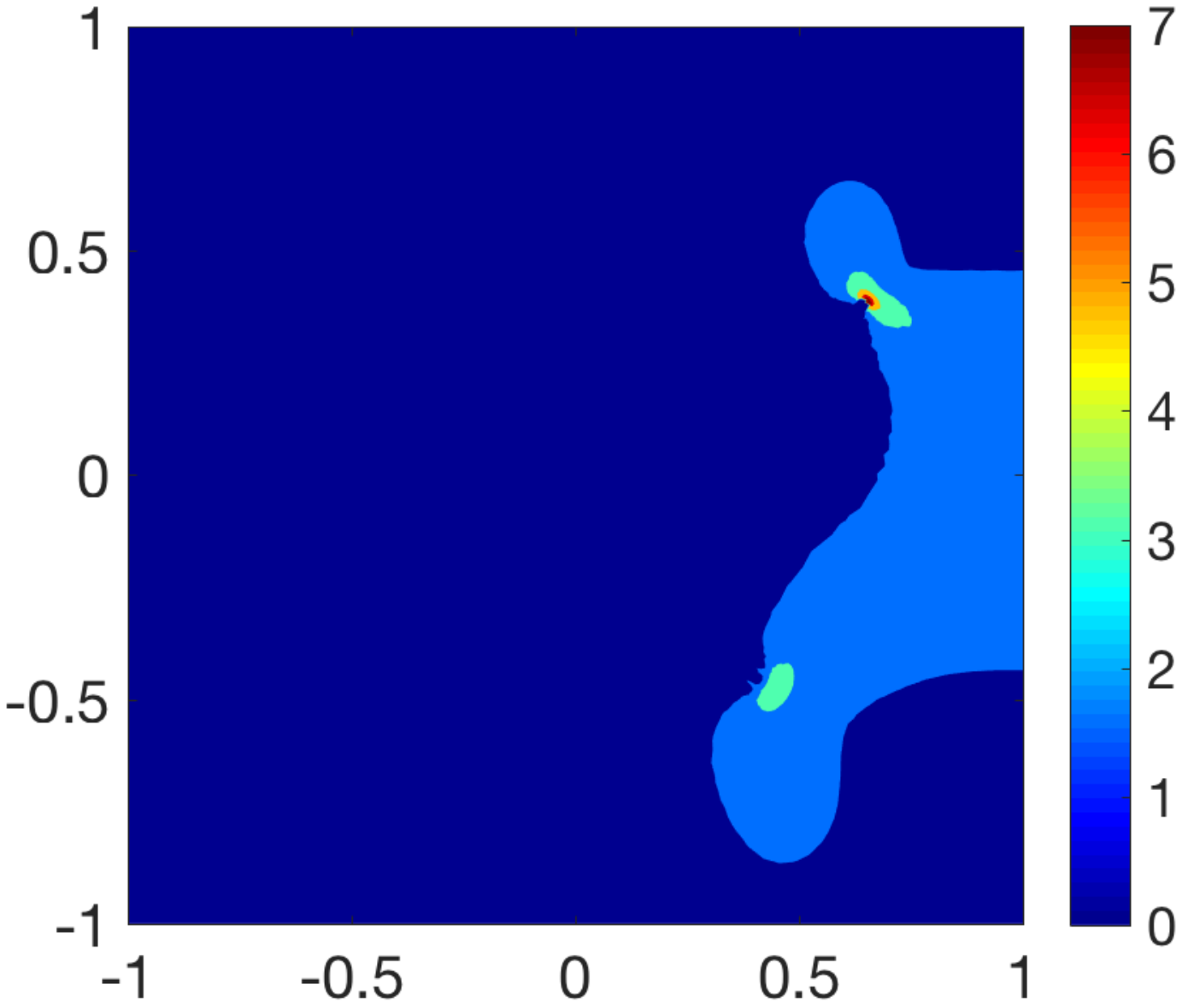}}
\subfigure[$U = 4.2 \times 10^{-2}$~mm]{\label{fig:subfig:TMNS_U420}
\includegraphics[width=0.225\linewidth]{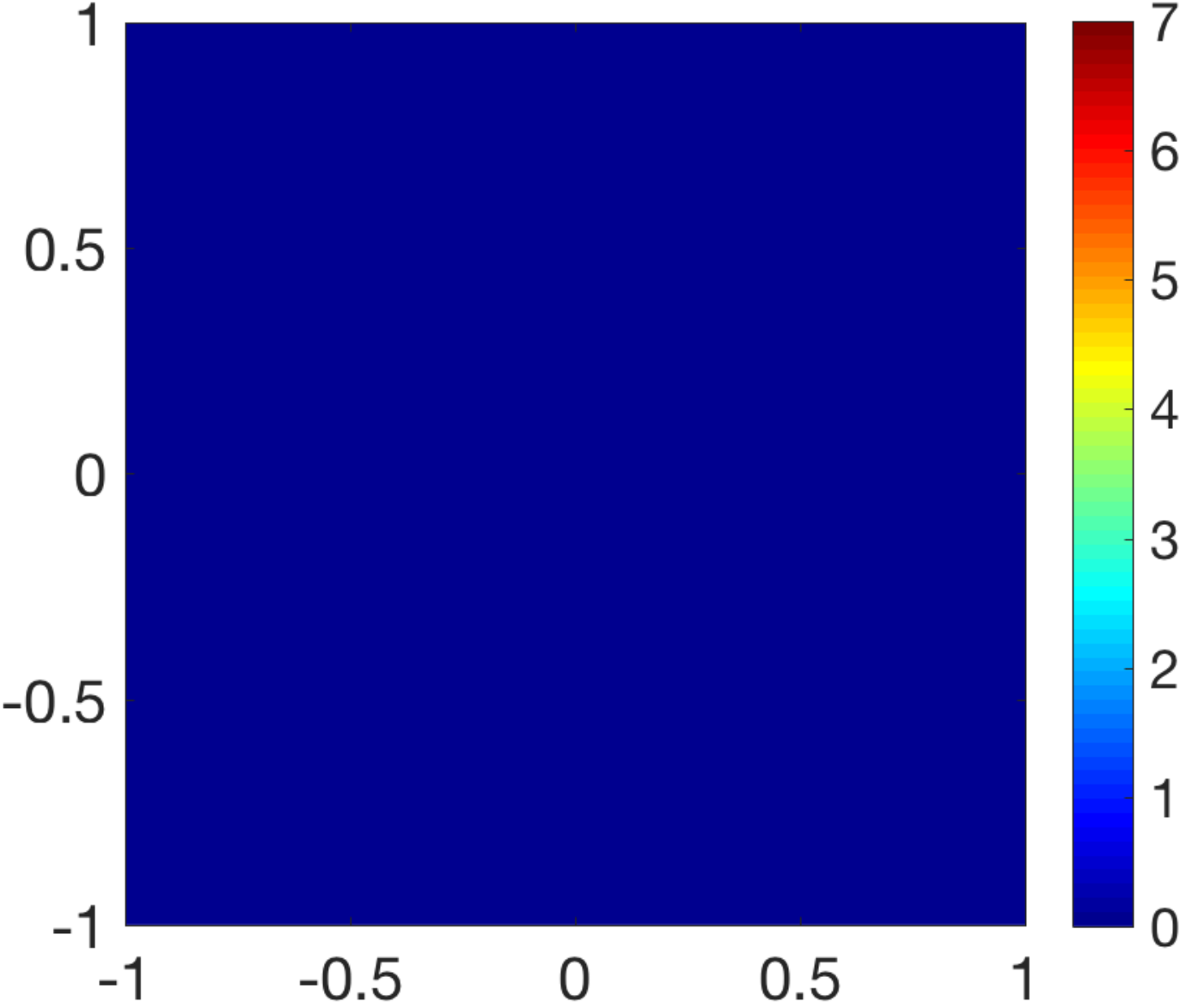}}
\caption{Example 3. The mesh and contours of the phase-field and von Mises stress distribution during
crack evolution for the two-crack shear test with $l = 0.00375$~mm, $N = 10,000\; (51\times51)$.
(improved v-d split with ItCBC)}
\label{fig:SNMBC}
\end{figure}

Fig. \ref{fig:SSOBC} shows the adaptive mesh and phase-field and von Mises stress distribution for spectral decomposition with original crack boundary conditions. The results for all three decomposition methods with ItCBC ($d_{cr} = 0.4$) are shown in Figs. \ref{fig:SSMBC}, \ref{fig:SAMBC}, and \ref{fig:SNMBC}, respectively. The results using improved v-d split in Fig. \ref{fig:SNMBC} agree well with those using spectral decomposition as shown in Fig. \ref{fig:SSMBC} whereas the v-d split (even with ItCBC) gives unphysical crack propagation.
This finding is consistent with the observations made from the previous examples.

\begin{figure} 
\centering 
\subfigure[$U = 0.179$~mm]{\label{fig:subfig:FOMD_U1790}
\includegraphics[width=0.18\linewidth]{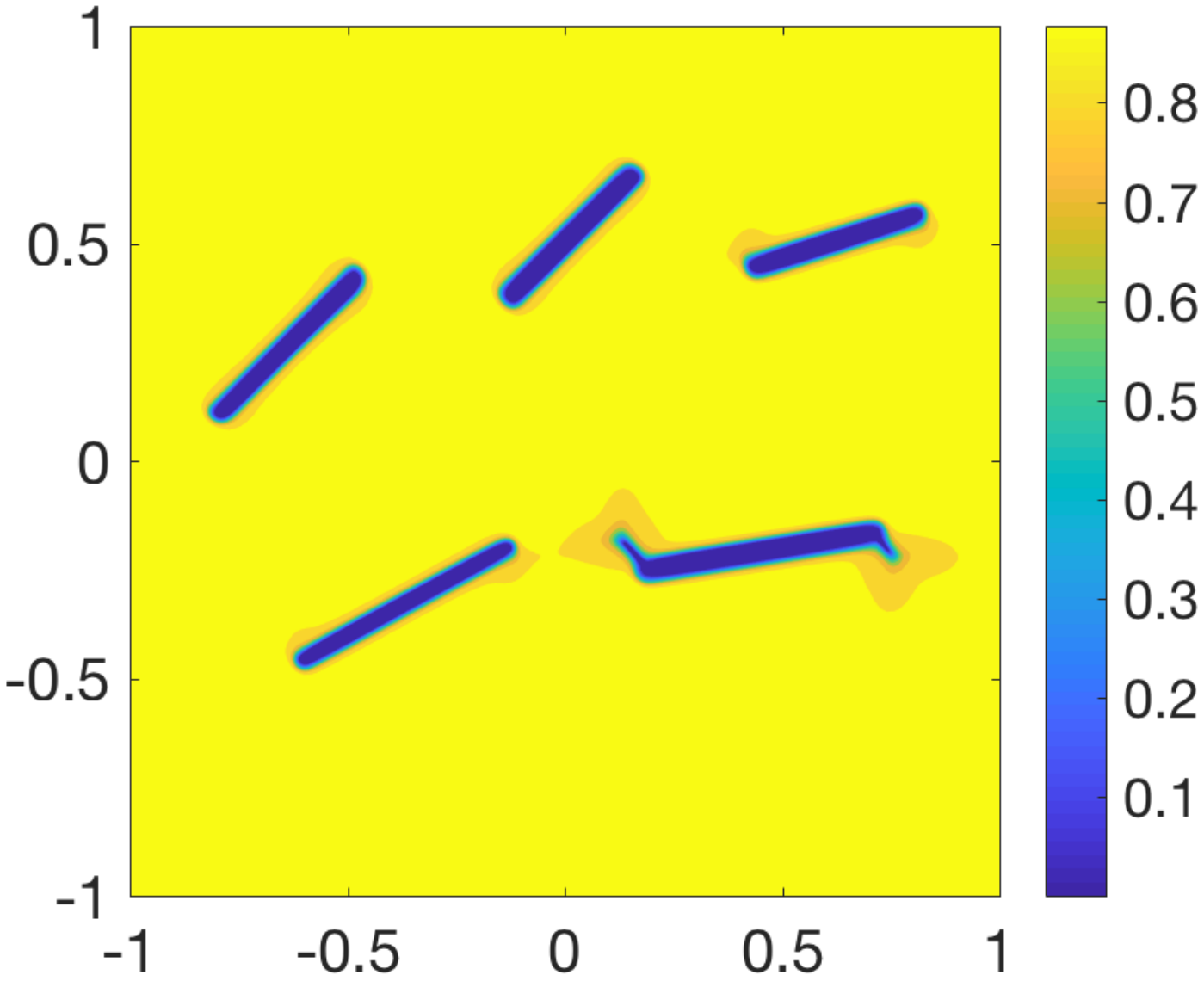}}
\subfigure[$U = 0.182$~mm]{\label{fig:subfig:FOMD_U1820}
\includegraphics[width=0.18\linewidth]{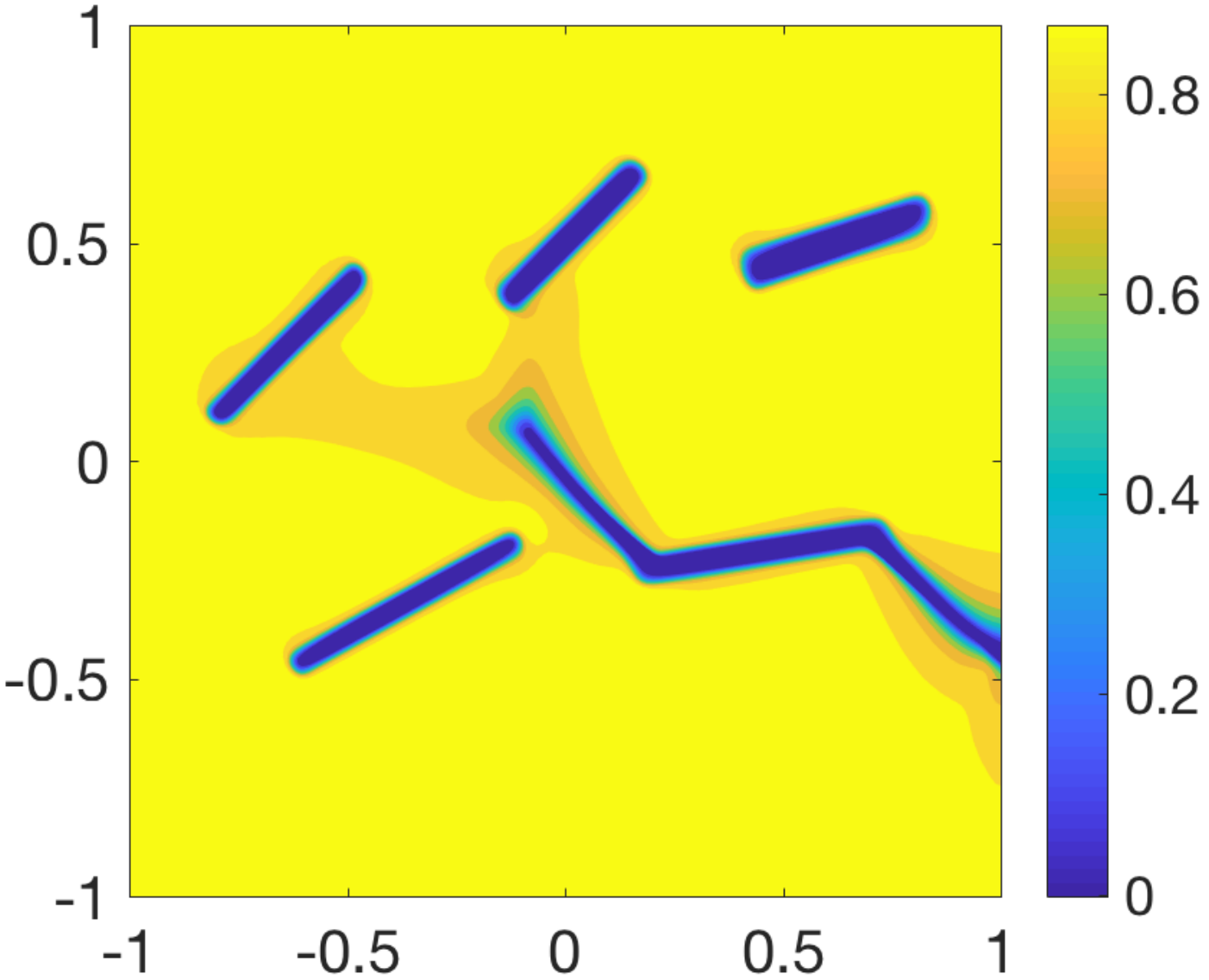}}
\subfigure[$U = 0.184$~mm]{\label{fig:subfig:FOMD_U1840}
\includegraphics[width=0.18\linewidth]{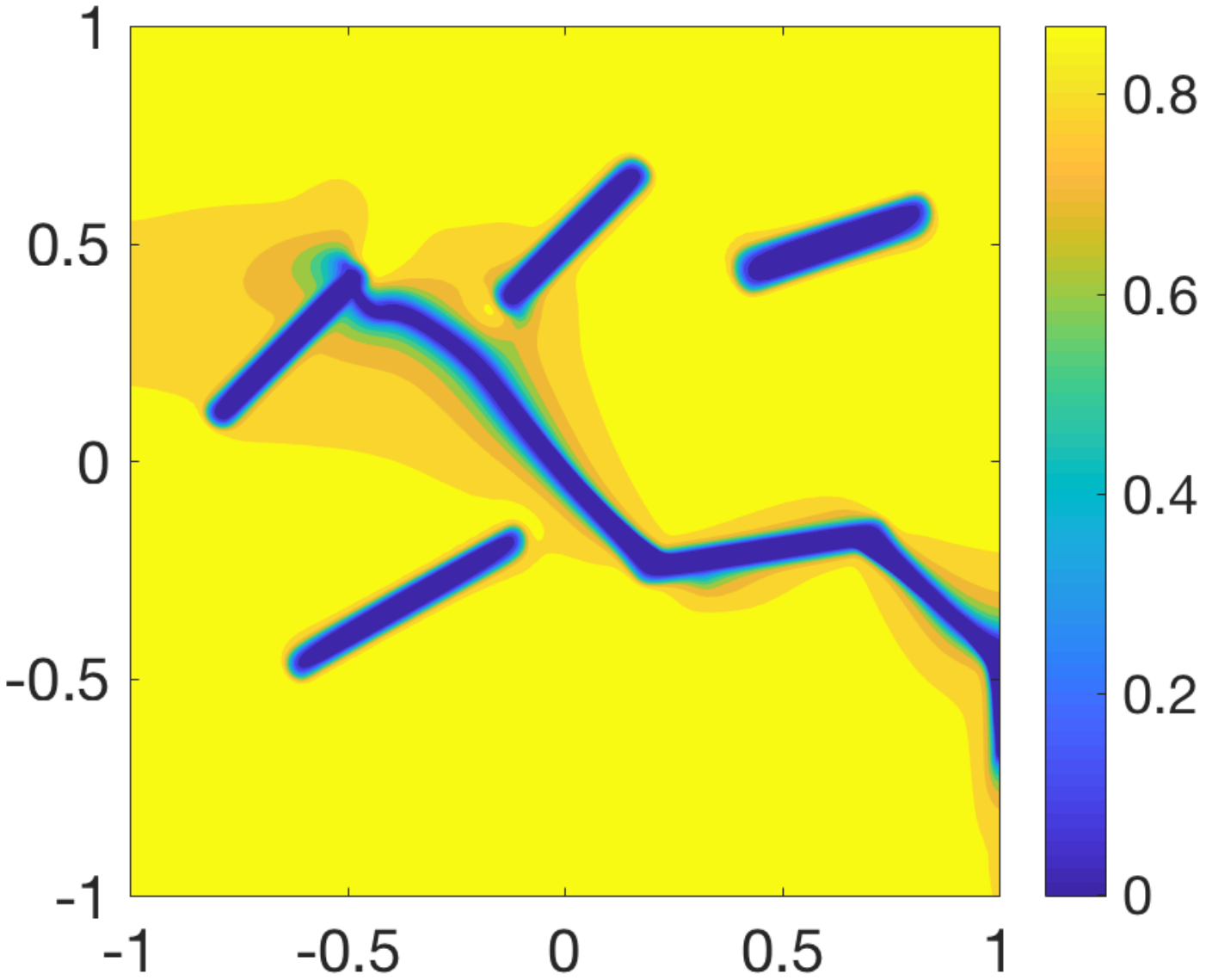}}
\subfigure[$U = 0.186$~mm]{\label{fig:subfig:FOMD_U1860}
\includegraphics[width=0.18\linewidth]{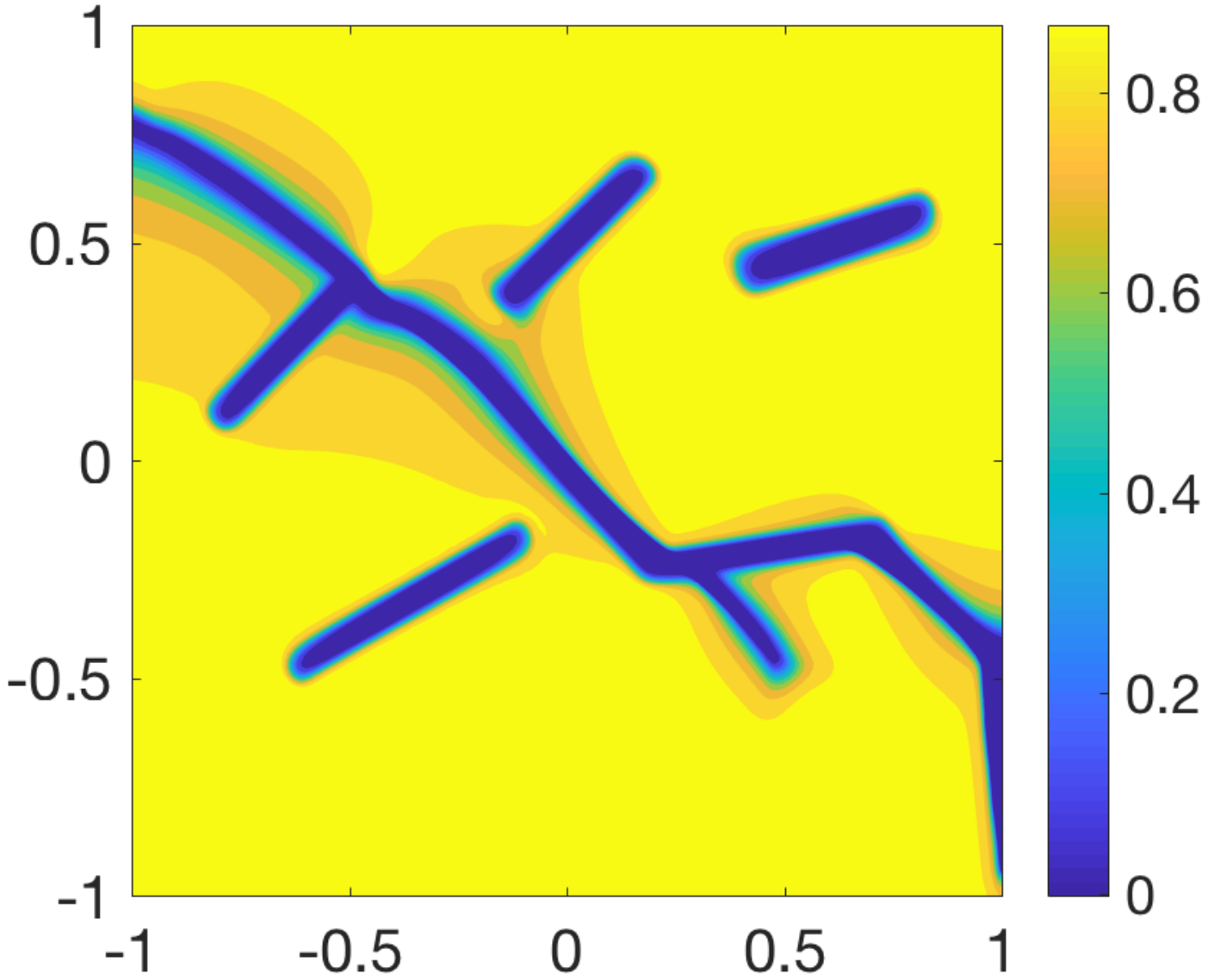}}
\subfigure[$U = 0.2$~mm]{\label{fig:subfig:FOMD_U2000}
\includegraphics[width=0.18\linewidth]{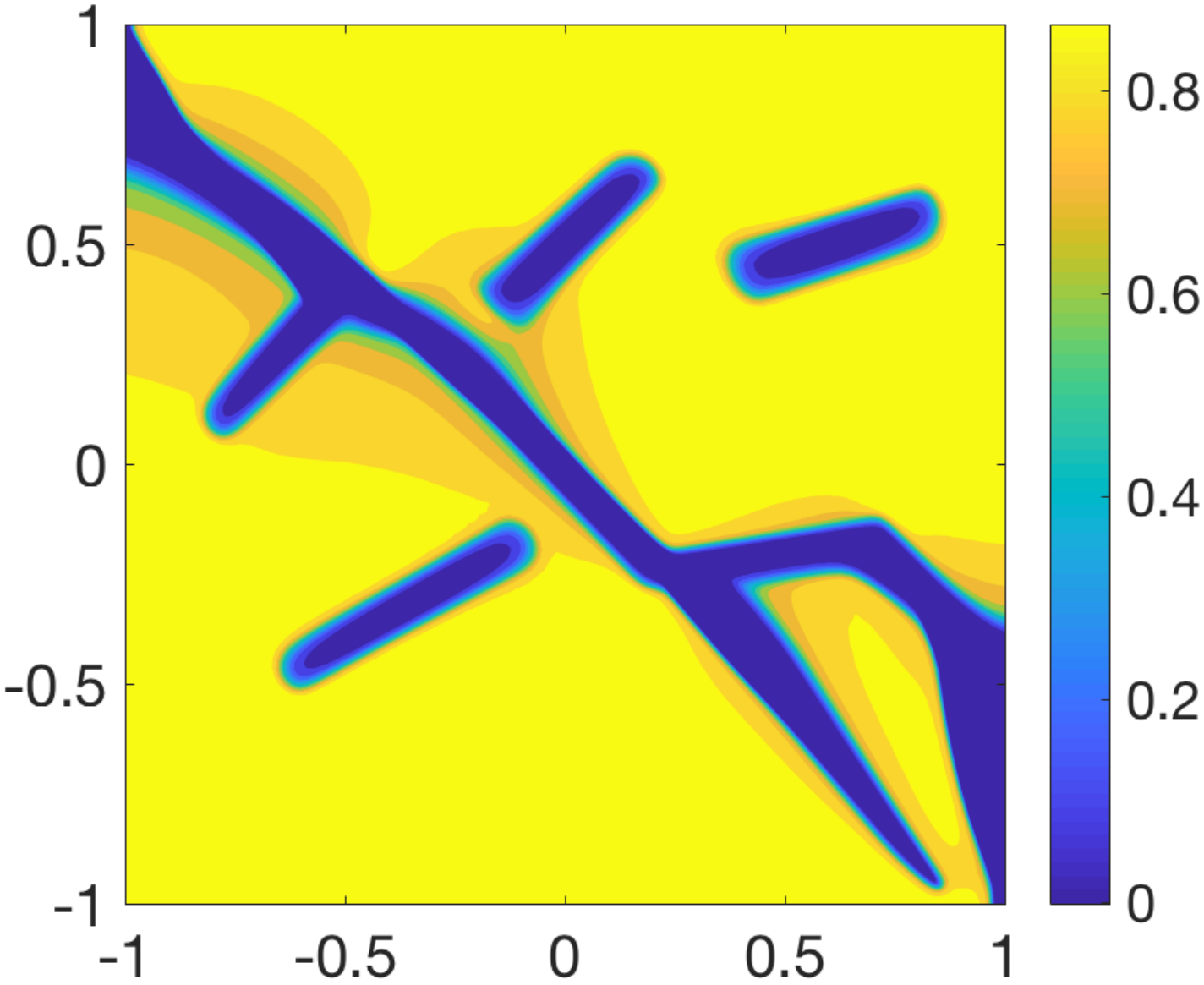}}
\vfill./
\subfigure[$U = 0.179$~mm]{\label{fig:subfig:FOMS_U1790}
\includegraphics[width=0.18\linewidth]{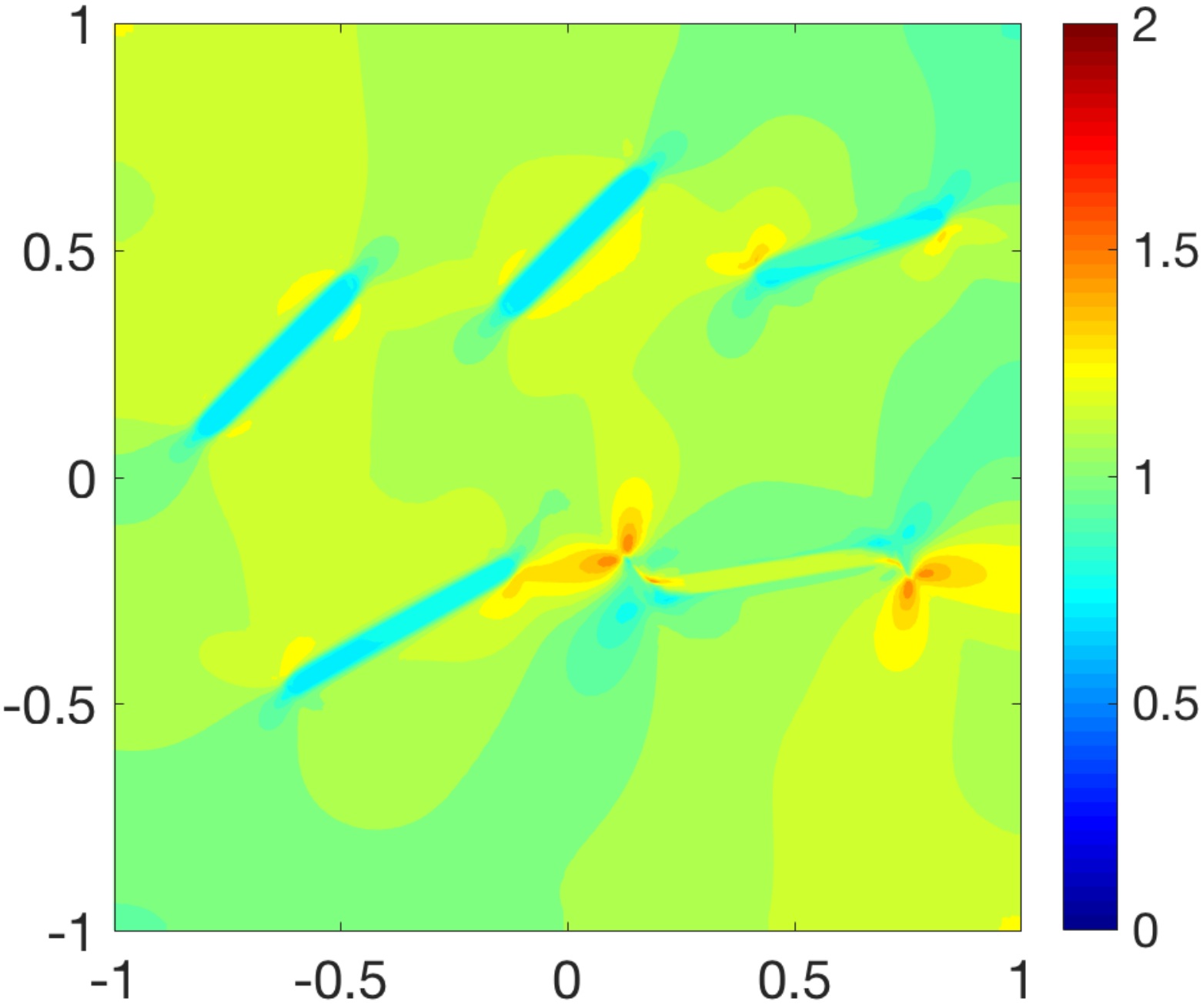}}
\subfigure[$U = 0.182$~mm]{\label{fig:subfig:FOMS_U1820}
\includegraphics[width=0.18\linewidth]{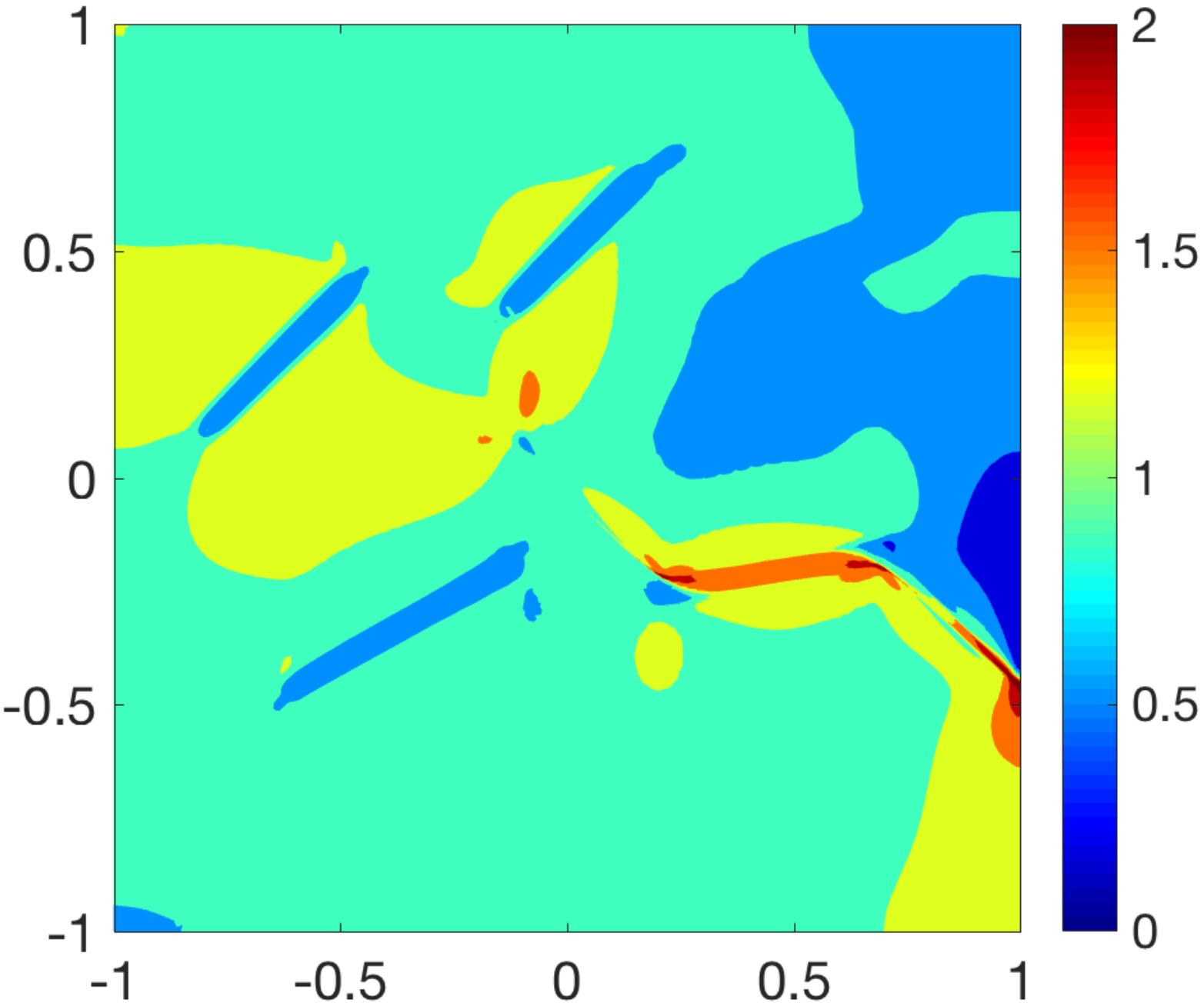}}
\subfigure[$U = 0.184$~mm]{\label{fig:subfig:FOMS_U1840}
\includegraphics[width=0.18\linewidth]{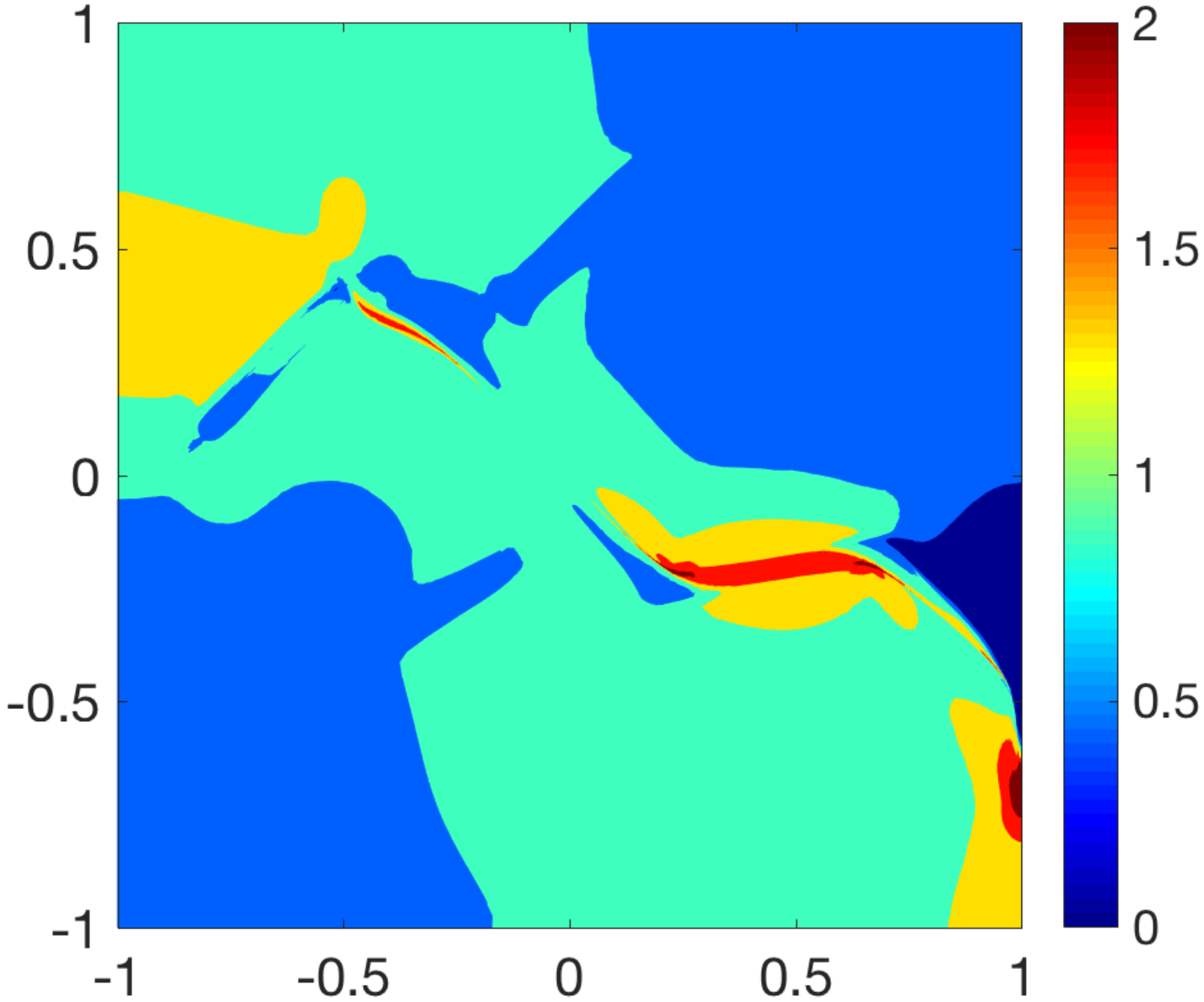}}
\subfigure[$U = 0.186$~mm]{\label{fig:subfig:FOMS_U1860}
\includegraphics[width=0.18\linewidth]{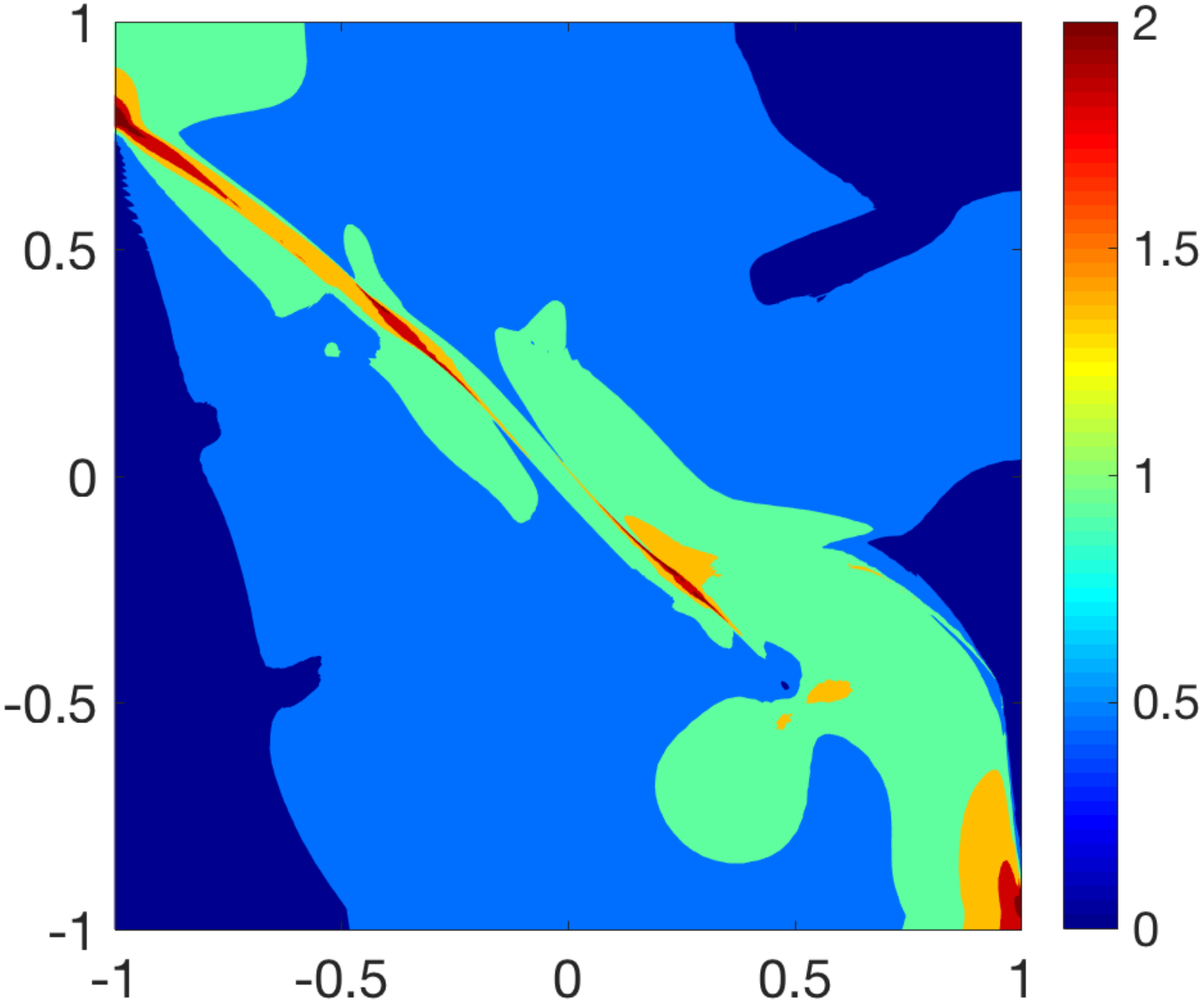}}
\subfigure[$U = 0.2$~mm]{\label{fig:subfig:FOMS_U2000}
\includegraphics[width=0.18\linewidth]{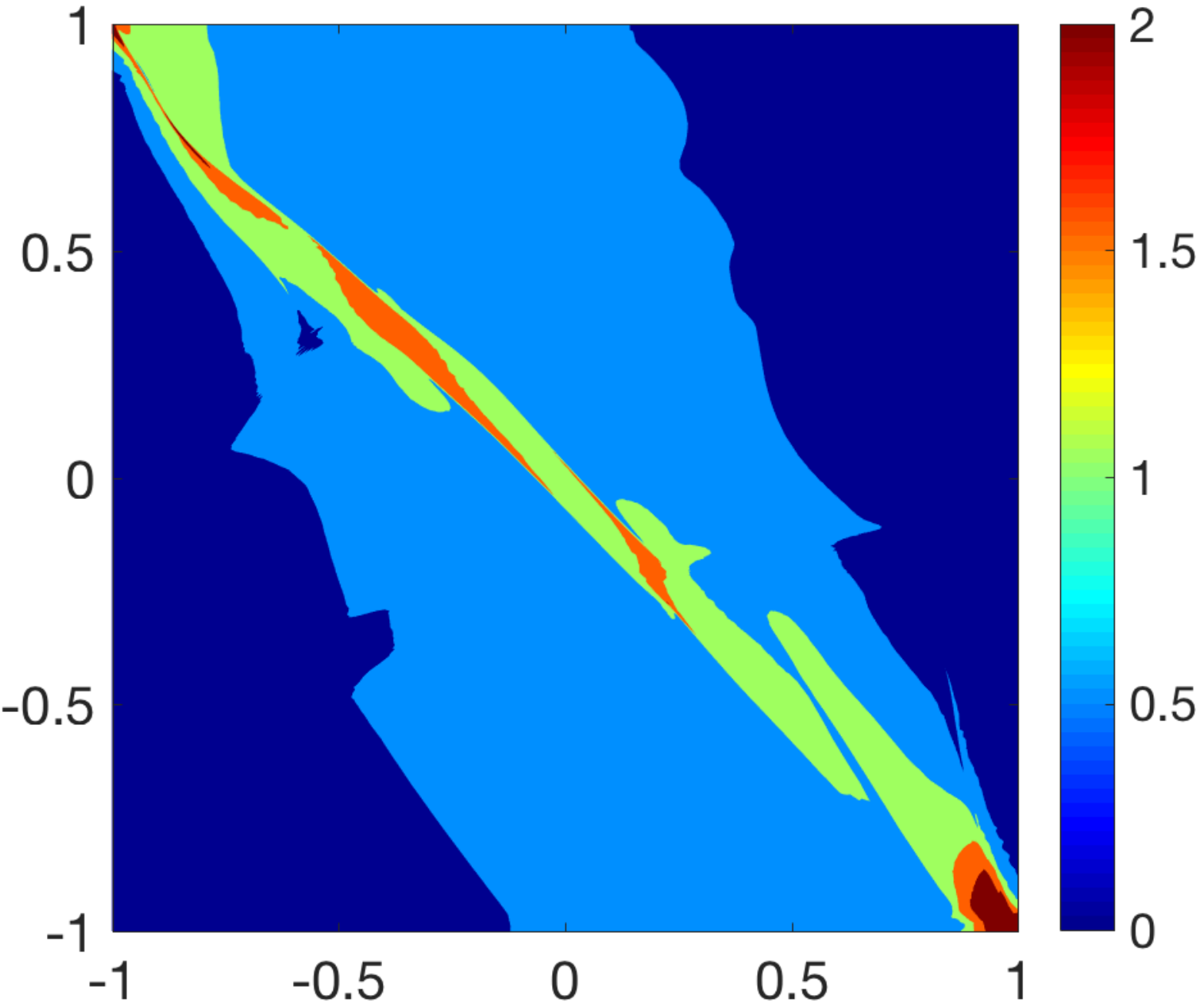}}
\caption{Example 3. The contour of the phase-field during crack evolution for the five-crack shear test with $l = 0.00375$~mm, $N = 25,600\; (81\times81)$.
(spectral decomposition with original crack boundary conditions)}
\label{fig:FSOBC}
\end{figure}

\begin{figure} 
\centering 
\subfigure[$U = 0.092$~mm]{\label{fig:subfig:FMMD_U920}
\includegraphics[width=0.18\linewidth]{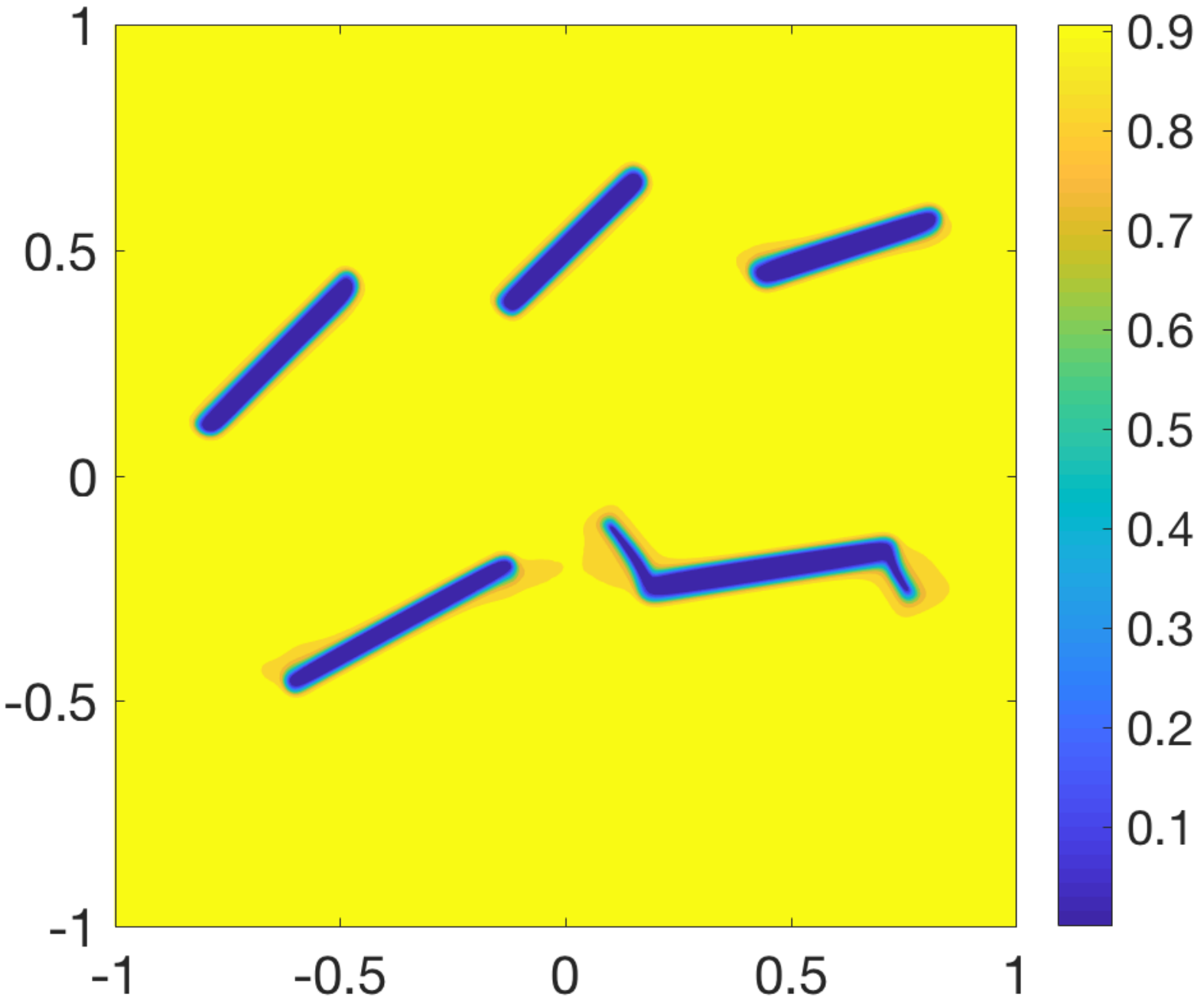}}
\subfigure[$U = 0.097$~mm]{\label{fig:subfig:FMMD_U970}
\includegraphics[width=0.18\linewidth]{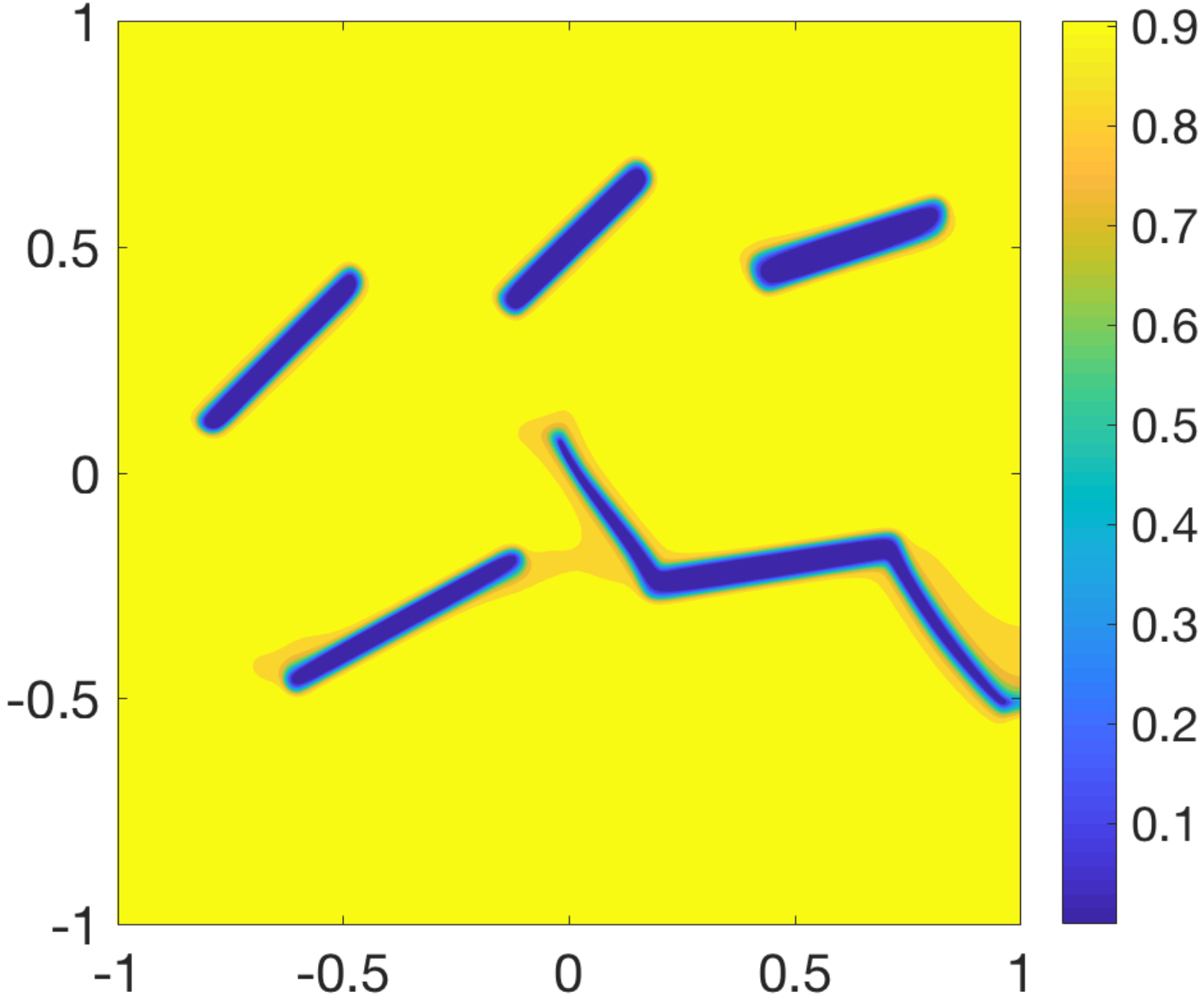}}
\subfigure[$U = 0.1$~mm]{\label{fig:subfig:FMMD_U1000}
\includegraphics[width=0.18\linewidth]{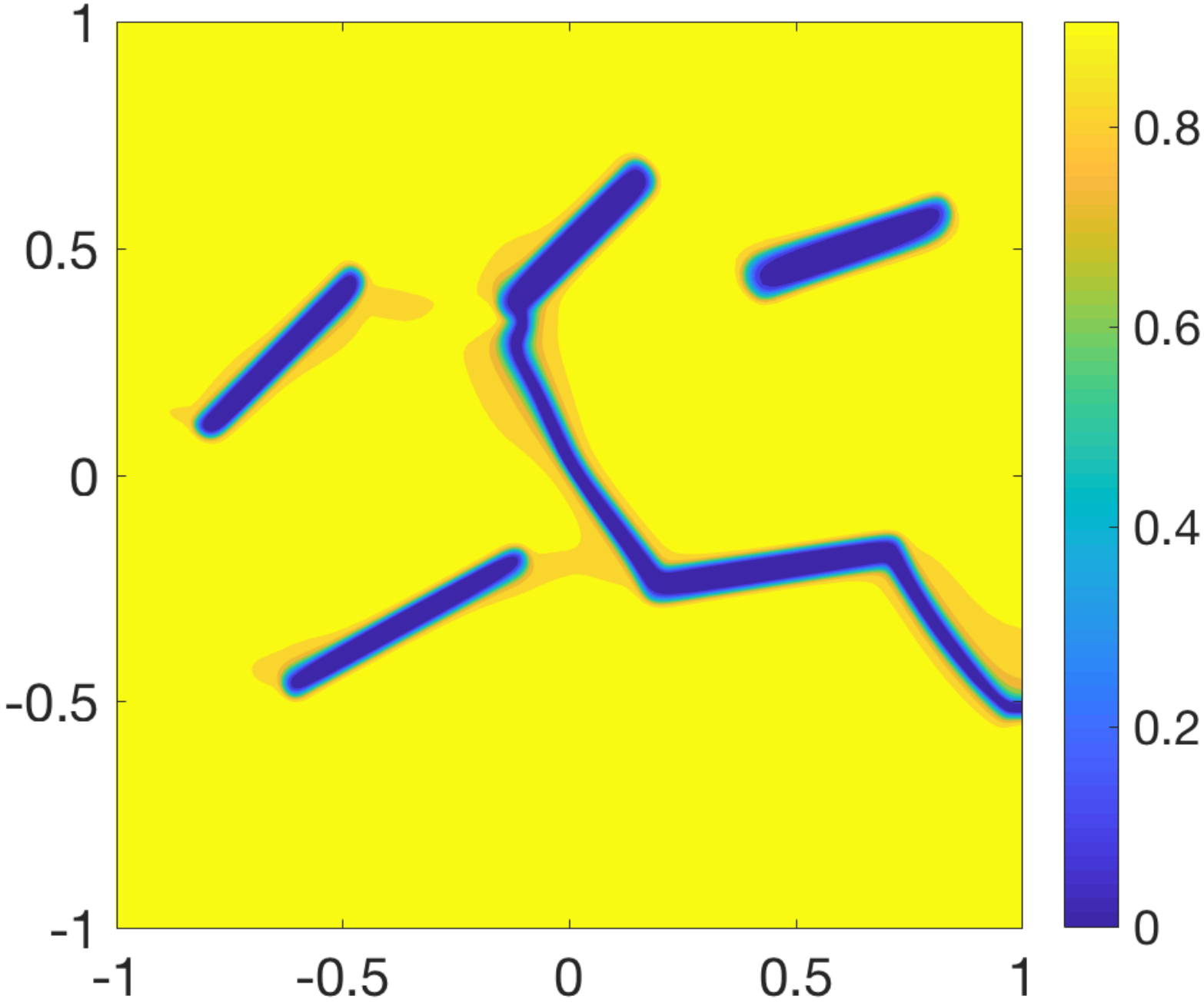}}
\subfigure[$U = 0.127$~mm]{\label{fig:subfig:FMMD_U1270}
\includegraphics[width=0.18\linewidth]{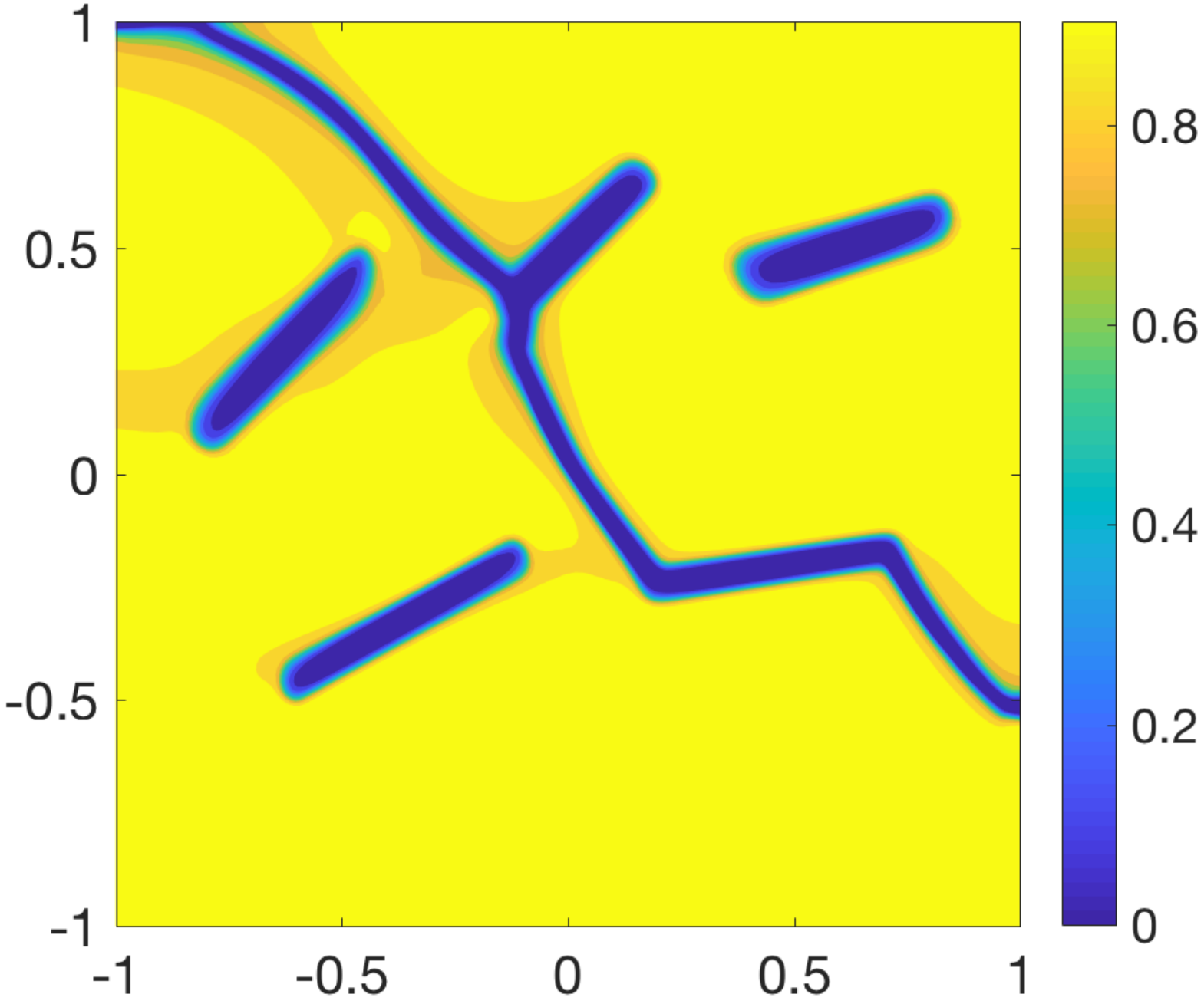}}
\subfigure[$U = 0.2$~mm]{\label{fig:subfig:FMMD_U2000}
\includegraphics[width=0.18\linewidth]{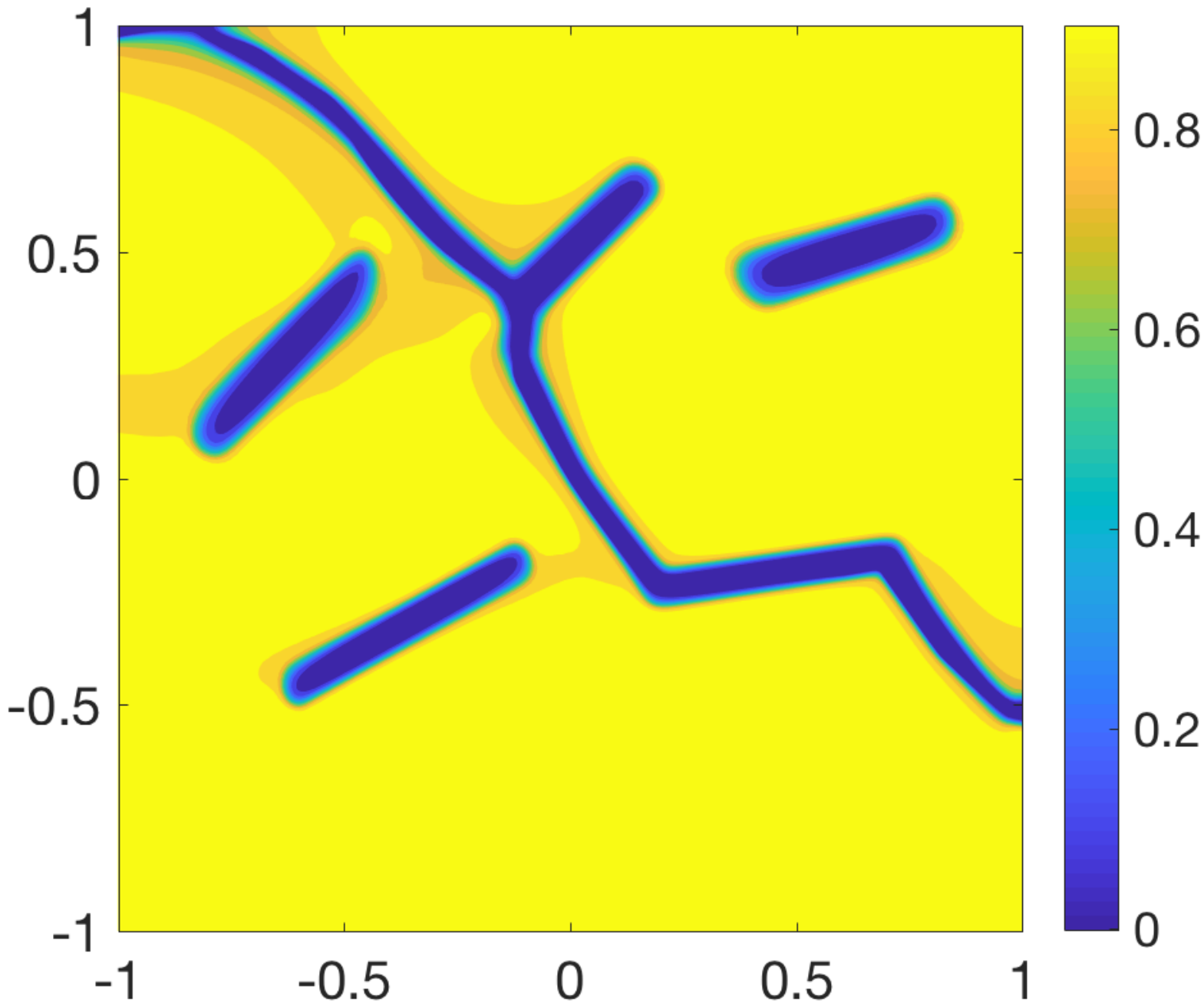}}
\vfill
\subfigure[$U = 0.092$~mm]{\label{fig:subfig:FMMS_U920}
\includegraphics[width=0.18\linewidth]{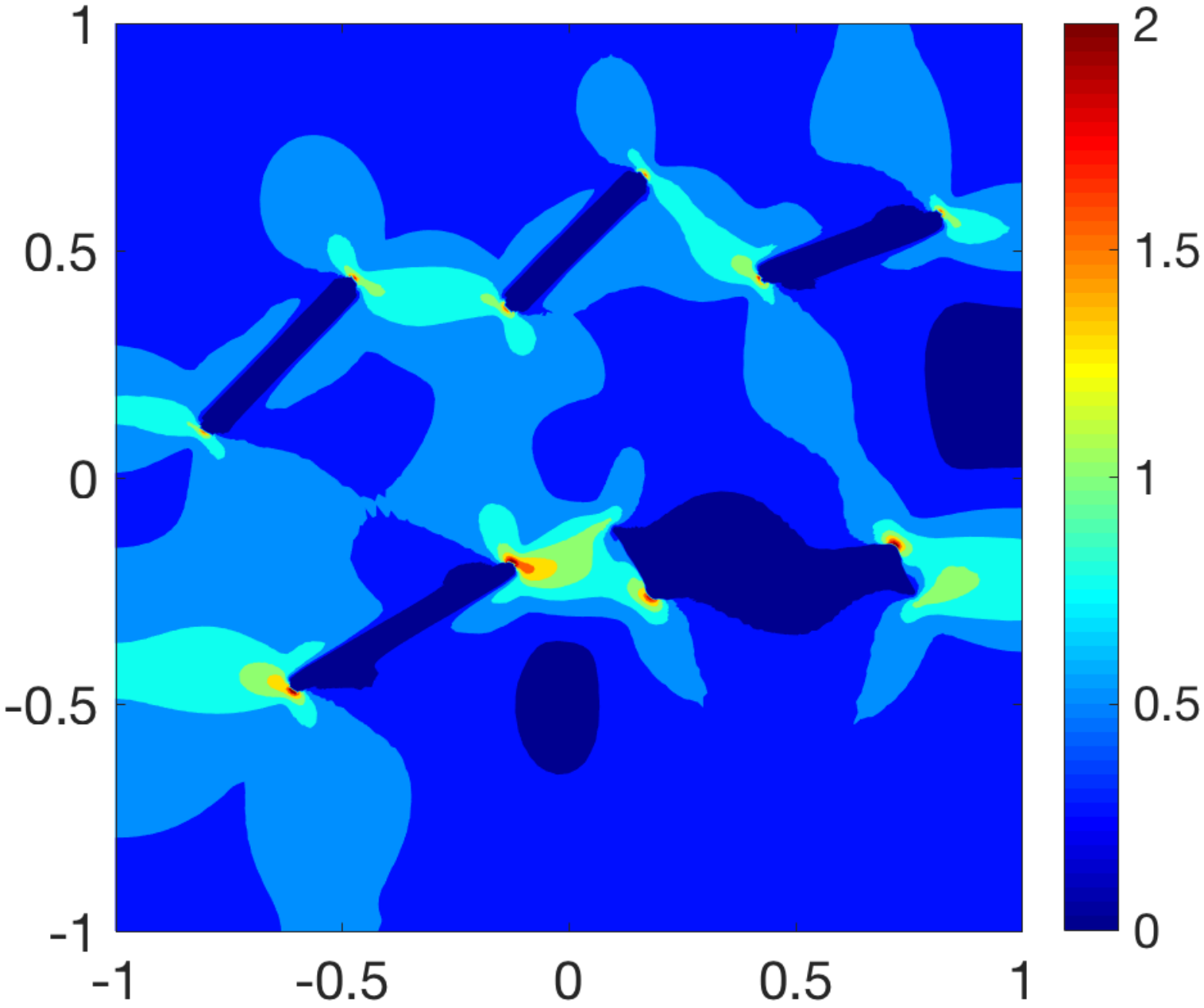}}
\subfigure[$U = 0.097$~mm]{\label{fig:subfig:FMMS_U970}
\includegraphics[width=0.18\linewidth]{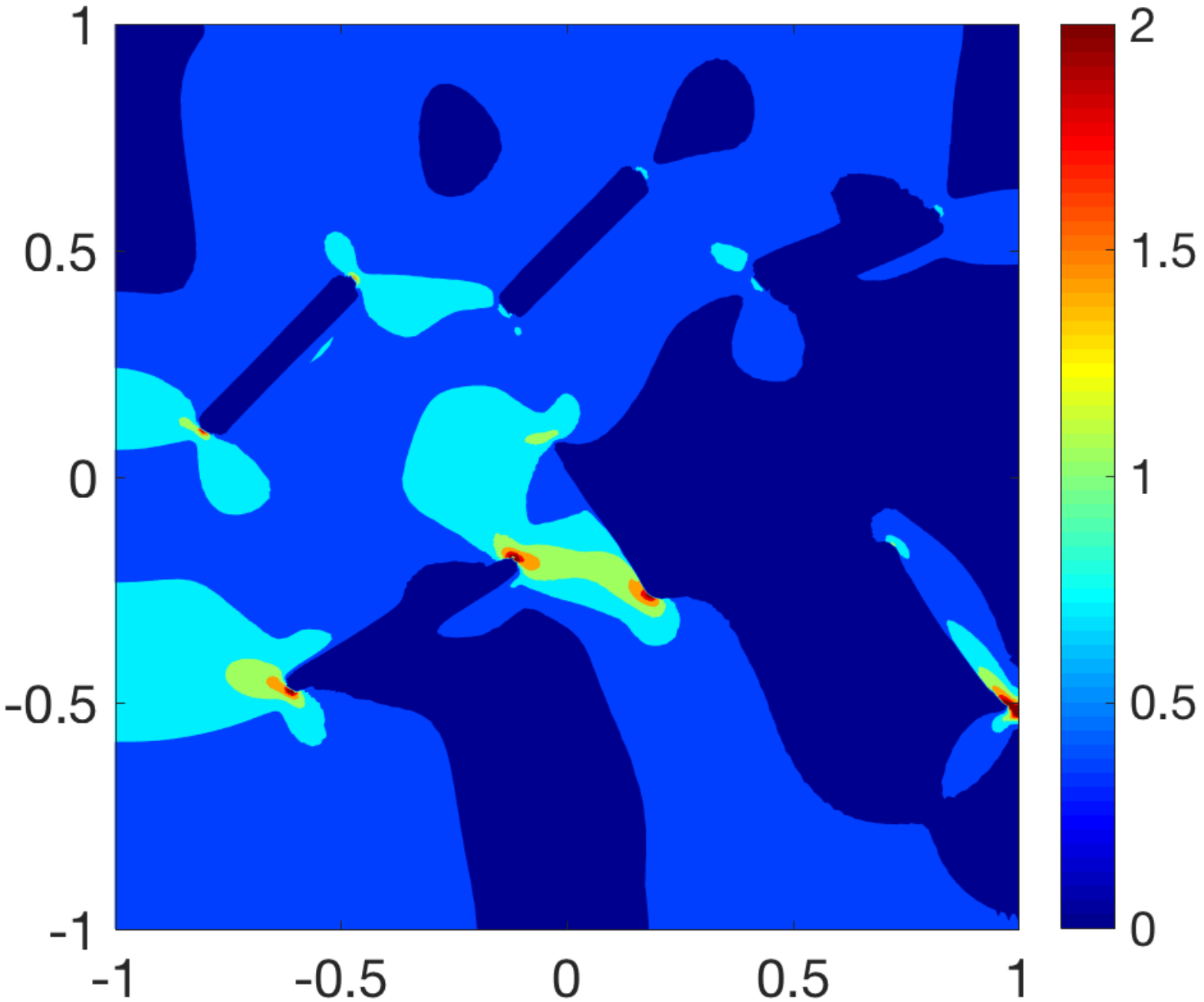}}
\subfigure[$U = 0.1$~mm]{\label{fig:subfig:FMMS_U1000}
\includegraphics[width=0.18\linewidth]{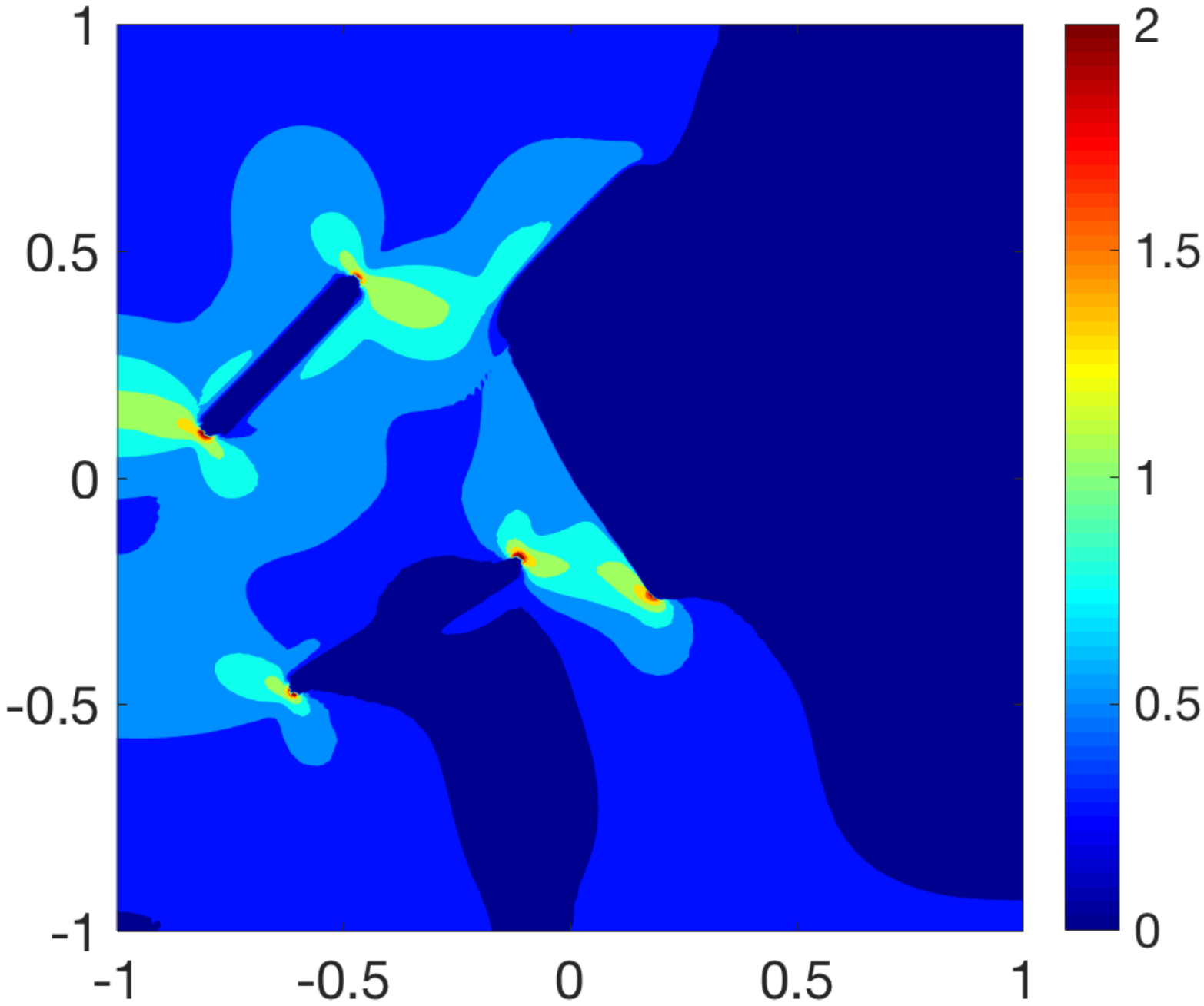}}
\subfigure[$U = 0.127$~mm]{\label{fig:subfig:FMMS_U1270}
\includegraphics[width=0.18\linewidth]{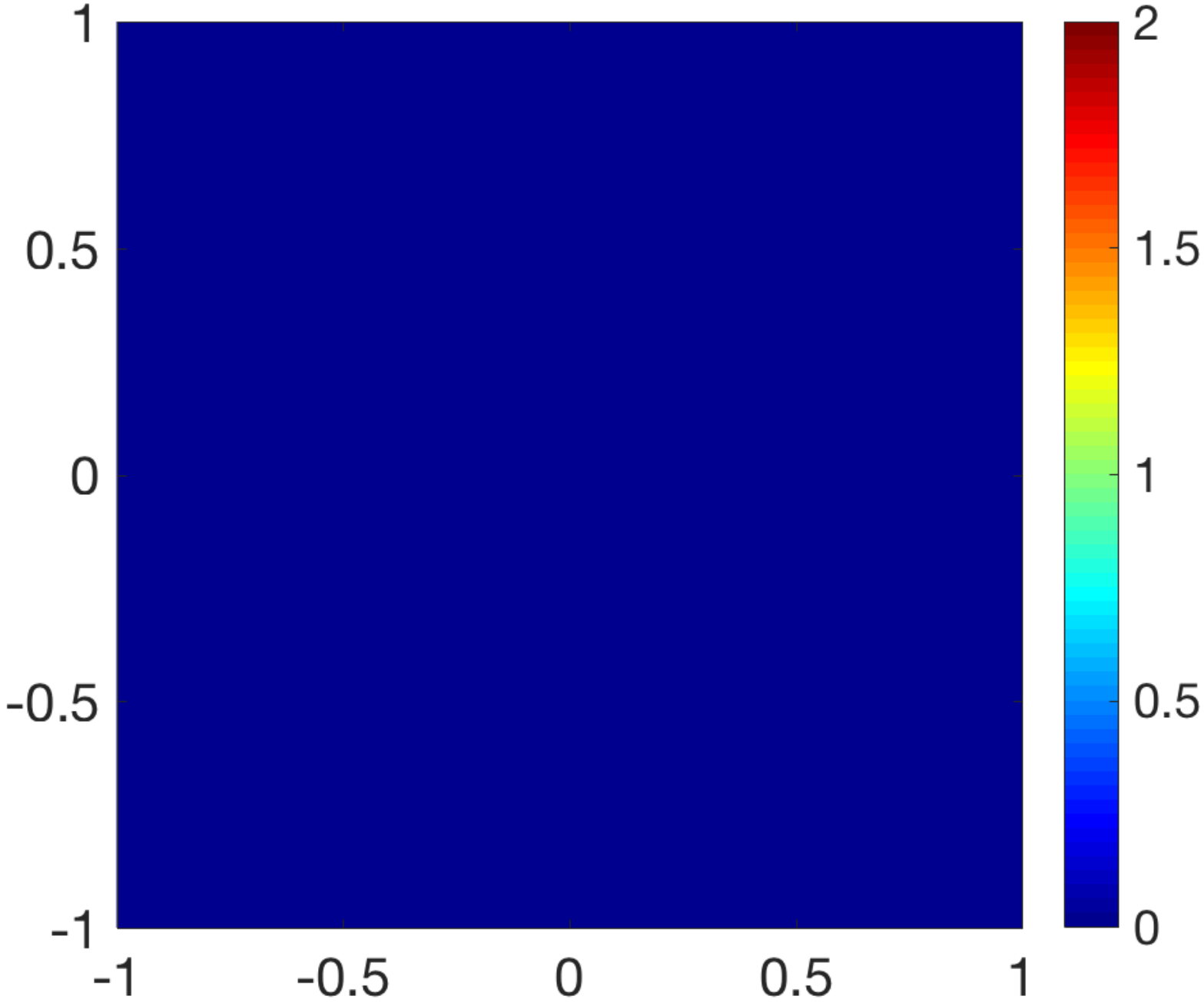}}
\subfigure[$U = 0.2$~mm]{\label{fig:subfig:FMMS_U2000}
\includegraphics[width=0.18\linewidth]{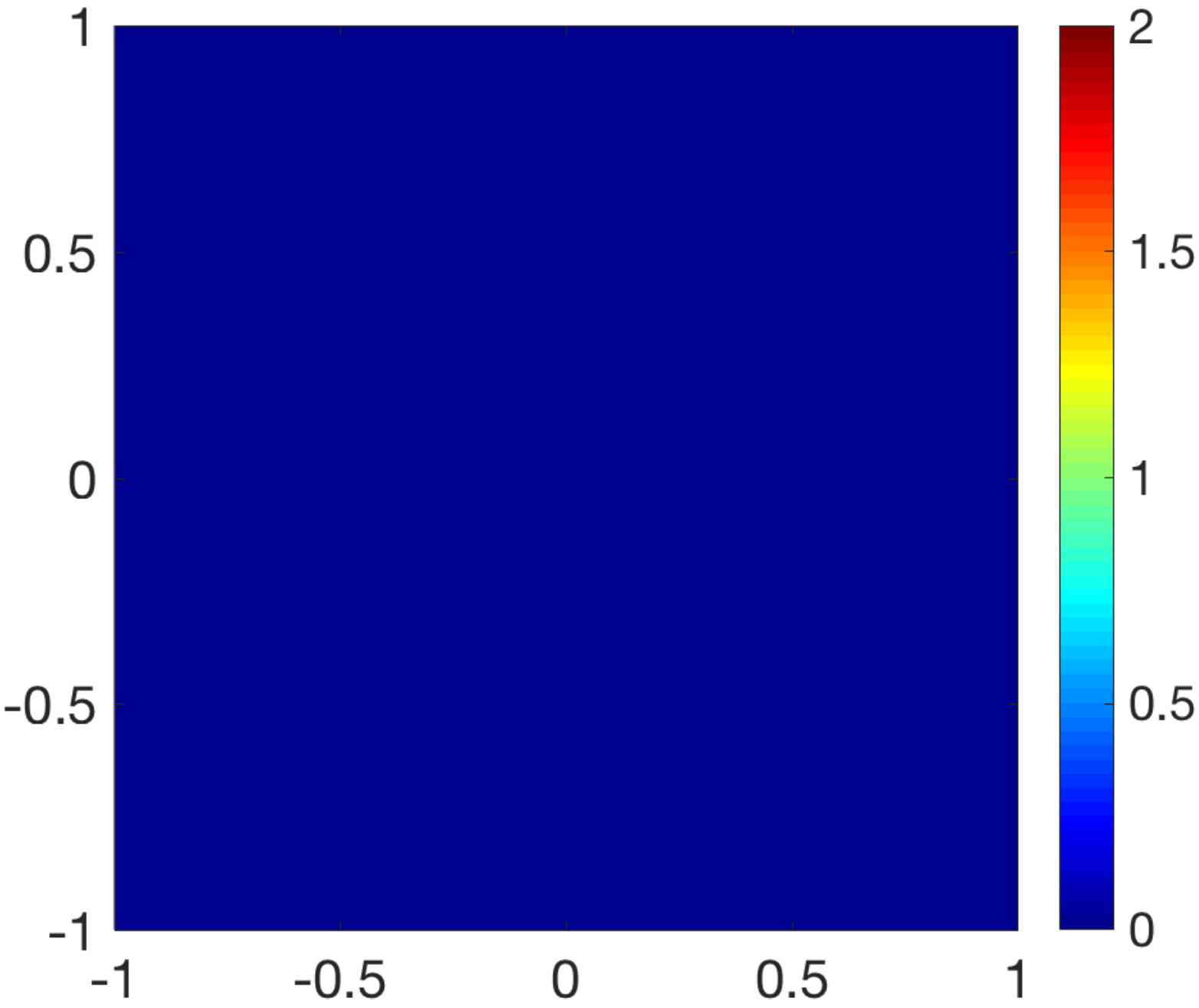}}
\caption{Example 3. The contour of the phase-field during
crack evolution for the five-crack shear test with $l = 0.00375$~mm, $N = 25,600\; (81\times81)$.
(spectral decomposition with ItCBC)}
\label{fig:FSMBC}
\end{figure}

For the five-crack problem, the lengths of Crack 1 to Crack 5 are $0.3$~mm, $0.35$~mm, $0.35$~mm, $0.5$~mm and $0.5$~mm, with polar angle
$30^{\circ}$, $45^{\circ}$, $17.2^{\circ}$, $28.6^{\circ}$ and $9^{\circ}$, respectively. The phase-field and stress distribution for the spectral decomposition method without and with ItCBC
($d_{cr} = 0.4$) are shown in Fig. \ref{fig:FSOBC} and \ref{fig:FSMBC}.
Again, unphysical stress concentration in the damaged zone is observed for the original treatment
of crack boundary conditions but vanishes with ItCBC.

\subsection{Example 4. A single edge notched shear test based on experiment}

\begin{figure} 
\centering 
\subfigure[experiment setup]{\label{fig:subfig:Experimental_shear}
\includegraphics[width=0.45\linewidth]{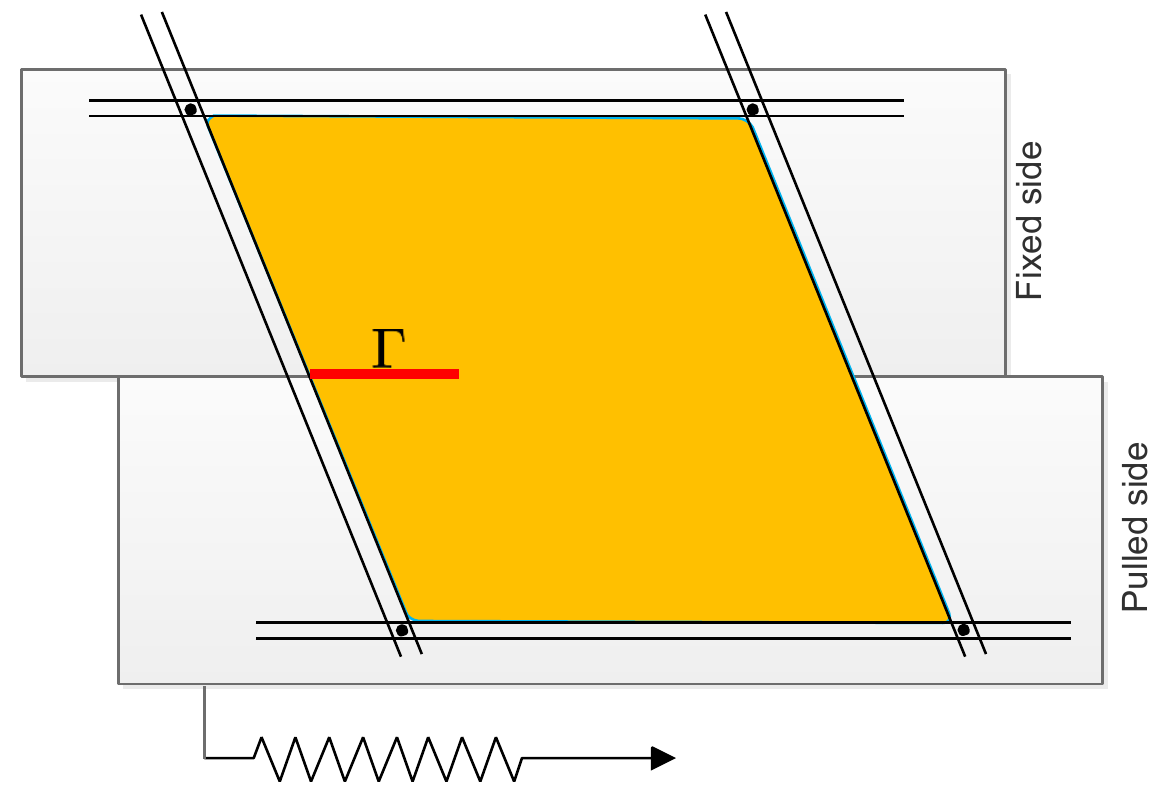}}
\subfigure[computation setup]{\label{fig:subfig:Computational_shear}
\includegraphics[width=0.5\linewidth]{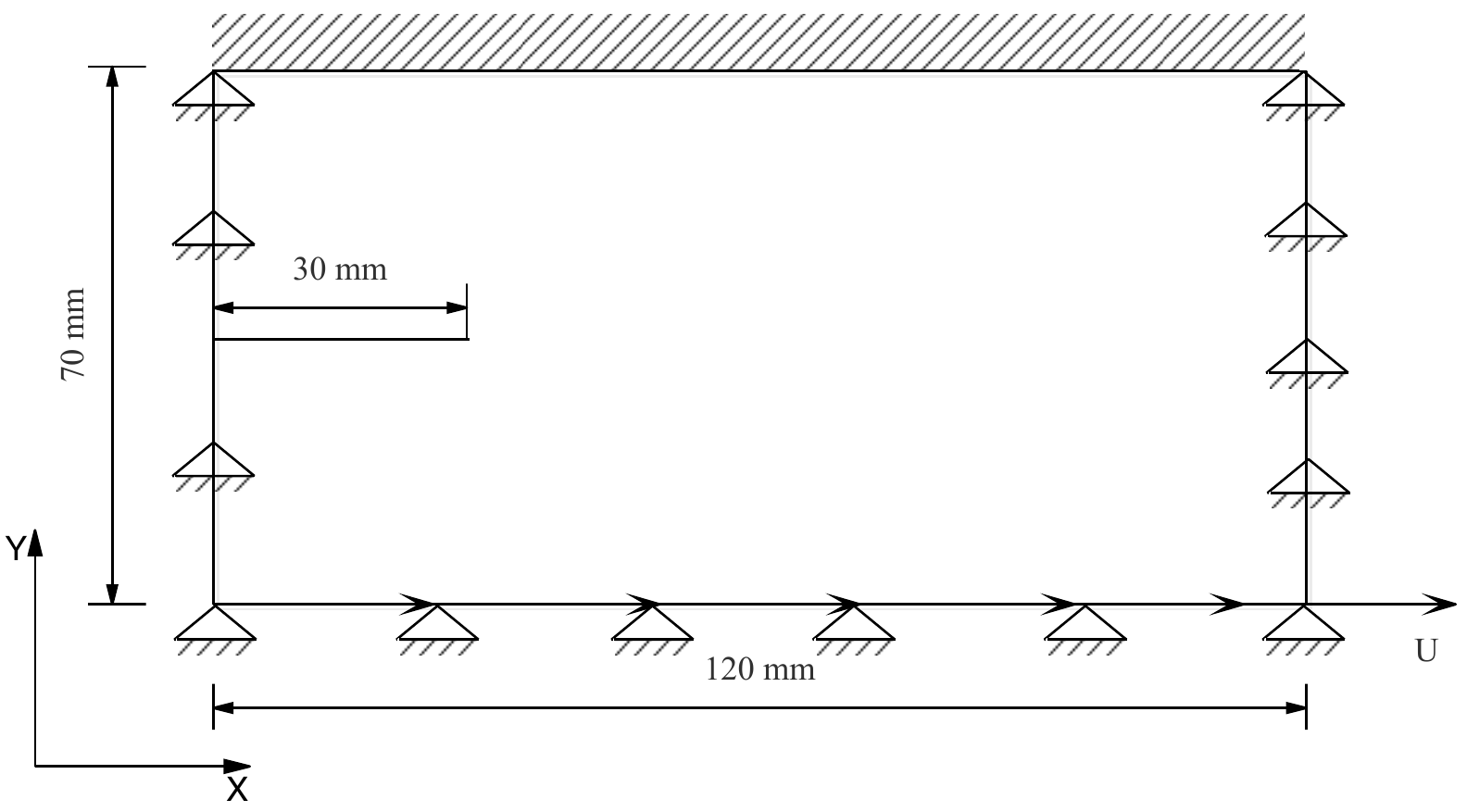}}
\caption{Example 4. Sketch of the experiment and computation setup.}
\label{fig:sketch_EC}
\end{figure}

\begin{figure} 
\centering
Physical experiment:\\ 
\subfigure[$U = 3.5$~mm]{\label{fig:subfig:Lee-1}
\includegraphics[width=0.28\linewidth]{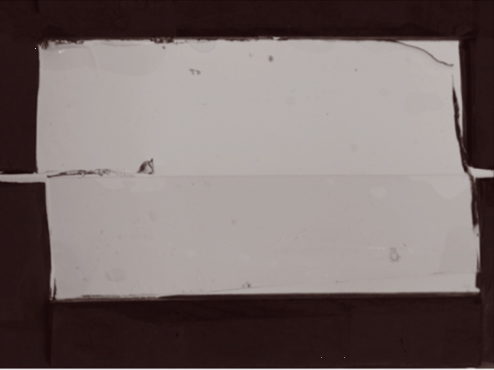}}
\subfigure[$U = 5.5$~mm]{\label{fig:subfig:Lee-2}
\includegraphics[width=0.28\linewidth]{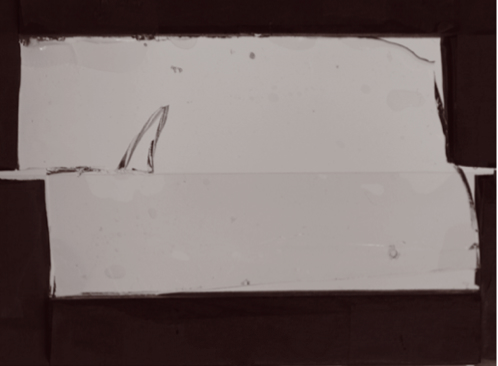}}
\subfigure[$U = 9$~mm]{\label{fig:subfig:Lee-3}
\includegraphics[width=0.28\linewidth]{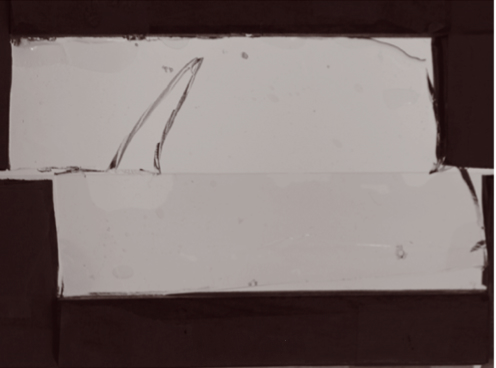}}
\vfill
Spectral decomposition:\\
\subfigure[$U = 3.8$~mm]{\label{fig:subfig:ExperM_D760}
\includegraphics[width=0.3\linewidth]{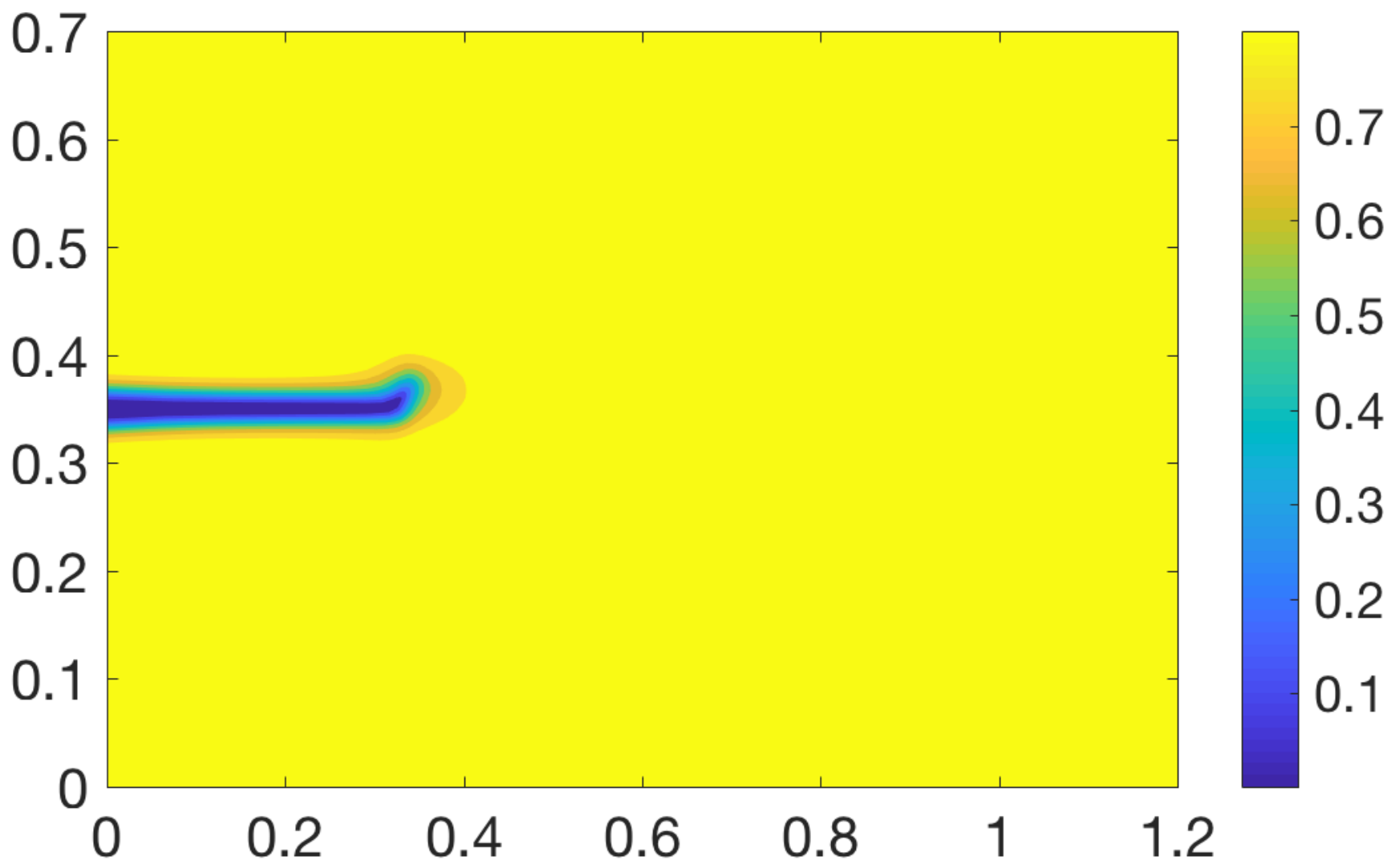}}
\subfigure[$U = 4$~mm]{\label{fig:subfig:ExperM_D800}
\includegraphics[width=0.3\linewidth]{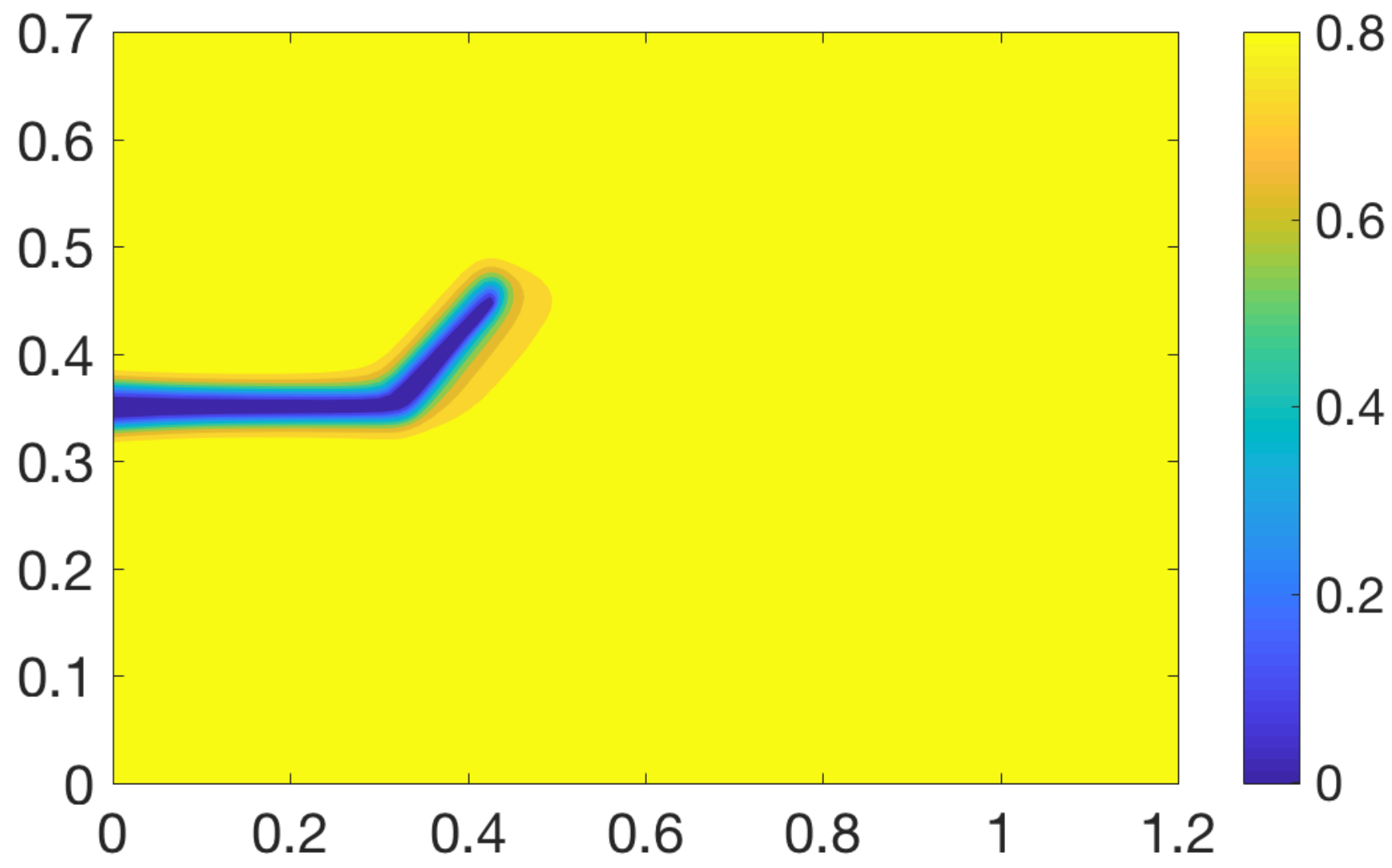}}
\subfigure[$U = 4.4$~mm]{\label{fig:subfig:ExperM_D880}
\includegraphics[width=0.3\linewidth]{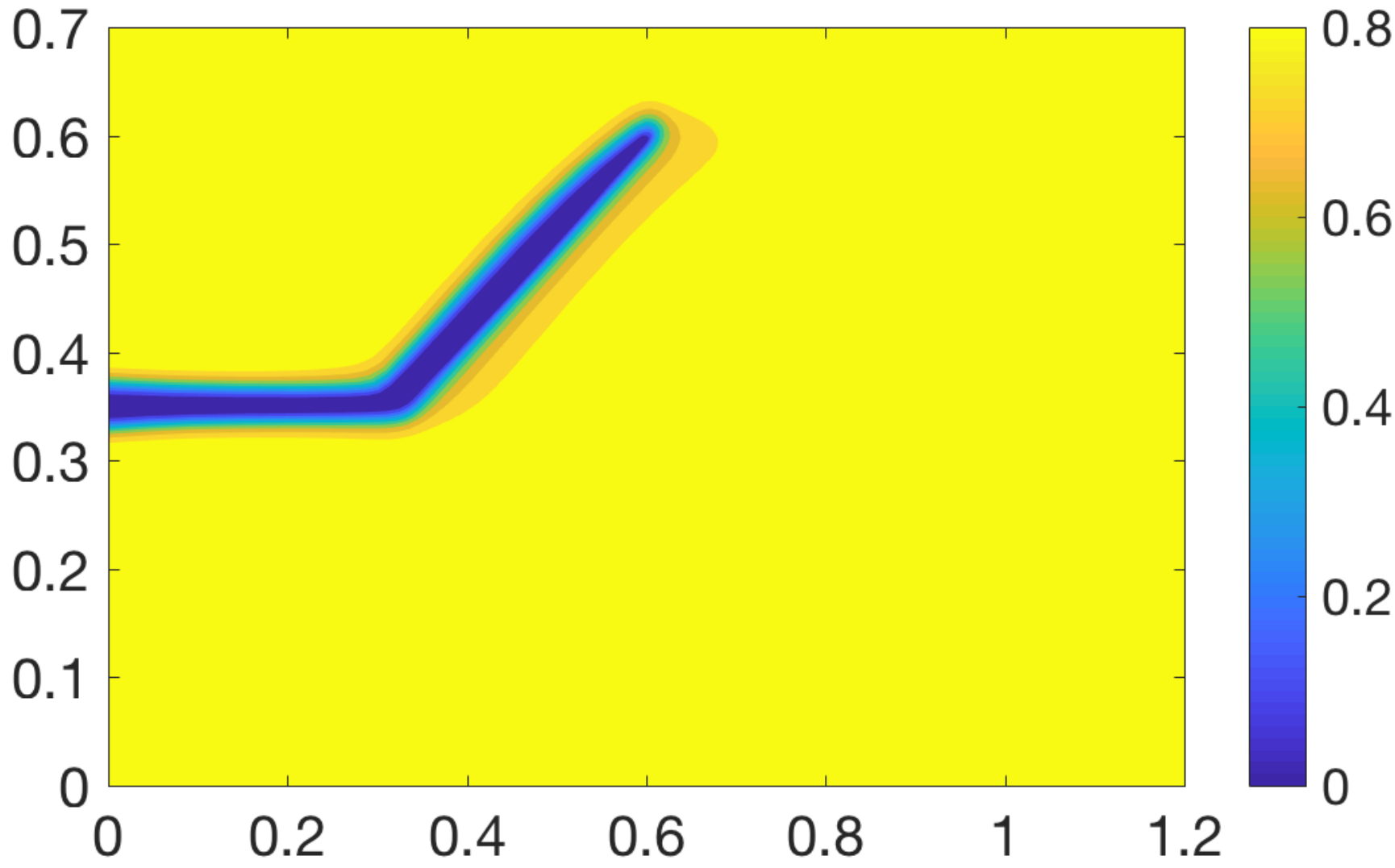}}
\vfill
v-d split:\\
\subfigure[$U = 3.1$~mm]{\label{fig:subfig:ExperA_D620}
\includegraphics[width=0.3\linewidth]{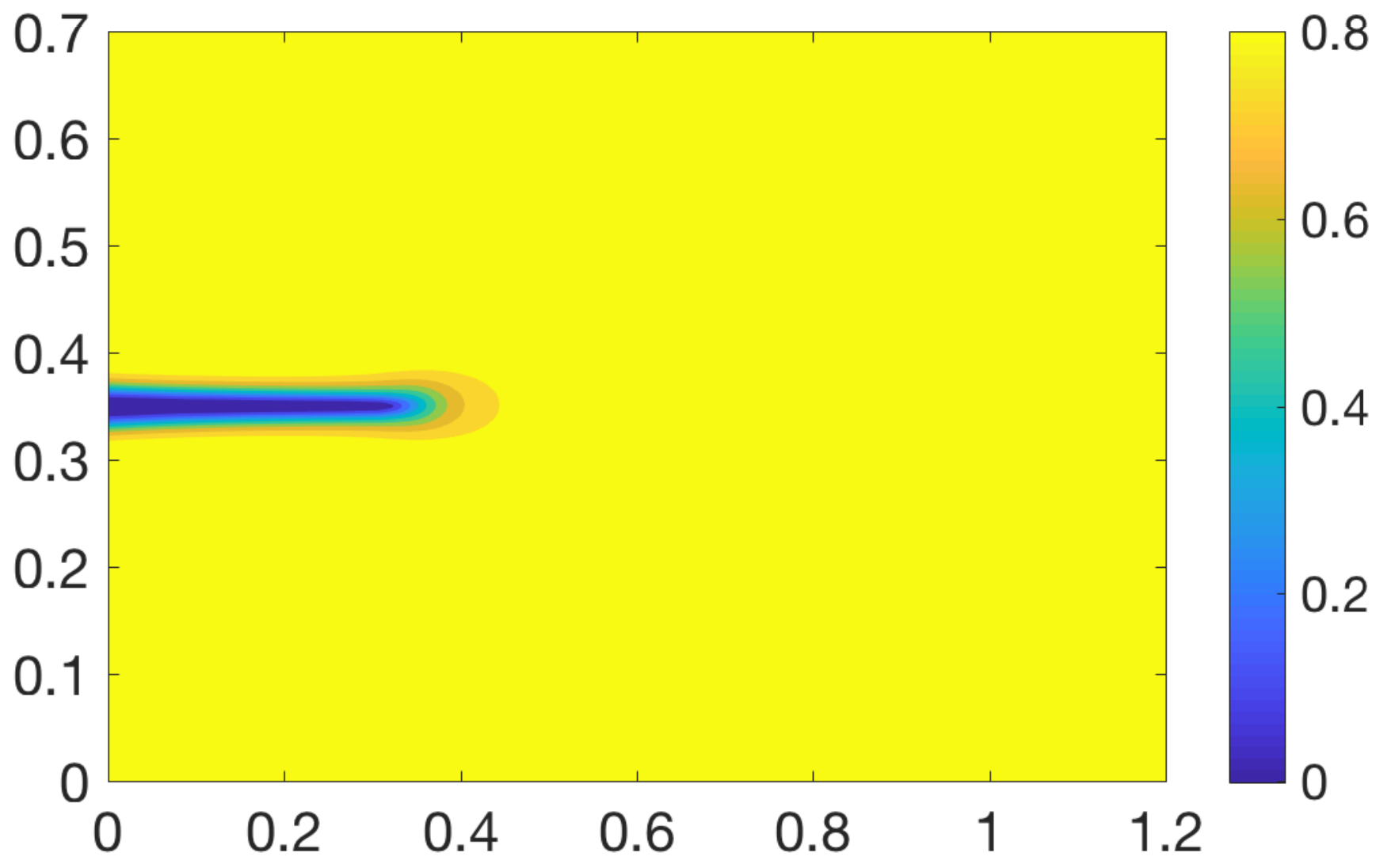}}
\subfigure[$U = 3.3$~mm]{\label{fig:subfig:ExperA_D660}
\includegraphics[width=0.3\linewidth]{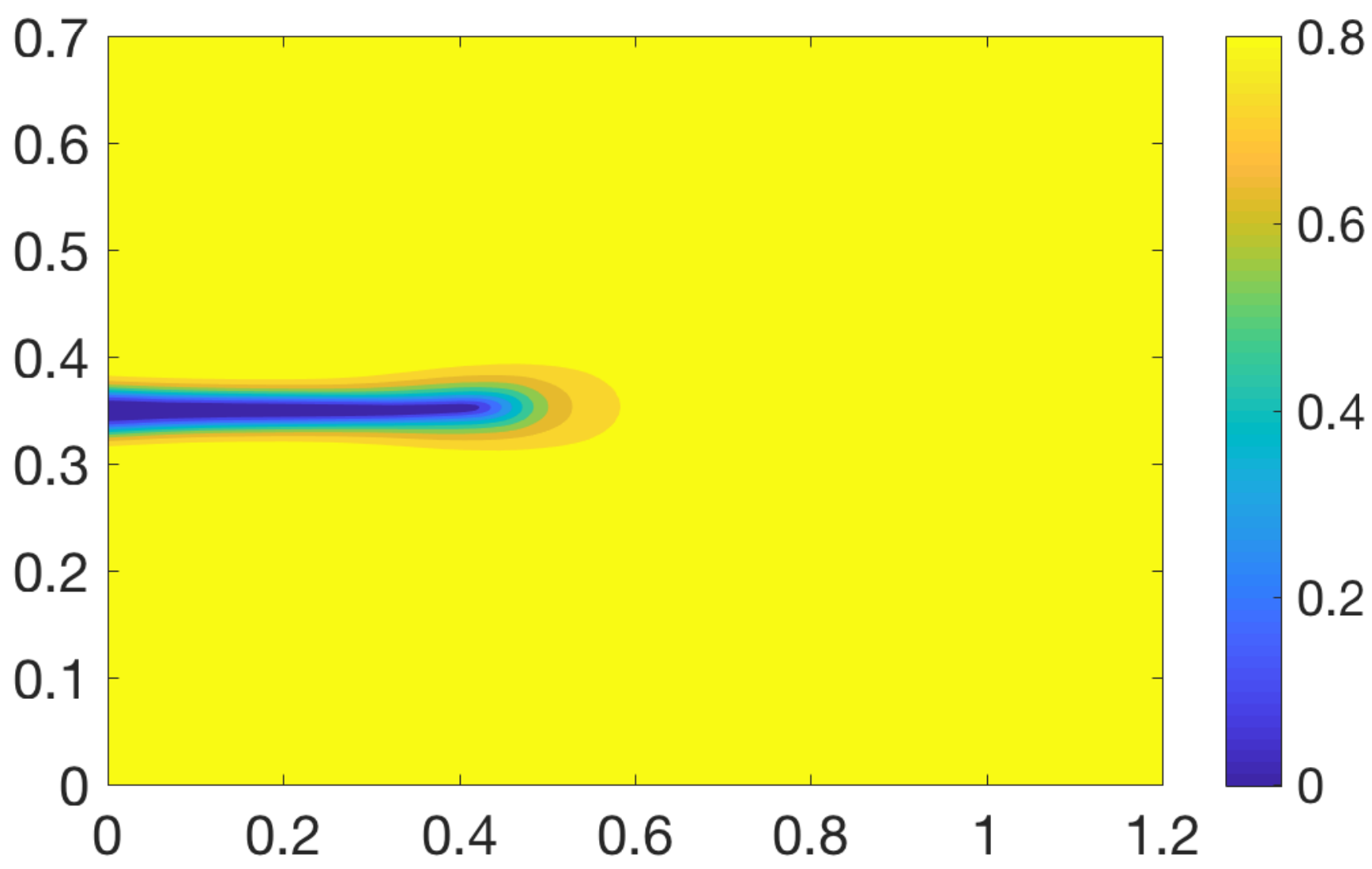}}
\subfigure[$U = 3.5$~mm]{\label{fig:subfig:ExperA_D700}
\includegraphics[width=0.3\linewidth]{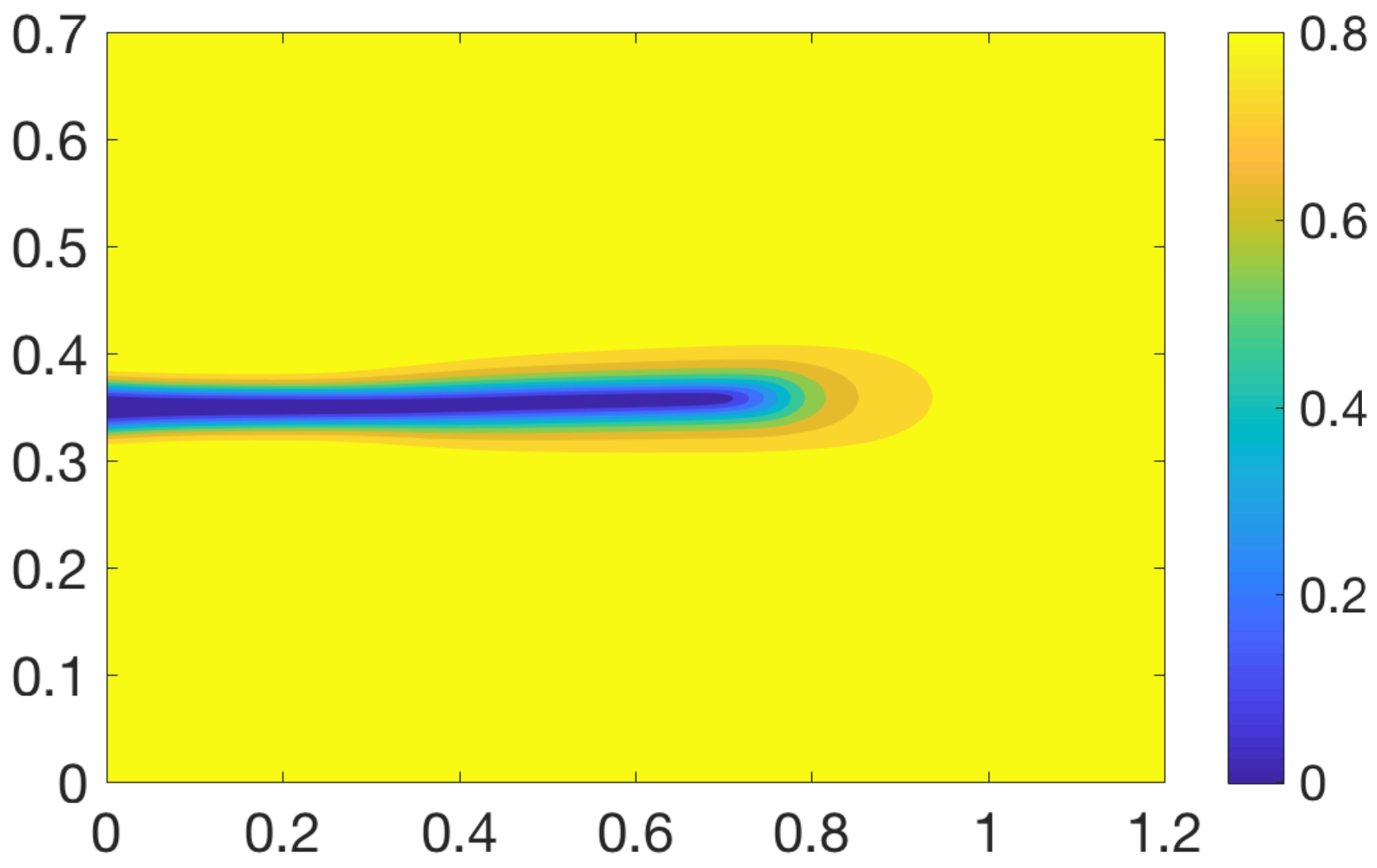}}\\
\vfill
Improved v-d split:\\
\subfigure[$U = 4.1$~mm]{\label{fig:subfig:ExperN_D820}
\includegraphics[width=0.3\linewidth]{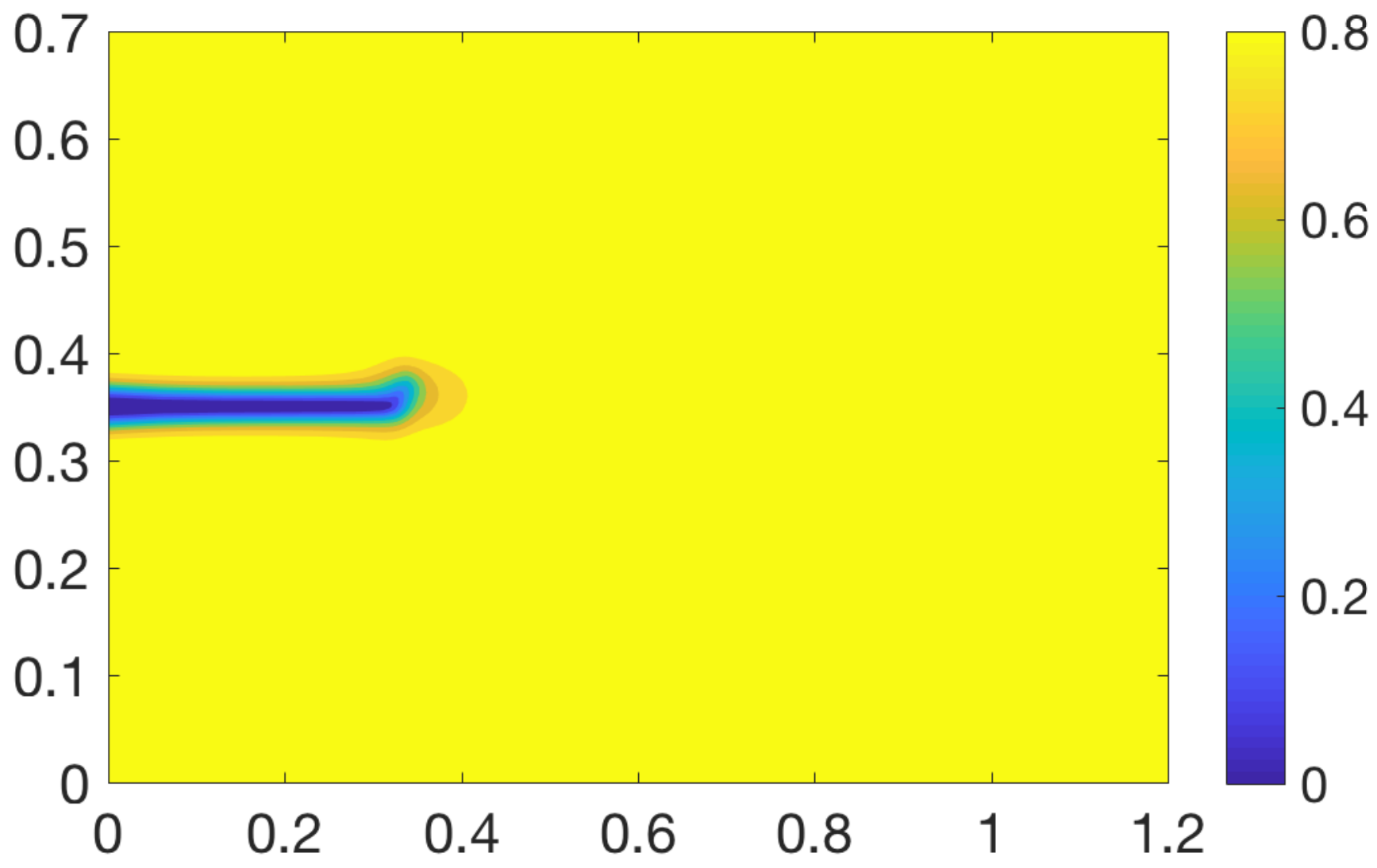}}
\subfigure[$U = 4.3$~mm]{\label{fig:subfig:ExperN_D860}
\includegraphics[width=0.3\linewidth]{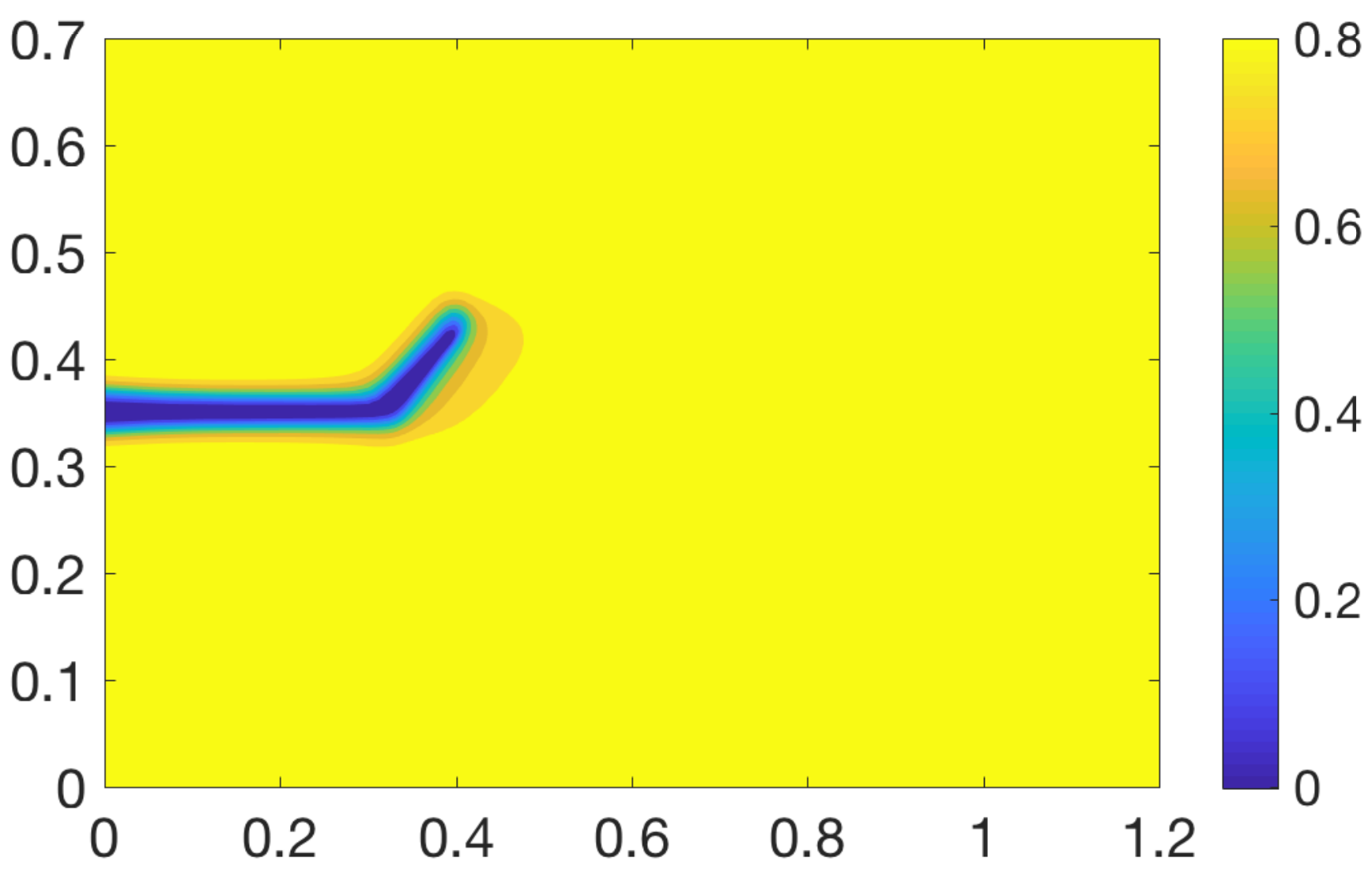}}
\subfigure[$U = 4.7$~mm]{\label{fig:subfig:ExperN_D940}
\includegraphics[width=0.3\linewidth]{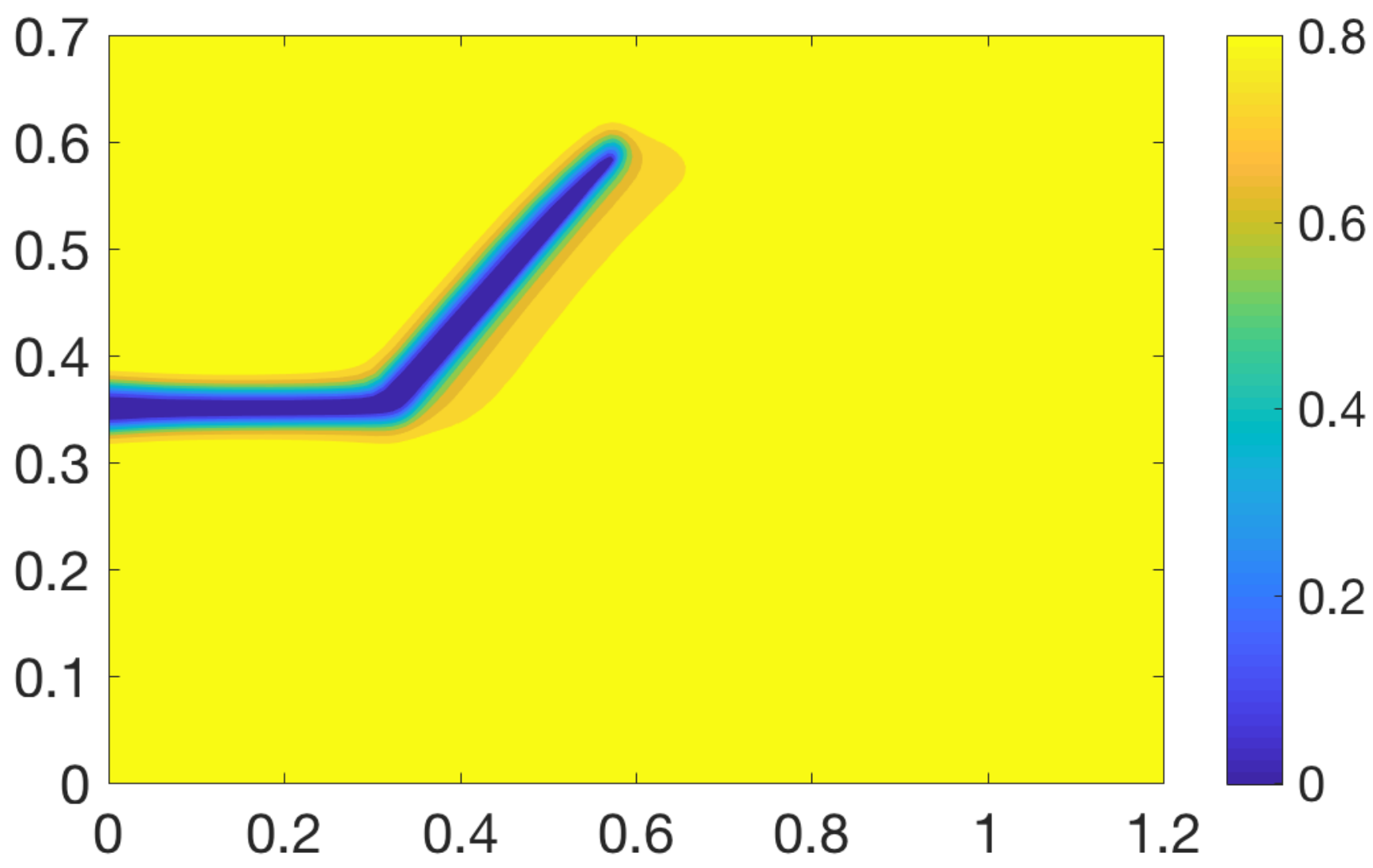}}
\caption{Example 4. Comparison of geometry of the fracture obtained from
the numerical computation with different decomposition models (with ItCBC $d_{cr} = 0.4$)
and the experiment \cite{LRHW16}.}
\label{fig:Comp-experimental}
\end{figure}

To further verify the decomposition models, we compare the simulation results with the 
experimental results obtained by Lee et al. \cite{LRHW16}.
The physical experiments were performed on a modified version of the shear apparatus
used in Reber et al. \cite{RLH15}. Gelatin was used as the material in the experiments.
The sketch of the experimental setup can be seen in Fig. \ref{fig:subfig:Experimental_shear}.
The spring black arrow represents the direction of the shear force.
The bottom side of the table moves under the constant velocity while the fixed side is stationary.
For the numerical computation, a rectangular plate with a length of $120$~mm and a width of $70$~mm is considered.
The initial fracture is located at the middle of the plate with the length of $30$~mm.
The top edge of the domain is fixed and the bottom edge is fixed along $y$-direction while a uniform 
$x$-displacement $U$ is increased with time ($\Delta U = 5 \times 10^{-3}$~mm).
We consider a mesh size as $41 \times 41$ ($N = 6,400$) and $l = 1.2$~mm.
The material properties are taken almost the same as the physical experiment:
the Young's modulus is $1.4 \times 10^{5}$~Pa, the Poisson ratio is $0.45$
and the fracture toughness is $1.96$~Pa $\cdot$ m.
(The Poisson ratio for gelatin is 0.499. We choose the Poisson ratio to be 0.45 to avoid
the locking effects in the finite element approximation of elasticity problems where
the performance of certain commonly used finite elements deteriorates when the Poisson ratio
is close to 0.5; e.g. see Babu$\check{s}$ka and Suri \cite{Babuska92}.)
These material properties correspond to
$\lambda = 4.345\times 10^{-4} $~kN/mm$^{2}$, $\mu = 4.829\times 10^{-5} $~kN/mm$^{2}$,
and $g_c = 1.96 \times 10^{-6}$~kN/mm.
The results from the computation with three decomposition models and
the experiment are comparable qualitatively in crack propagation.
As can be seen in Fig. \ref{fig:Comp-experimental}, with ItCBC ($d_{cr}$ = 0.4),
the spectral decomposition and the improved volumetric-deviatoric split lead to results
comparable with the physical experiment.
The displacement load is almost the same ($U = 3.5$~mm in the experiments and $U = 3.8$~mm in the models)
when the crack starts propagation.
However, the crack propagation speed in the numerical results is faster than the experiment's
(as can be seen in Fig. \ref{fig:subfig:ExperM_D800}, \ref{fig:subfig:ExperM_D880}, \ref{fig:subfig:ExperN_D860}
and \ref{fig:subfig:ExperN_D940}, the crack quickly extends to the top boundary after its breaking).
This discrepancy has also been observed by Lee et al. \cite{LRHW16}.
This may be attributed to the fact that the assumptions used in the mathematical model such as
the perfectly homogeneous model material, friction less deformation (no friction exist along the crack surface)
and perfect boundary conditions may not accurately describe the experimental setting.
Nevertheless, this will be an interesting topic for future investigations.
Moreover, as can be seen in Fig. \ref{fig:subfig:ExperA_D620}, \ref{fig:subfig:ExperA_D660} and \ref{fig:subfig:ExperA_D700}, the volumetric-deviatoric split model leads to incorrect crack propagation.

./

\section{Conclusions}
\label{SEC:conclusions}

A crucial component in the phase-field modeling of brittle fracture is the decomposition of the energy
into the active and passive components, with only the former contributing to the damage process
that results in fracture and thus being degraded. An improper decomposition can lead to
unphysical crack openings and propagation. In the previous sections, we have studied
two commonly used decomposition models, the spectral decomposition model \cite{MHW10}
and the volumetric-deviatoric (v-d) split model \cite{AMM09}. The active energy consists of those
related to the tensile strain component in the former model and the expansive volumetric strain energy
plus the total deviatoric strain energy in the latter model. Four examples have been used for those
studies. The first two are classical benchmark problems, single edge notched tension and shear tests.
The third example is designed to demonstrate the ability of the models to handle complex cracks.
The last example is also a single edge notched shear test but with the physical parameters
and domain geometry chosen based on an experiment setting.

Numerical results have shown that the v-d split model can lead to unphysical crack propagation
for some physical settings for which the spectral decomposition model gives correct crack propagation;
e.g., see Example 2 (a single edge notched shear test) and particularly Fig.~\ref{fig:shear a's OBC}.
An explanation for this is that the v-d split model includes the total deviatoric strain energy in the active energy
and thus accounts both the compressive and expansive deviatoric strain energies for contributing
to the damage process. To avoid this, we have proposed an improved v-d model (cf. Section~\ref{SEC:improved})
which degrades only the the expansive volumetric and deviatoric strain energies in the damage zone.
This improved v-d model improves the v-d model and leads to results comparable with those obtained with
the spectral decomposition model.

Numerical results have also shown that all of the models, including the spectral decomposition and
improved v-d models, can still lead to unphysical crack openings and propagation for various physical
settings. A close examination on this has revealed that they do not generally satisfy the vanishing stress condition
(\ref{CBC-1}) and there is stress remaining on the crack surface.
We have proposed a simple yet effective remedy (ItCBC or the Improved treatment of Crack Boundary
Conditions, cf. Section~\ref{SEC:ItCBC}):  define a critically damaged zone with $0 < d < d_{cr}$,
where $d_{cr}$ is a positive parameter, and then degrade both the active and passive components of the energy
in this zone. It has been shown that the spectral decomposition and improved
v-d models with ItCBC lead to correct crack propagation for all of the examples we have tested.
It should be emphasized that they include Example 3 with multiple cracks and Example 4 which is based on
an experimental setting. In the latter case, both the spectral decomposition and improved
v-d models with ItCBC yield comparable crack propagation results that also agree well qualitatively
with the experiment. Moreover, it has been shown that the numerical results are not very sensitive
to the choice of $d_{cr}$ although $d_{cr}$ in the range of [0.4, 0.6] seems to work best.

Finally, the numerical examples have demonstrated that the MMPDE moving mesh method
is able to dynamically concentrate the mesh elements around propagating cracks even for
complex crack systems.

\vspace{20pt}

\noindent
{\bf Acknowledgments.} 

F.Z. was supported by China Scholarship Council (CSC) and China University of Petroleum - Beijing (CUPB)
for his research visit to the University of Kansas from September of 2015 to September of 2017.
F.Z. is thankful to Department of Mathematics of the University of Kansas for the hospitality during his visit.


\begin{thebibliography}{10}

\bibitem{AMM09}
H.~Amor, J.-J.~Marigo, and C.~Maurini.
\newblock Regularized formulation of the variational brittle fracture with unilateral contact:
Numerical experiments.
\newblock {\em J. Mech. Phys. Solids.}, 57:1209--1229, 2009.

\bibitem{Babuska92}
I.~Babu$\check{s}$ka and M.~Suri.
\newblock Locking effects in the finite element approximation of elasticity problems.
\newblock {\em Numer. Math.}, 62:439-463, 1992.

\bibitem{Bor12}
M.~J. Borden.
\newblock {\em Isogeometric Analysis of Phase-field Models for Dynamic Brittle
  and Ductile Fracture}.
\newblock The University of Texas at Austin, 2012.
\newblock PhD thesis.

\bibitem{BHLV14}
M.~J. Borden, T.~J.~R. Hughes, C.~M. Landis, and C.~V. Verhoosel.
\newblock A higher-order phase-field model for brittle fracture: formulation
  and analysis within the isogeometric analysis framework.
\newblock {\em Comput. Methods Appl. Mech. Eng.}, 273:100--118, 2014.

\bibitem{BFM00}
B.~Bourdin, G.~A. Francfort, and J.~J. Marigo.
\newblock Numerical experiments in revisited brittle fracture.
\newblock {\em J. Mech. Phys. Solids}, 48:797--826, 2000.

\bibitem{BHR09}
C.~Budd, W.~Huang, and R.~Russell.
\newblock Adaptivity with moving grids.
\newblock {\em Acta Numerica}, 18:111--241, 2009.

\bibitem{FM98}
G.~A. Francfort and J.~J. Marigo.
\newblock Revisiting brittle fracture as an energy minimization problem.
\newblock {\em J. Mech. Phys. Solids}, 46:1319--1342, 1998.

\bibitem{Hua01}
W.~Huang.
\newblock Variational mesh adaptation: isotropy and equidistribution.
\newblock {\em J. Comput. Phys.}, 174:903--924, 2001.

\bibitem{Hua06}
W.~Huang.
\newblock Mathematical principles of anisotropic mesh adaptation.
\newblock {\em Comm. Comput. Phys.}, 1:276--310, 2006.

\bibitem{HK15a}
W.~Huang and L.~Kamenski.
\newblock A geometric discretization and a simple implementation for
  variational mesh generation and adaptation.
\newblock {\em J. Comput. Phys.}, 301:322--337, 2015.

\bibitem{HK15b}
W.~Huang and L.~Kamenski.
\newblock On the mesh nonsingularity of the moving mesh pde method.
\newblock {\em Math. Comp.}, 87:1887--1911, 2018.

\bibitem{HRR94a}
W.~Huang, Y.~Ren, and R.~D. Russell.
\newblock Moving mesh methods based on moving mesh partial differential
  equations.
\newblock {\em J. Comput. Phys.}, 113:279--290, 1994.

\bibitem{HRR94b}
W.~Huang, Y.~Ren, and R.~D. Russell.
\newblock Moving mesh partial differential equations (MMPDEs) based upon the
  equidistribution principle.
\newblock {\em SIAM J. Numer. Anal.}, 31:709--730, 1994.

\bibitem{HR11}
W.~Huang and R.~D. Russell.
\newblock {\em Adaptive Moving Mesh Methods}.
\newblock Springer, New York, 2011.
\newblock Applied Mathematical Sciences Series, Vol. 174.

\bibitem{KM10}
C.~Kuhn and R.~Muller.
\newblock A continuum phase field model for fracture.
\newblock {\em Eng. Frac. Mech.}, 77:3625--3634, 2010.

\bibitem{LRHW16}
S.~Lee, J. E.~Reber, N. W.~Hayman, and M. F.~Wheeler.
\newblock Investigation of wing crack formation with a combined phase-field
and experimental approach.
\newblock {\em Geophysical Research Letters.}, 43:7946--7952, 2016.

\bibitem{MVB15}
S.~May, J.~Vignollet, and R.~de~Borst.
\newblock A numerical assessment of phase-field models for brittle and cohesive
  fracture: $\Gamma$-convergence and stress oscillations.
\newblock {\em European J. Mech. A/Solids.}, 52:72--84, 2015.

\bibitem{MHW10}
C.~Miehe, M.~Hofacker, and F.~Welschinger.
\newblock A phase field model for rate-independent crack propagation: {R}obust
  algorithmic implementation based on operator splits.
\newblock {\em Comput. Methods Appl. Mech. Eng.}, 199:2765--2778, 2010.

\bibitem{MWH10}
C.~Miehe, F.~Welschinger, and M.~Hofacker.
\newblock Thermodynamically consistent phase-field models of fracture:
  Variational principles and multi-field {FE} implementations.
\newblock {\em Int. J. Numer. Meth. Eng.}, 83:1273--1311, 2010.

\bibitem{RLH15}
J.E.~Reber, L.L.~Lavier, and N.W.~Hayman.
\newblock Experimental demonstration of a semi-brittle origin for crustal strain transients.
\newblock {\em Nat. Geosci.}, 8:712--715, 2015.

\bibitem{SBS15}
D.~Schillinger, M. J.~Bordenr, and H. K.~Stolarski.
\newblock Isogeometric collocation for phase-field fracture models.
\newblock {\em Comput. Methods Appl. Mech. Eng.}, 284:583--610, 2015.

\bibitem{SS15}
M.~Stronl and T.~Seelig.
\newblock A novel treatment of crack boundary conditions in phase field models of fracture.
\newblock {\em Appl. Math. Mech.}, 15:155--156, 2015.

\bibitem{SS16}
M.~Stronl and T.~Seelig.
\newblock On constitutive assumptions in phase field approaches to brittle fracture.
\newblock {\em 21st European Conference on Fracture}, pages 3705--3712, 2016.
\newblock Catania, Italy.

\bibitem{VMBV14}
J.~Vignollet, S.~May, R.~de~Borst, and C.~V. Verhoosel.
\newblock Phase-field models for brittle and cohesive fracture.
\newblock {\em Meccanica}, 49:2587--2601, 2014.

\bibitem{ZHLZ17}
F.~Zhang, W.~Huang, X.~Li, and S. Zhang.
\newblock Moving mesh finite element simulation for phase-field modeling of brittle fracture
 and convergence of Newton's iteration.
\newblock {\em J. Comput. Phys.}, 356:127--149, 2018.

\end{thebibliography}

\end{document}